\documentclass{article}
\usepackage{graphicx}
\usepackage{amsmath}
\usepackage{amsfonts}
\usepackage{amsthm}
\usepackage{subfigure}
\usepackage{graphicx}

\newtheorem{theorem}{Theorem}

\newtheorem{corrolary}{Corrolary}

\newtheorem{remark}{Remark}

\newtheorem{definition}{Defining process}
\newtheorem{odefinition}{Definition}
\newtheorem{alemma}{Lemma}

\begin{document}

\title{On the multiresolution structure of Internet traffic traces\footnote{A preliminary version of this paper appeared in Proceedings of SPIE Vol. 4868 (2002) under the title: ``Some results on the multiresolution structure of Internet traffic traces''.}}
\author{Konstantinos Drakakis\footnote{The author is a scholar of the Lilian Boudouris Foundation.} \\
        Princeton University
        \and
        Dragan Radulovi\'c\footnote{The largest portion of the research for this paper was conducted 
        while the author was affiliated to Princeton University, and was partially supported 
        by AT\&T Research Center.}\\
        Yale University}         
\maketitle

\begin{abstract}
Internet traffic on a network link can be modeled as a stochastic process. After detecting and quantifying the 
properties of this process, using statistical tools, a series of mathematical models is developed, culminating in one 
that is able to generate ``traffic'' that exhibits --as a key feature--
the same difference in behavior for different time scales, as observed in real traffic, and is moreover indistinguishable from real traffic
by other statistical tests as well. Tools inspired from the models are then used to determine and calibrate
the type of activity taking place in each of the time scales. Surprisingly, the above procedure does \emph{not} require any detailed 
information originating from either the network dynamics, or the decomposition of the total traffic into its constituent 
user connections, but rather only the compliance of these connections to very weak conditions. 

\end{abstract}

\section{Introduction}

\subsection{The problem}

\label{problem}

The object of this paper is the study of the behavior of Internet traffic,
as observed on a network link (typically shared by a large number of
network users), aiming to the identification, quantification, and
justification of its salient features. Traffic observation here assumes the
form of \emph{data traces}, i.e. sequences of the form $\left\{ \left(
t_{i},d_{i}\right) \right\} _{i=1}^{N}$, where $d_{i}$ is a data quantity
(possibly measured in bits, bytes etc.), and $t_{i}$ is its \emph{time stamp}%
, the time when $d_{i}$ was observed. Throughout the paper, though, it will
be more convenient to use \emph{binned} data traces $X_{i}^{\Delta
},i=1,...,M=\left\lfloor \frac{t_{N}}{\Delta }\right\rfloor $, with bin size 
$\Delta $, defined as: $X_{i}^{\Delta }=\sum_{\left\{ j|i\Delta
<t_{j}<\left( i+1\right) \Delta \right\} }d_{j}$. Our ignorance of the exact
procedures taking place in the network, due either to their complexity,
or to lack of information, as well as of the behavior of its users, forces
us to view $X^{\Delta }$ as a stochastic process \cite{TTW1}.

This process has been the object of study of numerous
researchers. It has unusual and initially unexpected properties, such as
long range dependence \cite{ENW1,RE1}, nonstationarity \cite{ZPS1,CCLS1}, 
and different behavior on different time scales \cite{FGW1,FGWK1,RLV1,FGHW1}.
Moreover, even for the coarser time scales,
its marginal distributions may be neither Gaussian, nor p-Stable.
This is all the more surprising, since the Central Limit Theorems would lead one to expect the opposite, 
because a) the observed traffic is the aggregate result of
many users, possibly acting independently \cite{TWS1,LT1,RRCB1}
b) each $X_{i}^{\Delta }$ is actually the sum of a large number of $d_{i}$%
's, and c) for small time bins, traffic is uncorrelated.

\begin{figure}
 \centering
 \subfigure[]{\includegraphics[height=150pt, width=200pt]{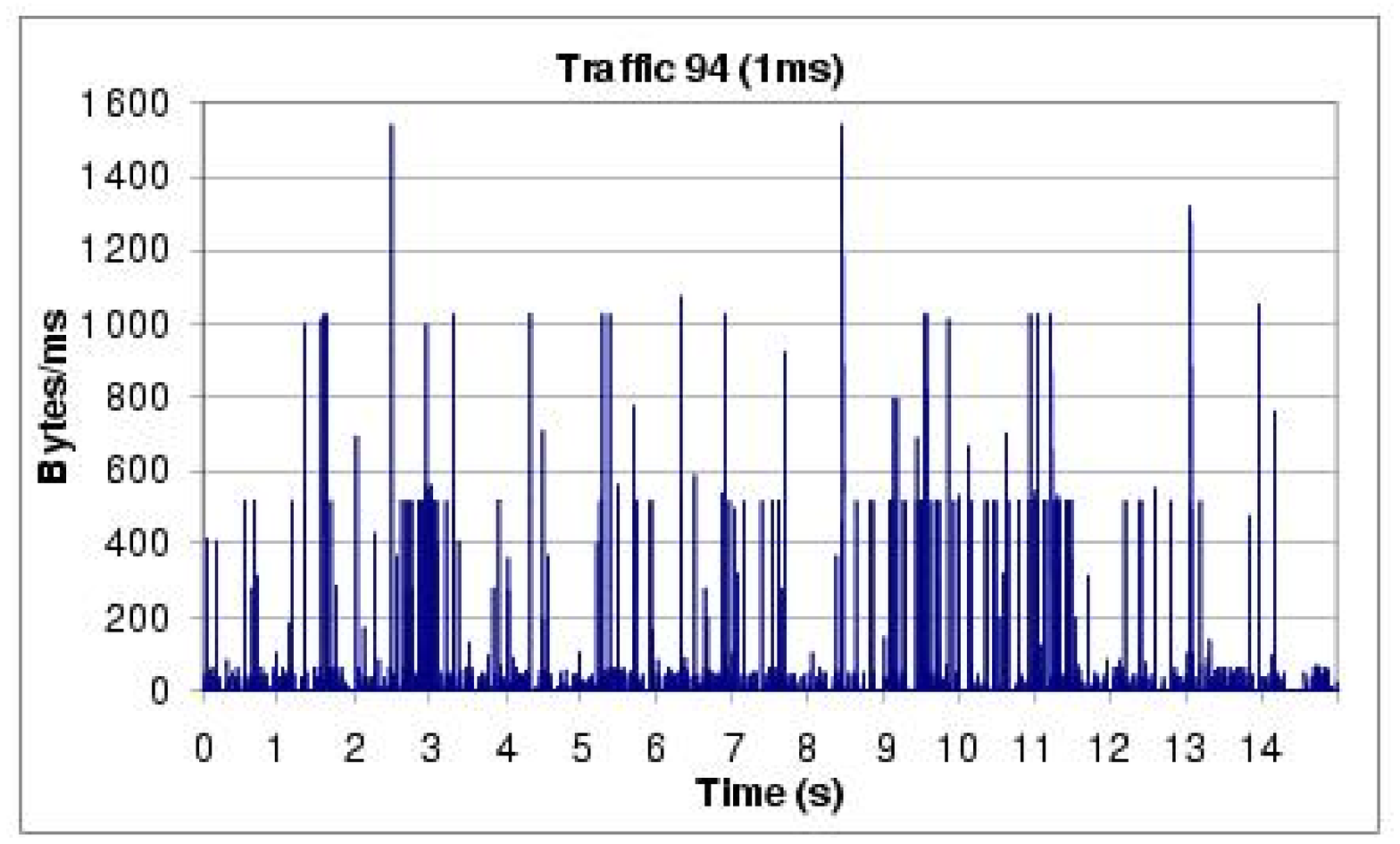}}
 \subfigure[]{\includegraphics[height=150pt, width=200pt]{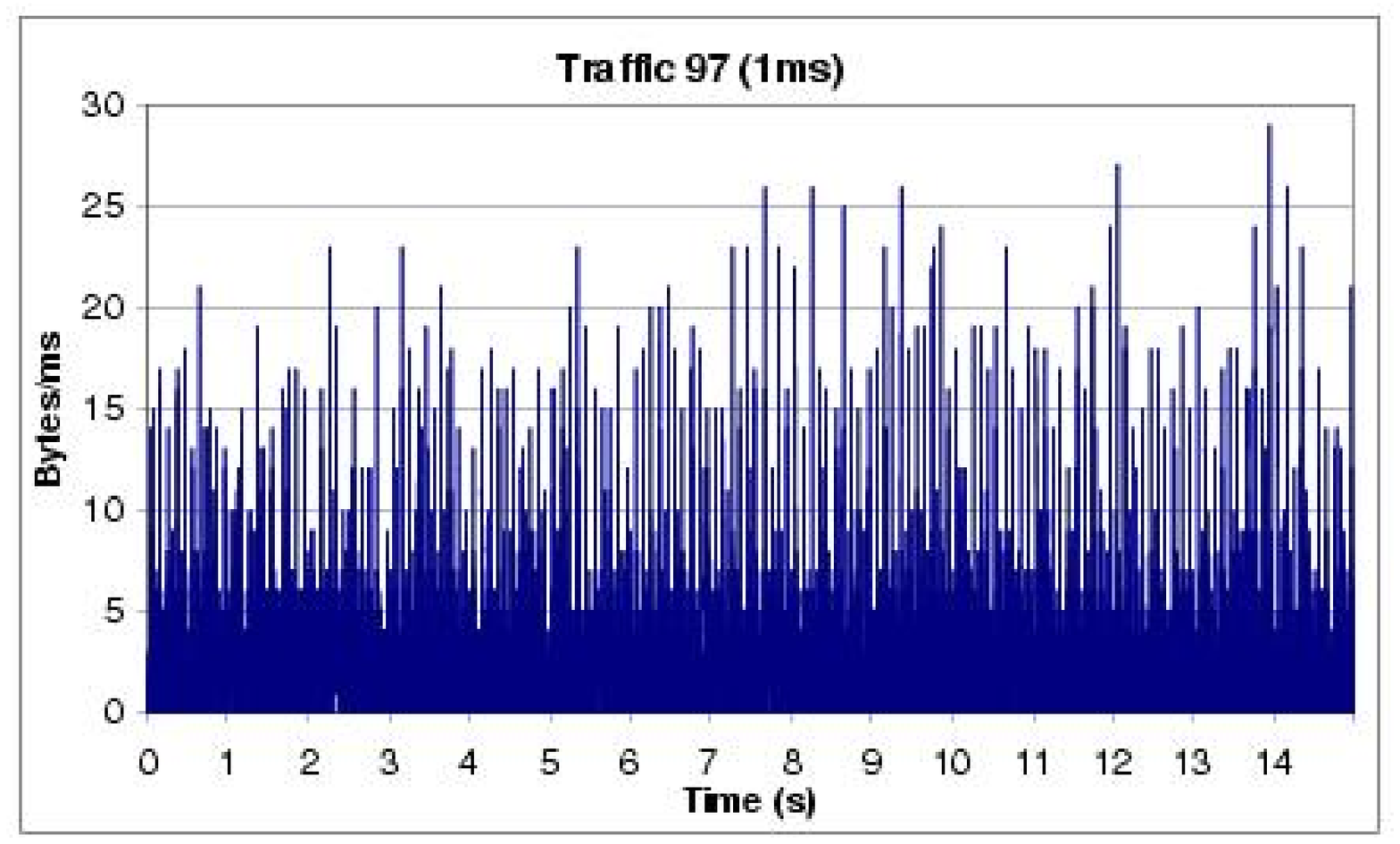}}
 \subfigure[]{\includegraphics[height=150pt, width=200pt]{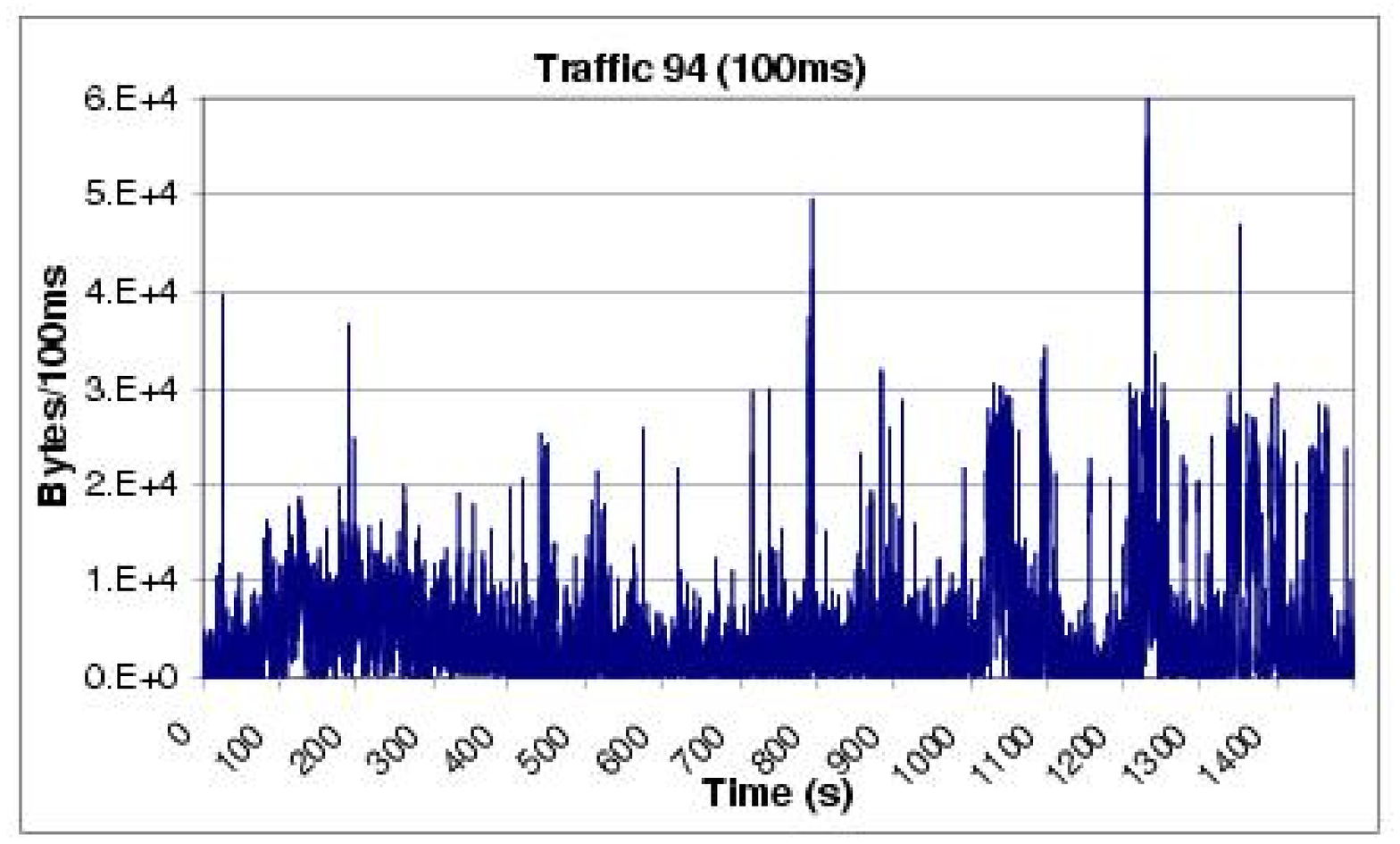}}
 \subfigure[]{\includegraphics[height=150pt, width=200pt]{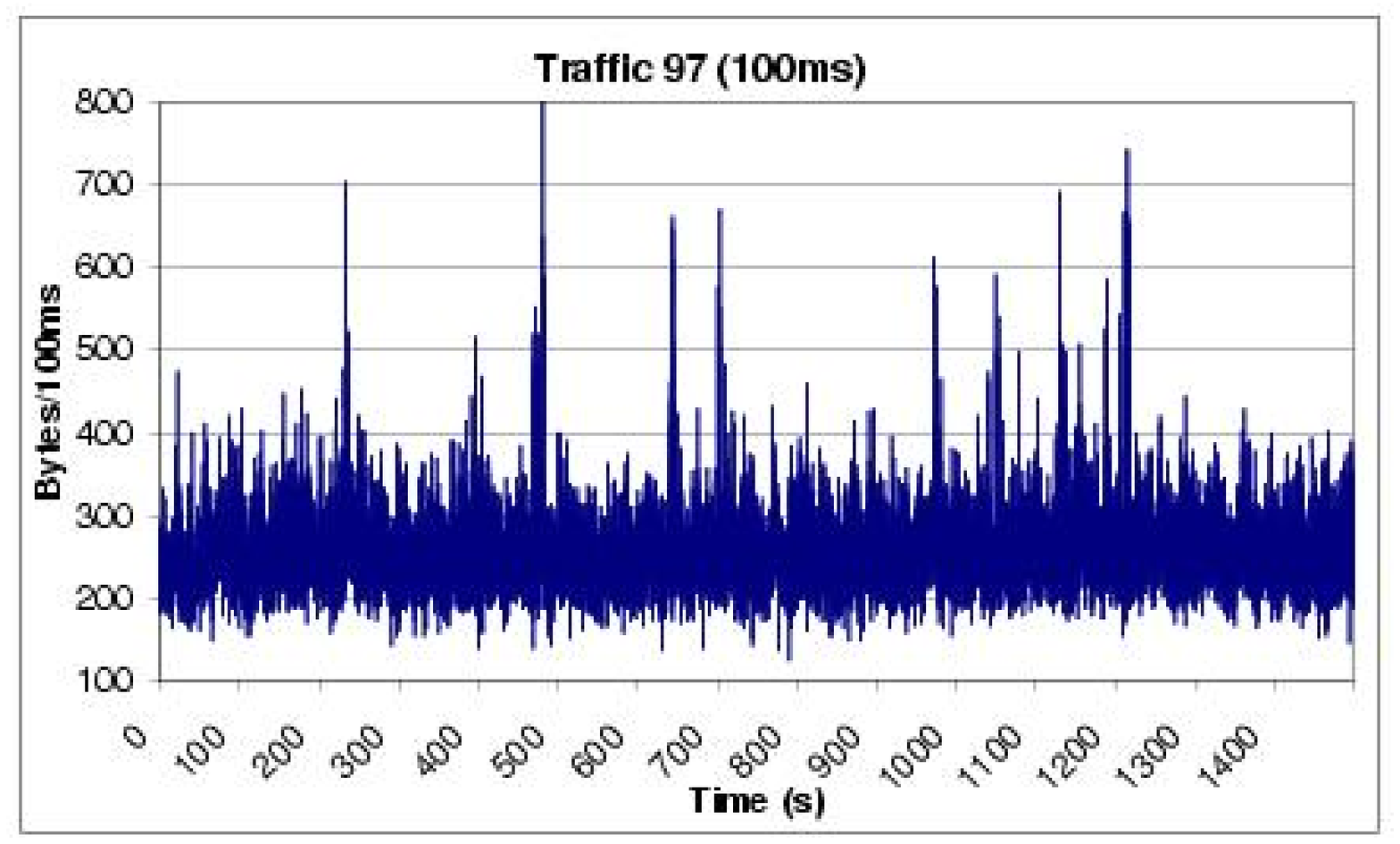}}
 \subfigure[]{\includegraphics[height=150pt, width=200pt]{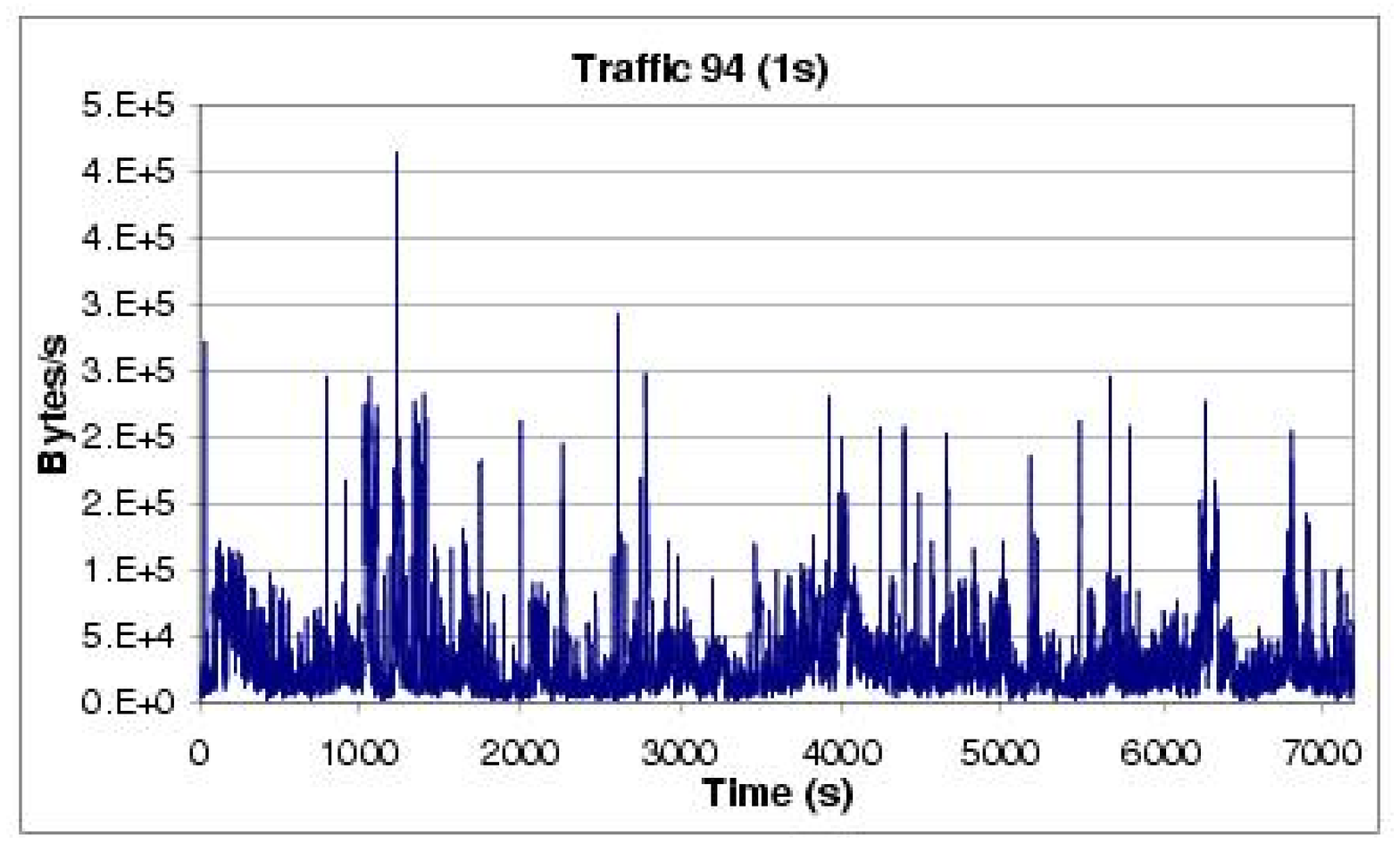}}
 \subfigure[]{\includegraphics[height=150pt, width=200pt]{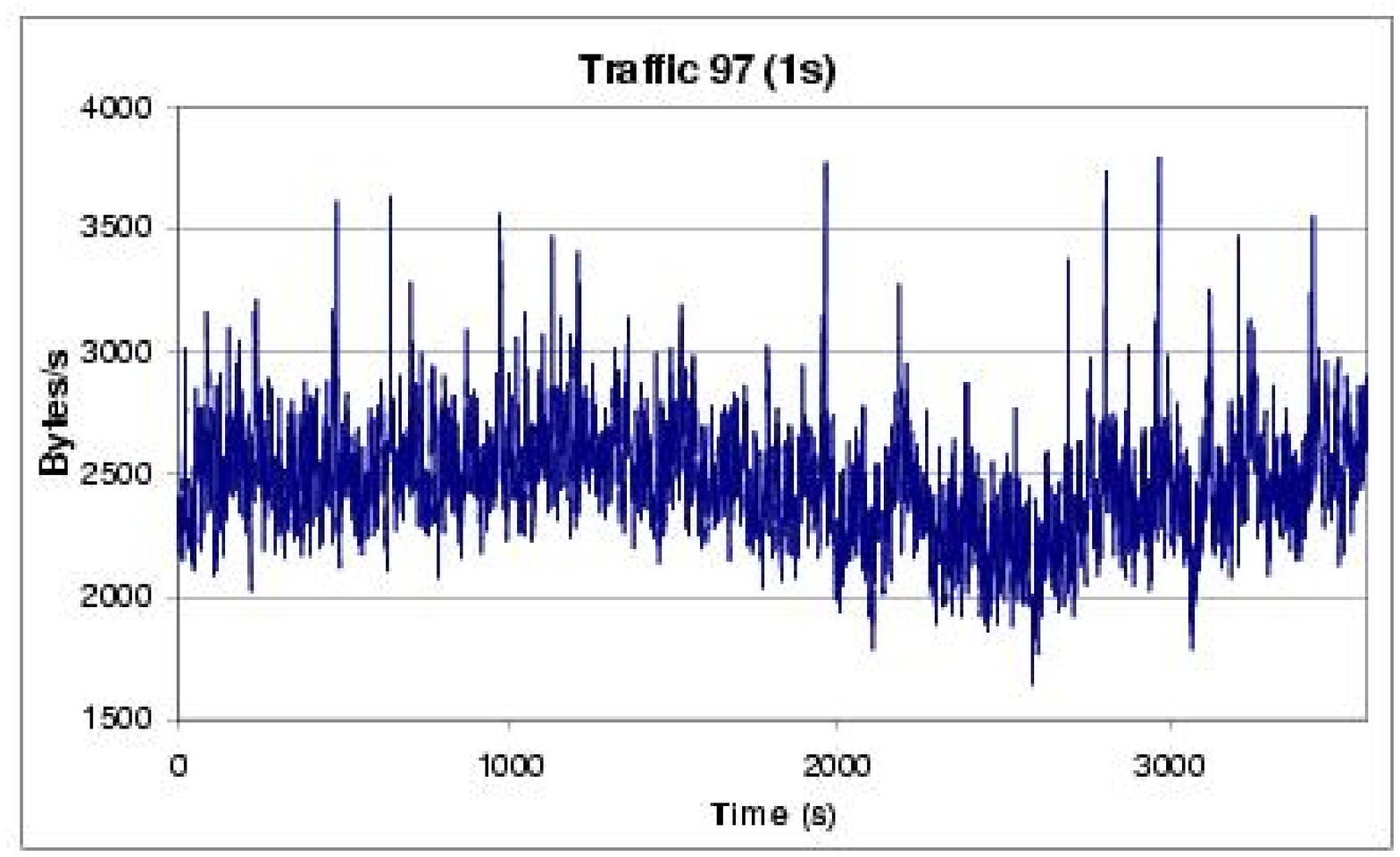}}
 \caption{
  Traces 94 -(a),(c),(e)- and 97 -(b),(d),(f)- binned with 3 different 
  time bins: 1ms,100ms,1s. The data sets 94 and 97 are two of the eight data sets used in this paper; for more details,
  see section 1.2 below.
  }
 \label{Traffic9497}
\end{figure}

Indeed, despite the fact that ``bursts'' and ``spikes'' seem to persist for
at least 4 orders of bin magnitude (Fig. \ref{Traffic9497} and \ref{TrafficM6U8}), which argues against
a Gaussian limit, it can be easily seen, using classical fitting
techniques, that the marginal distributions do not possess the characteristic heavy tails of p-Stable variables either;
instead, their tails appear to be exponential or Weibull (Fig. \ref{NoGauss}). Thus, the process seems
``unwilling'' to converge to any limit, in contrast to most of the processes observed
in nature or industry, which do not exhibit this type of ``wild'' behavior, but
rather get ``attracted'' to the (Gaussian or p-Stable) central limit much
sooner (this behavior will be quantified more precisely in section \ref{modelC} below). 

\addtocounter{figure}{-1}
\stepcounter{figure} 

\begin{figure}
 \centering
 \subfigure[]{\includegraphics[height=150pt, width=200pt]{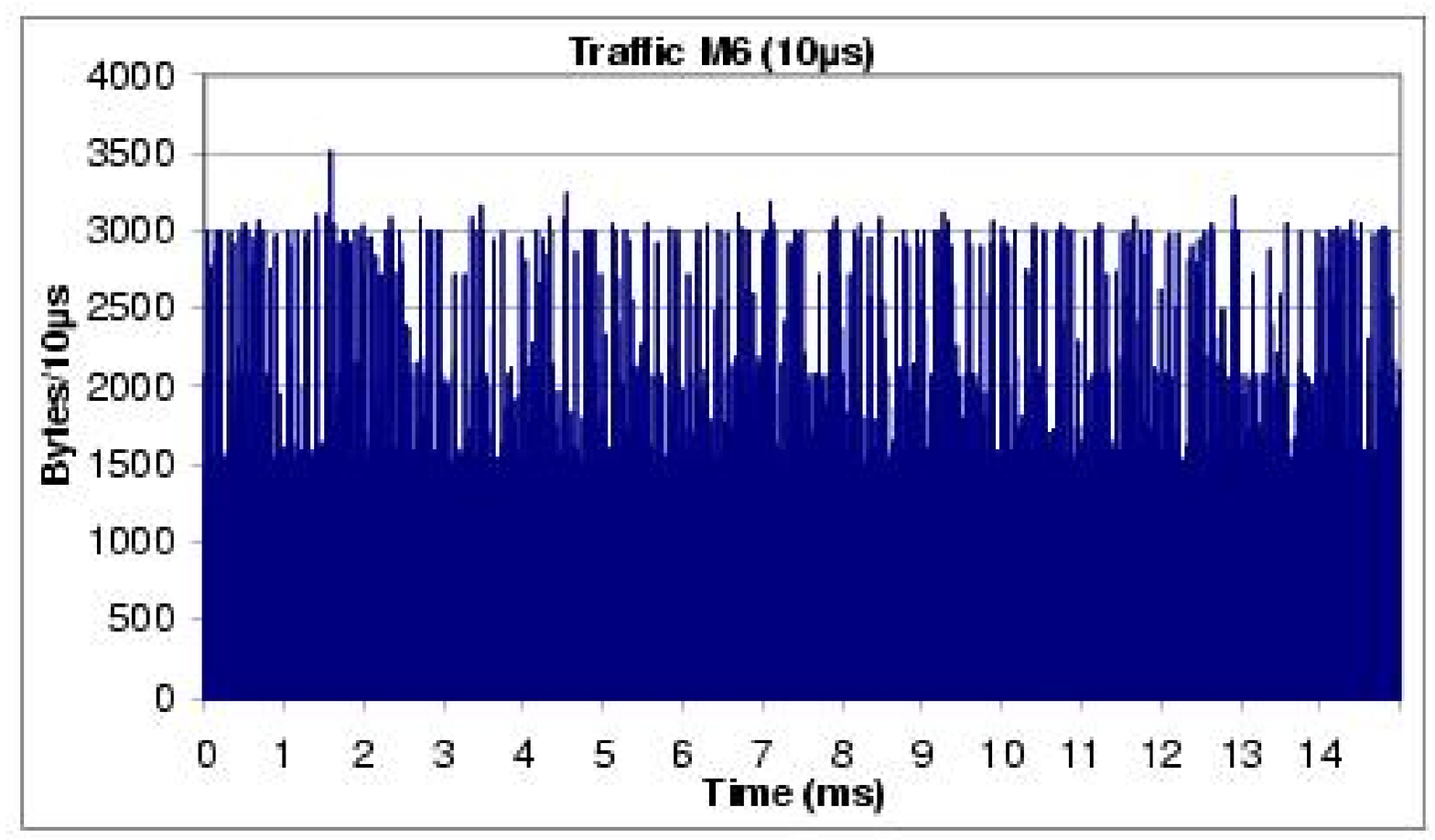}}
 \subfigure[]{\includegraphics[height=150pt, width=200pt]{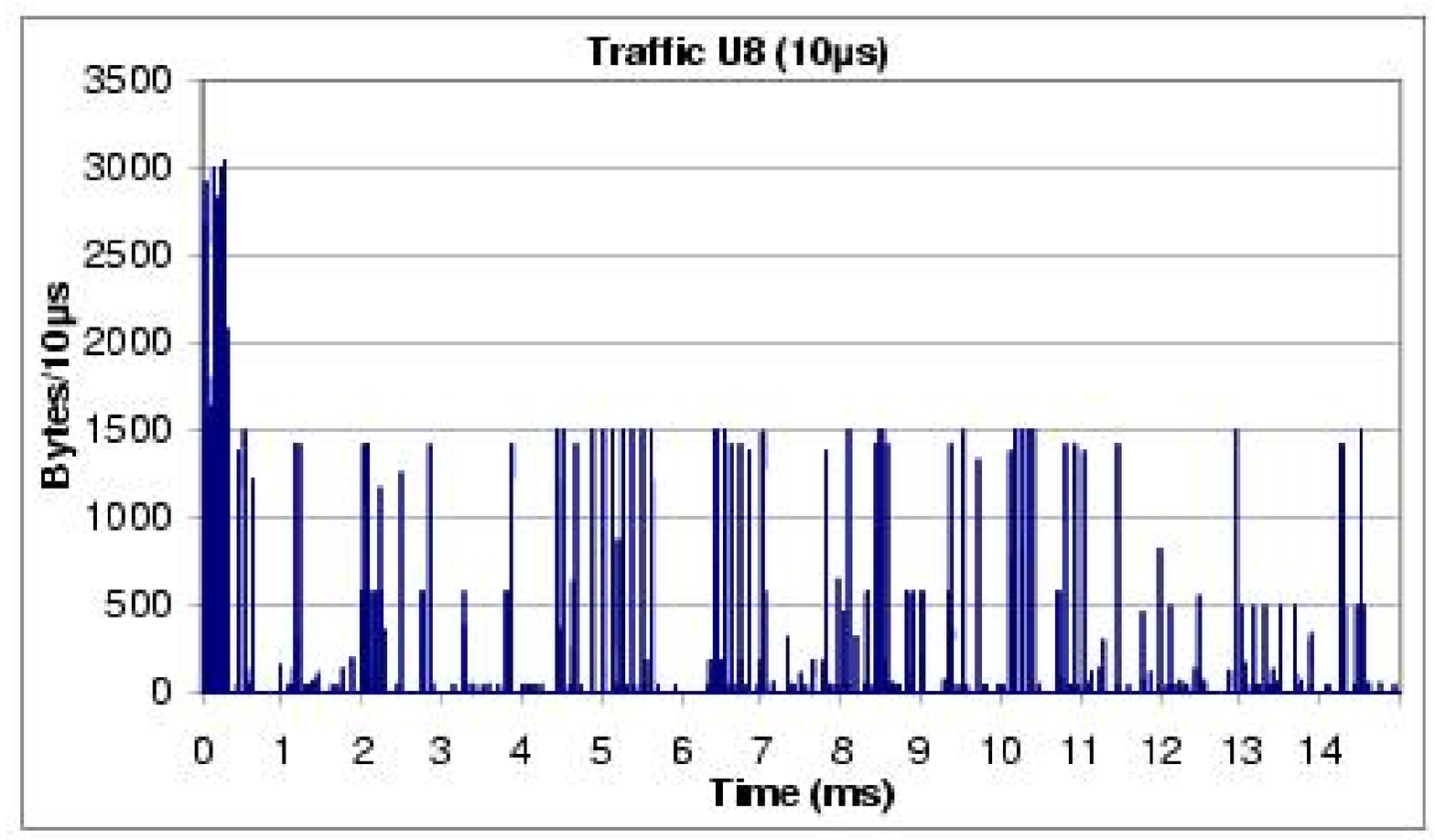}}
 \subfigure[]{\includegraphics[height=150pt, width=200pt]{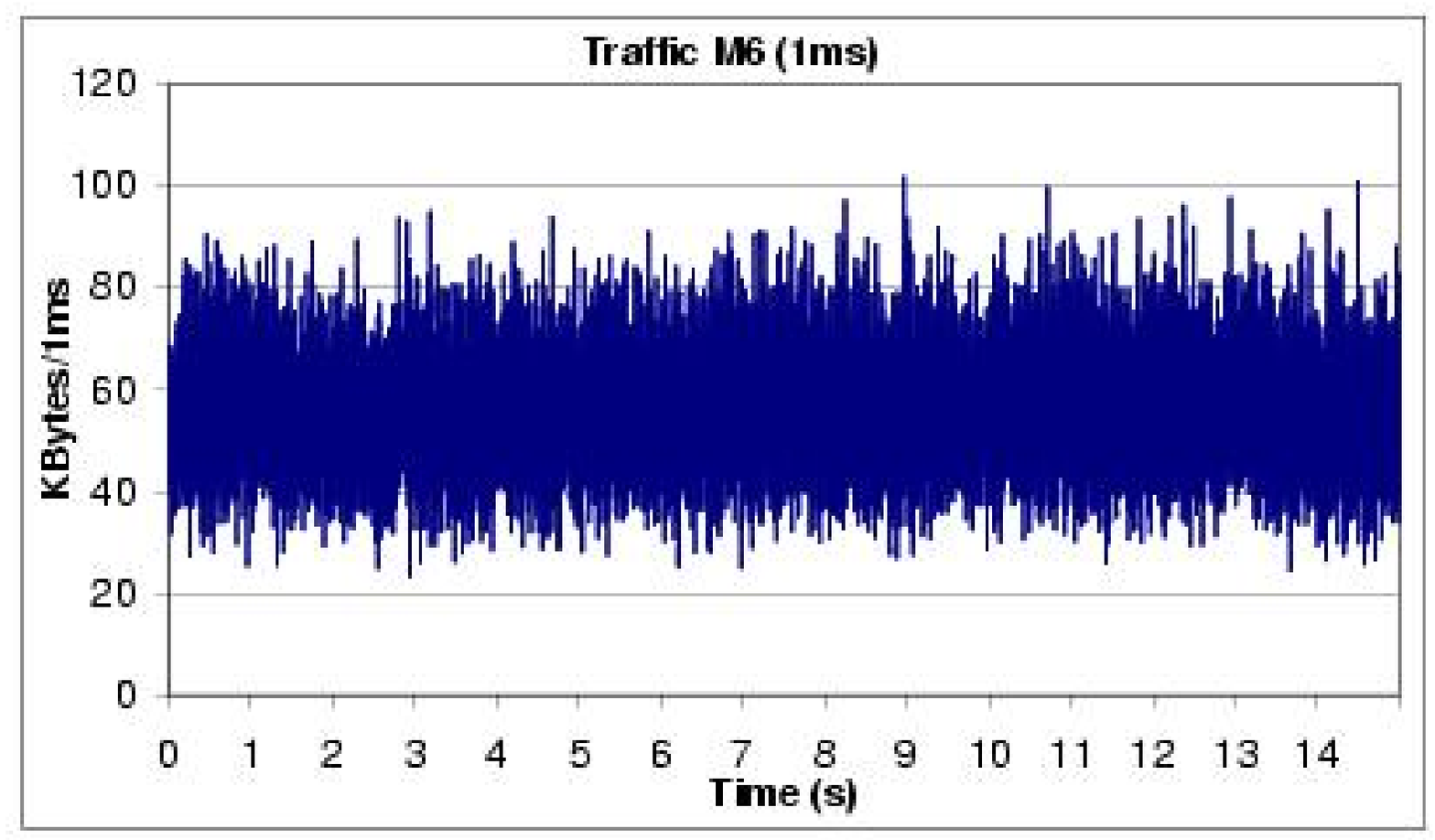}}
 \subfigure[]{\includegraphics[height=150pt, width=200pt]{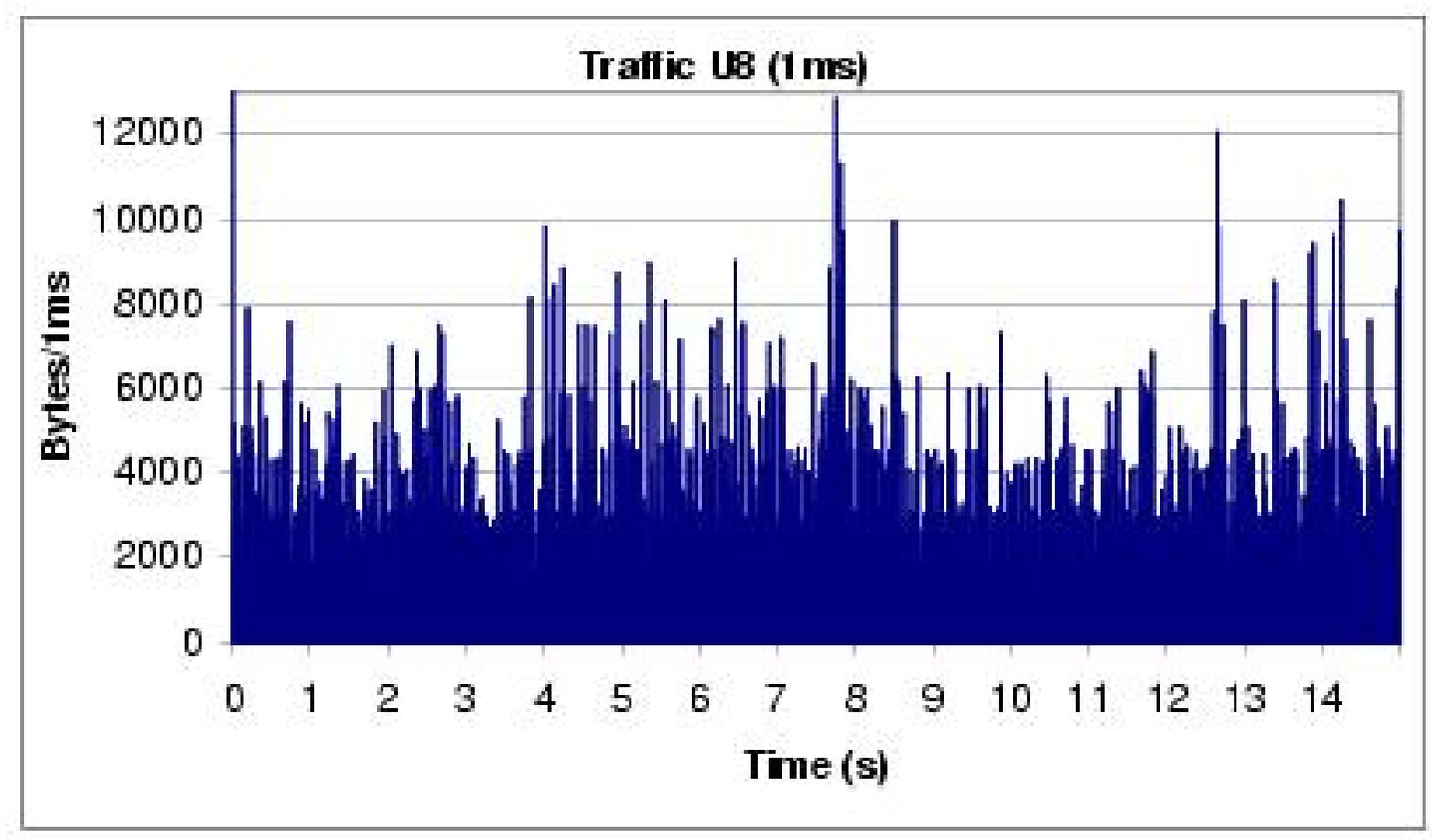}}
 \subfigure[]{\includegraphics[height=150pt, width=200pt]{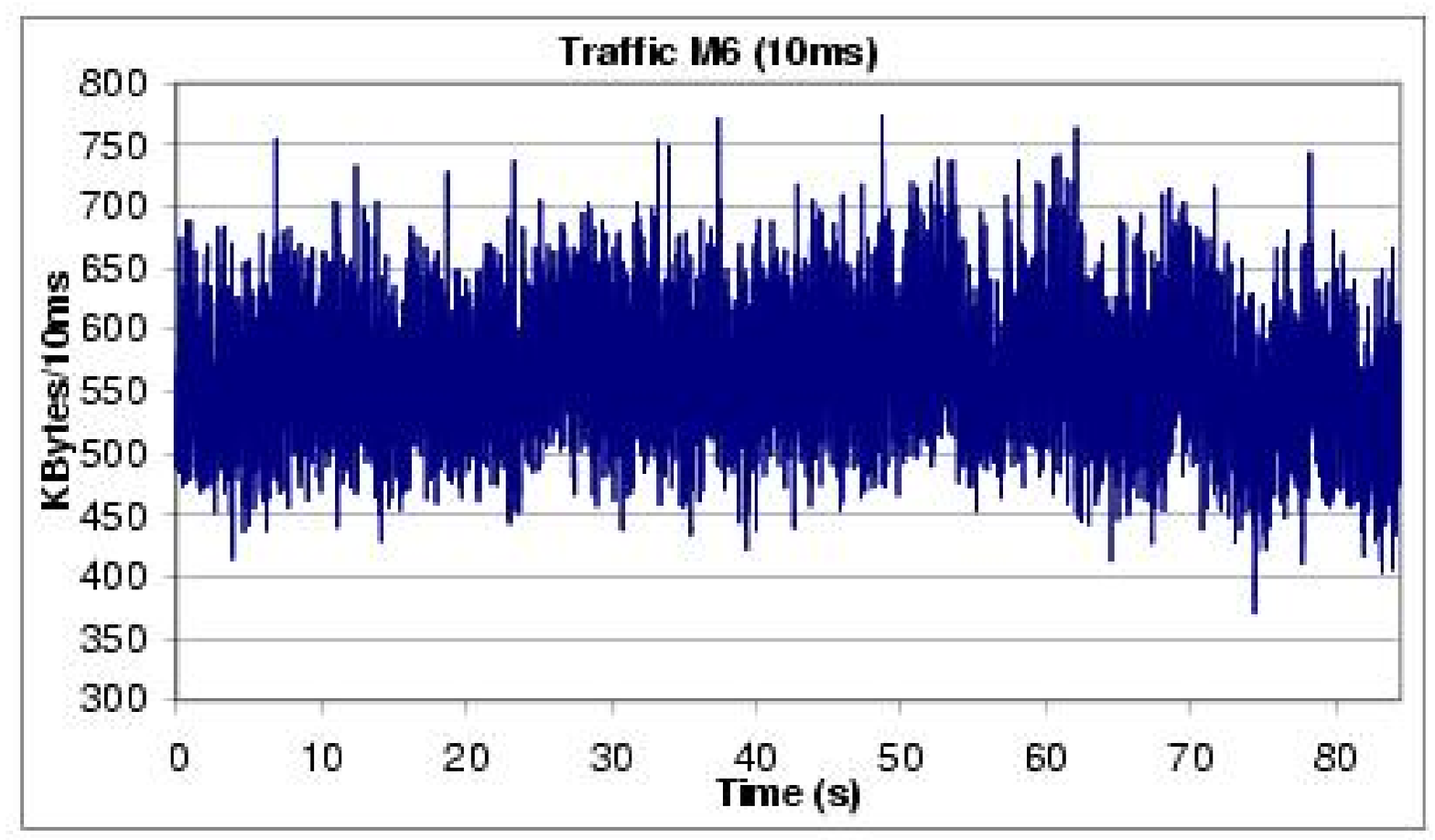}}
 \subfigure[]{\includegraphics[height=150pt, width=200pt]{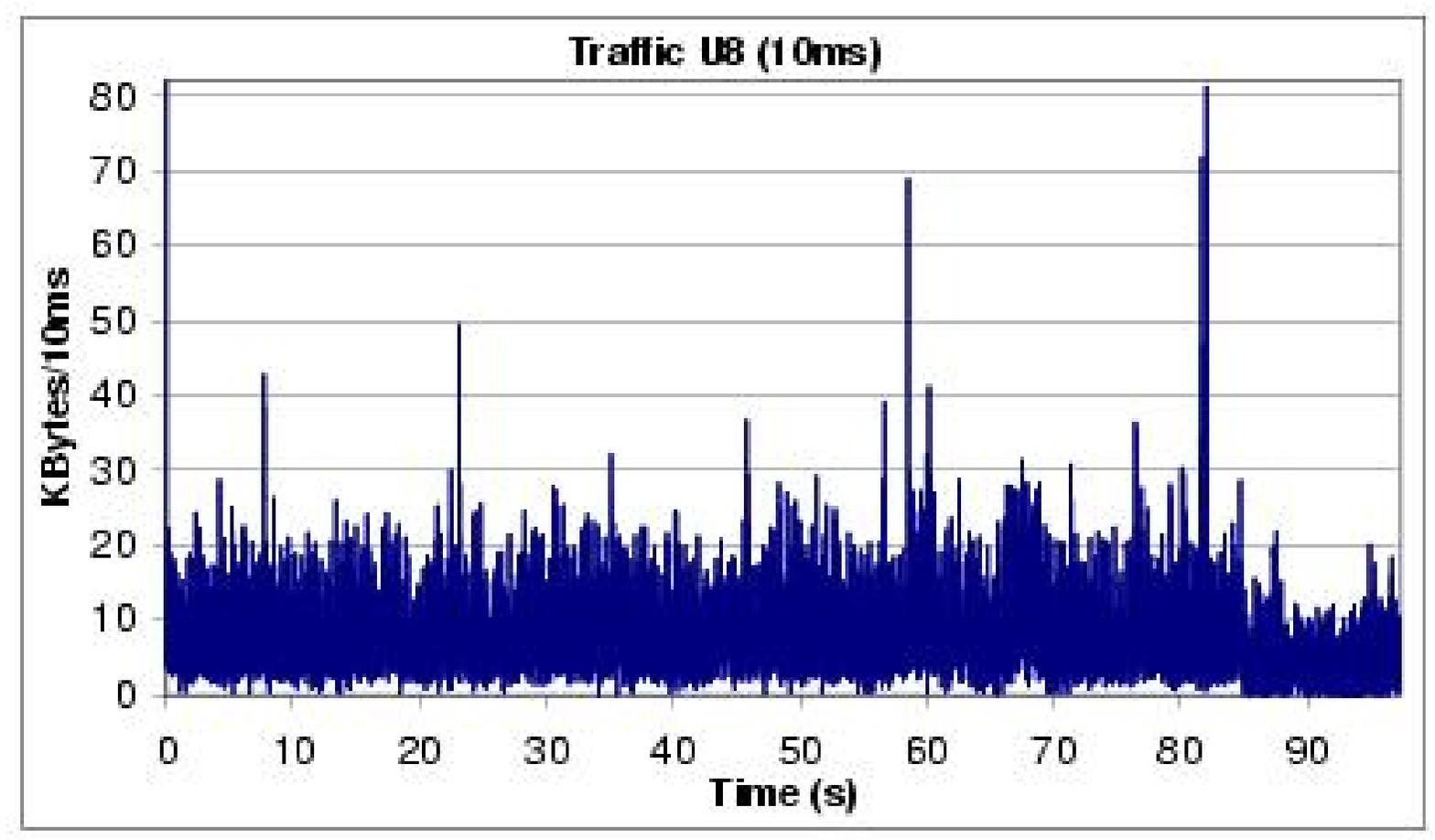}}
 \caption{
  Traces M6 -(a),(c),(e)- and U8 -(b),(d),(f)- binned with 3 different 
  time bins: 10$\mu$s,1ms,10ms. The data sets 94 and 97 are two of the eight data sets used in this paper; 
  for more details, see section 1.2 below.
  }
 \label{TrafficM6U8}
\end{figure}
\addtocounter{figure}{-1}
\stepcounter{figure} 

\begin{figure}
 \centering
 \subfigure[]{\includegraphics[height=150pt, width=200pt]{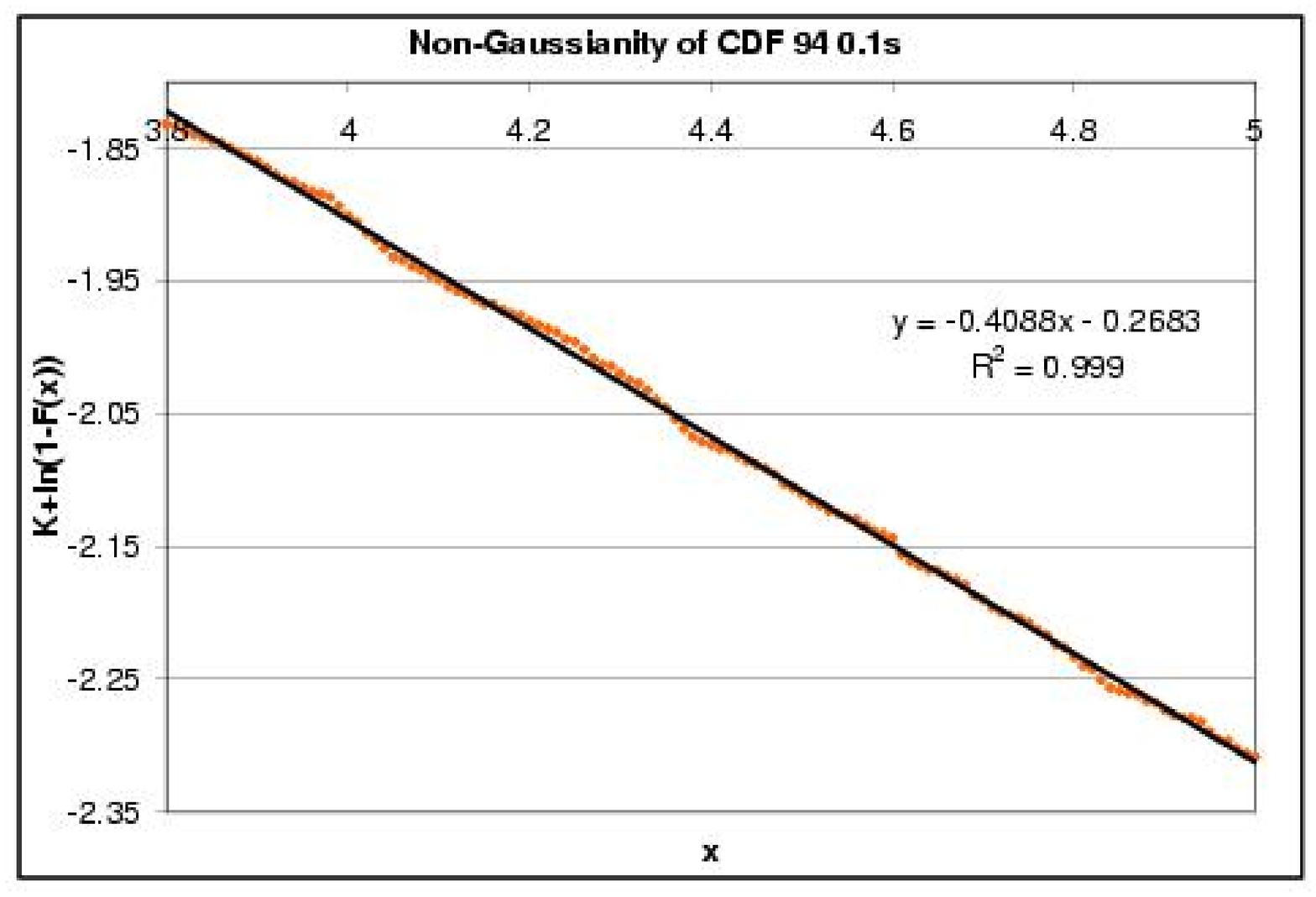}}
 \subfigure[]{\includegraphics[height=150pt, width=200pt]{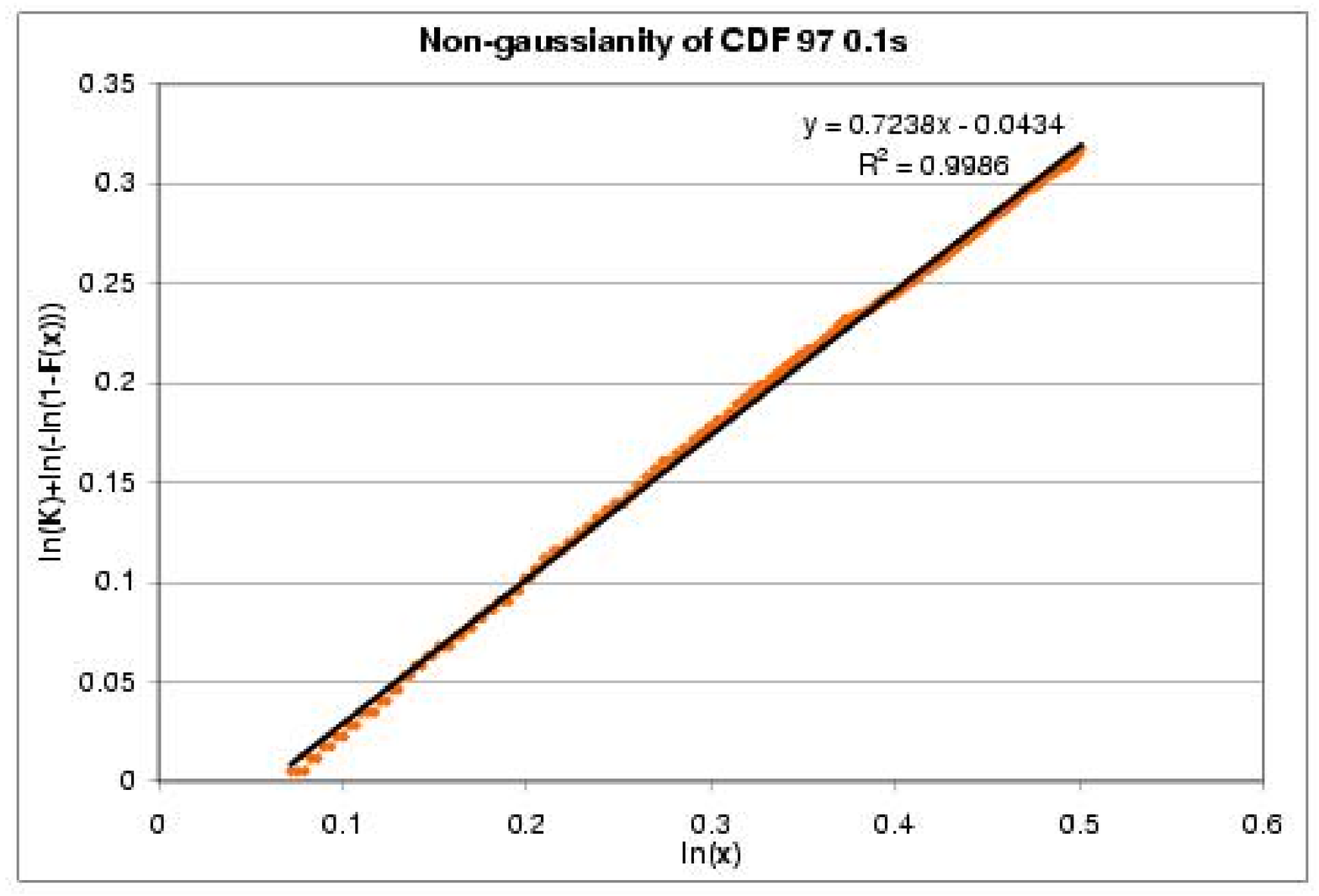}}
 \caption{
  Demonstration of the non-Gaussianity of the marginal CDFs of sets 94 -(a)-, and 97 -(b)-, through least squares  
  fitting: the 
  very high value of R in both cases shows that the curve is essentially linear, and thus that the tail in (a) is 
  exponential and in (b) Weibull. Both sets have been binned using a time bin of 0.1s. 
  }
 \label{NoGauss}
\end{figure}
\addtocounter{figure}{-1}
\stepcounter{figure}

The unusual properties of the traffic are generally attributed to the
diversity of user applications, the network protocols, and the complexity of
the network topology. Indeed, Internet traffic today is a result of numerous
applications: e-mail, WWW, FTP, music and video streams, etc. The duration
and the size of the sessions of these applications often vary within several
orders of magnitude. Furthermore, the details of the data transfer over the
network are regulated by  protocols, which control, among other things,
the transfer rate at each time instant, and the fragmentation of the data to
be transferred into packets. Finally, the transfer rate heavily depends on
the state of the network at a particular moment (link availability,
congestion, ISP contracts, routing tables,...), which in turn influences
(and is influenced by) network topology.

The traffic, though, is typically the result of the activity of a large
number of network users. It is expected that many of the features present in
each user's own subtraffic will just average out, leaving no trace in the
aggregate traffic, whereas others may persist; it is not obvious, however,
whether a particular feature will belong to the former, or the latter category.

To sum up, it appears there is a consensus about the factors that shape the
traffic, as well as about the inadequacy of the Gaussian distribution for the 
description of the shape of the traffic. This paper explores how different factors affect the traffic. 
This necessitates the development of tools that will calibrate
the traffic behavior, along with appropriate simulations, on which these
tools will be tested.

\subsection{Description of the data sets and the simulations}

\label{DatSim}

This section will provide a short description of both real and simulated
traces used. Some of the relevant definitions will be given later.

In what follows, eight data sets (traces), coded 89, 94, 97, L4, L5, M6, M7 and U8 will be
used systematically. Data sets 89, 94, L4, and L5 are available on the Web, at
http://ita.ee.lbl.gov/html/traces.html, and data sets M6, M7, and U8 are also available on the Web,
at http://pma.nlanr.net. These Web pages also contain extensive information about the data sets. More precisely:

\begin{itemize}
	\item 89 is the trace BCP-Aug89, and it contains ``one million packet arrivals seen on 
        an Ethernet at the Bellcore Morristown Research and Engineering facility''. It was collected 
        on AUG 29, 1989, starting at 11.25, and lasts approximately 3,143 seconds. It is by far the 
        oldest trace used here. The volume exchanged is approximately 434 MB.
  \item 94 is the trace LBL-TCP-3, collected on JAN 20, 1994, between 14.10 and 16.10. The volume 
        exchanged is approximately 244 MB. 
  \item L4 is the trace LBL-PKT-4, collected on Jan 28, 1994, between 14.00 and 15.00. The volume 
        exchanged is approximately 131 MB. 
  \item L5 is the trace LBL-PKT-5, collected on Jan 28, 1994, between 15.00 and 16.00. The volume 
        exchanged is approximately 94 MB.
  \item 97 contains 1 hour of the total traffic of users subscribing to
        WorldNet via modems, measured at the WorldNet gateway, on JUL 22, 1997,
        starting at 22.00, and using a time bin of 1ms. The volume 
        exchanged is approximately 9 MB.
  \item M6 is the trace MRA-1015686621-1, collected on MAR 09, 2002. Approximately, it lasts for 81 seconds, and the 
        volume exchanged is 4.8 GB.  
  \item M7 is the trace MRA-1015764522-1, collected on MAR 10, 2002. Approximately, it lasts for 82 seconds, and the 
        volume exchanged is 3.7 GB.   
  \item U8 is the trace ODU-1015816310, collected on MAR 11, 2002. Approximately, it lasts for 89 seconds, and the 
        volume exchanged is 76 MB.       
\end{itemize}

For data sets 89 and 97, information on their individual connections is not available. 

The sizes and dates of the traces described span quite a broad range: some of them contain only a few MB of 
information, whereas others contain 1000 times as much; some are extrememly recent, whereas others are more than a 
decade old. This variety will allow for the examination of the properties of the traffic, and for their 
validation, in general or specific network conditions.  

A set of simulated traces will also be used: 0-1 is a self-similar process
with continuous ON/OFF intervals (as in \cite{TWS1}) and will be the outcome of model A below; ARR
(standing for ARRival process) takes into account the
fragmentation of the information to be transmitted into packets (model B below); RH  (standing for Random Heights)
is the reward model, and RH HT (standing for Random Heights + Heavy Tails) is the $\alpha-\beta$ 
traffic model described in \cite{SRB1}; ARRRH (standing for ARRival process with 
Random Heights) is a combination of RH and ARR; EXP IID is a sequence of i.i.d.\ exponential
variables, and, similarly, HT IID is a sequence of heavy-tailed i.i.d.\
variables. Finally, a model that combines ARR with the new idea of ``levels'', introduced 
in section \ref{modelD}, will also be used; its simulations will be labeled by the time
scales of the levels, and an S will be added if the levels are ``sharp'', e.g. 8, 12S, 7S/12/17 etc
(see section \ref{modelD}).

\subsection{The State of the Art}

As a result of the intensive research on Internet traffic, many models
have been proposed. A requirement imposed was that they be \emph{%
parsimonious }\cite{ENW1}: they should use a very limited set of parameters,
as opposed to the millions of details present in a real network. In what
follows, some of the most common ones will be presented and discussed.

\bigskip

\noindent \textbf{The Self-Similar Model (SSM)}: This is perhaps the simplest model
conceptually. Its main assumption is that each of the $n$ users
generates a data stream or connection $W_{i}(t)$ , $i=1,...,n,$ the traffic being their
sum: $X\left( t\right) =\sum_{i=1}^{n}W_{i}(t)$, $t\in \mathbf{R}^{+}$. The
users are supposed to act independently, and the data streams are composed
of intervals of value $0$ (where the user is silent), alternating with intervals of value 
$1$ (where the user sends data), whose lengths are random,
independent, and cross-independent. The degrees of freedom of this model
then are just two: the distributions of the lengths of the aforementioned
intervals. If chosen appropriately, e.g. \emph{heavy-tailed}\footnote{%
A distribution will be henceforth called \emph{heavy-tailed} iff
its variance is infinite, otherwise it will be called \emph{light-tailed}.
This definition differs slightly from the one implied in \cite{NSSW1}.} 
distributions are allowed, they can lead to \emph{self-similar} simulated traffic \cite{TWS1}
For a more precise statement, see Theorem A (Section \ref{modelA}).

\begin{figure}
 \centering
 \subfigure[]{\includegraphics[height=150pt, width=200pt]{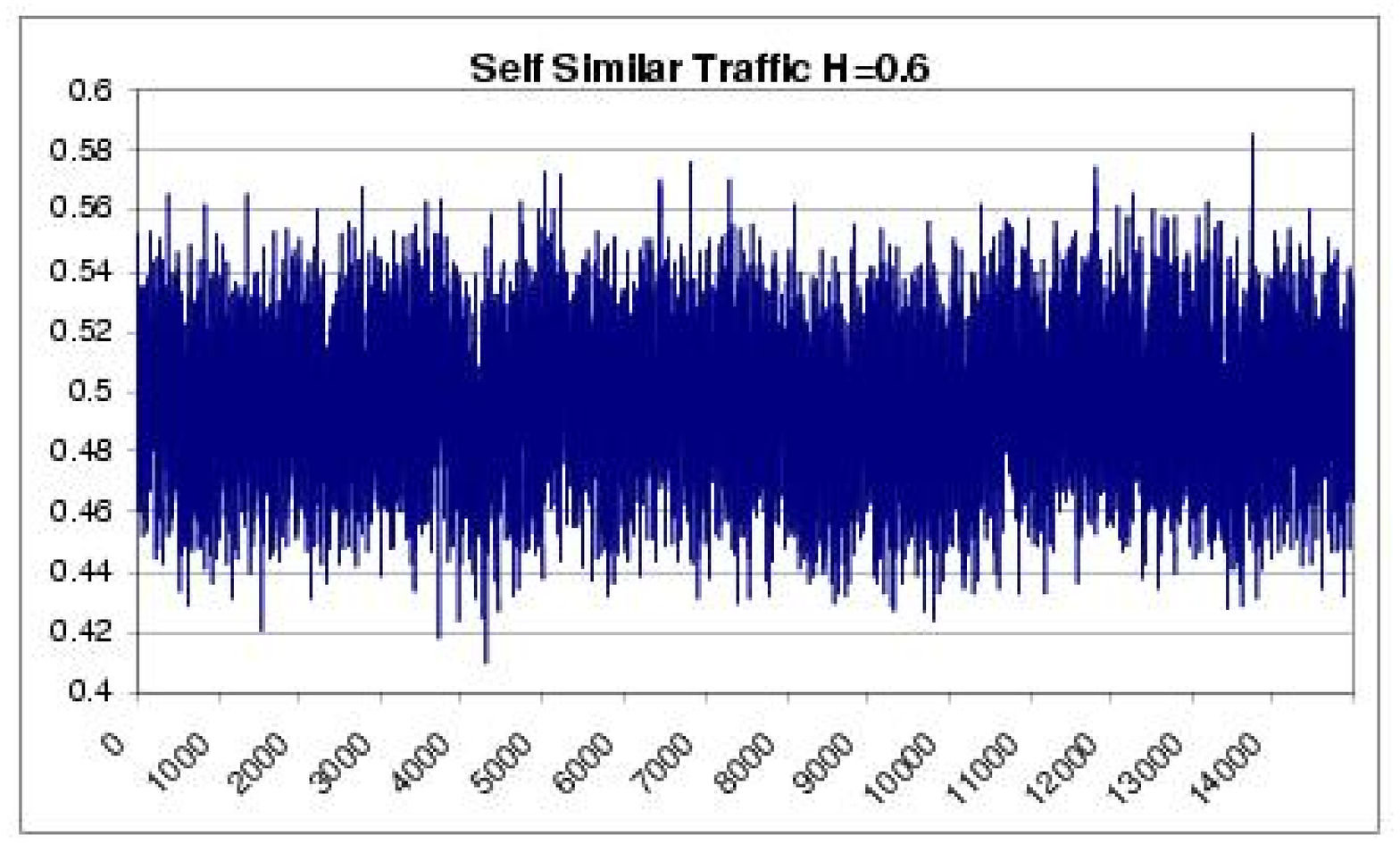}}
 \subfigure[]{\includegraphics[height=150pt, width=200pt]{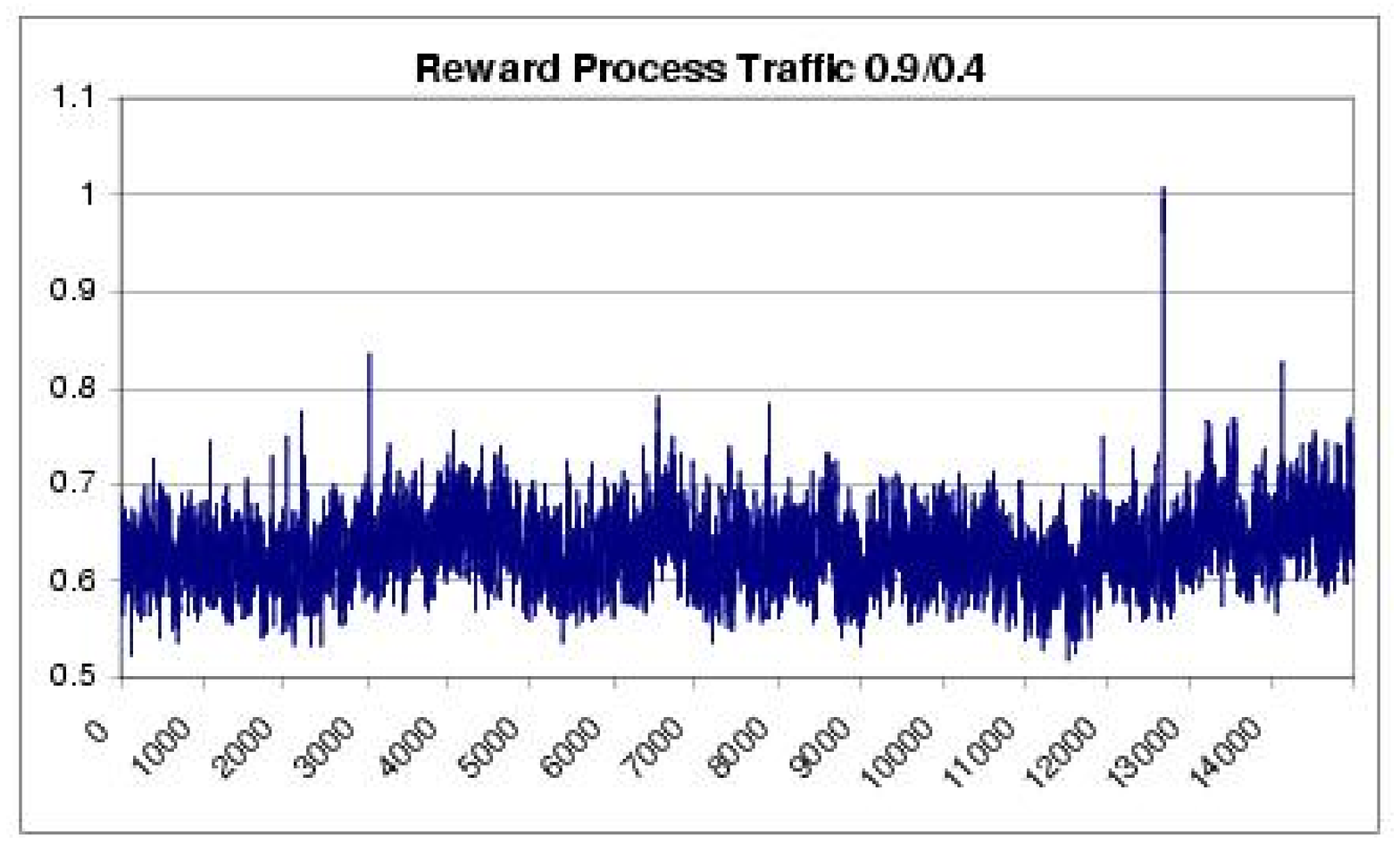}}
 \subfigure[]{\includegraphics[height=150pt, width=200pt]{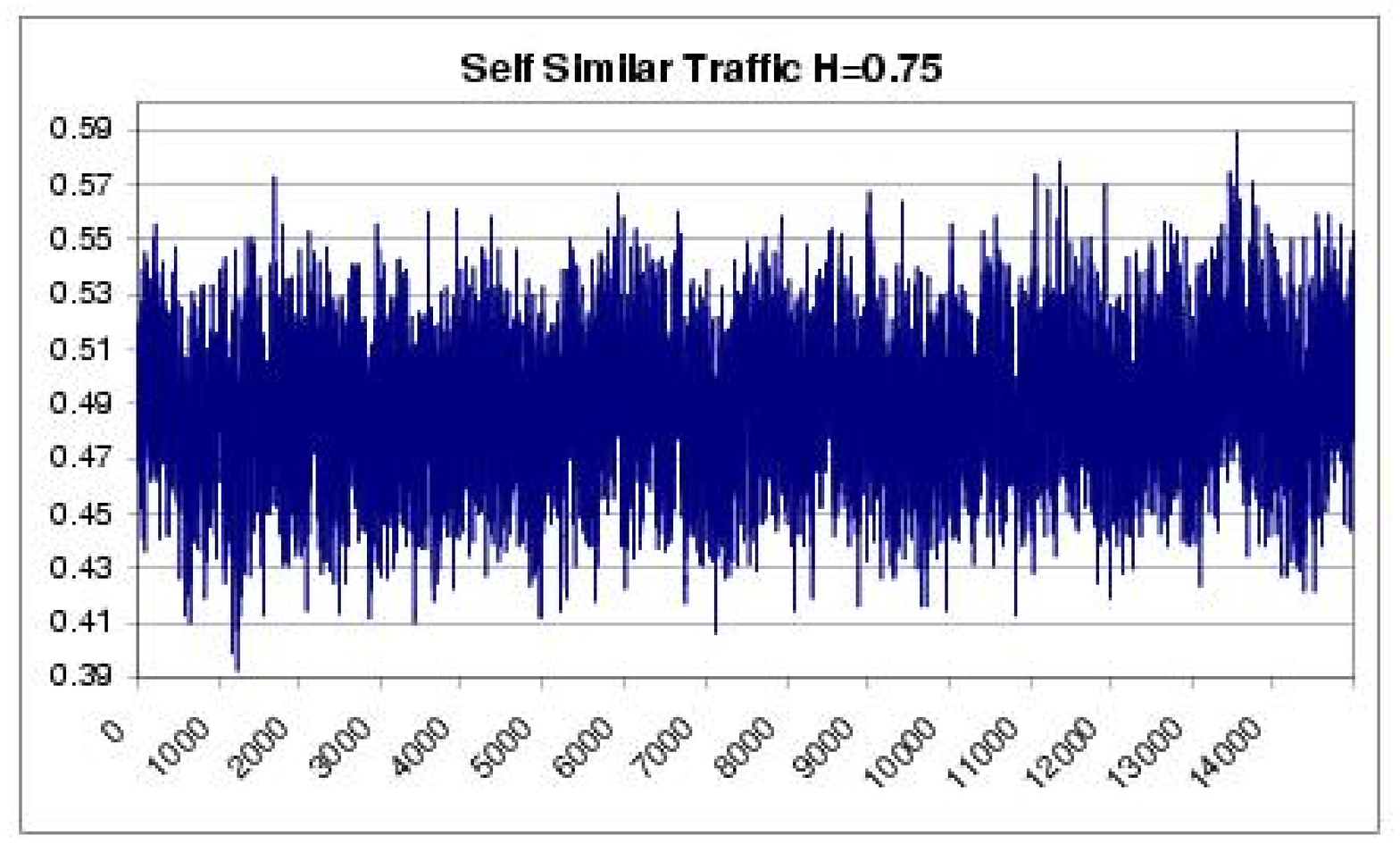}}
 \subfigure[]{\includegraphics[height=150pt, width=200pt]{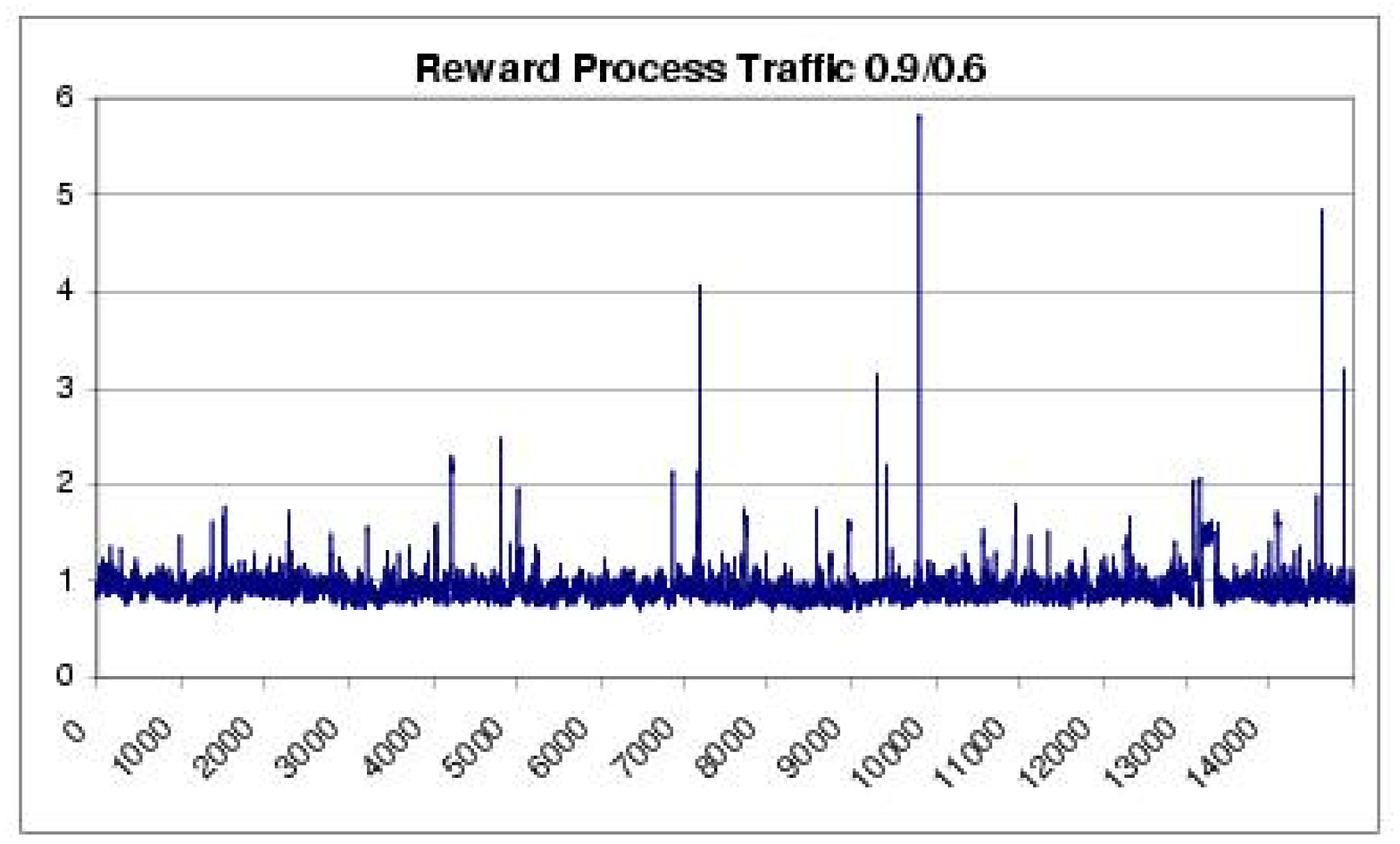}}
 \subfigure[]{\includegraphics[height=150pt, width=200pt]{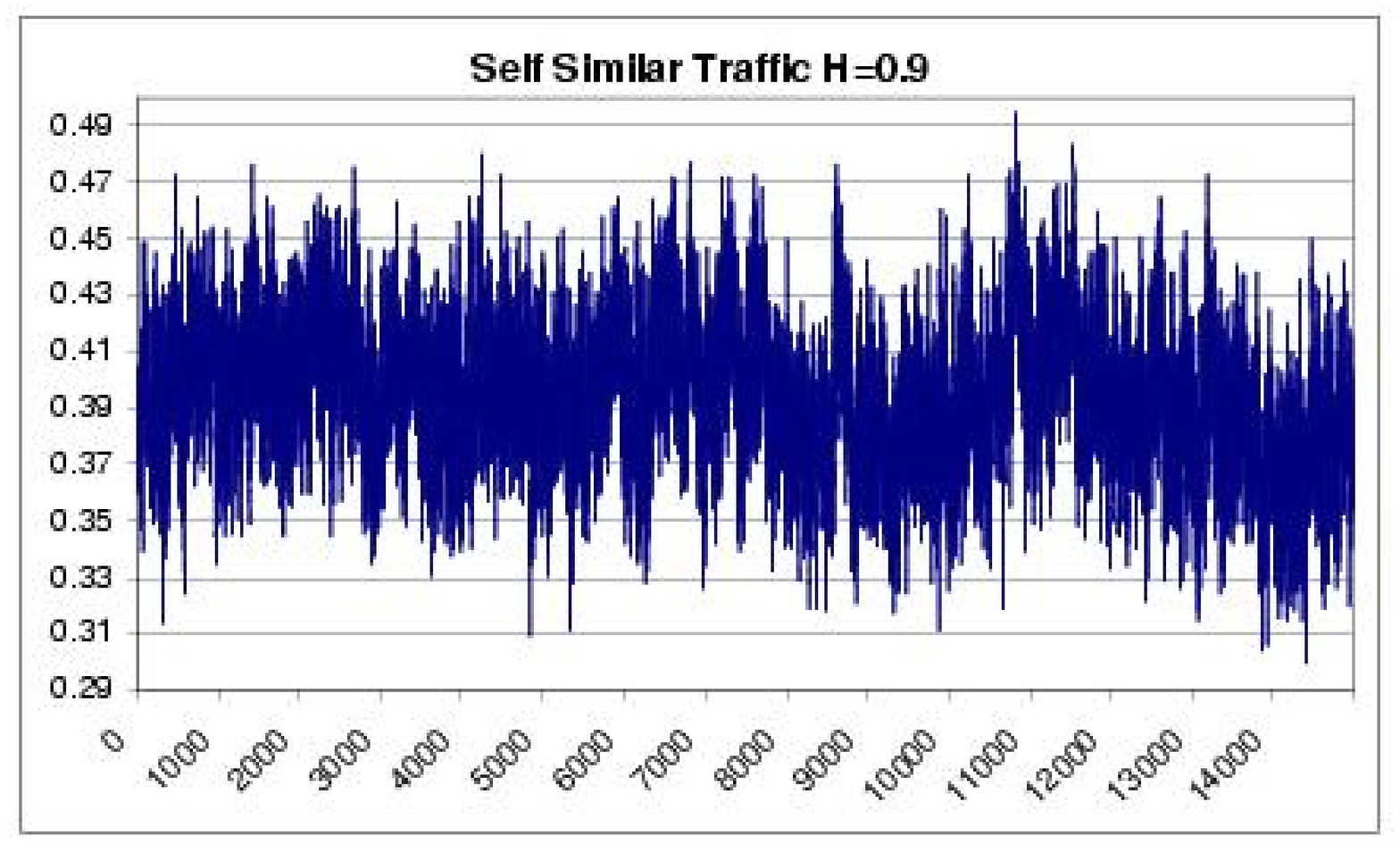}}
 \subfigure[]{\includegraphics[height=150pt, width=200pt]{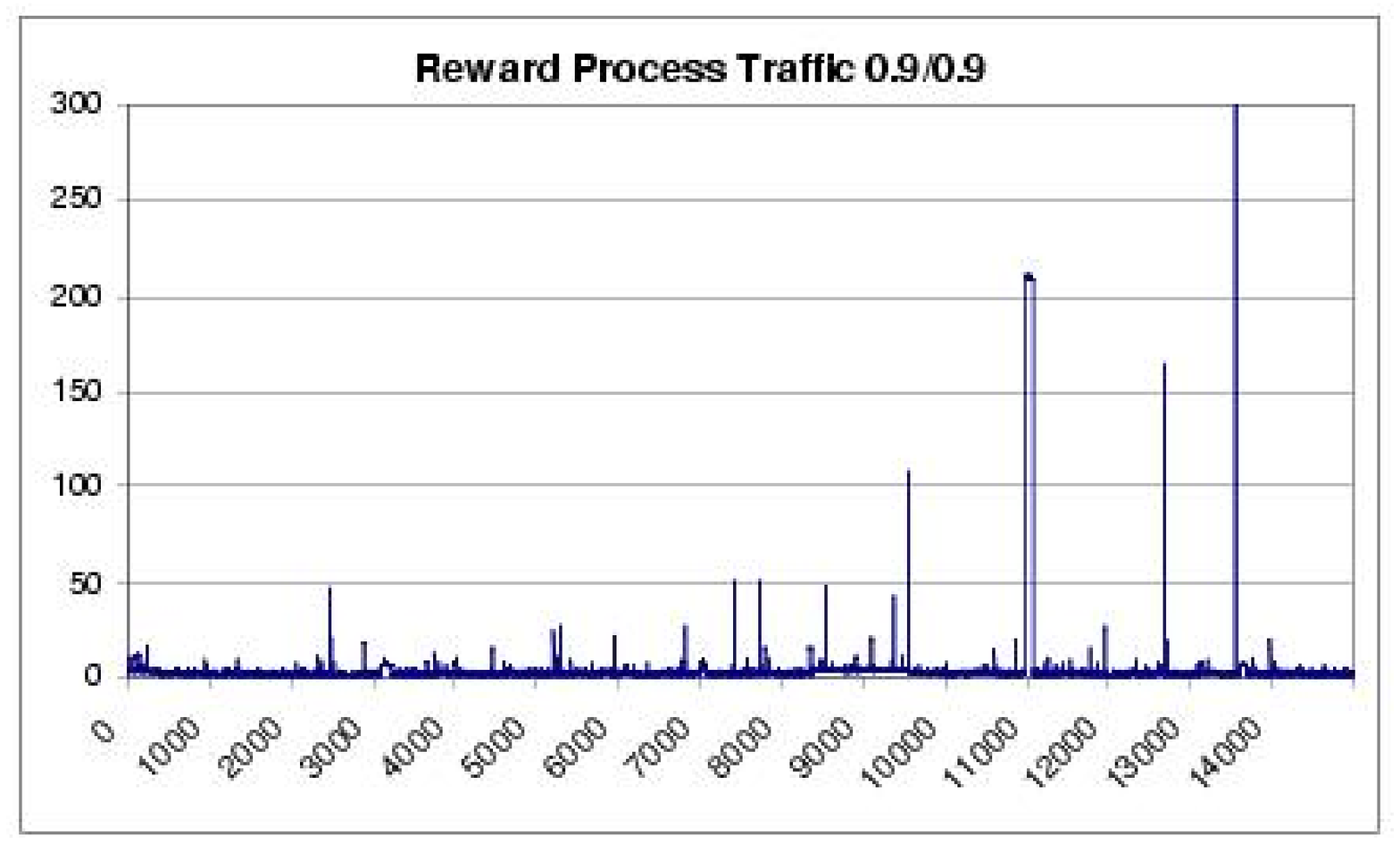}}
 \caption{ 
 Three instances of a Self-similar and a Reward
 process: the durations of the ON intervals and, in the latter case, the height of the ON intervals, 
 are distributed according to $\frac{1}{x^{H}},\ x \sim U\left(0,1\right) $.
 (a),(c),(e) features H = .6, .75, and .9, respectively, whereas in (b),(d),(f) H = .9 fixed for the 
 ON-interval lengths, and H = .4, .6, .9 for the heights, respectively. OFF-intervals are exponentially 
 distributed with mean equal to their corresponding ON-intervals.  
 }
 \label{SSRW}
\end{figure}

The self-similarity of the traffic, and the heavy-tailed distributions of
the intervals of silence and/or activity for each user, have been
verified in real traces. The violation of the traditional Poisson model \cite
{PF1,NSSW1}, and the long range dependence, i.e. slowly decaying
autocorrelation (illustrated in Fig. \ref{AC} and \ref{ACND} in section \ref{modelA}), which 
are both consequences of the heavy tail
distributions, are perhaps the most striking results. However, this model is
not entirely successful: the traffic it produces is not as ``spiky'' as real traffic 
(Fig. \ref{SSRW}). Moreover, the \emph{Energy function}, which gives, for each time scale, the contribution 
of the behavior at that scale to the total energy (more details in Section \ref{EnAvDef}, or \cite{FGW1,FGHW1}),
can clearly distinguish between the simulation and
the real traffic. For the simulation, the Energy function gives a straight line, whereas the Energy 
function of the real traffic curves (see Fig. \ref{EnOld} below, summarizing the Energy for data as well as simulations according to several models). 

The ``curving of the Energy function'' for Internet traffic traces has played an important role as a motivation for certain models of Internet data traffic, ever since its first observation in \cite{FGW1}. It will be used as one of the cornerstones of this work as well. As Section \ref{EnAvDef} will explain, though, a slightly different definition for the Energy function than the one in \ref{EnAvDef} will be used; the Energy function herein used is the average of the translation non-invariant concept in \ref{EnAvDef} over all possible choices of the ``origin'' ($t=0$ plays a special role in \ref{EnAvDef}). This is similar to averaging procedures in many other applications of wavelet expansions. The result is that the Energy function, as defined in this paper, remains reliable even for very large scales (whereas the standard Energy function produces large oscillations and becomes unreliable for scales above $15$ in Fig. \ref{EnOld}), making it possible to observe and quantify the ``curving'' or ``bumpiness'' of the Energy function even at these large scales. The issue will be further analyzed in Section \ref{modelD}.   

In this paper, simulations according to this model (our model A) will bear the label ``0-1''.

\addtocounter{figure}{-1}
\stepcounter{figure}

\bigskip

\noindent \textbf{Reward Model (RM)}: In an attempt to make the SSM ``spikier'', it
was proposed \cite{LT1,PT1} to change the values of ON-intervals from $%
1$ to some positive random value, which would be kept constant throughout
the interval, values of different intervals being independent. This
would reflect the non-equivalence of the users, in terms of the
amount of data they can send: some may have faster connections, may run
applications generating higher volumes of data etc. Needless to say, the
distribution of values for the ON-intervals can be itself heavy-tailed. A
result proved for this model states that, in the limit of infinite users,
the (properly normalized) total traffic can be either a self-similar or a 
\emph{p-Stable} process \cite{R2}.

This model resembles real traffic a bit better (compare Fig. \ref{SSRW} and \ref{Traffic9497}), however the
two types of traffic can still be discriminated visually. Moreover, the marginal
distribution of a part of the traffic generated by this model will always be
either p-Stable or Gaussian, whereas in the case of the real traffic
different parts may have different marginals (see Fig. \ref{KolDist} in section \ref{modelB}). Also, heavy-tailed
value distributions lead to very high spikes, while distributions of finite
second mean average out in the limit, leading to the same results as SSM
(Fig. \ref{SSRW} (b),(d), and (f)). Thus, the model does not allow for very fine control of the form of
the spikes. In addition, its Energy function is once more a
straight line, but an additional risk is that it can be \emph{meaningless},
in the case of a heavy-tailed value distribution (see Section \ref{EnAvDef}). Finally,
as noted earlier, real traffic \emph{does not have heavy-tailed spikes}.

In this paper, simulations according to this model (our model B) will bear the label ``RH''.

\bigskip

\noindent $\alpha$ \textbf{and} $\beta$ \textbf{Traffic model}: This model,
proposed and studied in \cite{SRB1}, adopts a different approach,
decoupling the spikes and the main body of the traffic, and dealing with each of
the two separately. It suggests that Internet traffic can be viewed as a
superposition of two independent processes: the first ($\beta $ traffic) is
self-similar Gaussian, and has no spikes, the second ($\alpha $ traffic) 
is a train of isolated spikes. The inter-spike intervals are
independent and approximately exponential, i.e. spike arrivals are
approximately Poisson.

\begin{figure}[t]
 \centering
 \subfigure[]{\includegraphics[height=150pt, width=200pt]{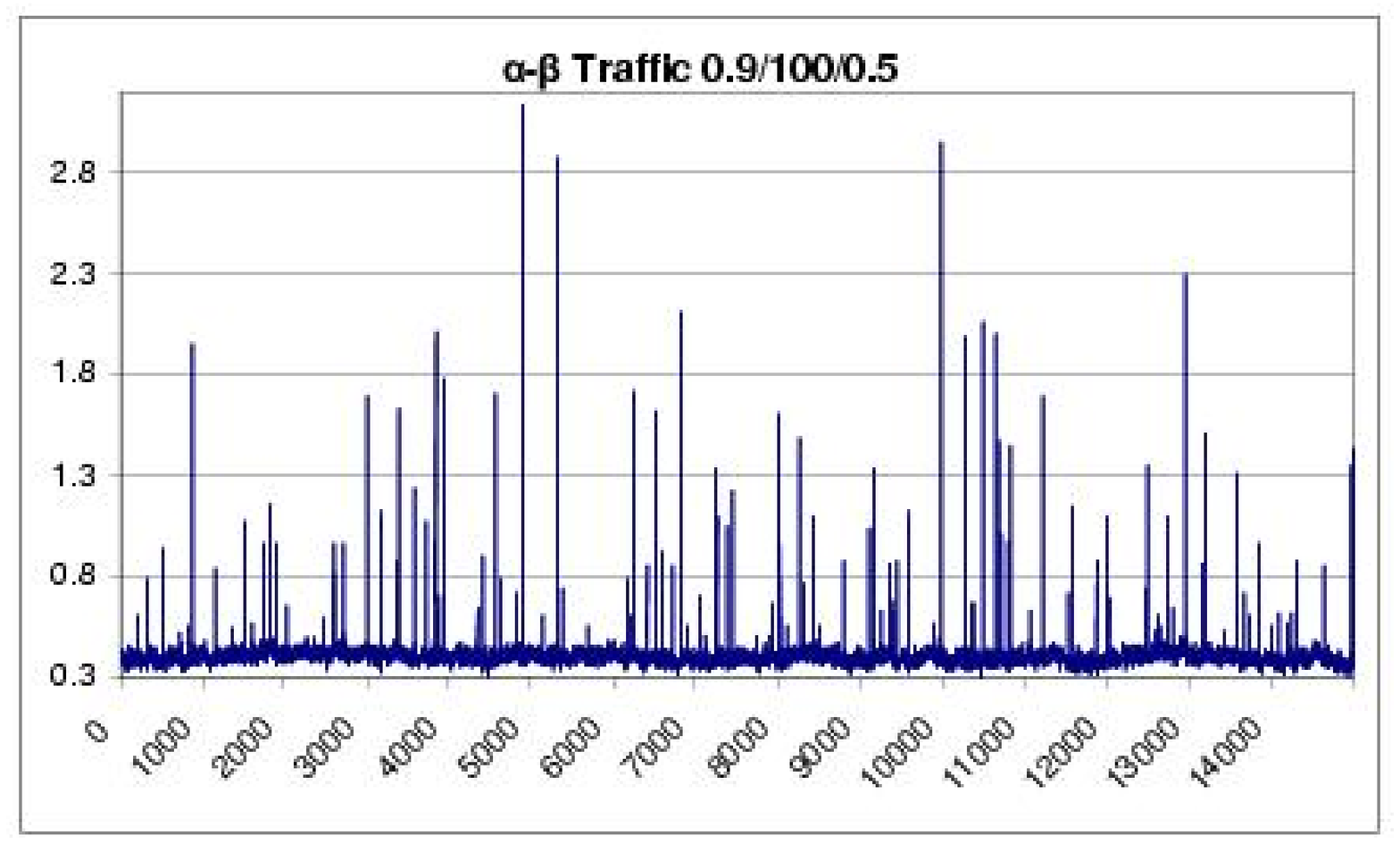}}
 \subfigure[]{\includegraphics[height=150pt, width=200pt]{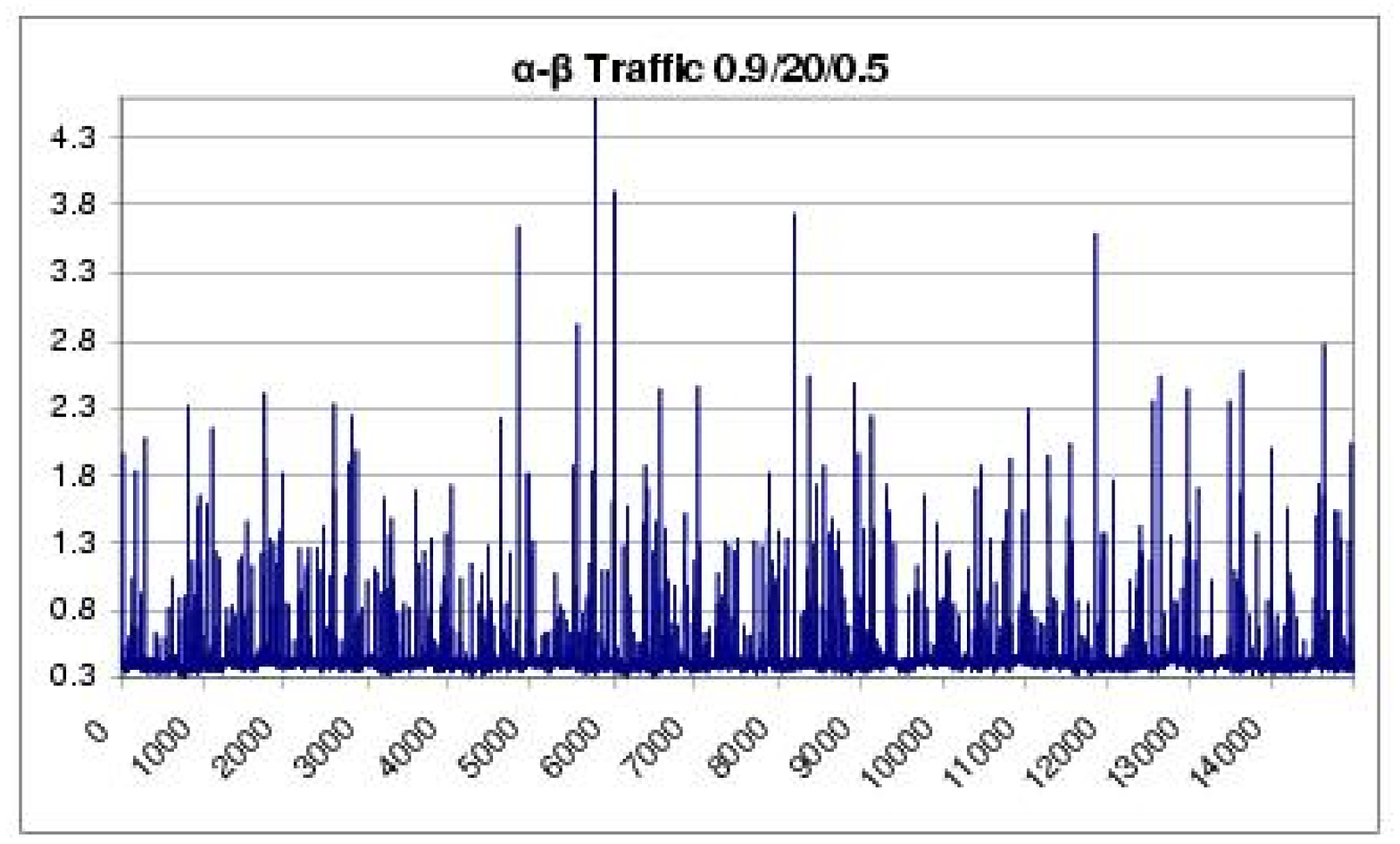}}
 \subfigure[]{\includegraphics[height=150pt, width=200pt]{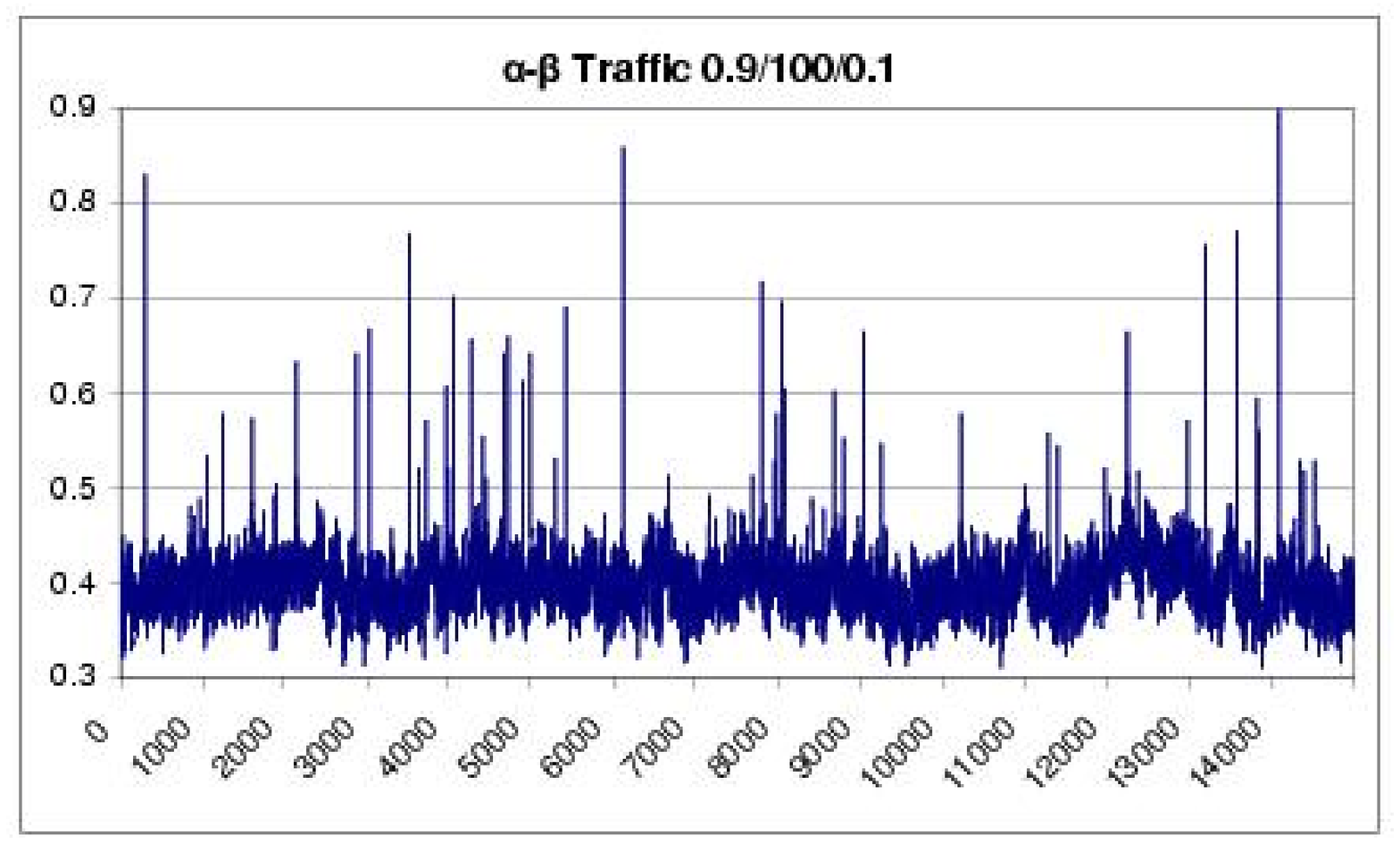}}
 \subfigure[]{\includegraphics[height=150pt, width=200pt]{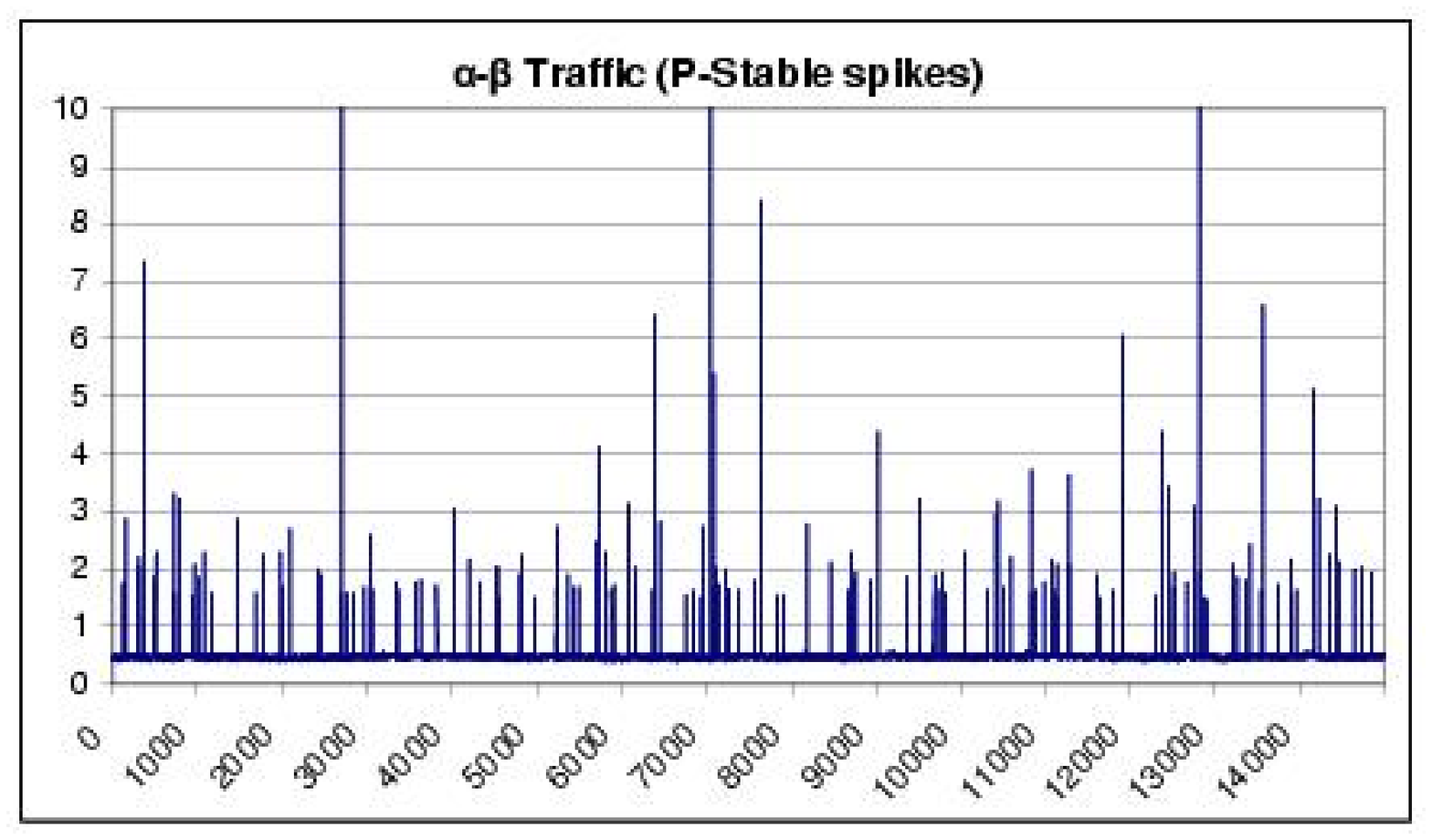}}
 \caption{
 Traffic generated by the $\alpha $ and $\beta $ model: (a), (b), and (c) feature 
 spikes with exponentially distributed heights (with mean 0.5, 0.5, and 0.1 respectively)
 and exponentially distributed interarrival times (with
 means 100, 20, and 100 respectively). (d) features p-Stable spikes having exponential interarrival
 times with mean 100. For the $\beta $ traffic, the same model was used as in Fig. \ref{SSRW}, with H = .9.}
 \label{abtraffic}
\end{figure}

This model was derived from, and calibrated on real traffic traces, and all of
the assumptions above were shown to hold. No distribution for the
spike heights was proposed in \cite{SRB1}; motivated, however, by the observed exponential marginals of
the traces considered here, we will be using an exponential one in the simulations. The result of the
simulation and the real traffic are almost indistinguishable visually
(compare Fig. \ref{abtraffic}, \ref{Traffic9497} and \ref{TrafficM6U8}). An example with a p-Stable 
distribution is also provided.

To sum up, $\beta$-traffic accounts for the long range dependence, and $%
\alpha $-traffic for the spikes. Thus, the model matches, to some extent,
the appearance, and the slowly decaying autocorrelation of the real traffic (Fig. \ref{AC} and \ref{ACND}).
Even better, the marginal distributions of parts of the traffic oscillate
between Gaussianity and p-Stability, a feature shared by real traffic \cite{SRB1}. Unfortunately, the Energy
function of the simulation is again oversimplifying reality (Fig. \ref{EnOld}).

Simulations according to this model will bear the label ``RH HT''.

\bigskip

\noindent \textbf{Infinite Source Poisson Model (ISPM)}: This is a variant of the SSM. According to this model, the
traffic is the sum of independent piecewise continuous connections, whose range is $\{0,1\}$, and whose support is 
compact and connected. The beginning points of their supports are given by a Poisson process of rate $\lambda$, and 
their sizes by a heavy tailed distribution of infinite variance (see \cite{MRRS1} and references therein for more
details). 

Obviously, the support of the total traffic consisting of a finite number of connections is also finite; thus, 
contrary to the SSM, this model evolves in time, as the number of connections $n$ increases. In order to capture this
evolution, the traffic can be normalized in time by a parameter $T$, which denotes time resolution, i.e. the
time scale at which traffic is viewed. If the two limits $n\rightarrow\infty$ and $T\rightarrow\infty$ are taken
simultaneously, the limit process can be either Gaussian or p-Stable, depending on their relative rates (see the
discussion in \cite{MRRS1}). This behavior will appear again in Section \ref{modelC}.

Unfortunately, it seems that this limit process cannot be linked back to real traffic in an obvious way.
  
\addtocounter{figure}{-1}
\stepcounter{figure}

\bigskip

It is very striking that the Energy function for these four models are very linear (as shown in Fig. \ref{EnOld}, at least up to scale $15$; at coarser scales, Fig. \ref{EnOld} is no longer accurate, but the modified definition of the Energy function given in Section \ref{EnAvDef}, which remains accurate at large scales, the linearity is seen to persist until even the coarsest scales). Indeed, as many authors (including us) have found, ``local'' models (in which the traffic is generated as an average over simulated sessions ``unrolling'' in time) typically produce such nice, linear, but therefore unrealistic Energy functions. 

Although each of the models is represented by a single simulation only in Fig. \ref{EnOld}(b), we did in fact check many values for the parameters of these models; for other parameter values, one finds a change in the slope of the corresponding Energy function but not in its essentially linear character. 

In simulations for the models listed above, the features of the model seem to control only the slope, or, at best,
introduce only one slope change from one scale onwards (Fig. \ref{EnOld}). Moreover, the Energy
function appears to be insensitive to the ``spikiness'' of the process
(compare EXP IID with HT IID in Fig. \ref{EnOld}(b)), although it seems to depend on the autocorrelation. 

The failure of these models to produce curved Energy functions stimulated researchers 
to seek more ``exotic'' models, such as cascades, and multifractal processes in
general \cite{R1,R2,FGW1,RRCB1}. These models, however, are not linked in any way to a description of 
user or connection behavior.  

\bigskip

\noindent \textbf{Multifractal models}: Multifractal models, unlike all of the preceding 
models, produce curving Energy functions \cite{FGW1}. They also produce nontrivial \emph{multifractal spectra} \cite{R2}. The simplest representatives of this category are the \emph{%
random cascades} \cite{FGW1,R2}. However, in our view, these models are ``global'' rather than ``local'': the stochastic process that generates the simulated sessions is not inspired by network mechanisms or user characteristics, or the structure of these in time. To the contrary, these processes typically treat the whole observation time interval as an entity, and link probabilities of certain behavior in a subinterval to the probability of other behavior elsewhere, without any motivation, other than that it produces Energy functions similar to those observed in real traces. Since they do not take into account network mechanisms, at least in a straightforward manner, they offer little help towards understanding these mechanisms. For this reason they will be considered no further here.

\bigskip

\noindent \textbf{Network simulation models}: These models adopt a completely
different point of view, trying to incorporate into the model every single
detail of the routers, the links, the protocols etc. Examples are
the Network Simulator and  (NS)(see \cite{FP1} and references therein), and the Scalable Simulation Framework
(SSF) (see http://www.ssfnet.org/homePage.html for more information). As expected, they produce very good
results, but the overwhelming number of parameters violates the parsimony
principle, and makes it hard to estimate the impact of each of them on the
traffic. For these reasons, this approach falls outside the scope of the present paper. 

\bigskip

\noindent \textbf{Summary}: In the list of models above, some follow an empirical approach, which attempts to capture properties of the traffic, but ignores the underlying mechanisms ($\alpha-\beta$ Traffic, Random Cascades); others present mathematical constructions that are motivated by network and user behavior, but produce traffic that only partially corresponds to reality (SSM, RH). As far as we know, there exists no model that produces traces a) that have the spikiness and long-range dependence of real traffic traces, b) whose marginals match the behavior of the observed marginals, which are sometimes close to and at other times deviate far from Gaussianity, and c) with a curving Energy function. Obtaining a curving Energy function via a ``local'' model seems particularly hard: extensive experimentation by the authors with simulations of the models described above and in Section \ref{DatSim} always produced linear, or at best essentially piecewise linear graphs (as in Fig. \ref{EnOld}(b) and \ref{EnAvCur}(a)), except for our model D described in Section \ref{modelD}. 

Random Cascade models do produce curving Energy functions; they are however not satisfactory, not only because they offer no explanation of the underlying mechanism, relating to network protocols or user behavior, but also because the resulting traffic still contradicts real traces in various other aspects (auto-correlation, marginal 
distributions etc.).

\begin{figure}
 \centering
 \subfigure[]{\includegraphics[height=150pt, width=200pt]{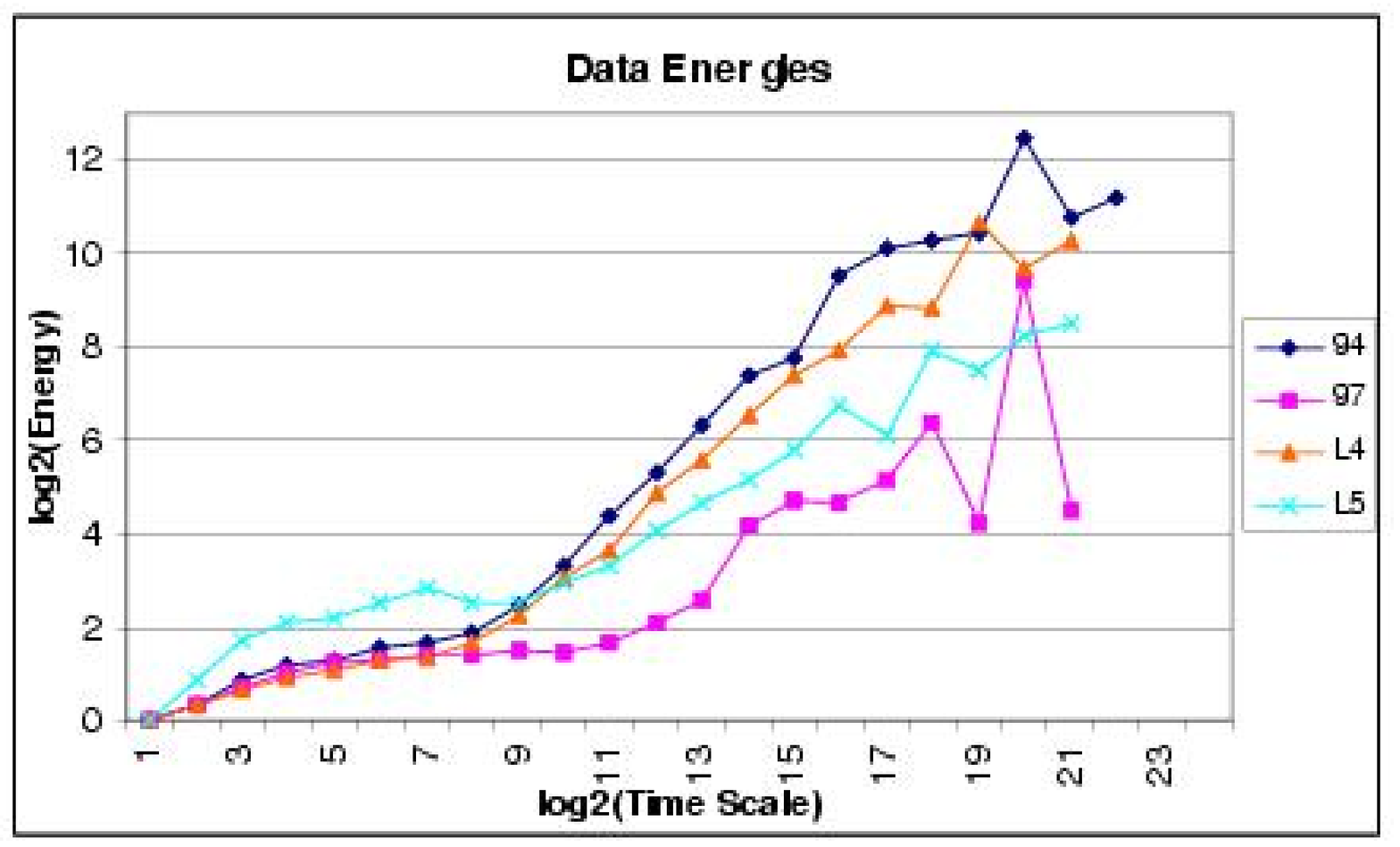}}
 \subfigure[]{\includegraphics[height=150pt, width=200pt]{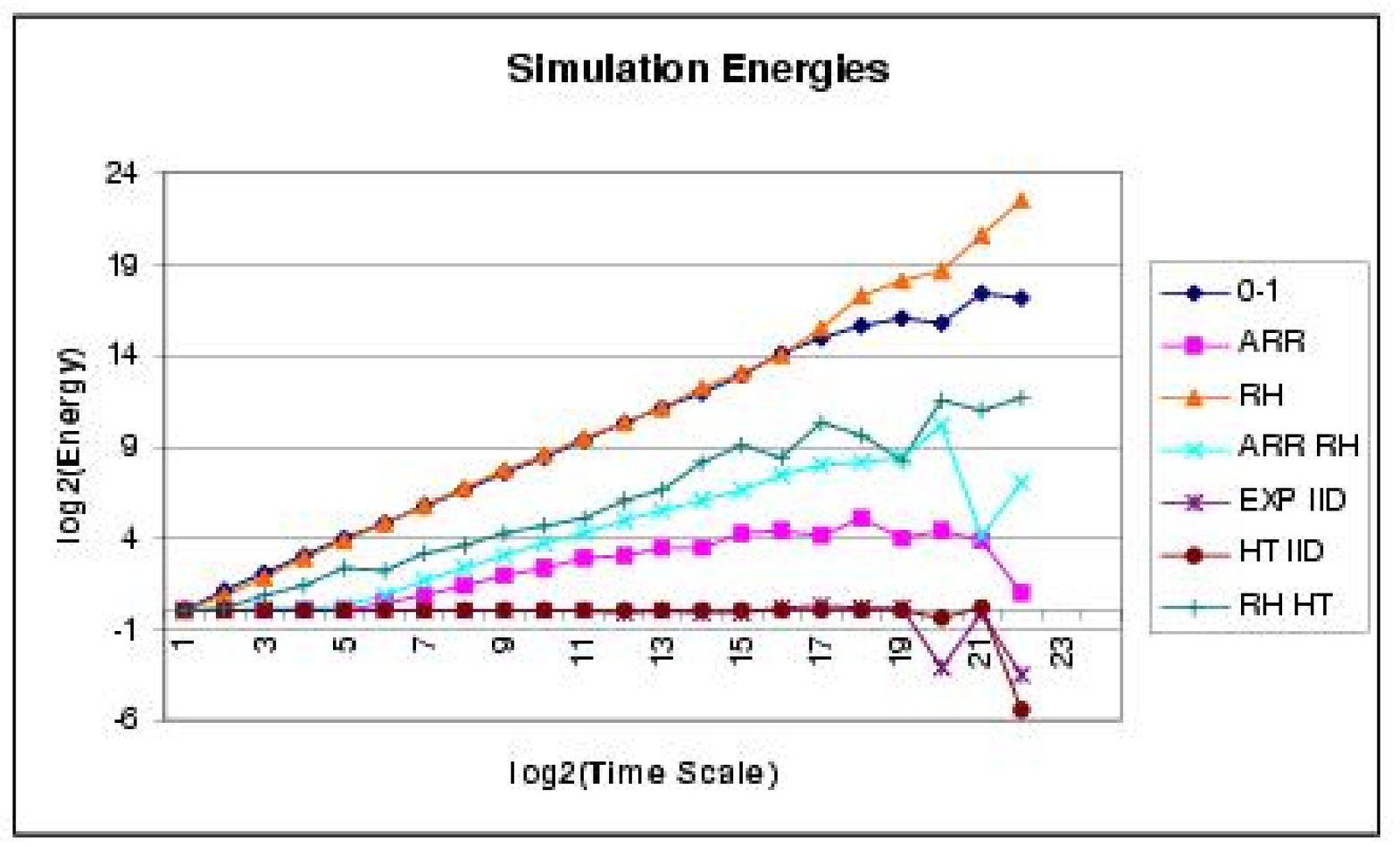}}
 \caption{Energy functions, according to the standard definition (Definition 1 in Section \ref{EnAvDef}, for (a) data 
          and (b) simulations. They are normalized 
          up to an additive constant, so that they all start at 0. Notice that the Energy functions in (b),
          corresponding to the simulations described in Section \ref{DatSim}, are approximately linear. Notice also 
          the large oscillations at coarse time scales, caused by a lack of accuracy at coarse scales; Definition 2 
          will take care of this. Finally, note that HT IID and 0-1 produce indistinguishable results, although they
          have completely different properties.}
 \label{EnOld}
\end{figure}
\addtocounter{figure}{-1}
\stepcounter{figure}

\subsection{Summary of our results}

The purpose of this paper is to present a model that describes Internet traffic on a network link; the main guideline in the construction of this model will be the incorporation, into the SSM, of additional network mechanisms and aspects of user behavior that will be identified as responsible for shaping the traffic. It will turn out that enriching the SSM with only a few extra features leads to simulations exhibiting spikiness (i.e. marginals), long-range dependence, curving of the Energy function, and wide oscillations in the deviation of its marginals from Gaussianity that are qualitatively and quantitatively indistinguishable from the same features of real traces.  
The impact of these features on the traffic will be rigorously determined by means of
theorems, and the predictions of these theorems will be subsequently validated against real traffic. The success of 
this validation will suggest that the factors that have been identified as important for shaping the traffic are indeed important, and may explain the behavior of the real traffic singled out by the four criteria above.   

In order to increase clarity, the construction will proceed incrementally in 4 steps, labeled A, B, C, and D.
At each step, a theorem will state and prove the properties of the model built so far.  
Step A will be the plain old SSM; step B will take into account that information transmission does not take place continuously, but rather in packets spaced apart by variable time intervals, which is turn are determined by the Round Trip Times; step C will then introduce the Slow Start feature of the TCP protocol into the model, arguing that the small volume of information exchanged by the majority of the connections is such as 
to allow Slow Start to impact the traffic heavily; and, 
finally, step D will enrich the model with other, longer time scales, characteristic of network activity, originating from either user behavior (e.g. switch of application), or protocol actions (e.g. congestion control).

\bigskip

\noindent \textbf{Step A}: The basic model, which will serve as the founding ground
for the rest, will be the SSM, introduced and studied in \cite{TWS1}. Theorem A states and proves its
properties.

\bigskip

\noindent \textbf{Step B}: Here, a finer structure will be imposed upon the ON-intervals of
the SSM: each value of 1 will be followed by a RTT, and
RTTs will be modeled by means of random variables. In general, little is known about their
distribution \cite{VB1, FP1}, but in large capacity links of the Internet backbone, packet
interarrival times tend to be exponentially distributed and independent~\cite{CCLS1}. Fortunately, this is not a 
serious issue, since Theorem B below asserts that the simulated traffic is insensitive to this distribution,
as long as it has finite variance. Moreover, Theorem B states that the traffic is actually Gaussian, and fairly uncorrelated at small time scales (Brownian motion), but that the
correlation grows stronger in larger time scales (causing the process to become
self-similar); this prediction is borne out by real data (Fig. \ref{AC} and \ref{ACND}). 

\bigskip

\noindent \textbf{Step C:} This step incorporates features of the 
TCP/IP combination of protocols, responsible for the bulk
of data transfer over the Internet today, into the model. One of the characteristics of TCP,
the \emph{Slow Start}, roughly stipulates that the number of packets sent
at each emission time increases \emph{geometrically} by a factor of 2 (1,2,4,8,...),
until either loss occurs, a preset maximum is reached, or all data is sent
(for more details, see e.g. \cite{PD1}). In detailed traces, which show individual connections, 
the Slow Start is easily discernible; this suggested its incorporation in the model, and the exploration 
of the consequences of its inclusion. 
Theorem C states that the traffic will tend to be either self-similar or
p-stable, contingent upon the Slow Start maximum and the number of users. In
practice, p-Stability cannot be distinguished from other very large deviations from Gaussianity, since observations can only last for a finite amount of time. Nevertheless, the fluctuation of the number of active users in time will allow
for both self-similar regimes and regimes that deviate far from this self-similar behavior within the same trace, as is also observed in real traces. An increase of the number of users steers the process towards self-similarity.

\bigskip

\noindent \textbf{Step D}: User connections over the Internet seem to be governed by a
hierarchy of time scales (\emph{levels}); for example, a particular session of data 94 
seems to live in six different time scales (Fig. \ref{Lev}), at approximately 2000, 500, 100, 40, 5, and 
0.1s, and a session of trace M6 in three different time 
scales (Fig. \ref{LevND}), at approximately 2, 0.4, and 0.01s. To the best of our knowledge, there has been no
effort yet towards the exact quantification and explanation of these levels at different scales, but it is
commonly accepted that they exist \cite{PF1,RRCB1,FGWK1}. A possible explanation could lie in the recursive structure within Internet applications (e.g. Web browsers opening multiple sessions, FTP applications establishing multiple
connections, etc., TCP, which sends and waits for receipt acknowledgement,
etc.), and the human work habits, according to which tasks are broken
into subtasks recursively. Analysis of 50 typical (long enough) sessions of
trace 94 reveals that levels are present in all of them, and are more or less similar. Moreover, their statistical properties (mean duration, variance, and range) are similar for other traces of the same time period. 

So far, model A had just 1 level, and models B and C each had 2 levels. This step introduces more
levels, at larger time scales. Theorem D then proves that such levels cause
the Energy function to curve.

\bigskip

The paper is organized as follows: Section \ref{EnAvDef} describes the statistical tools
used (including the modified definition of the Energy function), and 
sections \ref{modelA}, \ref{modelB}, \ref{modelC}, and \ref{modelD} 
deal with models A,B,C and D respectively. Section 7 concludes the paper with a discussion 
and some applications; in particular, several simple tools are introduced that detect the presence and number 
of levels, and that quantify the behavior of the traffic.  

\section{Statistical Tools}

\label{EnAvDef}

Let $\left\{ X_{k}\right\} _{k=1}^{n}$ represent a time series. In what
follows, extensive use of three standard statistical tools will be made:
\begin{enumerate}
	\item Empirical marginal distributions: 
	\[ F_{n}(t)=n^{-1}\sum_{i=1}^{n}1_{X_{i}\leq t} \]
	
	\item Autocorrelations: 	  
	    \[Corr(k)=\frac{\sum_{i=1}^{n-k}(X_{i}-\overline{X}_{n})(X_{i+k}-\overline{X}_{n})}
	    {\sum_{i=1}^{n}(X_{i}-\overline{X}_{n})^{2}} \] where $\overline{X}_{n}=\frac{1}{n}\sum_{i=1}^{n}X_{i}$

  \item The \emph{Kolmogorov distance} between two 
   Cumulative Density Functions (CDFs), which is defined as the maximum of the absolute value of 
   their difference:
   \begin{equation}
   d_K(F_1,F_2):=\underset{x \in \mathbb{R}}{\sup}|F_1(x)-F_2(x)|
   \label{KolD}
   \end{equation}
\end{enumerate}
In addition, the
Energy function will be used, as well as a more intuitively appealing (in our opinion) variant,
the \emph{p-Averaging function} (similar to the Normalized Variance appearing in 
\cite{PF1}), introduced below.

The definition of the $p$-Averaging function will be completed in two steps:\
initially, a definition that ensures compatibility with the well established
Energy function \cite{FGW1,FGHW1} (in a sense to be made clear below)
will be proposed. Then, a simple observation will lead to the redefinition
of both Averaging and Energy functions.

\bigskip

\begin{odefinition}
\upshape
Given a stationary sequence $\{X_{i}\}_{i=1}^{M}$, for some $M=2^{m}$, define 
\begin{equation*}
\begin{array}{llcll}
X_{k}^{j} &=&2^{-j}\displaystyle\sum_{i=(k-1)2^{j}+1}^{k2^{j}}X_{i} & j=0,...,m & k=1,...,2^{m-j} \\
m_{k}^{j} &=&\left|X_{2k}^{j}-X_{2k-1}^{j}\right| & j=0,...,m-1 & k=1,...,2^{m-j-1} \\
\widetilde{X}_{k}^{j} &=&2^{-j/2}\displaystyle\sum_{i=(k-1)2^{j}+1}^{k2^{j}}X_{i} & j=0,...,m & k=1,...,2^{m-j} \\
d_{k}^{j} &=&\dfrac{1}{\sqrt{2}}\left|\widetilde{X}_{2k}^{j}-\widetilde{X}_{2k-1}^{j}\right| & j=0,...,m-1 & k=1,...,2^{m-j-1}
\end{array}
\end{equation*}
For $M_{j}=M/2^{j+1}$ and $p>0$, define the \textbf{p-Averaging} \[A_{j+1}^{[p]}
=\left( \frac{1}{M_{j}}\sum_{k=1}^{M_{j}}(m_{j,k})^{p}\right) ^{1/p}\] and 
the \textbf{Energy}
 \[ E_{j+1}=\frac{1}{M_{j}}\sum_{k=1}^{M_{j}}(d_{j,k})^{2},\
j=0,...,m-1\] Notice the connection of this definition of the Energy function to the Haar 
wavelet coefficients as in \cite{FGW1,FGHW1}. If $p=2$, the relation between the two functions is elementary: 

\begin{equation*}
\log_{2}(E_{j})=j-2+2\log _{2}A_{j}^{[2]},\text{ } j=1,...,m
\end{equation*}
\end{odefinition}

In the sequel, $A_{j}^{[2]}$ will be used almost exclusively, because it is easier to 
analyze, and simply connected to the commonly adopted $E_{j}$; it will be indicated
just by $A_{j}$ and referred to simply as \emph{Averaging (function)}. When $p\neq 2$ is considered, 
it will be stated explicitly. 

This set of definitions has been motivated by the following heuristic:

\bigskip

\noindent \textbf{Heuristic}: Let $X_{i}$ represent the value of the traffic for the $%
i $th time bin, taken to be the \emph{initial} time bin, and assume the sequence 
is ergodic, (i.e. $\displaystyle \frac{1}{n}\sum_{i=1}^{n}X_{i}\rightarrow EX_{1}$). What is the 
rate of this convergence?  
The way to determine the rate is to find a function $\psi(n)$ such that 
$\displaystyle \psi(n)\sqrt{\bold{E}\left(\left|\frac{1}{n}\sum_{i=1}^{n}X_{i}-EX_{1}\right|^2\right)}\underset{n\rightarrow \infty}{\longrightarrow} 1$. Typically, the function $\psi$ is increasing in $n$, at least for $n$ sufficiently large, expressing the fact that the average $\displaystyle \frac{1}{n}\sum_{i=1}^{n}X_i$ gets closer to $\bold{E}(X_1)$ when $n$ increases. In practice, because $\bold{E}(X_1)$ is unknown, one can still get a handle on $\psi$ by comparing two different averages: 

\begin{multline}
\label{heur}
\bold{E}\left(\psi^2(n) \left|\frac{1}{n}\sum_{i=1}^{n}X_i-\frac{1}{n}\sum_{i=1}^{n}X_{i+n}\right|^2\right)
=\psi^2(n)\bold{E} \left(\left|\left(\frac{1}{n}\sum_{i=1}^{n}X_i-\bold{E}(X_1)\right)-\left(\frac{1}{n}\sum_{i=1}^{n}X_{i+n}-\bold{E}(X_1)\right)\right|^2\right)= \\
=\psi^2(n)\bold{E}\left( \left\{\left|\frac{1}{n}\sum_{i=1}^{n}Y_i\right|^2+\left|\frac{1}{n}\sum_{i=1}^{n}Y_{i+n}\right|^2-2\left(\frac{1}{n}\sum_{i=1}^{n}Y_i\right)\left(\frac{1}{n}\sum_{i=1}^{n}Y_{i+n}\right)\right\}\right)\overset{*}{=} \\
=\frac{\psi^2(n)}{n^2} \left\{\sum_{i=1}^{n}\sum_{j=1}^{n}\bold{E}(Y_iY_j)+\sum_{i=1}^{n}\sum_{j=1}^{n}\bold{E}(Y_{i+n}Y_{j+n})-2\bold{E}\left(\left(\sum_{i=1}^{n}Y_i\right)\left(\sum_{i=1}^{n}Y_{i+n}\right)\right)\right\}\overset{**}{=} \\
=2\frac{\psi^2(n)}{n^2} \left\{\sum_{i=1}^{n}\sum_{j=1}^{n}\bold{E}(Y_iY_j)-\sum_{i=1}^{n}\sum_{j=1}^{n}\bold{E}(Y_iY_{j+n})\right\}=
2\frac{\psi^2(n)}{n^2} \left\{\sum_{i=1}^{n}\sum_{j=1}^{n}[\bold{E}(Y_iY_j)-\bold{E}(Y_iY_{j+n})]\right\}
\end{multline}
where $Y_i=X_i-\bold{E}(X_1)$ was set. In ($**$) we used stationarity. Notice also that ($*$) proves that $0\leq\lim$(\ref{heur})$\leq2$.  

We claim that, if the $X_i$ are \emph{asymptotically independent}, i.e. $\forall i,\bold{E}(X_iX_{i+j})-\bold{E}(X_1)^2)=\bold{E}(Y_iY_{i+j})\rightarrow 0$, as $j\rightarrow \infty$ (irrespectively of the rate), then (\ref{heur}) has a \emph{strictly} positive limit, and, therefore, that it can be used to estimate $\psi$ (within a multiplicative constant). Because asymptotical independence is satisfied by most processes arising in applications (it is extremely unlikely that two random variables arbitrarily long apart will exhibit strong dependence), this method can be applied to a wide range of processes, and is, for this reason, very powerful.    

Indeed, as $n\rightarrow \infty$, for fixed $i,j$, $\bold{E}(Y_iY_{n+j})\rightarrow 0$, and therefore it becomes infinitely smaller than $\bold{E}(Y_iY_j)$. The RHS of (\ref{heur}) shows then that the limit has to be strictly positive. 

Thus, if a random variable $\epsilon_k^j$ is defined by:
\[\epsilon^j_k=\psi(2^j)\frac{1}{2^j}\left|\sum_{k2^j+1}^{(k+1)2^j}(X_i-\bold{E}(X_1))\right|=\psi(2^j)|m_{j,k}|\]
then 
\[A_{j+1}=\left(\frac{1}{M_j}\sum_{k=1}^{M_j}|m_{j,k}|^2\right)^{1/2}=\frac{1}{\psi(2^j)}\left(\frac{1}{M_j}\sum_{k=1}^{M_j}|\epsilon_k^j|^2\right)^{1/2}\approx \frac{C}{\psi(2^j)}\]

If $M_j$ is large, $A_{j+1}$ gives us a way to estimate $\psi(2^j)$. In particular, if $\psi (n)=n^{\alpha}$, then:

\begin{equation*}
\log _{2}A_{j+1}^{[p]}=K-\alpha j\text{ and }\log _{2}E_{j+1}=2K-1+j(1-2\alpha)
\end{equation*}

Therefore, the binary logarithms of the Averaging resp. Energy functions of the processes
that ``average'' with a rate of $n^{-a}$ will be straight lines with a slope
of $-\alpha $ resp. $1-2\alpha$. Throughout
the paper, the ``curvature of the Energy or Averaging function'' will refer to 
the curvature of their binary logarithms, plotted as functions of $j$. 

\bigskip

\begin{remark}
\upshape
In the case of processes with a $p$-Stable component (heavy
tailed marginals), the Energy function, as well as the $A_{j}^{[p]},\  p\geq 2,$ are
not well defined: their expected value is infinite, because the $\varepsilon _{k}^{j}$ do not have a second moment.
This is the reason why the Energy function fails to distinguish spiky processes (with heavy
tailed marginals) from ordinary white noise processes, as shown in Fig. \ref{EnOld}(b), where the EXP IID simulation, which is very close to white noise, and the heavily spiky HT IID simulation are seen to have the same Energy function.  The few experiments related to distributions with heavy-tailed marginals will be performed with $p=1$. 
\end{remark}

\bigskip

\begin{remark}
\upshape
$A_{j}$ has another interpretation that will
prove essential later. Namely: 
\begin{equation*}
A_{j+1}=\sqrt{\frac{1}{M_{j}}\sum_{k=1}^{M_{j}}\left(m_{k}^{j}\right)^{2}}\approx \sqrt{
\textbf{E}\left(m_{1}^{j}\right)^{2}}
\end{equation*}
\begin{equation*}
=\frac{1}{2^{j}}\sqrt{\textbf{E}\left(X_{2}^{j}-\textbf{E}X_{2}^{j}+
\textbf{E}X_{1}^{j}-X_{1}^{j}\right)^{2}}=\frac{\sqrt{2%
}}{2^{j}}\sqrt{Var\left(X_{1}^{j}\right)-Cov\left(X_{1}^{j},X_{2}^{j}\right)}
\end{equation*}
Applied to the process $\displaystyle X_{i}:=\frac{1}{\Delta }\int_{i\Delta }^{(i+1)\Delta }%
\frac{1}{\sqrt{n}}\sum_{l=1}^{n}W_{l}(t)dt,$ the above computation yields 
\begin{equation*}
A_{j+1}\approx K\frac{\sqrt{Var\left(\int_{0}^{2^{j}\Delta
}W_{1}(t)dt\right)-Cov\left(\int_{0}^{2^{j}\Delta }W_{1}(t)dt,\int_{2^{j}\Delta
}^{2^{j+1}\Delta }W_{1}(t)dt\right)}}{2^{j}\Delta }.
\end{equation*}
for some positive constant $K$, which also absorbs $n$, as it is now fixed. 
This formula connects the Energy/p-Averaging function of the
total traffic to the behavior of one user, who will be ``typical'' in the mathematical sense.
\end{remark}

\begin{remark}
\upshape
Although the name ``$p$-Averaging'' may appear to be inappropriately chosen at first sight, it actually describes well the fundamental principle of the function, which may be buried under its complicated definition.

The $p$-Averaging function endeavors to measure the rate of convergence of averages of $n$ consecutive points of the process to its true mean, as $n$ increases, where ``rate of convergence'' denotes the deviation between the $n$ point mean value estimator and the true mean. Unfortunately, the true mean is unknown. One way to circumvent this problem is to measure the rate by which two $n$ point estimators converge to each other; the two aforementioned convergence rates should be essentially the same, their relative difference being at most $2$ at all times (see the related discussion in the Heuristic above). 

This basic idea will then subsequently be elaborated in Definition 2 below, which ensures that the rate gets measured accurately, by comparing a large number $N$ of pairs of $n$ point estimators, $n$ increasing geometrically, and, in addition, that it gets measured equally accurately for all $n$, by making $N$ independent of $n$. These details aside, though, it is clear that the $p$-Averaging function ``averages'' the traces in order to estimate their mean, or rather how difficult it is to estimate their mean, hence the name ``Averaging function''.  
\end{remark}

It has been implicitly assumed, of course, that $X_{i}$ is stationary. Although there are many arguments against
Internet traffic being stationary, it can be assumed to be so over reasonably short time intervals. This 
could mean a couple of hours, according to \cite{CCLS1}), so the older data sets (before 1997), with a duration of 
at most 2 hours, can be considered to be stationary. The more recent ones only last for less than 2 minutes, but they
are also the products of a much faster network. Are they also stationary? We give two reasons why we believe they are: 
a) Non-stationarities, such as the diurnal day cycle, that are an artifact of human activity, do not become 
faster, as the network gets faster. b) In the older traces, consecutive packets in the link seem to be 0.1ms apart, whereas in the newer traces they seem to be 1$\mu$s apart, in order of magnitude,  i.e. the network has 
become faster by approximately 100 times. So, what needed an hour before needs now approximately a minute, and 
thus the duration of the old and new traces, normalized by the network speed, is the same.

Under the definitions above, the computation at the
scale $2^{j}$, for either the Energy or the Averaging functions, involves a mean value over $%
2^{m-j-1}$ numbers. Consequently, for $j$ close to $m$, the sum will comprise insufficiently many terms to produce 
a reliable mean, assuming that such a mean is indeed the
ultimate goal of the computation, and undesirable oscillations may be observed (see Fig \ref{EnOld}).
Instinctively, it is expected that averaging over more numbers will lead to smoother curves. 
Thus, a modification of the definition of 
the Energy and p-Averaging functions is in order; the definition that follows will allow for averaging 
over the maximum number of terms possible,
i.e. $2^{m-1}$, by using a standard technique of statistics, known as \emph{overlapping blocks}. 

From the point of view of wavelets, the definition change amounts to averaging the Energy function computations over all possible choices of the ``origin'' for the Haar basis underlying the definition of the $d_{j,k}$.; this type of averaging is now commonly done in situations where the use of a wavelet basis introduces a translation non-invariance that was not present earlier.

\begin{odefinition}
\upshape
The operator $\displaystyle \frac{1}{M_{j}}\sum_{\cdot=1}^{M_{j}}$ in the above
definitions will be formally replaced by $\textbf{E}$, the expectation operator.
Moreover, for reasons of clarity, the domain of $A$ will be shifted left by
1, i.e. it will start with 0. The same name and notation will be kept
for the new functions emerging out of these definitions, which, explicitly,
take the form: 
\[A_{j}^{[p]}=\left( \mathbf{E}(m_{j,1})^{p}\right) ^{\frac{1}{p}}\] and 
\[E_{j}=\mathbf{E}(d_{j,1})^{2}\] 
Remember here the stationarity assumption, which
implies $m_{j,k}=m_{j,1}$, and $d_{j,k}=d_{j,1}$, the equality holding in distribution.
\end{odefinition}

How different is the new definition from the old? Stationarity can now
be fully exploited, in the sense that there is no ``beginning of time'' for
the sequences: namely, it can be assumed that they are ordered circularly
(on a ring), not linearly (on an interval). Thus, for example, computing $%
A\left( 0\right) $, in addition to considering the differences used before,
i.e. $X_{0,1}-X_{0,2},X_{0,3}-X_{0,4},...,X_{0,2^{m}-1}-X_{0,2^{m}}$, one
will also take into account $X_{0,2}-X_{0,3},X_{0,4}-X_{0,5},...,X_{0,2^{m}-2}-X_{0,2^{m}-1},X_{0,2^{m}}-X_{0,1}$.
The circular folding of the indices is justified by the stationarity assumption. In practice, this means that the computations will be done according to the formula:

\bigskip 

\begin{eqnarray*}
A_{j}^{[p]}&\approx&\left( \frac{1}{2^{j+1}}\sum_{l=0}^{2^{j+1}-1}\left( 
\frac{1}{2^{m-j-1}}\sum_{k=1}^{2^{m-j-1}}\left| F_{kl}^{1}(j)-F_{kl}^{2}(j)\right|^{p}\right) \right) ^{1/p} \\
F_{kl}^{1}(j)&=&\frac{1}{2^{j}} \sum_{i=1}^{2^{j}}X_{1+\left( \left( k-1\right)2^{j+1}+l+i-1\right) \bmod 2^{m}} \\
F_{kl}^{2}(j)&=&\frac{1}{2^{j}} \sum_{i=1}^{2^{j}}X_{1+\left( \left( k-1\right) 2^{j+1}+2^{j}+l+i-1\right) \bmod 2^{m}}
\end{eqnarray*}

which can be simplified into:

\begin{eqnarray*}
A_{j}^{[p]}&\approx&\left( \frac{1}{2^{m+jp}}\sum_{l=0}^{2^{j+1}-1}%
\sum_{k=1}^{2^{m-j-1}}\left| \sum_{i=1}^{2^{j}}\left( G_{ikl}^{1}(j)-G_{ikl}^{2}(j)\right)
\right| ^{p}\right) ^{1/p} \\
G_{ikl}^{1}(j) &=&X_{1+\left( \left( k-1\right) 2^{j+1}+l+i-1\right) \bmod 2^{m}} \\ 
G_{ikl}^{2}(j) &=&X_{1+\left( \left( k-1\right) 2^{j+1}+2^{j}+l+i-1\right) \bmod 2^{m}}
\end{eqnarray*}
and similarly for the Energy function. Observe that the inner mean value (on $%
k$) is the previous definition. The outer (on $l$) is the improvement: now,
for every $j$, the calculation involves a mean value taken on $2^{m}$ points.

\begin{figure}
 \centering
 \subfigure[]{\includegraphics[height=150pt, width=200pt]{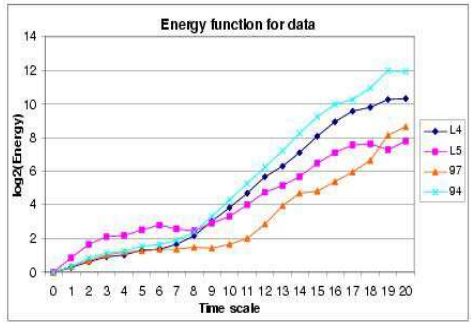}}
 \subfigure[]{\includegraphics[height=150pt, width=200pt]{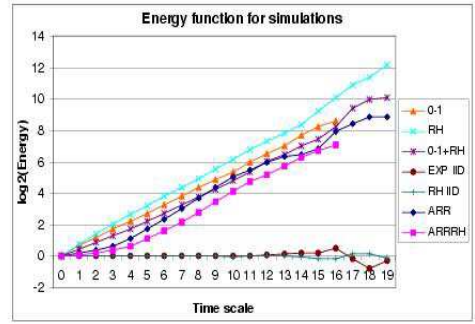}}
 \subfigure[]{\includegraphics[height=150pt, width=200pt]{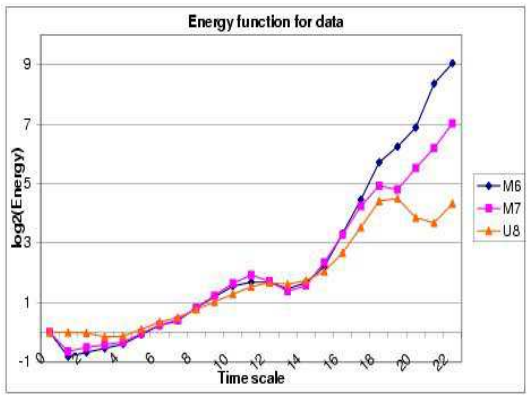}}
 \subfigure[]{\includegraphics[height=150pt, width=200pt]{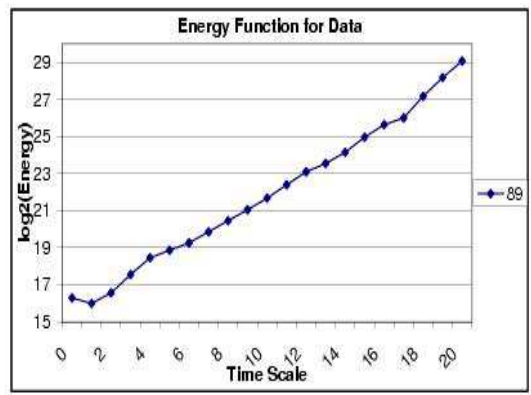}}
 \caption{Energy functions, according to Definition 2, for  data -(a), (c), and (d)- and simulations -(b). Compare 
          to Fig. \ref{EnOld}.}
 \label{En}
\end{figure}

In the case $p=2$, an alternative representation exists:

\begin{equation}
\label{ACAv}
A_{j}^{[2]}=\frac{Var\left( X\right) }{2^{j-1}}\left\{ 1-R\left( 2^{j}\right)
+2\sum_{i=1}^{2^{j}-1}\left[ \left( 1-\frac{i}{2^{j}}\right) \left[ 2R\left(
i\right) -R\left( 2^{j}+i\right) -R\left( 2^{j}-i\right) \right] \right]\right\}
\end{equation}
where $R$ is the autocorrelation of the traffic. This representation proves
that $A_{j}$ is really a linear transformation of the autocorrelation, so
one should expect it to be insensitive to the ``spikiness'' of the traffic.
This is indeed the case (Fig. \ref{En}). Note that the $A_j^{[p]}$ can be computed via a fast algorithm of complexity
$\Theta \left( m2^{m}\right)$; for $A_j$ this can be seen directly from (\ref{ACAv}), since $R$ as well as the convolution in the last term of (\ref{ACAv}) can be computed via the FFT. For general $p$, one can use the multiresolution character of what is essentially a computation with Haar wavelet coefficients to obtain a $\Theta \left( m2^{m}\right)$ algorithm (this also works for $p=2$).

Under the new definition, data and simulation Energy functions assume a 
slightly different form (Fig.~\ref{En}, compare to Fig.~\ref{EnOld}).

\addtocounter{figure}{-1}
\stepcounter{figure}

\section{Model A}

\label {modelA}

This model assumes that the traffic on a link consists of the data that $n$ users are sending, i.e. it is the sum 
of their data streams $W_{i}\left(t\right),\ i=1,...,n$. These are modeled as bi-valuate stochastic processes (at a given $t$, $W\left(t\right)=0\text{ or}\ 1$), whose paths are piecewise constant. The intervals of value $1$ and $0$ will be strictly 
alternating, and their durations will be assumed to be i.i.d.\ for each of the two interval types and for 
each user, cross-independent between interval type and cross-independent between users. 

\begin{definition}
\upshape
Let $\{O_{n}^{a}(i)\}_{i=1}^{\infty }$
and $\{O_{f}^{a}(i)\}_{i=1}^{\infty }$ be i.i.d.\ sequences of positive random
variables with the following properties:
 
\begin{equation}
t^{p} \bold{P}\left(O_{n}^{a}(1)>t\right) \overset{t \rightarrow \infty}{\longrightarrow} C\ 
\text{for}\ p\in (1,2)\ \text{and}\ 0<C<\infty 
\end{equation}

\noindent and 

\begin{equation*}
EO_{f}^{a}(i)^{2}<\infty
\end{equation*}

\noindent Define a process $\widetilde{W}^{a}(t)$ as follows:
 
\begin{equation*}
\widetilde{W}^{a}(t)=\left\{ 
\begin{array}{cc}
1 & \text{ if }\sum_{i=1}^{k-1}\left(O_{f}^{a}(i)+O_{n}^{a}(i)\right)+O_{f}(k)\leq
t\leq \sum_{i=1}^{k}\left(O_{f}^{a}(i)+O_{n}^{a}(i)\right) \\ 
0 & \text{Otherwise}
\end{array}
\right.
\end{equation*}

\noindent and consider $W^{a}(t)$ to be a ``stationary'' version of $\widetilde{W}^{a}(t)$,
i.e. let $u$ be uniformly distributed on the interval 
$\left[0,O_{f}^{a}(1)+O_{n}^{a}(1)\right]$ and define $W^{a}(t)=\widetilde{W}^{a}(u+t)\ $\footnote{
This technique was also used in \cite{MRRS1}}.
Finally, let $\{W_{i}^{a}(t)\}$ be i.i.d.\ such that $W_{i}^{a}(t)\stackrel{\mathcal{L}}{=}W^{a}(t)$,
where $\stackrel{\mathcal{L}}{=}$ denotes equality in the sense of all 
finite dimensional distributions of any order, i.e.
\begin{equation*}
\forall n \in N,\ \forall (t_{1},\ldots, t_{n}),\ \bold{p}(W(t_{1}),\ldots,W(t_{n}))=\bold{p}(W_{i}(t_{1}),\ldots,W_{i}(t_{n}))
\end{equation*}
We will also use the notions of \emph{weak convergence}   
in the sense of the finite dimensional distributions, i.e. \emph{finite dimensional convergence},
denoted by $W_{m}(t)\stackrel{W}{\rightarrow}W(t)$:  
\begin{equation*}
\forall n \in N,\ \forall (t_{1},\ldots, t_{n}),\ \lim_{m\rightarrow\infty} \bold{p}(W_{m}(t_{1}),\ldots,W_{m}(t_{n}))= \bold{p}(W(t_{1}),\ldots,W(t_{n}))
\end{equation*}
and \emph{stochastic equicontinuity}:
\begin{equation*}
\underset{\delta \rightarrow 0}{\lim }\underset{n\rightarrow \infty }{\lim
\sup }\bold{P}\left(\underset{|t-s|\leq \delta }{\sup }\left|W_{n}(t)-W_{n}(s)\right|>%
\varepsilon \right)=0 
\end{equation*}

Finally, a stochastic process $W\left(t\right)$ will be called \emph{self-similar} iff
$\exists\ H>0$: $\forall t>0$, $W\left(t\right)\stackrel{\mathcal{L}}{=}t^{H}W\left(1\right)$. This is certainly
not the only possible definition (see \cite{TTW1} for some alternatives), but it is the most usual. 
 
\end{definition}

This setting yielded a seminal mathematical result about the properties 
of Internet traffic, explaining its self-similarity in a LAN \cite{TWS1}. 
Theorem A is closely connected to this. The partial proof offered spells out in more detail some steps skipped 
in \cite{TWS1}:
 
\begin{theorem}
Let $\displaystyle T_{n}^{a}(t)=\int_{0}^{t}\frac{1}{\sqrt{n}}\sum_{i=1}^{n}W_{i}^{a}(s)ds$. Then:  
\begin{equation*}
Z_{n}^{a}(t)=T_{n}^{a}(t)-\bold{E}T_{n}^{a}(t)\overset{W}{\rightarrow }G^{a}(t),%
\ t\in \lbrack 0,T]
\end{equation*}
where $G(t)$ is a centered Gaussian self-similar process with covariance structure: 
\begin{eqnarray*}
\bold{Cov}\left(G^{a}(s),G^{a}(t)\right)
=\bold{E}\left(\int_{0}^{s}W_{1}^{a}(x)dx\int_{0}^{t}W_{1}^{a}(x)dx\right)-\bold{E}\left(%
\int_{0}^{s}W_{1}^{a}(x)dx\right)\bold{E}\left(\int_{0}^{t}W_{1}^{a}(x)dx\right)
\end{eqnarray*}
\end{theorem}

\begin{proof}

First of all, it is straightforward that:
\begin{equation*}
\bold{Cov}(Z_{n}^{a}(s),Z_{n}^{a}(t))=\frac{1}{n}\sum_{i=1}^{n}%
\sum_{j=1}^{n}\mathbf{E}\left((X_{i}^{s}-\mathbf{E}X_{i}^{s})(X_{j}^{t}-\mathbf{E}X_{j}^{t})\right)
\end{equation*}
\begin{equation*}
=\mathbf{E}\left(\int_{0}^{s}W_{1}^{a}(x)dx\int_{0}^{t}W_{1}^{a}(x)dx\right)-\mathbf{E}\left(%
\int_{0}^{s}W_{1}^{a}(x)dx\right)\mathbf{E}\left(\int_{0}^{t}W_{1}^{a}(x)dx\right)
\end{equation*}

Now, fix $t\in [0,T]$ and define $\displaystyle X_{i}^{t}=\int_{0}^{t}W_{i}^{a}(x)dx$.
With this notation $\displaystyle Z_{n}^{a}(t)=\frac{1}{\sqrt{n}}\sum_{i=1}^{n}\left(X_{i}^{t}-\mathbf{E}X_{i}^{t}\right)$. 

To show finite dimensional convergence, it is sufficient to prove that: $\forall k\in \mathbf{N},\forall t_{1},\ldots,t_{k}\in [0,T],\linebreak[0] \forall \alpha_{1},\ldots,\alpha _{k}\in \mathbb{R}$
$\sum_{j=1}^{k}\alpha _{j}Z_{n}^{a}(t_{j})\rightarrow N(0,\sigma ^{2})$,
where $\sigma ^{2}$ depends on $\alpha ^{\prime }s$ and $t^{\prime }s$.
Easily 
\begin{equation*}
\sum_{j=1}^{k}\alpha _{j}Z_{n}^{a}(t_{j})=\frac{1}{\sqrt{n}}%
\sum_{i=1}^{n}\sum_{j=1}^{k}\left(\alpha _{j}X_{i}^{t_{j}}-\bold{E}\left(\alpha
_{j}X_{i}^{t_{j}}\right)\right)=\frac{1}{\sqrt{n}}\sum_{i=1}^{n}\left(Y_{i}^{a}-\bold{E}Y_{i}^{a}\right)
\end{equation*}
where $\displaystyle Y_{i}^{a}=\sum_{j=1}^{k}\alpha _{j}\int_{0}^{t_{j}}W_{i}^{a}(x)dx$.
Clearly $Y_{i}^{a}$'s are i.i.d.\ and bounded so classical CLT can be invoked.

The next step is to establish tightness. For this it is sufficient (and necessary)
to show \emph{stochastic equicontinuity}. That is, $\forall \varepsilon >0$%
\begin{equation}
\underset{\delta \rightarrow 0}{\lim }\underset{n\rightarrow \infty }{\lim
\sup }\bold{P}\left(\underset{|t-s|\leq \delta }{\sup }\left|Z_{n}^{a}(t)-Z_{n}^{a}(s)\right|>%
\varepsilon \right)=0  \label{StochEqui}
\end{equation}

The following estimate will be needed:

\begin{equation*}
\mathbf{E}(Z_{n}^{a}(t)-Z_{n}^{a}(s))^{2}=Var(Z_{n}^{a}(t)-Z_{n}^{a}(s))
\end{equation*}
\begin{equation*}
=Var(X_{i}^{t}-X_{i}^{s})\leq
\mathbf{E}(X_{i}^{t}-X_{i}^{s})^{2}=\mathbf{E}\left(\int_{s}^{t}W_{i}^{a}(x)dx\right)^{2}\leq (t-s)^{2}
\end{equation*}

A classical ``chaining'' result (see \cite{Tallagrand} Thm. 11.1) yields the following:

\begin{alemma}
Let $\xi _{t}$ be a real-valued
stochastic process $\xi _{t}$ such that $\forall s,t \in 
\mathbb{R}:\ (\mathbf{E}(\xi _{t}-\xi _{s})^{p})^{1/p}\leq d(s,t)$. Then
$\forall A\subset \mathbb{R}:\ Card(A)<\infty $

\begin{equation*}
\mathbf{E}\left(\underset{d(s,t)\leq \delta :s,t\in A}{\sup }|\xi _{t}-\xi _{s}|\right)\leq
8\int_{0}^{\delta }\sqrt[p]{N(A,d,x)}dx
\end{equation*}

\noindent where $d(\cdot,\cdot)$ is any pseudo distance
on $\mathbb{R}$ and $N(A,d,x)$ is minimum number of ``balls'' of
radius $x$ (under distance $d$) needed to cover $A$.
\end{alemma}

\bigskip

In order to use this lemma we first observe that in our case $p=2,\ d(s,t)=|t-s|$ and 

\begin{equation*}
N(A,d,x)\leq N([0,T],d,x)=Tx^{-1}
\end{equation*}

Since the right hand side does not depend on $A$ or $n$, the
restriction $Card(A)<\infty$ can be removed (see Theorem 11.6 in \cite{Tallagrand}). Thus:
 
\begin{equation*}
\bold{P}\left(\underset{|t-s|\leq \delta }{\sup }|Z_{n}(t)-Z_{n}(s)|>\varepsilon \right)\leq
\varepsilon ^{-1}\mathbf{E}\left(\underset{|t-s|\leq \delta }{\sup }|Z_{n}(t)-Z_{n}(s)|\right)
\leq \frac{8}{\varepsilon }\int_{0}^{\delta }\sqrt{N([0,T],d,x)}dx=\frac{16}{%
\varepsilon }T^{1/2}\delta ^{1/2}
\end{equation*}
This establishes stochastic equicontinuity. 

The proof of the self-similarity of the process $G^a$ has already appeared in \cite{TWS1}, and 
therefore it will not be reproduced here. 

\end{proof}

\section{Model B}

\label{modelB}

The first criticism to the usual bi-valuate SSM (model A) is the continuity
of the ON-intervals. The data transfer between computers is quantized
(in \emph{packets}), and although a continuous approximation may be acceptable, if
only time scales much larger than the typical RTT are of interest, for
studies of traffic in the sub-second ranges a more detailed model has to be
considered.

Thus, the first modification is the discretization of ON-intervals. 
Consider the usual SSM, with $W_{i}^{b}(t)$ independent (for $%
i=1,...,n).$ For a fixed connection (fixed $i$), the OFF-intervals will
have lengths ${O_{f}^{b}(j)},\ j\in N $, where the $O_{f}$'s are i.i.d.\ random
variables of finite variance. At the end of an OFF-interval, a \emph{%
heavy tail distributed random integer} $L$ will be assigned, representing
the ``load'' (number of packets) the user wishes to send in this particular
session. More precisely, 

\begin{equation}
\exists p\in (1,2),\exists C>0:\ t^{p} \bold{P}\left(L>t\right) \overset{t \rightarrow \infty}{\longrightarrow} C
\label{PP}
\end{equation}

At this point, an ON-interval begins, but, in this model, packet emissions assume the form of delta 
functions, expressing mathematically the fact that information is transmitted in individual packets 
and in particular time instants, and the $j$-th packet emission of the $i$-th 
ON-interval is followed by the idle time interval $R_{i,j}>0$: these are random variables
i.i.d.\ in $i$ and $j$ (with $\bold{E}(R_{i,j}^{2})<\infty $), independent of the OFF-interval lenghts and the session
loads $L$. After all of the $L$
packets have been sent, an OFF-interval follows, and the cycle restarts. Notice that
ON-intervals (representing sessions) are heavy-tailed, as consisting of a heavy-tailed number
of light tailed time intervals.

The random inter-arrival times $R_{i,j}$'s are designed to address the
uncertainty of the length of RTT's. Since exponential distribution has the
highest entropy (for any positive random variable), it would be a reasonable
choice.  

The total traffic from time $0$ until time $t$ is defined by $\displaystyle
T_{n}^{b}(t)=\int_{0}^{t}\sum_{i=1}^{n}W_{i}^{b}(s)ds.$ Theorem B below asserts the
convergence of properly normalized total traffic, under very general
conditions on the distributions of the $L$ and $R_{i,j}$, to a centered Gaussian
process with continuous sample paths.

\begin{definition}
\upshape
Let $\{O_{f}^{b}(i)\}_{i=1}^{\infty}$ be the same as 
$\{O_{f}^{a}(i)\}_{i=1}^{\infty }$ and let $\{L^{b}(i)\}_{i=1}^{\infty}$
be the same as $\{O_{n}^{a}(i)\}_{i=1}^{\infty }$, but with the extra requirement
that $L^{b}(i)$ be integer valued. Let also $\{R^{b}_{i,j}\}_{i,j=1}^{\infty }$ be
an i.i.d.\ array of random variables so that:
 
\begin{equation*}
R^b_{i,j}\geq 0\ a.s.,\ \bold{E}(R^{b}_{i,j})^{2}<\infty,\  R^{b}_{i,0}=0\ a.s.
\end{equation*}
and 
\begin{equation}
\underset{\delta \rightarrow 0}{\lim }\delta ^{-1}\bold{P}(R^{b}_{i,j}<\delta )\leq
M<\infty  \label{Alpha}
\end{equation}

Let also $\displaystyle O_{n}^{b}(i)=\sum_{j=1}^{L(i)}R^{b}_{i,j}.$ Finally, define formally the
``delta'' function $\widehat{\delta }_{x}(s)$ with the properties that $\widehat{%
\delta }_{x}(s)=0$ if $x\neq s$ and $\displaystyle \int_{x-h}^{x+h}\widehat{\delta }%
_{x}(s)ds=1$ for all $h>0$. Using these conventions, define further: 
\begin{equation*}
\widetilde{W}^{b}(t)=\left\{ 
\begin{array}{cc}
\widehat{\delta }_{t}(t) & \text{if }t=%
\sum_{i=1}^{k-1}\left(O_{f}^{b}(i)+O_{n}^{b}(i)\right)+O_{f}^{b}(k)+%
\sum_{i=1}^{l}R^{b}_{k,i}\text{ , }l\leq L^{b}(k) \\ 
0 & \text{Otherwise}
\end{array}
\right.
\end{equation*}
As was the case before, define $W^{b}(t)$ to be a stationary version of $%
\widetilde{W}^{b}(t)$, and let $\{W_{i}^{b}(t)\}$ be i.i.d.\ versions of $%
W^{b}(t)$. Finally, define a constant $\bold{E}W^{b}$ such that $%
E\int_{s}^{t}W^{b}(x)dx=\bold{E}W^{b}|t-s|$.
\end{definition}

For this process, a result analogous to Theorem A holds:

\begin{theorem}
Let $\displaystyle T_{n}^{b}(t)=\int_{0}^{t}\frac{1}{\sqrt{n}}\sum_{i=1}^{n}W_{i}^{b}(s)ds$. Then: 
 
\begin{equation*}
Z_{n}^{b}(t)=T_{n}^{b}(t)-ET_{n}^{b}(t)\overset{W}{\rightarrow }G^{b}(t),\ for\ t\in [0,T]
\end{equation*}

where $G^{b}(t)$ is a centered Gaussian process with covariance structure:
\begin{eqnarray*}
\bold{Cov}(G^{b}(s),G^{b}(t))
=\bold{E}\left(\int_{0}^{s}W_{1}^{b}(x)dx\int_{0}^{t}W_{1}^{b}(x)dx\right)-\bold{E}\left(%
\int_{0}^{s}W_{1}^{b}(x)dx\right)\bold{E}\left(\int_{0}^{t}W_{1}^{b}(x)dx\right)
\end{eqnarray*}
\end{theorem}

\begin{proof}
 
The assertion will be proved under a slightly stronger condition on the $R^{b}_{i,j}$ than (\ref{Alpha}); more precisely, it will be assumed here (as well as throughout the paper) that $\exists \beta >0$:

\begin{equation}
0<\beta \leq R^{b}_{i,j}\text{ a.s.}  \label{Beta}
\end{equation}
The proof under the less restrictive assumption (\ref{Alpha}) is more technical, and the differences will
be commented upon in the course of the proof; these comments will be indicated by square brackets [ ]. 

The finite dimensional convergence follows exactly in the same line as in Theorem
A. Namely, $\displaystyle Z_{n}^{b}(t) =\frac{1}{\sqrt{n}}\sum_{i=1}^{n}\left(X_{i}^{t}-\bold{E}\left(X_{i}^{t}\right)\right)$, where $\displaystyle X_{i}^{t}=\int_{0}^{t}W_{i}^{b}(x)dx.$ Under the assumption (\ref{Beta}), $X_{i}^{t}$ are
bounded by $\beta ^{-1}T$ and classical CLT can be used. [A short computation
reveals that, under the original assumption (\ref{Alpha}), $X_{i}^{t}$ has
exponential tails and one can still use the same argument.]

Stochastic equicontinuity (\ref{StochEqui}) will be proved in two steps:

\bigskip

\noindent \textit{Step 1}. It is sufficient to show that for every $\varepsilon >0$
and $\gamma >0$ 
\begin{equation}
\underset{\delta \rightarrow 0}{\lim }\underset{n\rightarrow \infty }{\lim
\sup }\bold{P}\left(\underset{\gamma n^{-1/2}\leq |t-s|\leq \delta }{\sup }%
\left|Z_{n}^{b}(t)-Z_{n}^{b}(s)\right|>\varepsilon \right)=0
\label{StEq}
\end{equation}
\noindent Proof: $\forall \varepsilon >0$ consider the partition $\displaystyle \Lambda
_{n}= \left\{t_{i}=\frac{i}{K_{\varepsilon }\sqrt{n}}\ : \ i=0,...,K_{\varepsilon
}\sqrt{n}T\right\} $ where $K_{\varepsilon }=\varepsilon ^{-1}8\bold{E}W_{1}^{b}$. For
every $t$, define functions $\phi _{1}^{n}$ and $\phi _{2}^{n}$: $%
[0,T]\rightarrow \Lambda _{n}$, such that $\phi _{1}^{n}(t)=t_{j}\leq t\leq
t_{j+1}=\phi _{2}^{n}(t)$ for some $j\leq K_{\varepsilon }\sqrt{n}T$.

In other words $\phi _{1}^{n}$ and $\phi _{2}^{n}$ are closest elements to $t$ from
the cover $\Lambda _{n}$ that enclose $t$.

Clearly 
\begin{equation*}
T_{n}^{b}(\phi _{1}^{n}(t))\leq T_{n}^{b}(t)\leq T_{n}^{b}(\phi _{2}^{n}(t))
\end{equation*}
which implies 
\begin{multline}
\label{Bcomput}
Z_{n}^{b}(t)-Z_{n}^{b}(s)=
T_{n}^{b}(t)-\bold{E}\left(T_{n}^{b}(t)\right)-T_{n}^{b}(s)+\bold{E}\left(T_{n}^{b}(s)\right) \\
\leq T_{n}^{b}(\phi _{2}^{n}(t))-\bold{E}\left(T_{n}^{b}(\phi _{1}^{n}(t))\right)-T_{n}^{b}(\phi
_{1}^{n}(s))+\bold{E}\left(T_{n}^{b}(\phi _{2}^{n}(t))\right)=\\
=T_{n}^{b}(\phi _{2}^{n}(t))-\bold{E}\left(T_{n}^{b}(\phi _{2}^{n}(t))\right)+\bold{E}\left(T_{n}^{b}(\phi
_{2}^{n}(t))\right)-\bold{E}\left(T_{n}^{b}(\phi _{1}^{n}(t))\right) \\
+\bold{E}\left(T_{n}^{b}(\phi _{1}^{n}(s))\right)-T_{n}^{b}(\phi _{1}^{n}(s))+\bold{E}\left(T_{n}^{b}(\phi
_{2}^{n}(t))\right)-\bold{E}\left(T_{n}^{b}(\phi _{1}^{n}(s))\right)
\end{multline}
Since $\displaystyle \bold{E}\left(T_{n}^{b}(\phi _{2}^{n}(t))\right)-\bold{E}\left(T_{n}^{b}(\phi _{1}^{n}(t))\right)=\bold{E}\left(\frac{1}{\sqrt{n}}\sum_{i=1}^{n}\int_{\phi _{1}^{n}(t)}^{\phi
_{2}^{n}(t)}W_{1}^{b}(x)dx\right)=\sqrt{n}(\phi _{2}^{n}(t)-\phi
_{1}^{n}(t))\bold{E}W_{1}^{b}= \bold{E}W_{1}^{b}/K_{\varepsilon }=\varepsilon /8$, one gets: 
\begin{equation*}
Z_{n}^{b}(t)-Z_{n}^{b}(s)\leq Z_{n}^{b}(\phi _{2}^{n}(t))-Z_{n}^{b}(\phi_{1}^{n}(s))+\varepsilon /4
\end{equation*}
Similarly 
\begin{equation*}
Z_{n}^{b}(s)-Z_{n}^{b}(t)\leq Z_{n}^{b}(\phi _{2}^{n}(s))-Z_{n}^{b}(\phi_{1}^{n}(t))+\varepsilon /4
\end{equation*}
Combining the two inequalities: 

\begin{equation*}
\left|Z_{n}^{b}(t)-Z_{n}^{b}(s)\right|\leq \left|Z_{n}^{b}(\phi _{2}^{n}(t))-Z_{n}^{b}(\phi
_{1}^{n}(s))\right|+\left|Z_{n}^{b}(\phi _{2}^{n}(s))-Z_{n}^{b}(\phi
_{1}^{n}(t))\right|+\varepsilon /2
\end{equation*}

and since $\underset{|t-s|\leq \delta }{\sup }\left|Z_{n}^{b}(\phi
_{2}^{n}(t))-Z_{n}^{b}(\phi _{1}^{n}(s))\right|=\underset{|t-s|\leq \delta }{\sup }%
\left|Z_{n}^{b}(\phi _{2}^{n}(s))-Z_{n}^{b}(\phi _{1}^{n}(t))\right|$: 

\begin{multline*}
\bold{P}\left(\underset{|t-s|\leq \delta }{\sup }|Z_{n}^{b}(t)-Z_{n}^{b}(s)|>\varepsilon
\right)\leq \bold{P}\left(\underset{|t-s|\leq \delta }{\sup }\left|Z_{n}^{b}(\phi
_{2}^{n}(t))-Z_{n}^{b}(\phi _{1}^{n}(s))\right|>\varepsilon /4\right) \overset{*}{\leq}\\
\bold{P}\left(\underset{K_{\varepsilon }^{-1}n^{-1/2}\leq |t-s|\leq \delta }{\sup }%
|Z_{n}^{b}(\phi _{2}^{n}(t))-Z_{n}^{b}(\phi _{1}^{n}(s))|>\varepsilon /4\right)
\end{multline*}

where in ($*$) we used that $|t-s|\leq \delta \Rightarrow \phi _{2}^{n}(t)-\phi _{1}^{n}(s)\geq
K_{\varepsilon }^{-1}n^{-1/2}$.

This proves Step 1.

\bigskip

\noindent \textit{Step 2}. The goal here is to show that (\ref{StEq}) holds for all $\varepsilon ,\gamma >0$ 

\begin{equation*}
\underset{\delta \rightarrow 0}{\lim }\underset{n\rightarrow \infty }{\lim
\sup }\bold{P}\left(\underset{\gamma n^{-1/2}\leq |t-s|\leq \delta }{\sup }%
\left|Z_{n}^{b}(t)-Z_{n}^{b}(s)\right|>\varepsilon \right)=0
\end{equation*}

Using the notation $\displaystyle X_{i}^{t}=\int_{0}^{t}W_{i}^{b}(x)dx,\ 
Z_{n}^{b}(t)=\frac{1}{\sqrt{n}}\sum_{i=1}^{n}(X_{i}^{t}-\bold{E}\left(X_{i}^{t}\right)$:

\begin{multline}
\bold{E}\left(Z_{n}^{b}(t)-Z_{n}^{b}(s)\right)^{4}=\bold{E}\left( \frac{1}{\sqrt{n}}%
\sum_{i=1}^{n}\left[ (X_{i}^{t}-X_{i}^{s})-\bold{E}(X_{i}^{t}-X_{i}^{s})\right]
\right)^{4} \\
\leq \frac{1}{n}\bold{E}\left( (X_{1}^{t}-X_{1}^{s})-\bold{E}(X_{1}^{t}-X_{1}^{s})\right)
^{4}+6\left( \bold{E}\left( (X_{1}^{t}-X_{1}^{s})-\bold{E}(X_{1}^{t}-X_{1}^{s})\right)
^{2}\right) ^{2} \\
\leq \frac{1}{n}\left(
\bold{E}(X_{1}^{t}-X_{1}^{s})^{4}+6\bold{E}(X_{1}^{t}-X_{1}^{s})^{2}\left(
\bold{E}(X_{1}^{t}-X_{1}^{s})\right) ^{2}\right) +6\left(
\bold{E}(X_{1}^{t}-X_{1}^{s})^{2}\right) ^{2}  \label{**}
\end{multline}

Let us examine the distribution of $Q_{s,t}=X_{1}^{t}-X_{1}^{s}=%
\int_{s}^{t}W_{1}^{b}(x)dx.$ The assumption (\ref{Beta}) coupled with the
definition of $W_{i}^{b}(x)$ implies that for $\delta $ small enough, there
exists a constant $C$ such that 

\begin{equation*}
Q_{s,t}=1\text{ with probability }C|t-s|
\end{equation*}

\begin{equation}
Q_{s,t}=0\text{ with probability }1-C|t-s|  \label{BetaDist}
\end{equation}

Thus, for $\delta $ small enough, the above expression (\ref{**}) is bounded
by $C(|t-s|^{2}+|t-s|/n)$ for some (possibly different) constant $C.$ \
Therefore, for $1>\delta \geq |t-s|>1/\sqrt{n},\ \exists C_{1}>0$ such that: 

\begin{equation}
\left\| Z_{n}^{b}(t)-Z_{n}^{b}(s)\right\| _{L_{4}}\leq C_{1}|t-s|^{1/2}
\label{L4}
\end{equation}

[It is not hard to show that, under the original assumption (\ref{Alpha}), there
exist constants $\nu _{1},\nu _{2}>0$ such that, for $\delta $ small enough, 

\begin{equation}
P(Q_{s,t}=k)\leq \nu _{1}|t-s|^{k},\ P(Q_{s,t}=0)\geq 1-\nu_{2}|t-s|.\text{]}  
\label{AlphaDist}
\end{equation}

With this estimate one can easily recover (\ref{L4}) under the original assumption (\ref{Alpha}).

Lemma A can be invoked now. Letting $d(s,t)=|t-s|^{1/2}$ we have (up to a
constant $C_{1}$) 
\begin{equation*}
\bold{E}\left(\underset{\gamma n^{-1/2}\leq |t-s|\leq \delta }{\sup }%
\left|Z_{n}^{b}(t)-Z_{n}^{b}(s)\right|\right)\leq 8\int_{0}^{\delta }\sqrt[4]{N([0,T],d,x)}dx
\leq 8\int_{0}^{\delta }T^{1/2}x^{-1/2}dx=16T^{1/2}\delta ^{1/2}
\end{equation*}

\end{proof}

Probably more interesting is Corollary B, which derives some properties of the limit process: 

\stepcounter{corrolary}
\begin{corrolary}
The Gaussian process $G^{b}(t)$ has the
following properties. There exist constants $0<C_{i}<\infty,\ i=1,2$, such that:

\newcounter{temp}
\begin{list}{\roman{temp})}{\usecounter{temp}}

\item $\displaystyle \frac{\bold{Var}(G^{b}(\Delta ))}{\Delta ^{2/p}}\underset{\Delta
\rightarrow \infty }{\longrightarrow }C_{1}$

\item $\displaystyle \frac{\bold{Var}(G^{b}(\Delta ))}{\Delta ^{1/2}}\underset{\Delta
\rightarrow 0}{\longrightarrow }C_{2}$

\item  
$\displaystyle \forall k\in N,\ \bold{Cov}\left(G^{b}(\Delta ),G^{b}((k+1)\Delta )-G^{b}(k\Delta )\right)
\underset{\Delta \rightarrow \infty }{\longrightarrow }k^{\frac{2-2p}{p}}$ (\upshape {$p$ was defined in {\upshape (\ref{PP})}})
 
\item $\displaystyle \forall k\in N,\ \bold{Corr}\left(G^{b}(\Delta ),G^{b}((k+1)\Delta )-G^{b}(k\Delta )\right)%
\underset{\Delta \rightarrow 0}{\longrightarrow }0$
\end{list}

\end{corrolary}

\begin{proof}

Arguments for $i)$ and $iii)$ are essentially the same as in
\cite{TWS1}. Namely, it is easy to show that for ``large''
$\Delta$ the processes $G^{a}$ and $G^{b}$ behave in a similar way. In particular, it
can be shown that there exists $0<K_{p}<\infty$, such that:
 
\begin{equation*}
\frac{\bold{E}\left(G^{b}(\Delta )\right)^{p}}{\bold{E}\left(G^{a}(\Delta )\right)^{p}}\underset{\Delta
\rightarrow \infty }{\longrightarrow }K_{p},\ \forall p>0
\end{equation*}

For \textit{ii)} and \textit{iv)} the same approach as in the proof
for Theorem B will be used. Namely, details will be presented only under assumption (\ref
{Beta}), and comments will be given on the extension of the results under the original assumption (\ref{Alpha}).

For \textit{ii)}, one can show easily (using (\ref{BetaDist})) that for $\Delta < \beta $%
\begin{equation*}
\bold{Var}\left(G^{b}(\Delta )\right)=\bold{Var}\left(\int_{0}^{\Delta }W_{1}^{b}(t)dt\right)=\bold{Var}\left(Q_{0,\Delta}\right)=C_{2}\Delta
\end{equation*}
[Under the original assumption (\ref{Alpha}), estimate (\ref{AlphaDist}) holds which also implies $ii)$.]

\bigskip

For \textit{iv)}, one gets:

\begin{multline*}
\bold{Corr}\left(G^{b}(\Delta ),G^{b}((k+1)\Delta )-G^{b}(k\Delta )\right)= \\
=\frac{\bold{E}\left( \int_{0}^{\Delta }W_{1}^{b}(t)dt\int_{k\Delta }^{(k+1)\Delta
}W_{1}^{b}(t)dt\right) -\bold{E}\left(\int_{0}^{\Delta }W_{1}^{b}(t)dt\right)\bold{E}\left(\int_{k\Delta
}^{(k+1)\Delta }W_{1}^{b}(t)dt\right)}{\bold{Var}\left(\int_{0}^{\Delta }W_{1}^{b}(t)dt\right)}
\end{multline*}

Invoking the distribution of $Q_{s,t}=$ $\int_{s}^{t}W_{1}^{b}(t)dt$ for
small $\Delta $, it follows that $Var(Q_{0,\Delta })\approx \Delta $ and $%
EQ_{0,\Delta }\approx \Delta ,$ as well as: 

\begin{equation*}
Q_{0,\Delta }Q_{k\Delta ,(k+1)\Delta }=0\ a.s.
\end{equation*}

The combination of these three estimates yields:
 
\begin{equation*}
\bold{Corr}\left(G^{b}(\Delta ),G^{b}((k+1)\Delta )-G^{b}(k\Delta )\right)\approx \Delta
\end{equation*}

[Again, under the original assumption (\ref{Alpha}), the proof turns out to be more involved. 
However, it is not hard to show (using the estimate (\ref{BetaDist}))
that in this case $\bold{E}\left(Q_{0,\Delta }Q_{k\Delta ,(k+1)\Delta }\right)=C\Delta
^{2}+o(\Delta ^{2})$ which implies $iv)$ as well.]

\end{proof}

\bigskip

The result suggests that the behavior of the process varies according to
time scale: for coarse ones, it behaves like the SSM, but for fine ones it
behaves more like a Brownian motion. This change of behavior is also
reflected in the autocorrelation.

\begin{figure}[t]
 \centering
 \subfigure[]{\includegraphics[height=150pt, width=200pt]{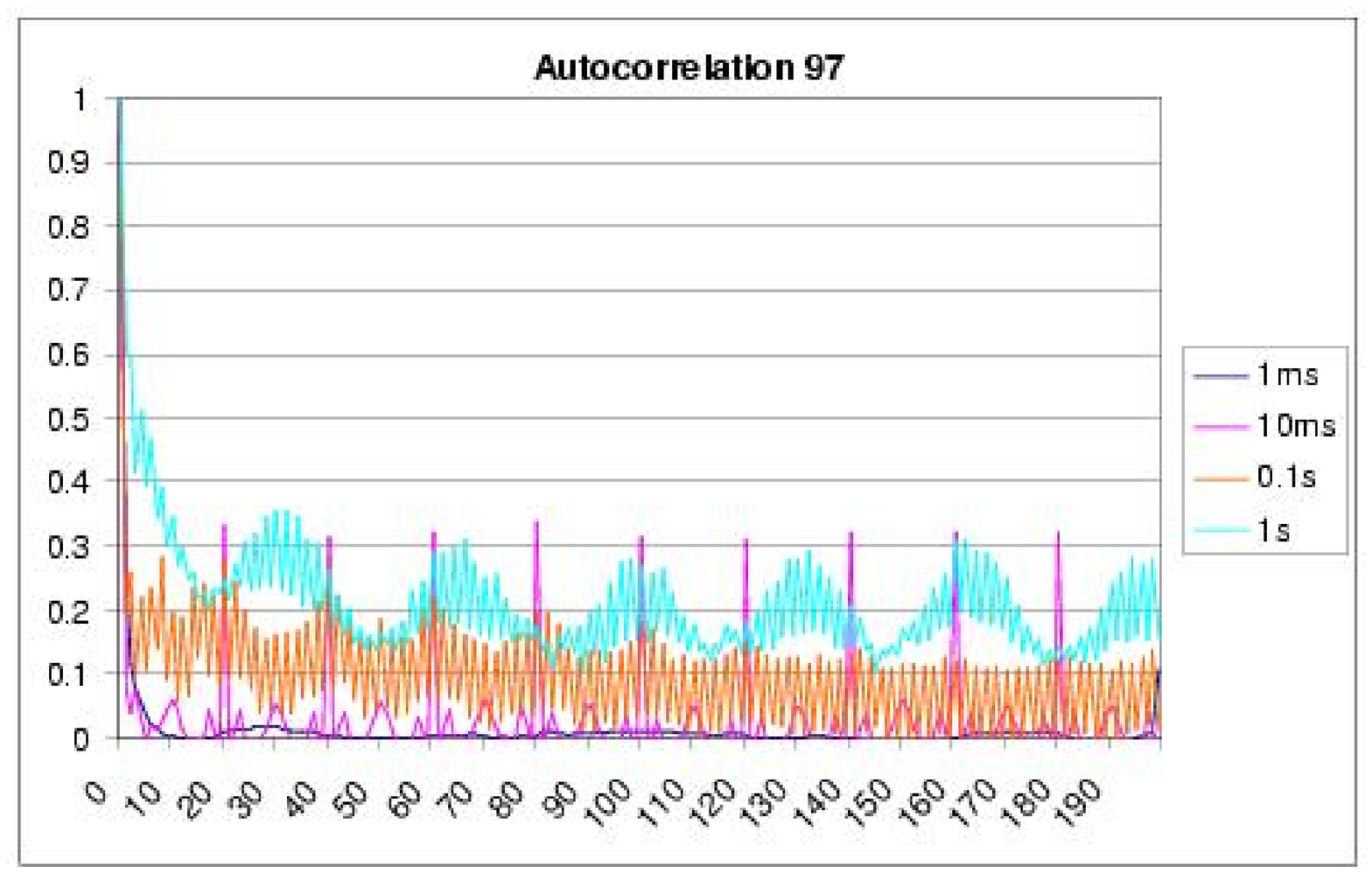}}
 \subfigure[]{\includegraphics[height=150pt, width=200pt]{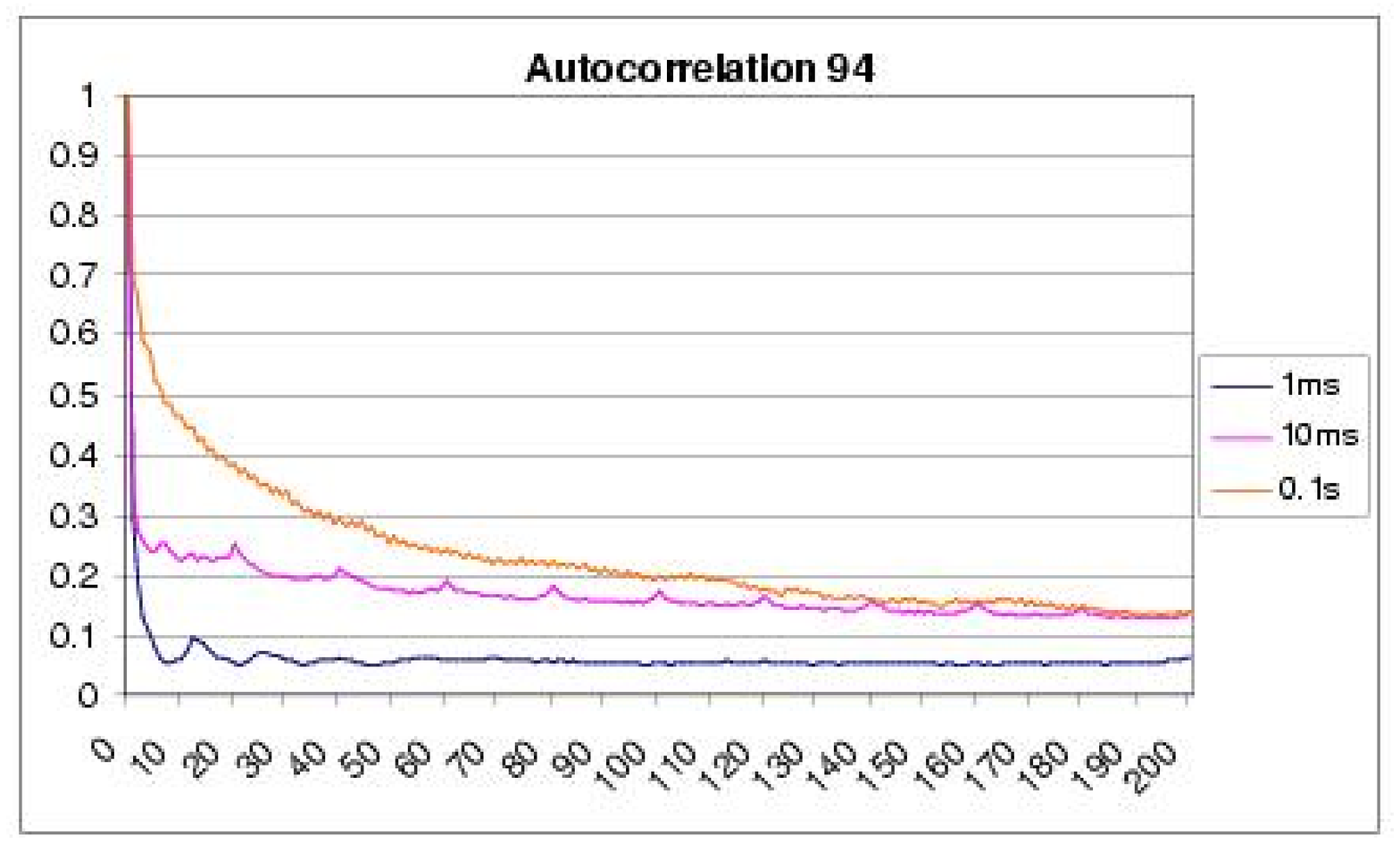}}
 \subfigure[]{\includegraphics[height=150pt, width=200pt]{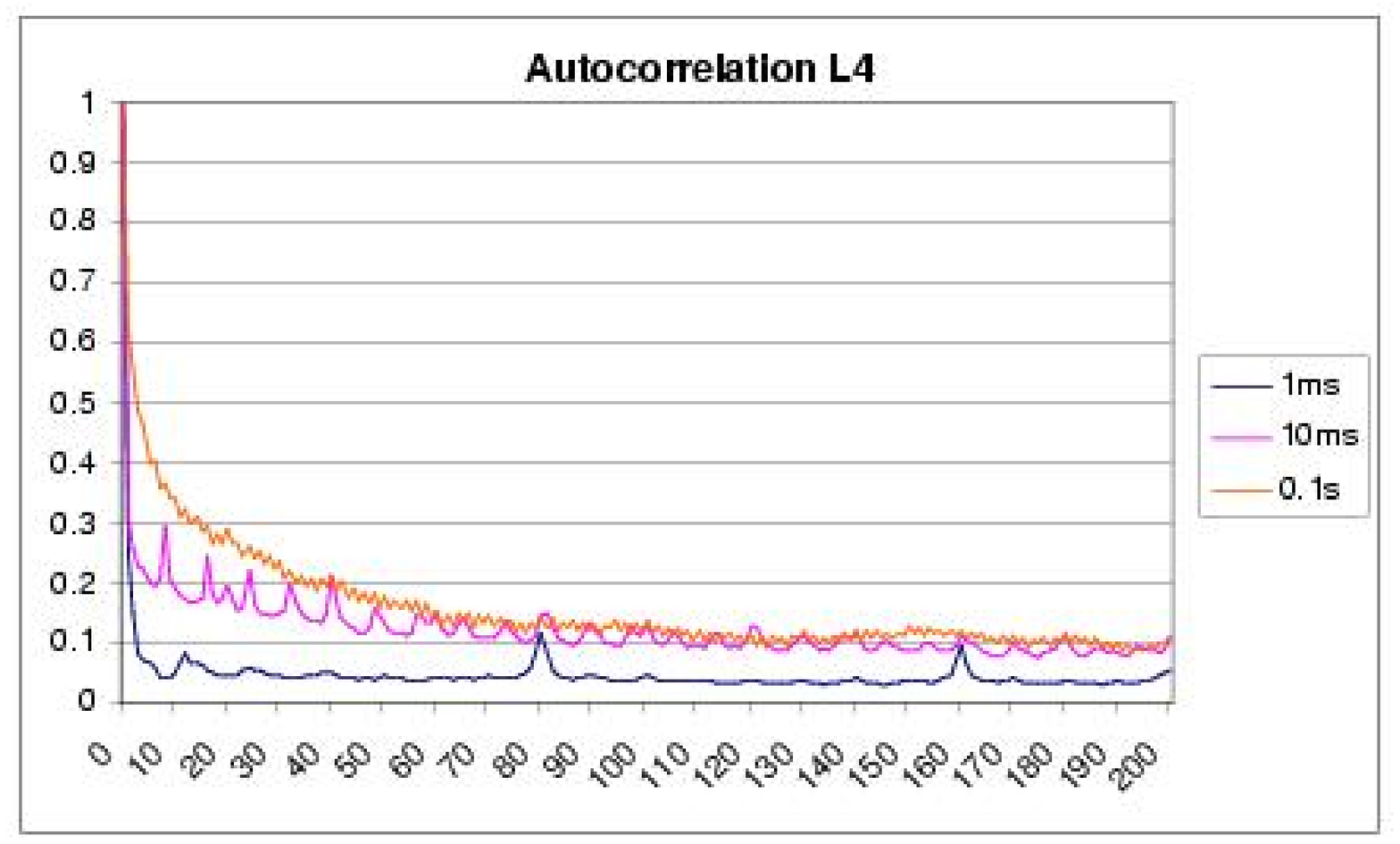}}
 \subfigure[]{\includegraphics[height=150pt, width=200pt]{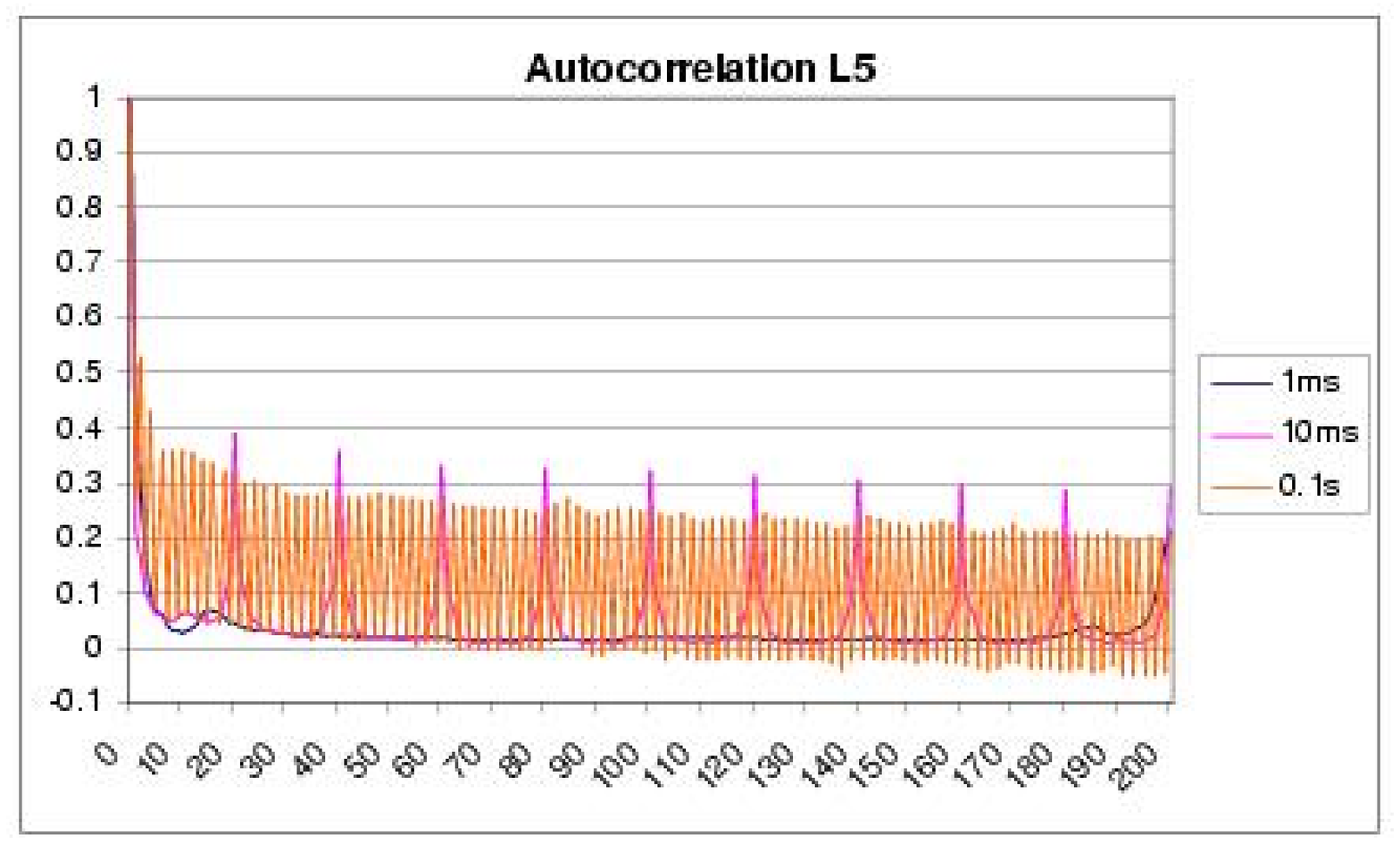}}
 \caption{Trace autocorrelations at four different time bins (1ms, 10ms, 100ms, 1s): the 
          long-range dependence in coarse time scales is obvious. 
          97 (a) and L5 (d) seem to have very strong periodic components, probably due to some monitoring protocol.}
 \label{AC}         
\end{figure}
\addtocounter{figure}{-1}
\stepcounter{figure}

\begin{figure}[t]
 \centering
 \subfigure[]{\includegraphics[height=150pt, width=200pt]{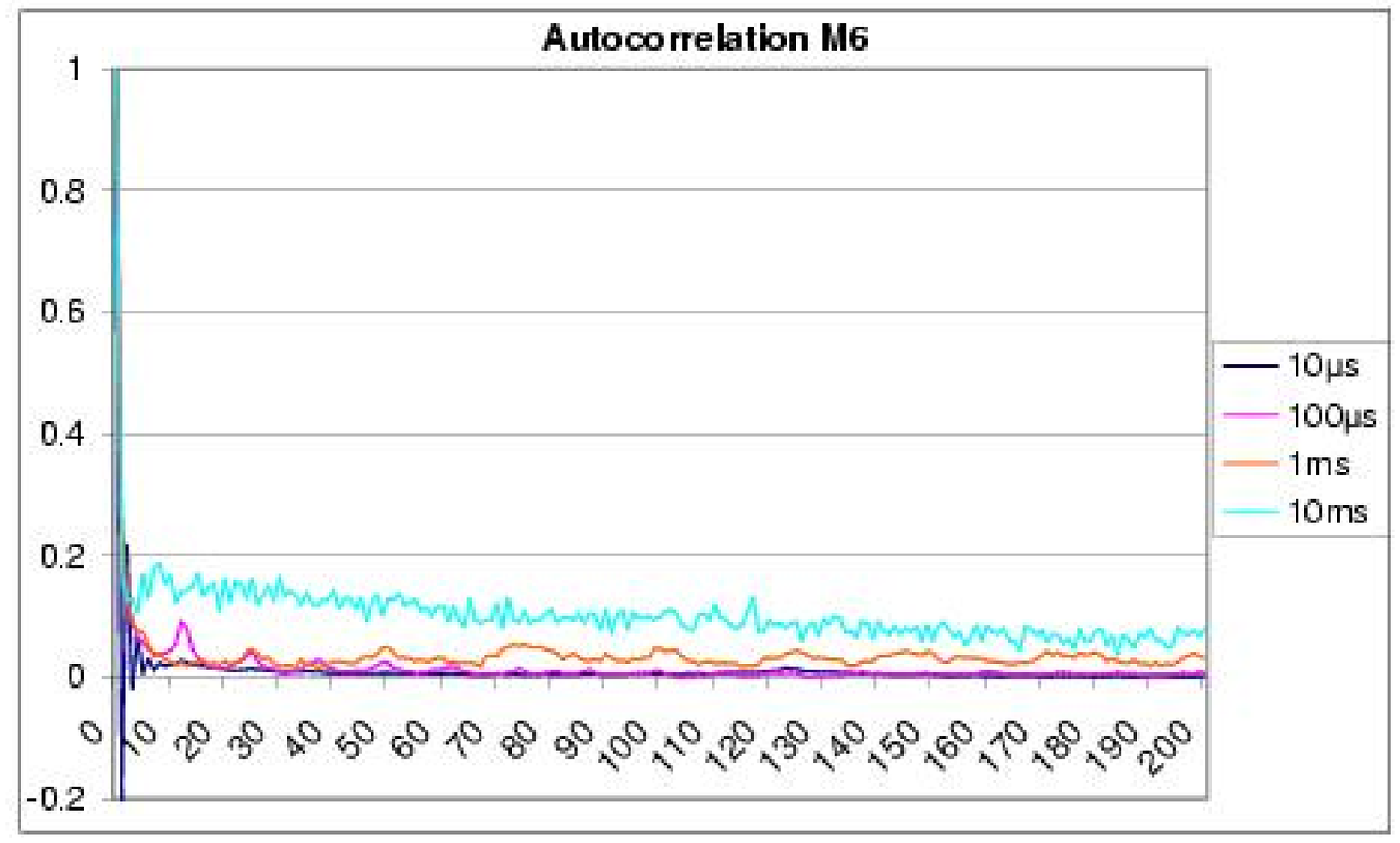}}
 \subfigure[]{\includegraphics[height=150pt, width=200pt]{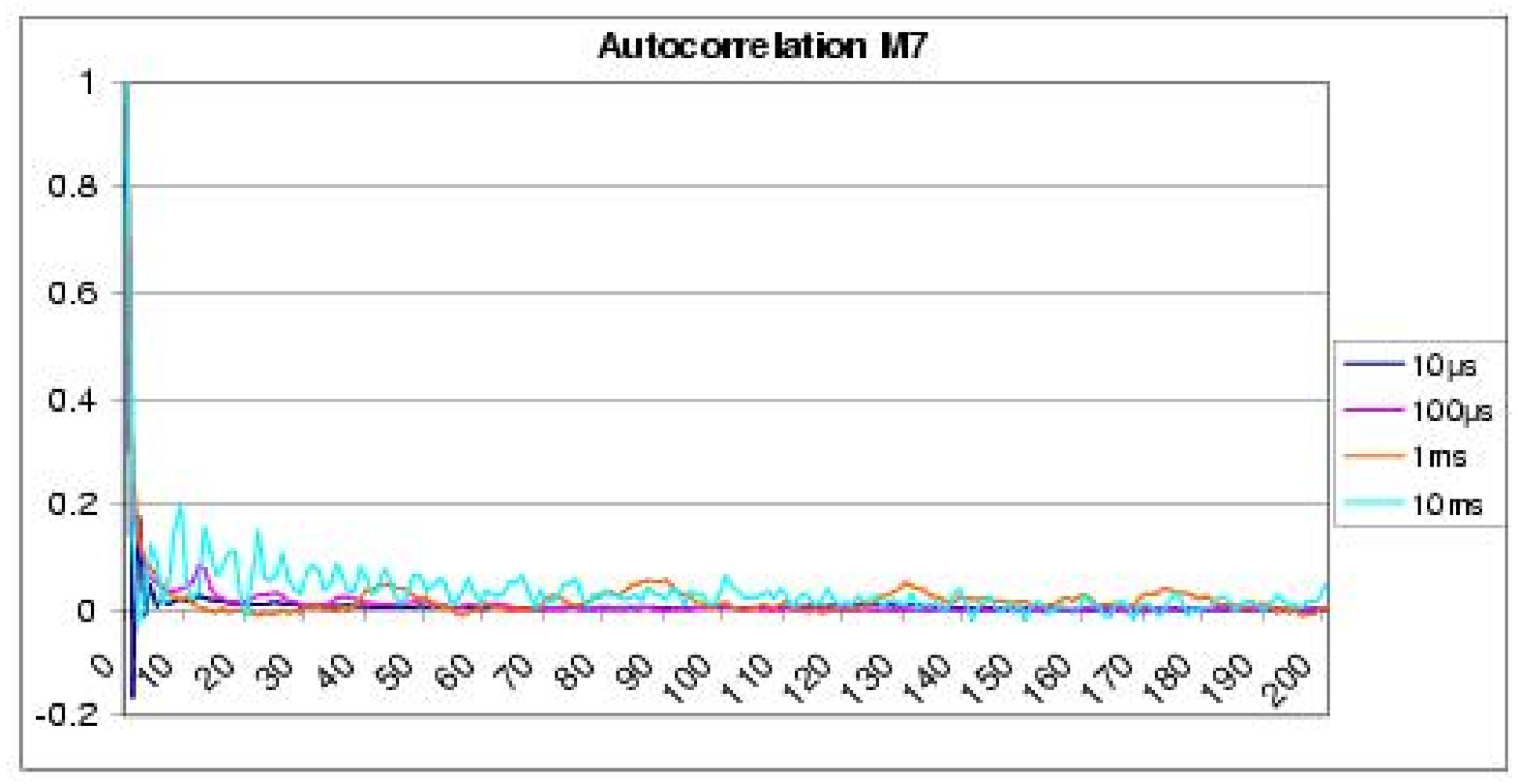}}
 \subfigure[]{\includegraphics[height=150pt, width=200pt]{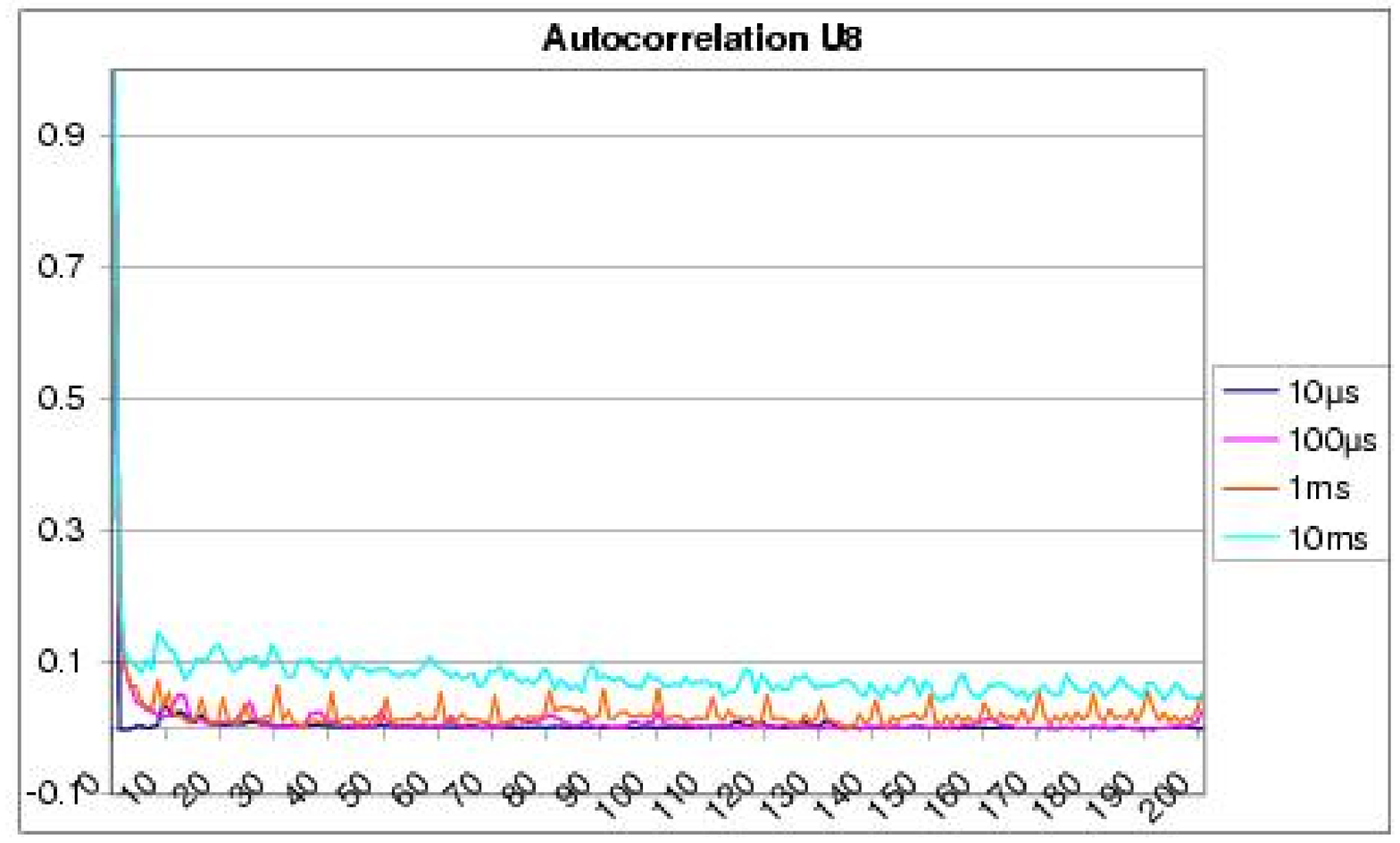}}
 \subfigure[]{\includegraphics[height=150pt, width=200pt]{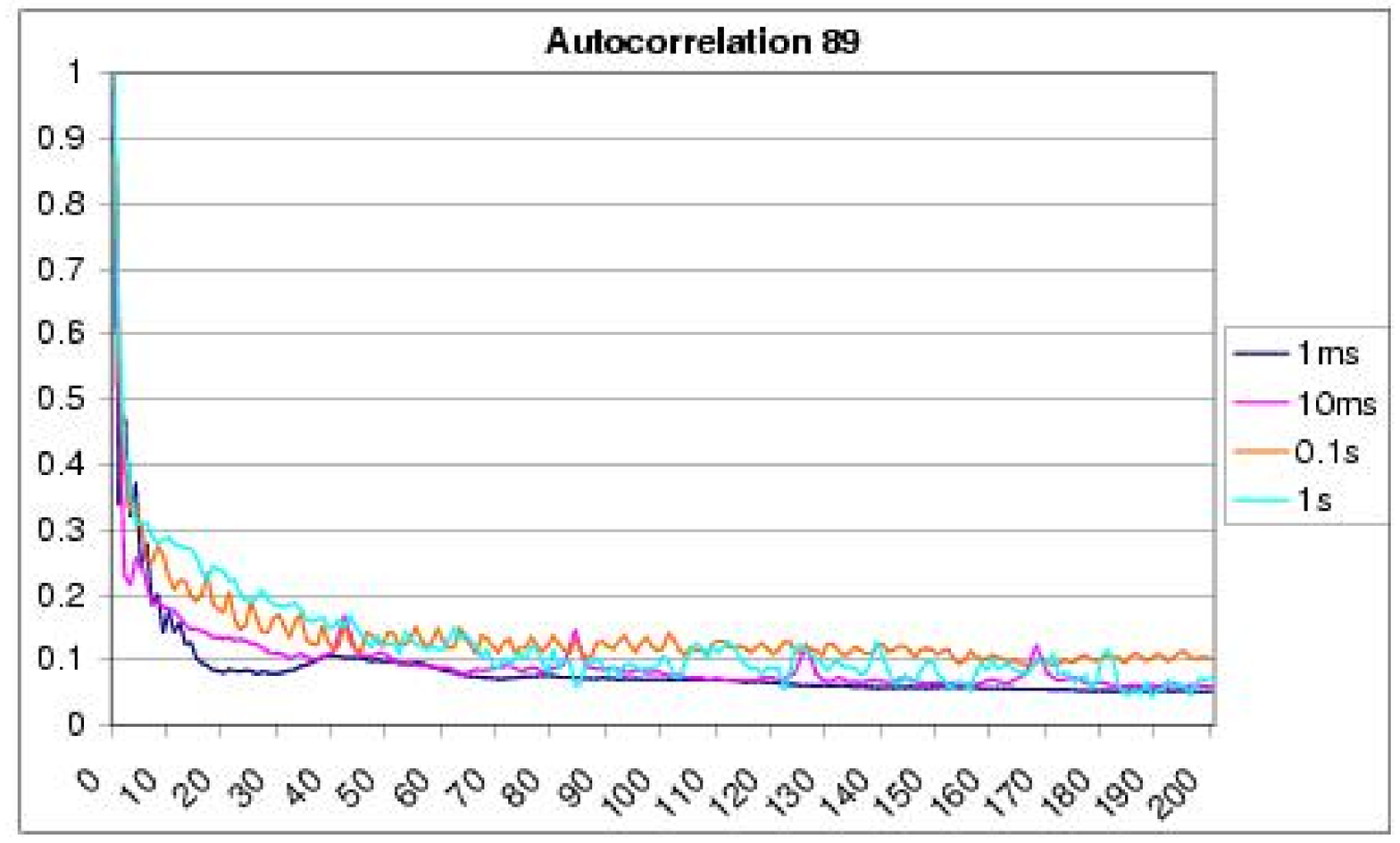}}
 \caption{Trace autocorrelations at four different time bins: 10$\mu$s, 100$\mu$s, 1ms, 10ms for M6, M7, and U8, and 
          1ms, 10ms, 100ms, 1s for 89. Trace 89 shows a relatively strong correlation irrespectively of the bin used.
          On the contrary, the newer traces seem to be uncorrelated, if binned with a bin other than 10ms. If, however, 
          the 10ms bin is used, M7 remains uncorrelated, but M6 and U8 show a weak correlation.}
 \label{ACND}         
\end{figure}
\addtocounter{figure}{-1}
\stepcounter{figure}

The autocorrelations of the data traces follow the same pattern (fast decay at fine time scales, long-range dependence at coarser scales;  see Fig. \ref{AC} and \ref{ACND}), consistent with the predictions of Corrolary B.
This argues against ``pure'' self-similarity of the true traffic, for it
behaves differently on different time scales. Moreover, Theorem B reveals
a mechanism that is responsible for this: the quantized form of data
transfer (packets and RTTs). Note that the result is independent of the
exact distribution of the RTTs, which makes it robust; it can be expected to 
persist even for the haphazard nature of RTT distribution behavior in the network.

It should be mentioned that similar behavior has been observed before. For
instance, in large capacity links of the Internet backbone, packet
interarrival times tend to be exponentially distributed and independent~\cite{CCLS1}.

Model B was successful in capturing the autocorrelation of real traffic, but what about spikiness? Actually, 
Theorem B states that model B traffic is only as spiky as the SSM. Model C will introduce a new feature that will add spikiness to simulated traffic. 

\section{Model C}

\label{modelC}

\subsection{Introduction}

\begin{figure}
 \centering
 \subfigure[]{\includegraphics[height=150pt, width=200pt]{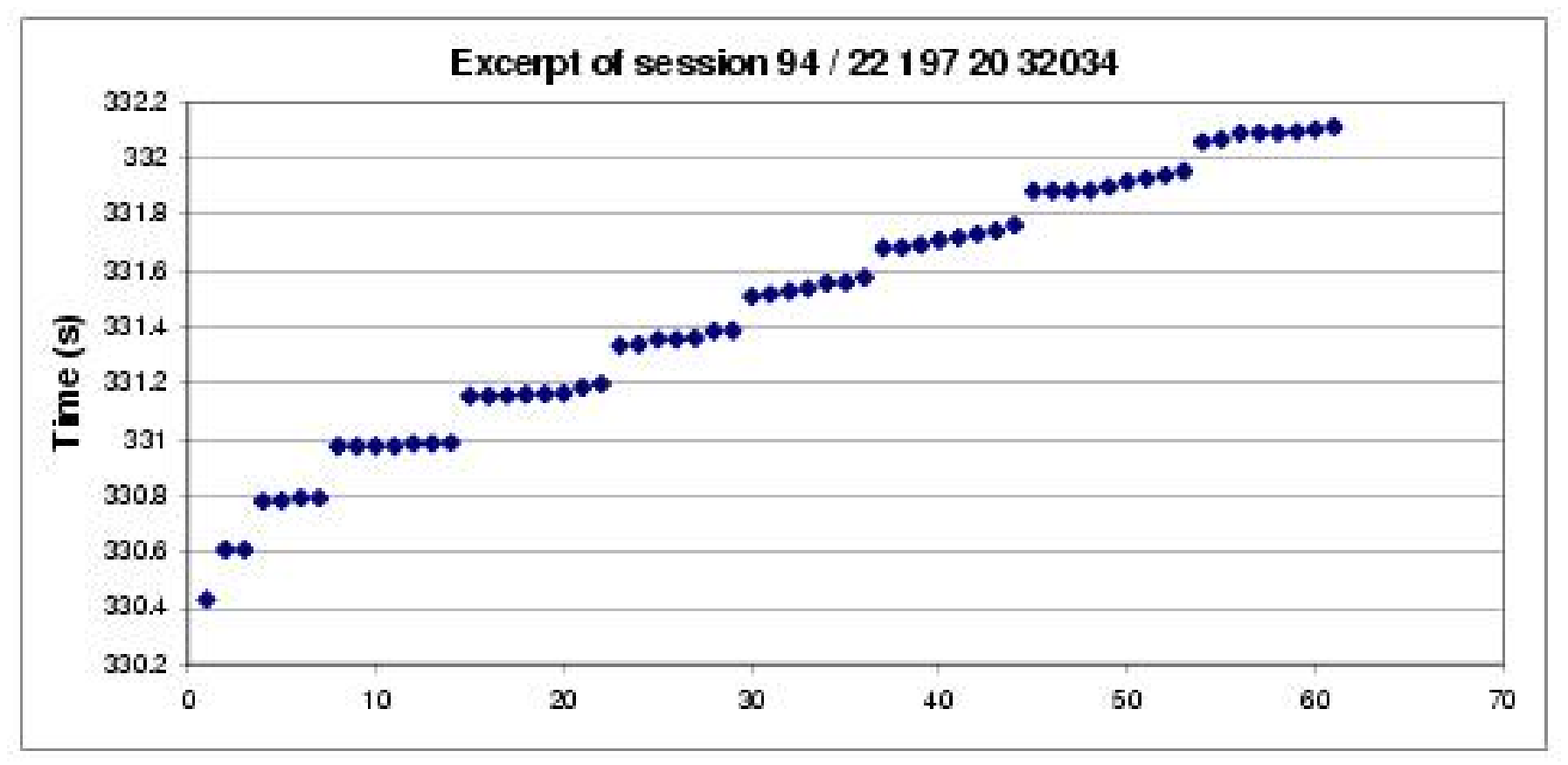}}
 \subfigure[]{\includegraphics[height=150pt, width=200pt]{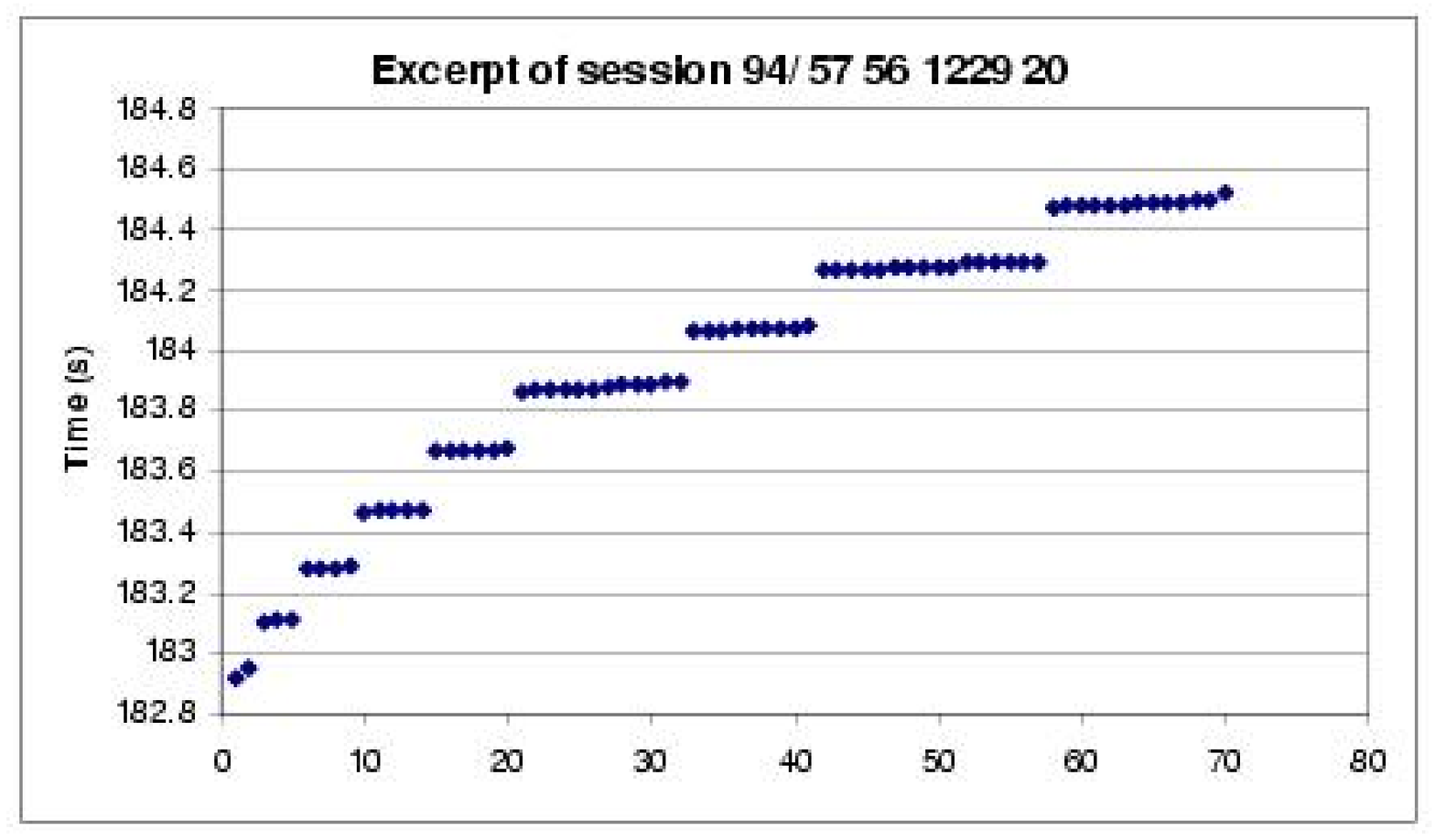}}
 \subfigure[]{\includegraphics[height=150pt, width=200pt]{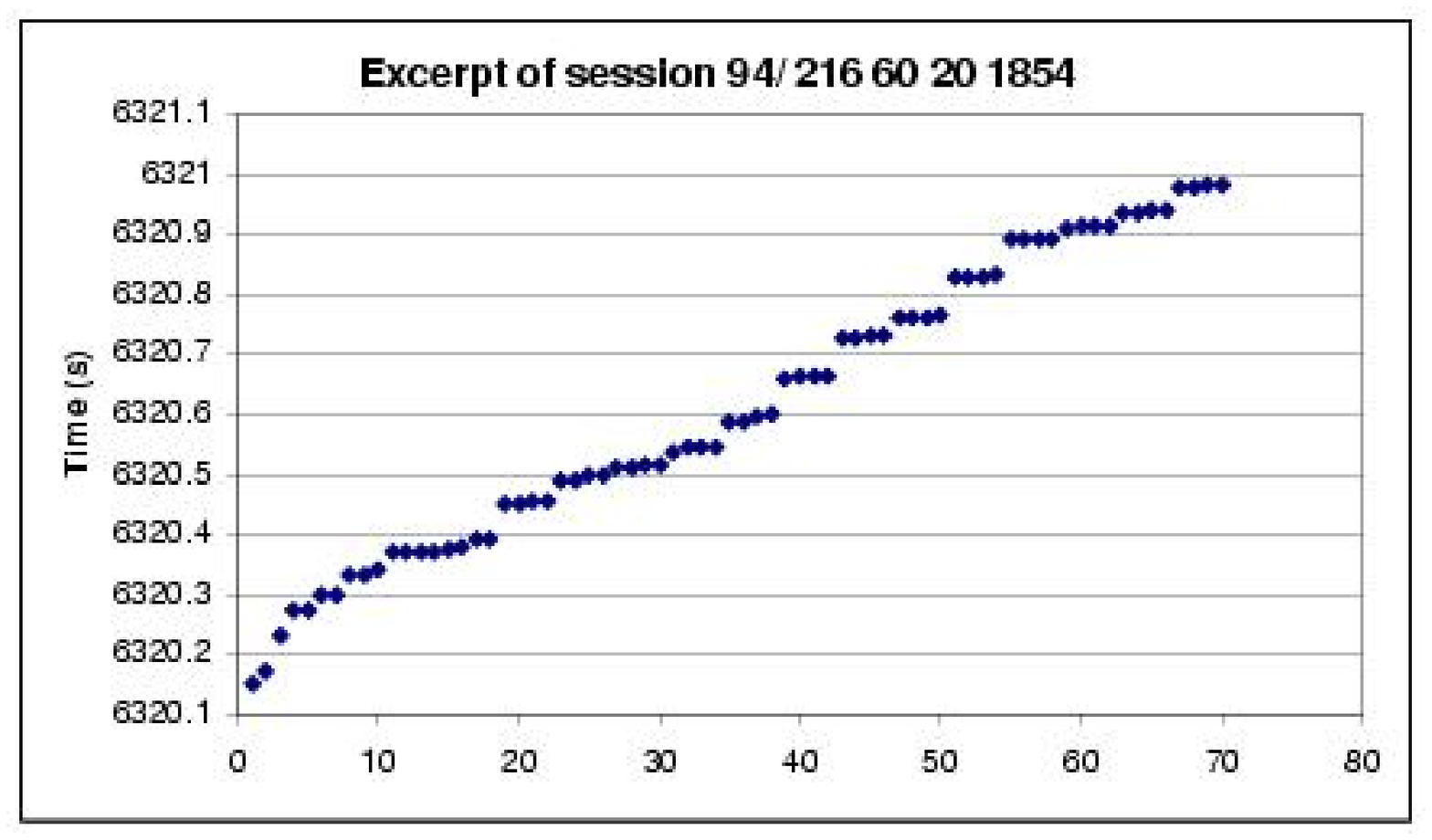}}
 \subfigure[]{\includegraphics[height=150pt, width=200pt]{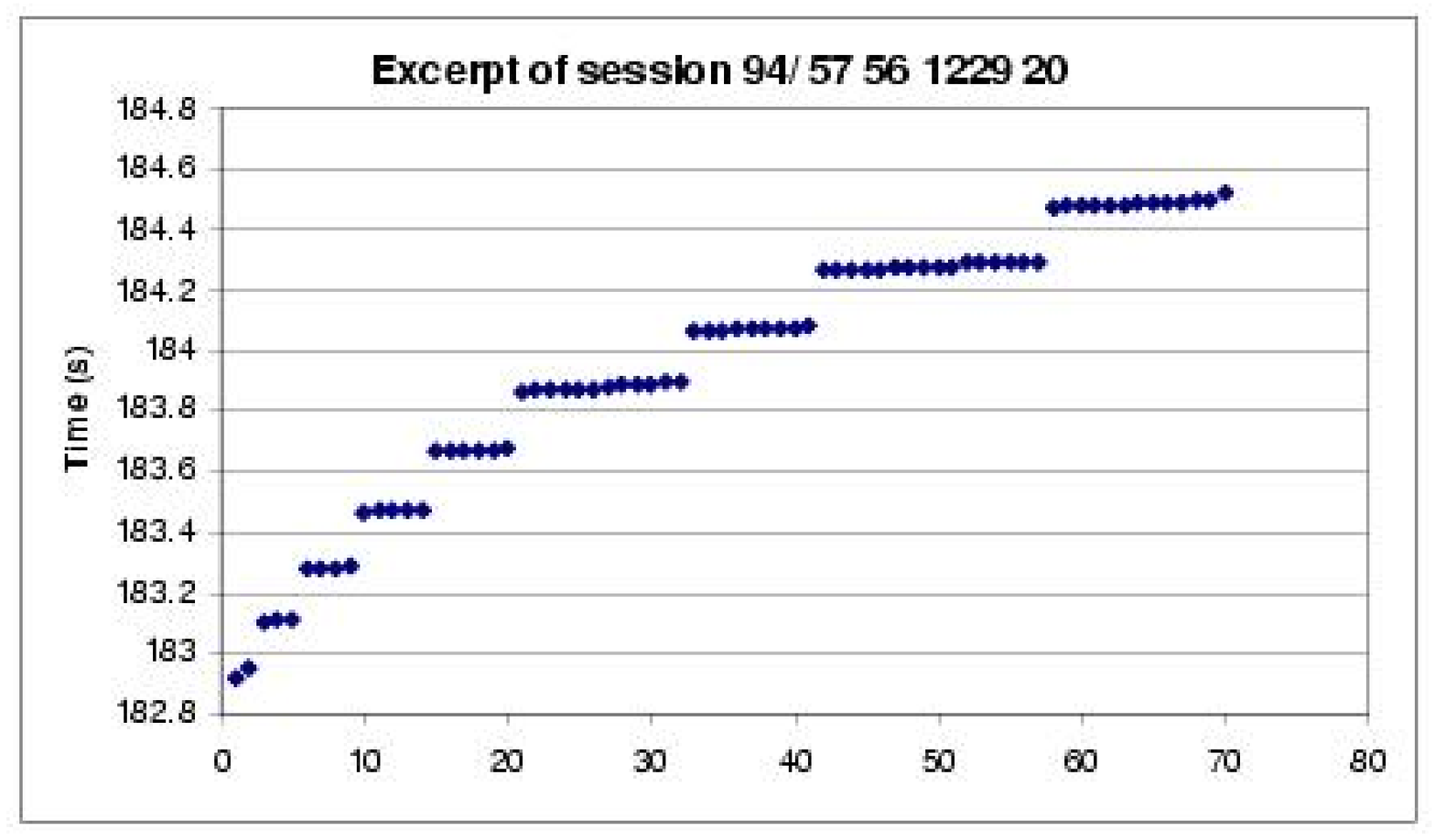}}
 \subfigure[]{\includegraphics[height=150pt, width=200pt]{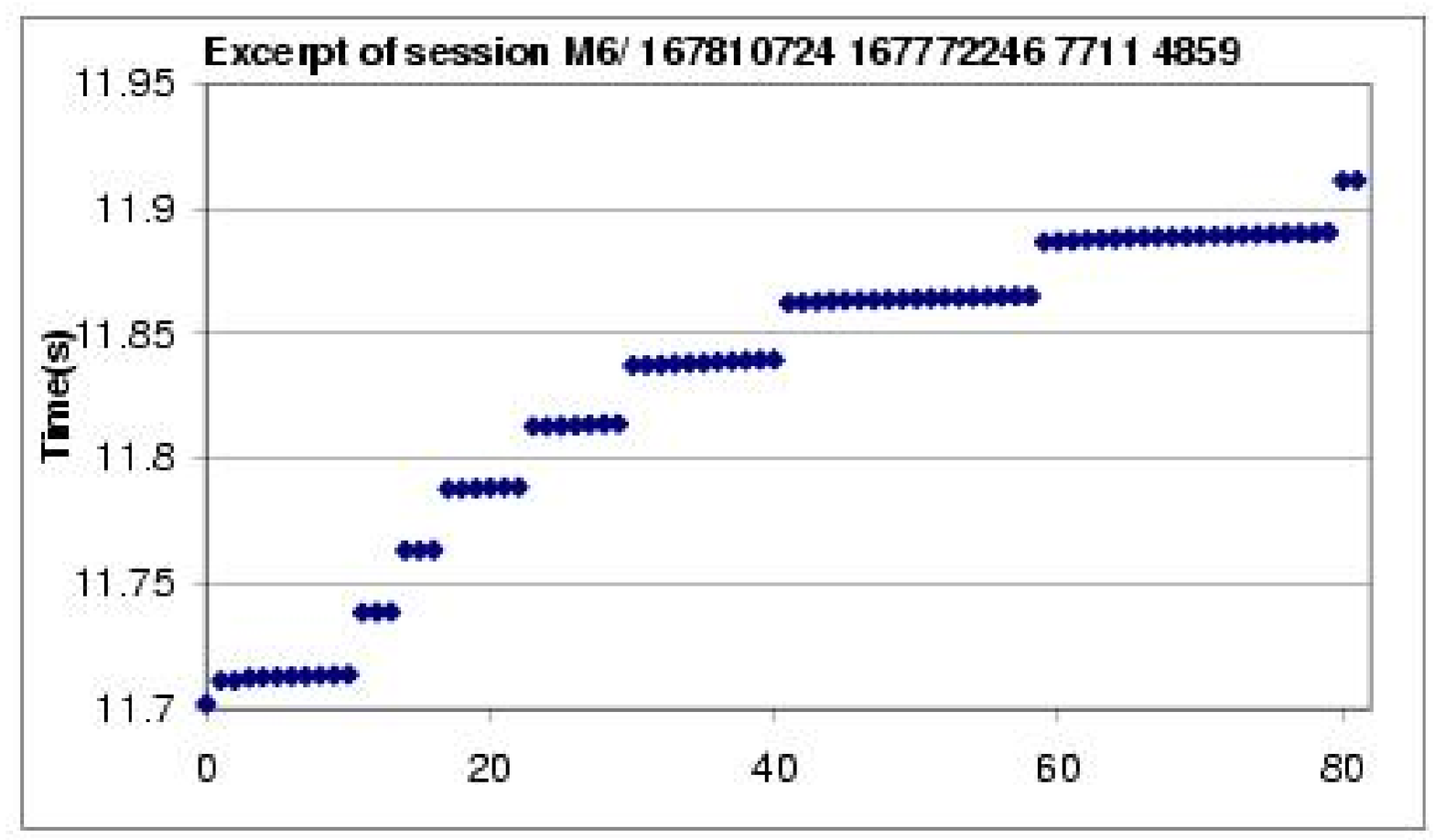}}
 \subfigure[]{\includegraphics[height=150pt, width=200pt]{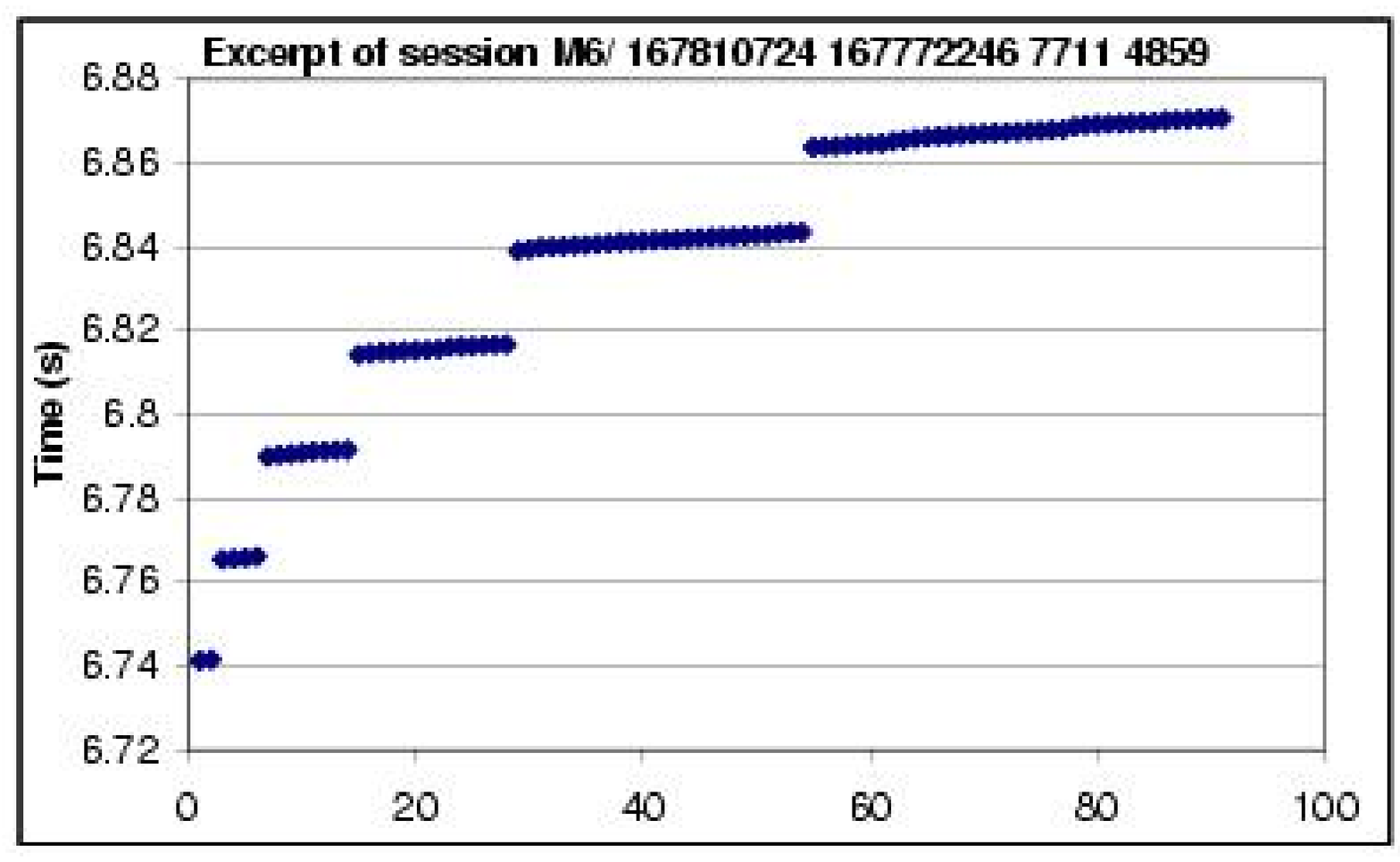}}
 \caption{Multiple packet emission according to TCP's Slow Start, as observed in 4 connections of trace 94 and 2 
          connections of trace M6 (remember here that a \emph{connection} is the subset of the traffic 
          corresponding to specific values of sender and receiver host and port numbers): in all cases except (c) Slow Start behaves \emph{exactly} according to the TCP 
          specifications! Slow Start may not always be as ``pure'' as in figures (a), (b), and (d), though; one also sees cases like (c), in which the increase of packet number is present but does not fit the protocol specifications as nicely. Nevertheless, our experimental observations correspond to our theoretical predictions.}
 \label{SSSess}         
\end{figure}

Model B oversimplifies the data transfer between two computers: although it
addresses the uncertainty in the network and the quantized nature of the
transfer, it fails to take into account protocol specifications and
connection speed variations. Indeed, the most widely adopted transport layer
protocol, TCP, has the Slow Start feature, according to which the number
of packets it sends each time increases geometrically by 2. This geometric
increase is readily observed in real traffic sessions (Fig. \ref{SSSess}). Also,
connection speeds can vary from 56K modems to extremely fast backbone links.

Therefore, within the framework of model B, a modification is made, so that
the number of packets sent each time increases by 2, until either all of the $L$
packets are sent, or a specified maximum $Max$ is reached, after which it remains
the same. Thus, the sequence of the number of packets sent, instead of $%
1,1,1,...,1$, is now going to be $1,2,4,8,...,Max,Max,...$, unless $L$ is
less than $1+2+4+...+Max$. In other words, the delta functions will no longer have equal weights, 
but, instead, they will be multiplied by geometrically increasing coefficients. The different connection speeds can be
modeled through different ``packet sizes'', so that the sequence of packets
sent from a particular user now becomes $\lambda ,2\lambda ,4\lambda
,...,\lambda Max,...$ \cite{SRB1}. Finally, packet losses can by
incorporated in the model by allowing $Max$ to be a random variable.

\addtocounter{figure}{-1}
\stepcounter{figure}

\subsection{$p$-Stable variables and processes}

\label{pSDef}

Theorem C will make use of \emph{$p$-Stable variables}. The purpose of this paragraph is to familiarize the reader with this quite ``exotic'' type of random variables, and the random processes they generate. More information can be found in \cite{R2}. 

A distribution is defined to be \emph{stable} iff it remains invariant ---within an affine transformation--- under convolutions. More precisely, choose $n$ and let $X_1,\ldots, X_n$ follow this distribution; it will be stable iff $\exists a_n, b_n:\ X_1+\ldots+X_n\overset{\mathcal{L}}{=}a_nX+b_n$. A $p$-Stable variable $V$ is represented through its characteristic function, as no closed form for its distribution function is generally known:

\[ \bold{E}[\exp(i\theta V)] = \exp(-\sigma^p|\theta|^p+i\mu\theta),\ 0<p\leq2\]

It is customary to write that $V\sim Sp S(\sigma)$. Closed forms for the distribution are known in the cases where $p=1$ (Cauchy) and $p=2$ (Normal). The tails of the distributions, however, are known to follow power laws:

\[\underset{\lambda\rightarrow \infty}{\lim} \lambda^p \bold{P}(V>\lambda)=\underset{\lambda\rightarrow \infty}{\lim} \lambda^p \bold{P}(V<-\lambda)= C_p \sigma^p\]   

A self-similar $p$-Stable random process with stationary increments and self-similarity factor $H$ is usually denoted by $L_{H,p}$ (see the definition of self-similarity in section \ref{modelA}). The stationary increments requirement, combined with self-similarity, imposes that:
\[L_{H,p}(t)-L_{H,p}(s)\sim Sp S(|t-s|^H\sigma_{V(1)})\]

The increments will also be independent iff $H=1/p$, in which case the process is called a \emph{L\'evy Stable Motion}.

An important difference between $p$-Stable and Gaussian processes is that, in general, the paths of a $p$-Stable process with $p<2$ are discontinuous with probability 1. Notice also the larger number of parameters needed to specify a $p$-Stable variable completely, compared to a Gaussian variable (mean, variance, and $p$ vs. mean and variance): it will become clear after the statement and proof of Theorem C, and the discussion afterwards, that estimating $p$ from a given data set is by no means easy, but, fortunately, for this paper's purposes, its value will not be needed. 

\subsection{Statement and proof of Theorem C}

Before the formal statement of Theorem C, a brief description of the
notation is in order: the stream of data corresponding to the $i$-th user,
with $Max=M$, will be $W_{i,M}^{c}(s)$. $T_{n}^{c}(t)$ will denote the total aggregated traffic until
time $t$. The function $C(s,t)$ is the covariance structure of the Gaussian
limit process, as given in (\ref{CovStr}). Finally, $p\/$ is the exponent in (\ref{PP}).

\begin{definition}
\upshape
\ Let $\{O_{f}^{c}(i)\}_{i=1}^{\infty
},\{R^{c}_{i,j}\}_{i,j=0}^{\infty }$ and $\{L^{c}(i)\}_{i=1}^{\infty }$ be the
same as $\{O_{f}^{b}(i)\}_{i=1}^{\infty },\{R^{b}_{i,j}\}_{i,j=0}^{\infty }$
and $\{L^{b}(i)\}_{i=1}^{\infty }$.

Given $M>0$, let $\displaystyle \phi (l)=\left\lfloor \log _{2}(l\wedge M)+\frac{%
(l-2M+1)_{+}}{M}\right\rfloor +1$, where $\left\lfloor X\right\rfloor $ stands
for the largest integer which is smaller than $X$ and $(X)_{+}$ stands for the positive
part of $X$, as usual. Let $\displaystyle O_{n,M}^{c}(i)=\sum_{j=1}^{\phi
(L^{c}(i))}R^{c}_{i,j}$. Recall also the definition of the ``delta'' function $b\widehat{\delta }%
_{x}(s)$ and its properties: $b\widehat{\delta }_{x}(s)=0$ if $x\neq s$ and $\displaystyle
\int_{x-h}^{x+h}b\widehat{\delta }_{x}(s)ds=b$ for all $h>0.$ With this convention, define:
 
\begin{equation*}
\widetilde{W}_{M}^{c}(t)=\left\{ 
\begin{array}{cc}
(2^{l}\wedge M)\widehat{\delta }_{t}(t) & \text{if }t=%
\displaystyle\sum_{i=1}^{k-1}\left(O_{f}^{c}(i)+O_{n}^{c}(i)\right)+O_{f}^{c}(k)+%
\displaystyle\sum_{i=1}^{l}R^{c}_{k,i},\text{ }l\leq \phi (L^{c}(k)) \\ 
0 & \text{Otherwise}
\end{array}
\right.
\end{equation*}

As before, define further $W_{M}^{c}(t)$ to be a stationary version of $\widetilde{%
W}_{M}^{c}(t).$ Finally, let $\{W_{i,M}^{c}(t)\}$ be i.i.d.\ versions of $%
W_{M}^{c}$, define a constant $\bold{E}W_{M}^{c}(x)$ such that
 
\begin{equation*}
\bold{E}\left(\int_{a}^{b}W_{M}^{c}(x)dx\right)=\bold{E}W_{M}^{c}|b-a|,
\end{equation*}

\noindent and let 

\begin{eqnarray}
\label{CovStr}
C(s,t)=\lim_{n\rightarrow \infty }\frac{1}{\sigma _{n}} \left(
\bold{E}\left(\int_{0}^{s}W_{1,M_{n}}^{c}(x)dx\int_{0}^{t}W_{1,M_{n}}^{c}(x)dx\right)- \bold{E}\left(%
\int_{0}^{s}W_{1,M_{n}}^{c}(x)dx\right)\bold{E}\left(\int_{0}^{t}W_{1,M_{n}}^{c}(x)dx\right)\right)
\end{eqnarray}
where $\displaystyle \sigma _{n}^{2}=\bold{Var}\left(\int_{0}^{T}W_{i,M_{n}}(s)ds\right)$.
Recall also that $\exists p\in (1,2):\ t^{p}\bold{P}\left(L^{c}(i)>t\right)%
=\Theta (1)$, as $t\rightarrow \infty $, which is another way to say that the LHS of the expression is bounded
above and below by a constant. 
\end{definition}

\bigskip
\begin{theorem}
Let $\displaystyle T_{n}^{c}(t)=\int_{0}^{t}\sum_{i=1}^{n}W_{i,M}^{c}(s)ds$. Then

\begin{list}{\roman{temp})}{\usecounter{temp} \setlength{\rightmargin}{\leftmargin}}
\item If $\ M=M_{n}\nearrow \infty $ and $M_{n}\leq n^{q}$
for $q<1/p$, then, for fixed $t$: 
\begin{equation*}
Z_{n}^{c}(t)=\frac{1}{\sigma _{n} \sqrt{n}}\left(T_{n}^{c}(t)-ET_{n}^{c}(t)\right)\overset{W}{%
\rightarrow }G^{c}(t),\ t\in [0,T]
\end{equation*}
where $\displaystyle \sigma _{n}^{2}=\bold{Var}\left(\int_{0}^{T}W_{i,M_{n}}(s)ds\right)$
and $G^{c}(t)$ is a centered Gaussian process with covariance structure $C(s,t)$.

\item If $M_{n}\geq n^{q}$ for $q>1/p$, then, for fixed $t$: 
\begin{equation*}
Z_{n}^{c}(t)=\frac{1}{n^{1/p}} \left[T_{n}^{c}(t)-\bold{E}\left(T_{n}^{c}(t)\right)\right]\overset{W}{%
\rightarrow }P^{c}(t),\ t\in [0,T]
\end{equation*}
where $P^{c}(t)$ is a p-Stable process.
\end{list}

\end{theorem}

Note: we will not pursue further the issues of what the exact value of $p$ for $P^c(t)$ is and how it can be estimated from a given data set (see the relevant remarks in section \ref{pSDef}).

\begin{proof}

First of all, recall the equality
\begin{equation}
\bold{E}\left(\left|X\right|\right)=\int_{0}^{\infty}P(X>t)dt
\label{MV}
\end{equation}
which will prove repeatedly useful later.  

For part \textit{i)}, finite dimensional
convergence has to be established first. An easy computation yields 
\begin{equation*}
\frac{\sigma _{t,n}^{2}}{\sigma _{n}^{2}}=\frac{\bold{Var}\left(%
\int_{0}^{t}W_{i,n}^{c}(s)ds\right)}{\bold{Var}\left(\int_{0}^{T}W_{i,n}^{c}(s)ds\right)}\underset{%
n\rightarrow \infty }{\longrightarrow }\gamma _{t}^{2}\text{ for }0<\gamma
_{t}^{2}\leq 1.
\end{equation*}
Using the notation $\displaystyle W_{1,n}^{c}=W_{i,M_{n}}^{c},\ X_{i,n}^{t}=%
\int_{0}^{t}W_{i,n}^{c}(x)dx,$ and $\displaystyle \widetilde{X}_{i,n}=\frac{1}{\sigma
_{t,n}}\left(X_{i,n}^{t}-\bold{E}\left(X_{i,n}^{t}\right)\right)$, it is sufficient to show that
\begin{equation*}
\frac{1}{\sqrt{n}}\sum_{i=1}^{n}\widetilde{X}_{i,n}=\frac{1}{\sigma _{t,n}%
\sqrt{n}}\sum_{i=1}^{n}\left[X_{i,n}^{t}-\bold{E}\left(X_{i,n}^{t}\right)\right]\overset{W}{\longrightarrow 
}N(0,1)
\end{equation*}
Since the $\widetilde{X}_{i,n}$ are obviously i.i.d., centered and with $%
\bold{Var}\left(\widetilde{X}_{i,n}\right)=1$, the above follows easily (using a classical
Lindeberg argument, see \cite{Pollard}), provided that 
\begin{equation*}
n^{-1/2}\bold{E}\left(\left|\widetilde{X}_{i,n}\right|^{3}\right)\rightarrow 0.
\end{equation*}
Clearly $\bold{E}\left(X_{i,n}^{t}\right)\rightarrow C<\infty $, hence the above is bounded (up to a
constant) by 
\begin{multline*}
\frac{\bold{E}\left(\left|X_{i,n}^{t}\right|^{3}\right)}{\sqrt{n}\left(\bold{E}\left(\left|X_{i,n}^{t}\right|^{2}\right)\right)^{3/2}}\leq C\frac{%
M_{n}\bold{E}\left(\left|X_{i,n}^{t}\right|^{2}\right)}{\sqrt{n}\left(\bold{E}\left(\left|X_{i,n}^{t}\right|^{2}\right)\right)^{3/2}}= C\frac{M_{n}}
{\sqrt{n \bold{E}\left(\left|X_{i,n}^{t}\right|^{2}\right)}}
\overset{*}{\leq} C'\frac{M_{n}^{p/2}}{\sqrt{n}}\leq C'n^{\frac{pq-1}{2}}=C'n^{\frac{p}{2}(q-\frac{1}{p})}
\overset{n\rightarrow\infty} {\longrightarrow} 0
\end{multline*}
where in ($*$) we used \ref{MV} and (\ref{PP}).

The stochastic equicontinuity follows after a computation very similar to the one
presented in the proof of Theorem B.

Namely, one needs to estimate 
\begin{multline*}
\bold{E}\left(Z_{n}^{c}(t)-Z_{n}^{c}(s)\right)^{4}=\bold{E}\left( \frac{1}{\sigma _{n}\sqrt{n}}%
\sum_{i=1}^{n}\left(
(X_{i,n}^{t}-X_{i,n}^{s})-\bold{E}\left(X_{i,n}^{t}-X_{i,n}^{s}\right)\right) \right) ^{4} \\
\leq \frac{1}{n\sigma _{n}^{4}}\bold{E}\left(
(X_{1,n}^{t}-X_{1,n}^{s})-\bold{E}\left(X_{1,n}^{t}-X_{1,n}^{s}\right)\right) ^{4}+\frac{6}{%
\sigma _{n}^{4}}\left( \bold{E}\left(
(X_{1,n}^{t}-X_{1,n}^{s})-\bold{E}\left(X_{1,n}^{t}-X_{1,n}^{s}\right)\right) ^{2}\right) ^{2} \\
\leq \frac{1}{n\sigma _{n}^{4}}\left(
\bold{E}\left(X_{1,n}^{t}-X_{1,n}^{s}\right)^{4}+6\bold{E}\left(X_{1,n}^{t}-X_{1,n}^{s}\right)^{2}\left(
\bold{E}\left(X_{1,n}^{t}-X_{1,n}^{s}\right)\right) ^{2}\right) +\frac{6}{\sigma _{n}^{4}}%
\left( \bold{E}\left(X_{1,n}^{t}-X_{1,n}^{s}\right)^{2}\right) ^{2}
\end{multline*}

It can be shown, using (\ref{MV}) that, for $\delta$ small enough and $\varepsilon=1/p-q$, $\exists C>0$: 
\begin{equation*}
\frac{\bold{E}\left(X_{1,n}^{t}-X_{1,n}^{s}\right)^{2}}{\sigma _{n}^{2}}\leq C|t-s|
\end{equation*}
and 
\begin{equation*}
\frac{\bold{E}\left(X_{1,n}^{t}-X_{1,n}^{s}\right)^{4}}{n\sigma _{n}^{4}}\leq Cn^{-\varepsilon}|t-s|
\end{equation*}

Therefore, for $|t-s|>n^{-1/2}$, $\exists C_1>0$: 
\begin{equation}
\left\| Z_{n}^{c}(t)-Z_{n}^{c}(s)\right\| _{L_{4}}\leq C_{1}|t-s|^{1/4+\varepsilon/2}
\label{L4C}
\end{equation}
As in the proof for Theorem B, this allows for the replacement of the standard requirement
for stochastic equicontinuity (\ref{StochEqui}) with: 
\begin{equation*}
\underset{\delta \rightarrow 0}{\lim }\underset{n\rightarrow \infty }{\lim
\sup }\bold{P}\left(\underset{\gamma n^{-1/2}\leq |t-s|\leq \delta }{\sup }%
\left|Z_{n}^{c}(t)-Z_{n}^{c}(s)\right|>\varepsilon \right)=0
\end{equation*}
In order to prove the above, one can proceed as in the proof of Theorem B:
construct the partition $\displaystyle \Lambda _{n}=\{t_{i}=\frac{i}{K_{\varepsilon }\sqrt{%
n}}:\ i=0,...,K_{\varepsilon }\sqrt{n}T\}$ where $K_{\varepsilon
}=\varepsilon ^{-1}8\bold{E}W_{M}^{c}$, and, for every $t$, define functions $%
\phi _{1}^{n}$ and $\phi _{2}^{n}$: $[0,T]\rightarrow \Lambda _{n}$ such
that $\phi _{1}^{n}(t)=t_{j}\leq t\leq t_{j+1}=\phi _{2}^{n}(t)$ for some $%
j\leq K_{\varepsilon }\sqrt{n}T$.

Clearly 
\begin{equation*}
T_{n}^{c}(\phi _{1}^{n}(t))\leq T_{n}^{c}(t)\leq T_{n}^{c}(\phi _{2}^{n}(t))
\end{equation*}
which implies (by essentially the same computation as the one presented in (\ref
{Bcomput})) 
\begin{equation*}
Z_{n}^{c}(t)-Z_{n}^{c}(s)\leq Z_{n}^{c}(\phi _{2}^{n}(t))-Z_{n}^{c}(\phi
_{1}^{n}(s))+\varepsilon /4
\end{equation*}
 Observe now that $\displaystyle \bold{E}\left(T_{n}^{c}(\phi _{2}^{n}(t))\right)-\bold{E}\left(T_{n}^{c}(\phi
_{1}^{n}(t))\right)=\bold{E}\left(\frac{1}{\sigma _{n}\sqrt{n}}\sum_{i=1}^{n}\int_{\phi
_{1}^{n}(t)}^{\phi _{2}^{n}(t)}W_{1,n}^{c}(x)dx\right)=\frac{\sqrt{n}}{\sigma _{n}}%
(\phi _{2}^{n}(t)-\phi _{1}^{n}(t))\bold{E}W_{M_n}^{c}\leq \bold{E}W_{M_n}^{c}/K_{\varepsilon
}=\varepsilon /8$.

This proves part \textit{i)} of Theorem C.

\bigskip

Proof for \textit{ii).} Let us begin by giving a summary of some relevant convergence results found in \cite{Gine}, section 2.6 (see also section \ref{pSDef} above): 

\stepcounter{alemma}
\stepcounter{alemma}

\begin{alemma}
Let $\{X_i\}$ be a sequence of i.i.d. random variables, and suppose there exist two sequences $\{a_n\}$ and $\{b_n\}$ such that $\displaystyle \frac{X_1+\ldots+X_n}{a_n}-b_n\overset{w}{\longrightarrow} Z$ (in which case we say that $X_i$ is \emph{in the domain of attraction} of $Z$). Then:
\begin{itemize}
	\item $Z$ can only be a \emph{$p$-Stable} r.v., where $0<p\leq 2$. 
	\item $\exists C>0: a_n=Cn^\frac{1}{p}$.
	\item $p=2\Leftrightarrow \bold{E}(X_1^2)<\infty$, or, equivalently, $p<2\Leftrightarrow \bold{E}(X_1^2)=\infty$.
	\item If $p=2$, $Z$ is Gaussian.
\end{itemize}
\end{alemma}

For fixed $t$, let $\displaystyle Y_{i}^{t}=\int_{0}^{t}W_{i,\infty }^{c}(x)dx$, where $W_{i,\infty }^{c}$ stands for $W_{i,M_{n}}^{c}$ with $M_{n}=\infty$. A straightforward calculation of $\mathcal{L}(Y_{i}^{t})$ shows that $\bold{E}[(Y_{i}^{t})^2]=\infty$; it follows then from Lemma D that $Y_{i}^{t}$ is in the domain of 
attraction of a $p$-Stable random variable, i.e. 
\begin{equation*}
n^{-1/p}\sum_{i=1}^{n}\left(Y_{i}^{t}-\bold{E}\left(Y_{i}^{t}\right)\right)\rightarrow p\text{-Stable r.v.}
\end{equation*}
It is also a classical result (see \cite{Gine}) that, if $Y_{i}^{t}$ is replaced with its truncated version $\displaystyle Y_{i,n}^{t}=Y_{i}^{t}\bold{1}_{|Y_{i}^{t}|<M_{n}}$, the above convergence still
holds. The result now follows since the tail behavior of $Y_{i,n}^{t}$ is the
same as for $\displaystyle X_{i,n}^{t}=\int_{0}^{t}W_{i,M_{n}}^{c}(x)dx$. 

\end{proof}

\subsection{An interesting prediction: oscillation between Gaussianity and p-Stability}

According to Theorem C, and contingent upon the relation between $n$ and the maximum window size $M_{n}$, 
the limit process can be either Gaussian or p-Stable. Notice, however,
that this condition does not depend on the length $T$ of the time interval
under consideration. Therefore, it can be conjectured that the result applies ``locally''.
This, in turn, suggests two possible mechanisms causing 
oscillations between Gaussianity and p-Stability.

On the one hand, although
formally $n$ users are active in $[0,T]$, their actual number may fluctuate,
if examined on smaller time intervals (Fig. \ref{Stat94}). Let then $N_{n}\left(
t\right) $ be the number of users, among the total $n$ users of the simulation, who are in
an ON-interval at time $t$. If $\left(N_{n}(t)\right)^{q}<M_{n}$, the local behavior around $t$
will tend to be p-Stable, if $\left(N_{n}( t)\right)^{q}>M_{n}$, it will
tend to be Gaussian. This oscillation is observed in practice
(Fig. \ref{KolDist}), and mere visual inspection reveals that the number of users and
the Gaussianity of the process are well correlated (Table \ref{CrossCorr}). The fluctuation of the number of active users in time can apparently be well modeled by model A. Finally, notice that the number of bins used to 
compute the marginals is \emph{not} responsible for this oscillation, as it is kept fixed throughout 
the computation. 

On the other hand, $M_{n}$ itself can change, too. This happens actually in real networks built on TCP/IP,
where, according to the protocol specifications, packet loss halves the
window size. In other words, according to Theorem C, the more congested the
network is, the more Gaussian the traffic should look like. This has indeed
been confirmed \cite{SRB1,CCLS1}.

As a consequence, if incorporation of packet losses into this model is
required, $M_{n}$ can be allowed to be different for each user, and randomly
selected, but constant.

\begin{figure}
 \centering
 \subfigure[]{\includegraphics[height=150pt, width=200pt]{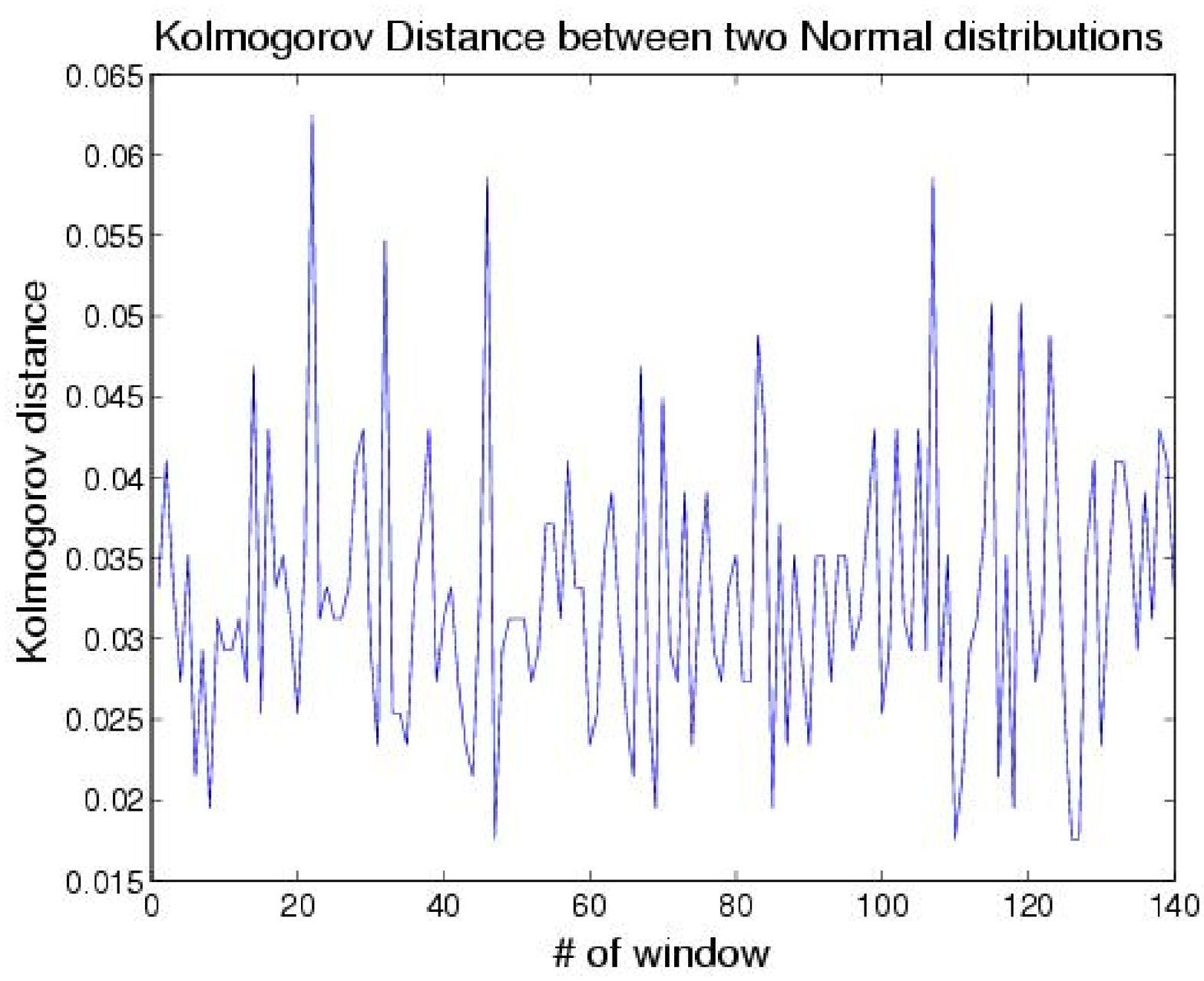}}
 \subfigure[]{\includegraphics[height=150pt, width=200pt]{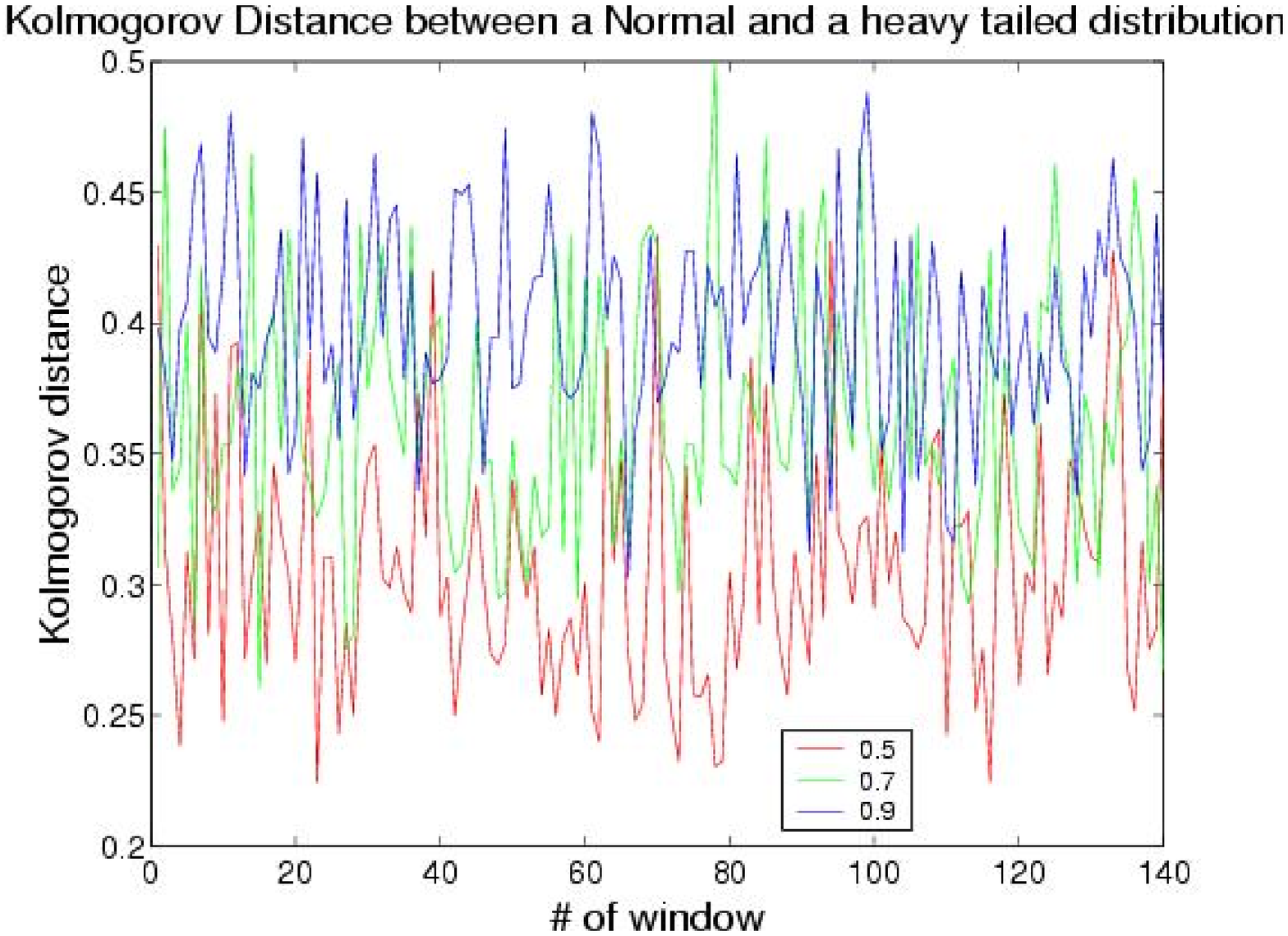}}
 \subfigure[]{\includegraphics[height=150pt, width=200pt]{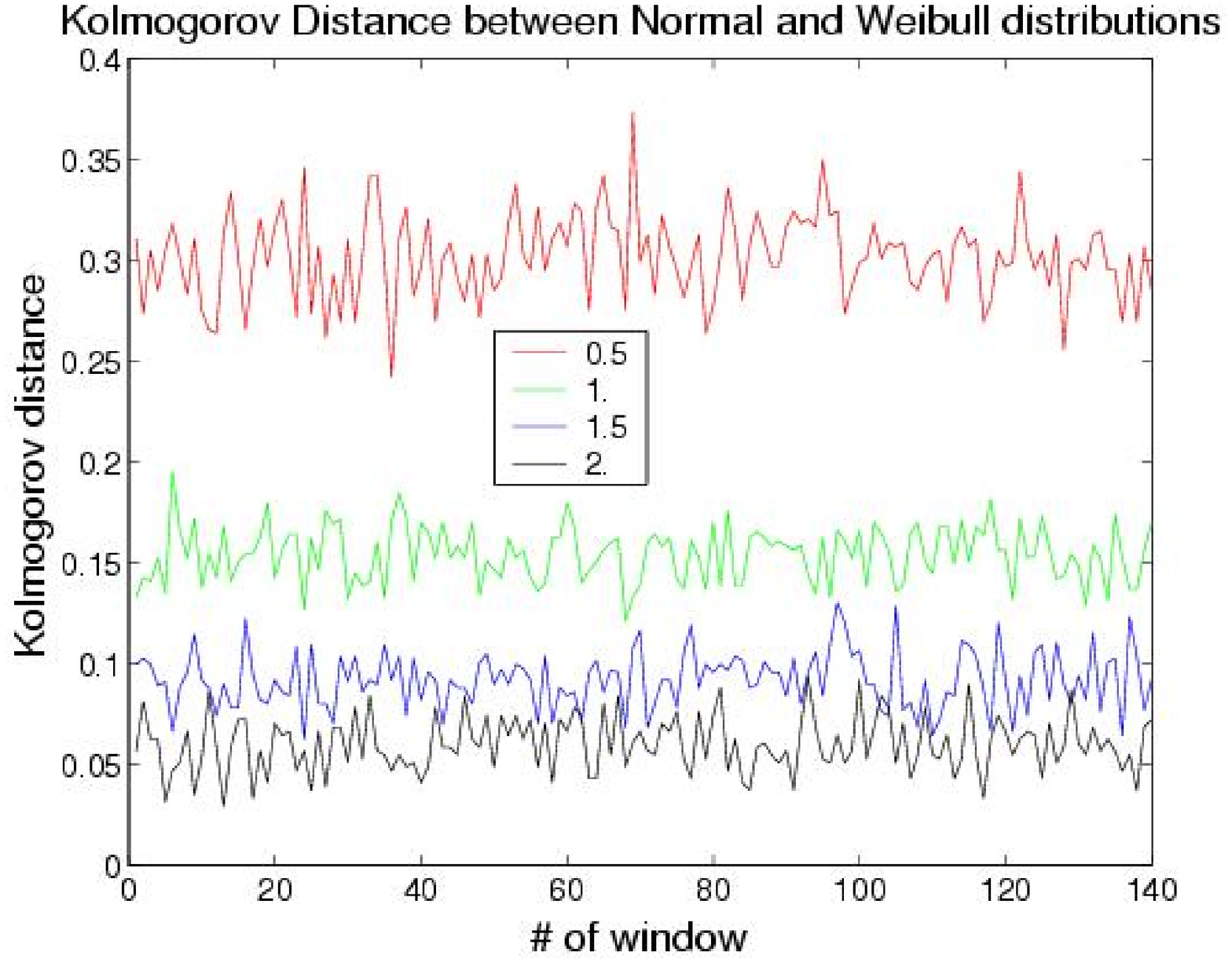}}
 \subfigure[]{\includegraphics[height=150pt, width=200pt]{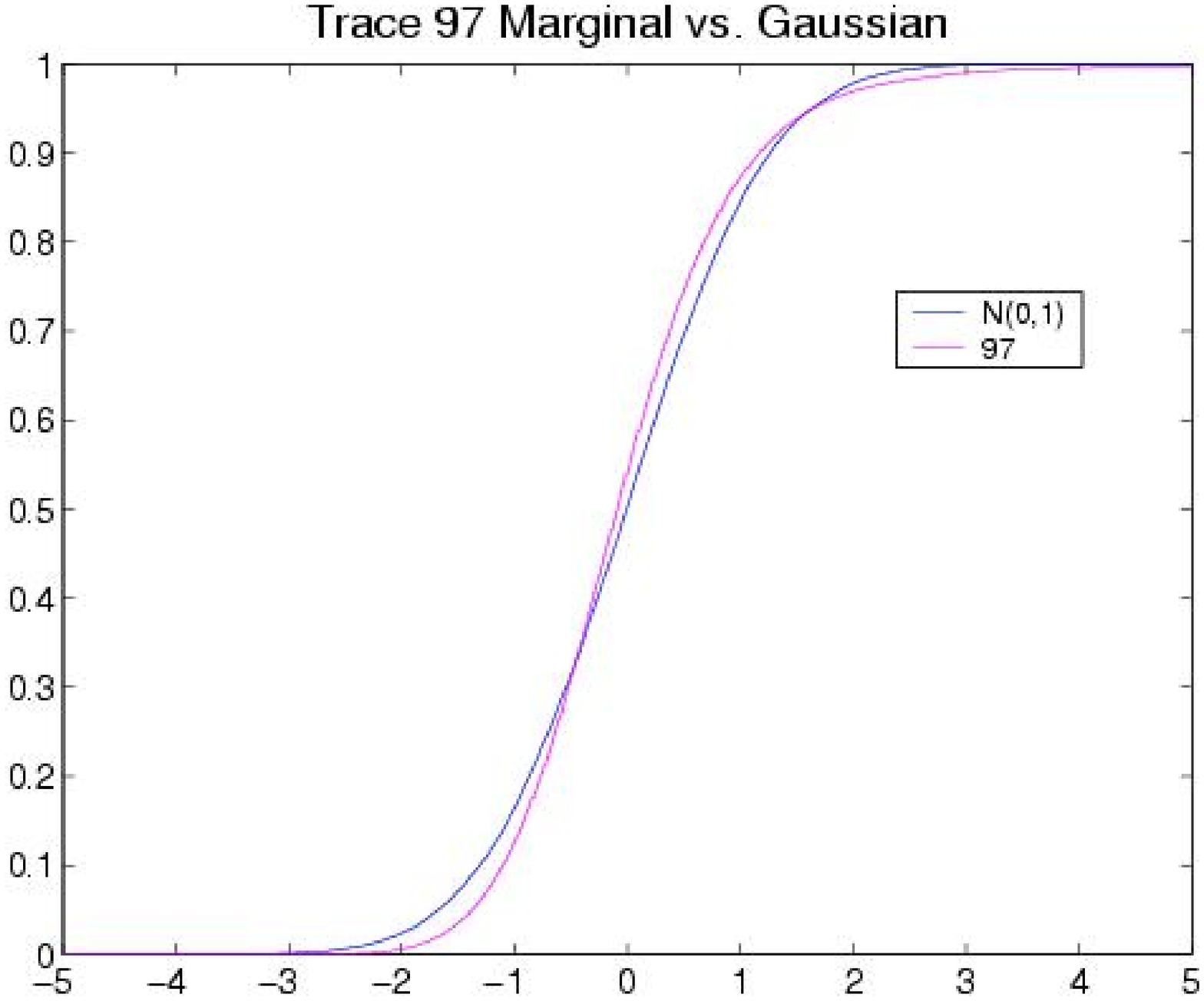}}
 \subfigure[]{\includegraphics[height=150pt, width=200pt]{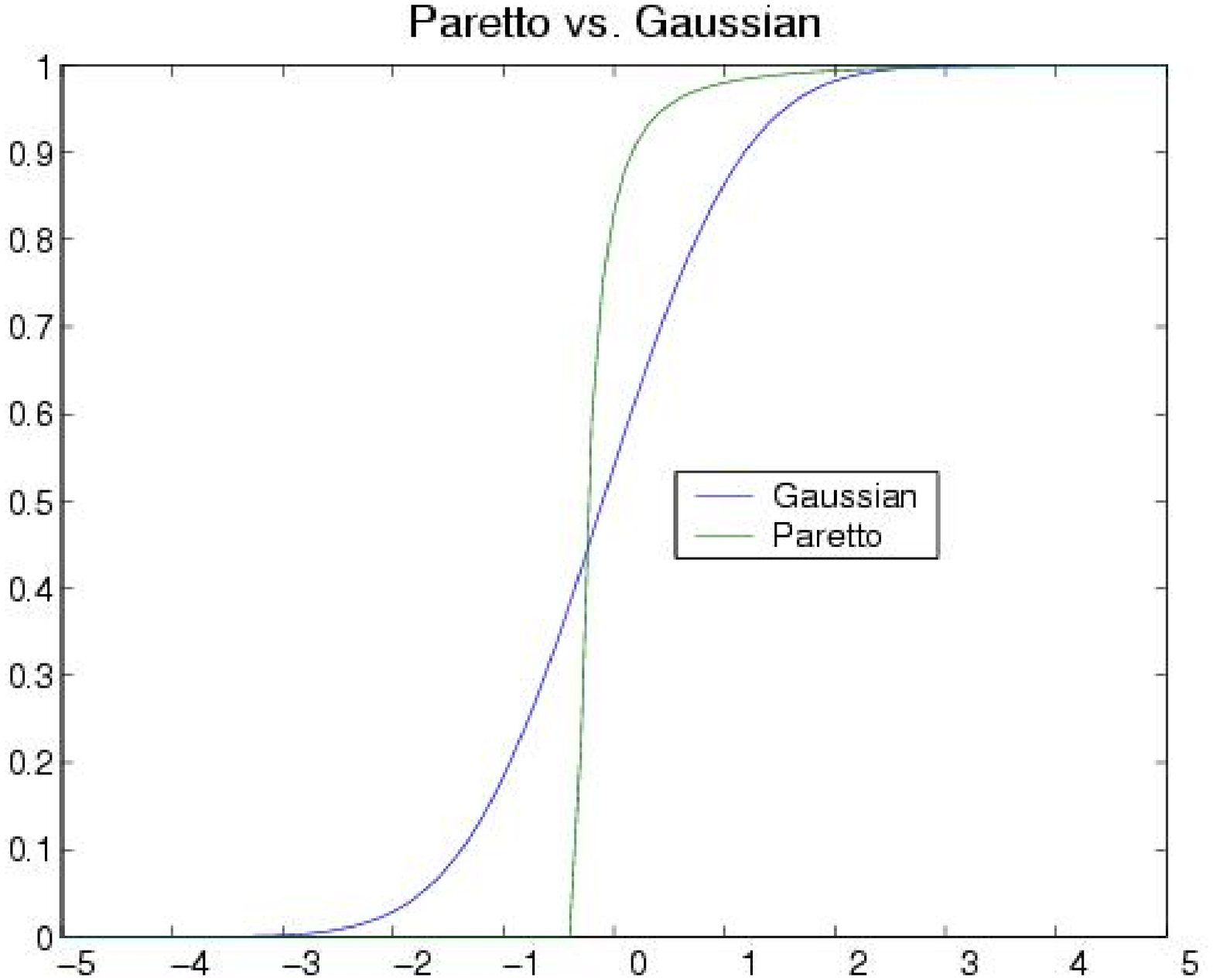}}
 \subfigure[]{\includegraphics[height=150pt, width=200pt]{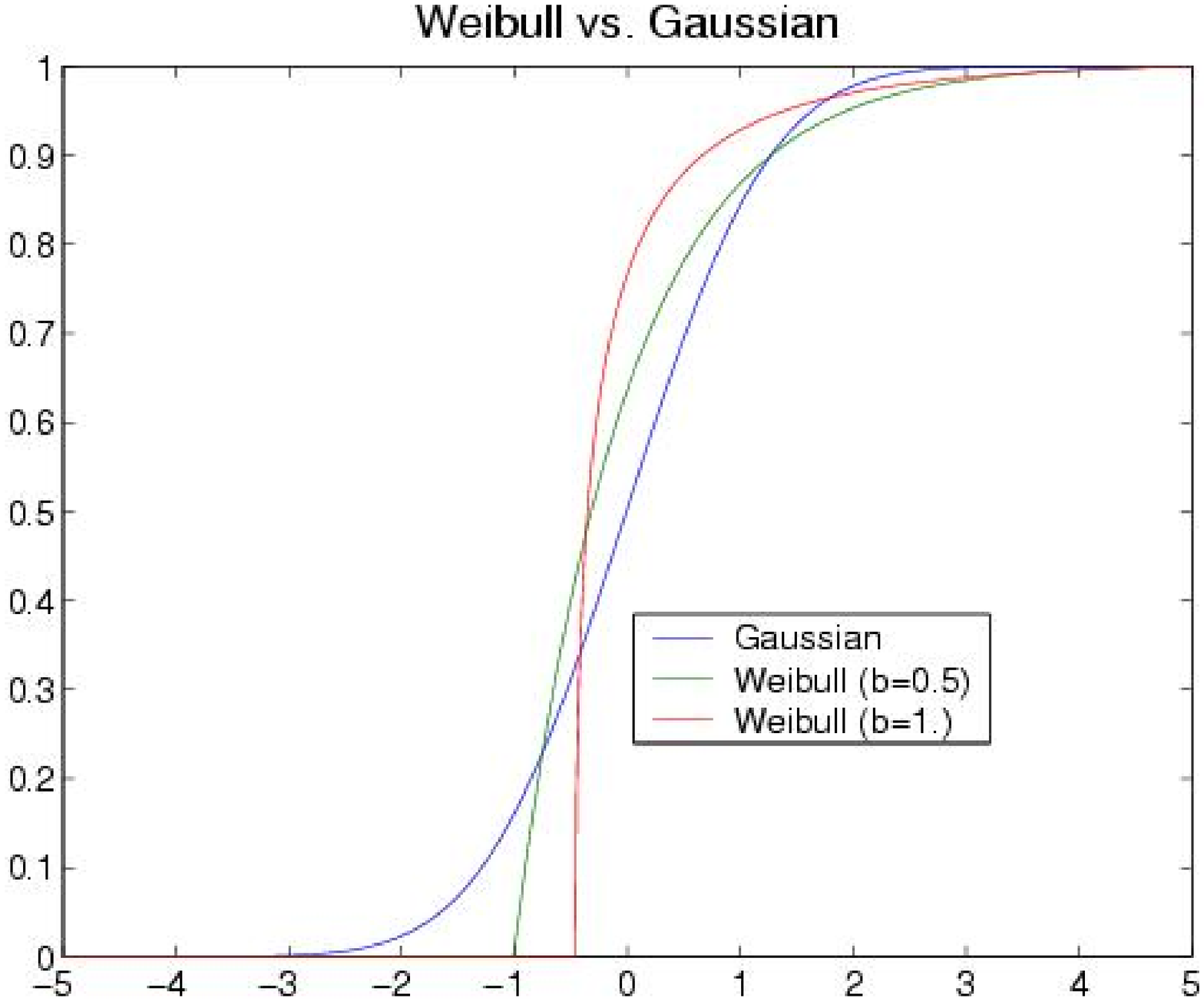}}
 \caption{Kolmogorov Distances of normalized (a) Gaussian, (b) heavy-tailed, and (c) Weibull
          marginals from the N(0,1) distribution. Marginals have been computed using non-overlapping consecutive
          windows of $2^{9}=512$ bins. The formula used for the Weibull is: $\displaystyle F(x)=1-e^{-ax^b}$. The 
          heavy tailed distribution corresponds to the pure power law (Paretto) $\displaystyle Y=\frac{a}{X^b},\ X\sim 
          U(0,1)$. In both cases, $a$ does not affect the results, due to the normalization, whereas the values for 
          $b$ used are shown in the legends. (d) shows that the trace 97 marginal is very close to a Gaussian, 
          although its tails are Weibull, as seen in Section \ref{problem}. Finally, (e) and (f) show the distance 
          between a Gaussian and a Paretto (with $b=0.8$), and between a Gaussian and a Weibull (one with $b=1.$, and 
          one with $b=0.5$), respectively. Notice the similarity between the Weibull with $b=0.5$ and the Paretto.}
 \label{KolDistTest}         
\end{figure}
\addtocounter{figure}{-1}
\stepcounter{figure}

\begin{figure}[t]
 \centering
 \subfigure[]{\includegraphics[height=150pt, width=200pt]{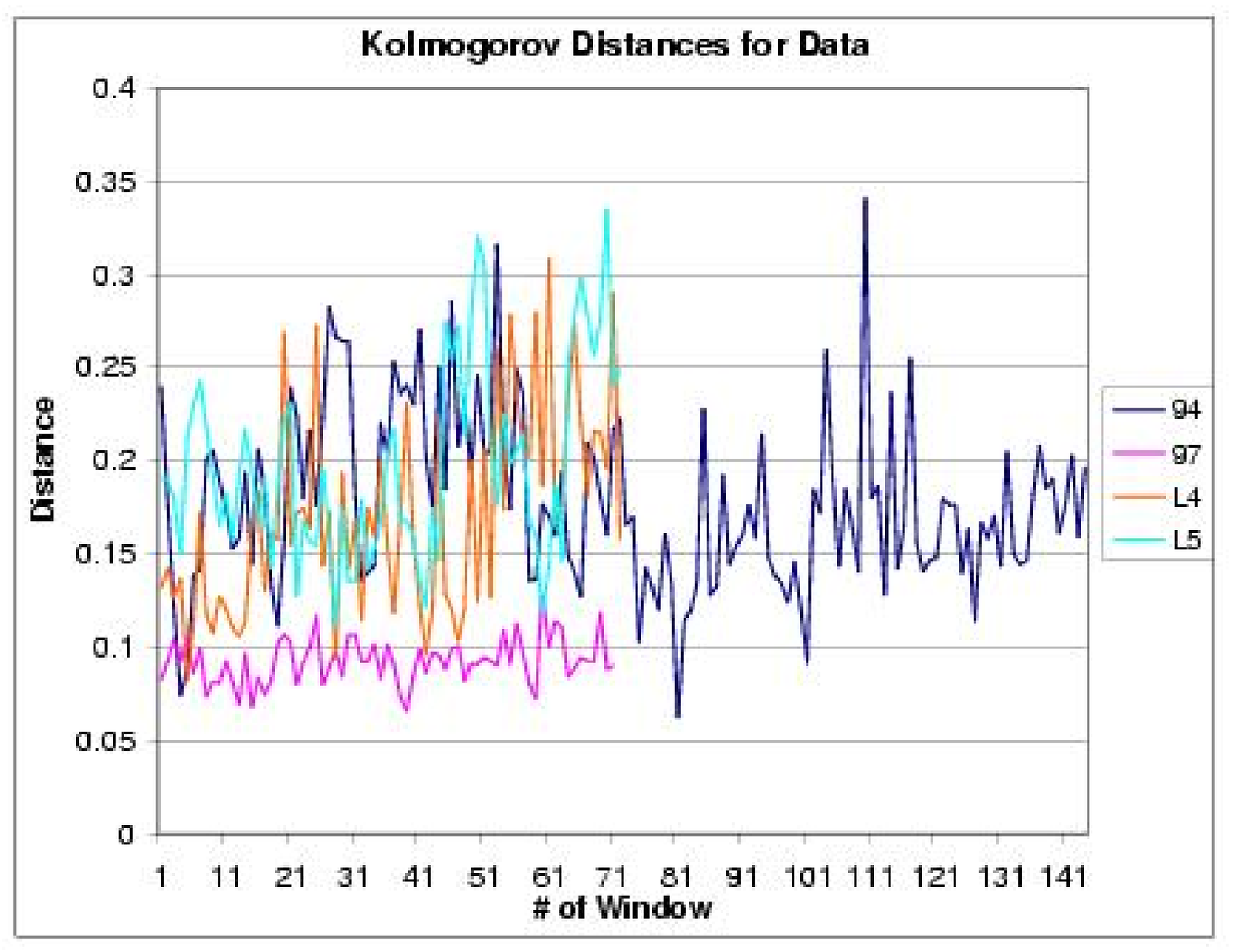}}
 \subfigure[]{\includegraphics[height=150pt, width=200pt]{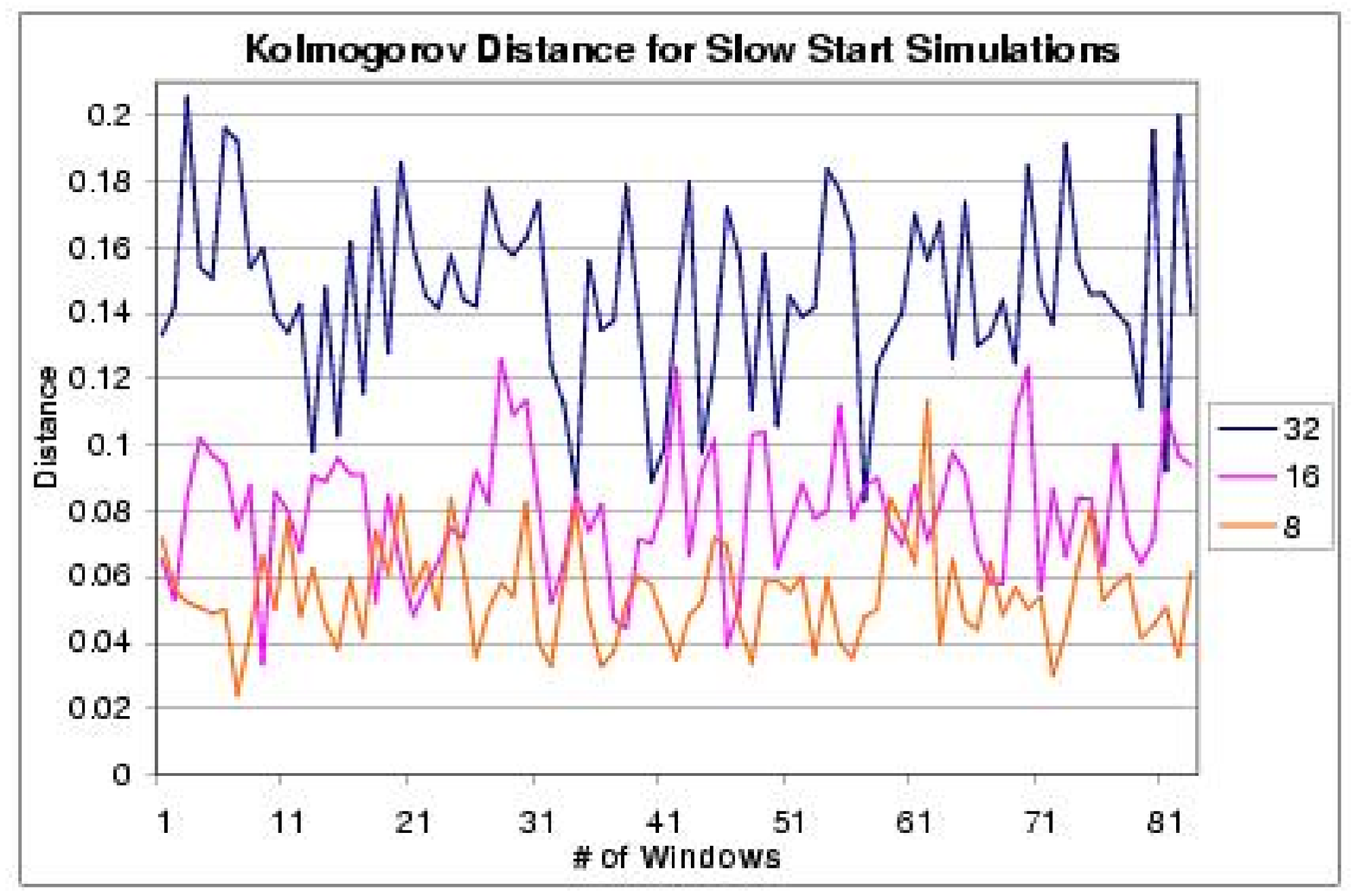}}
 \subfigure[]{\includegraphics[height=150pt, width=200pt]{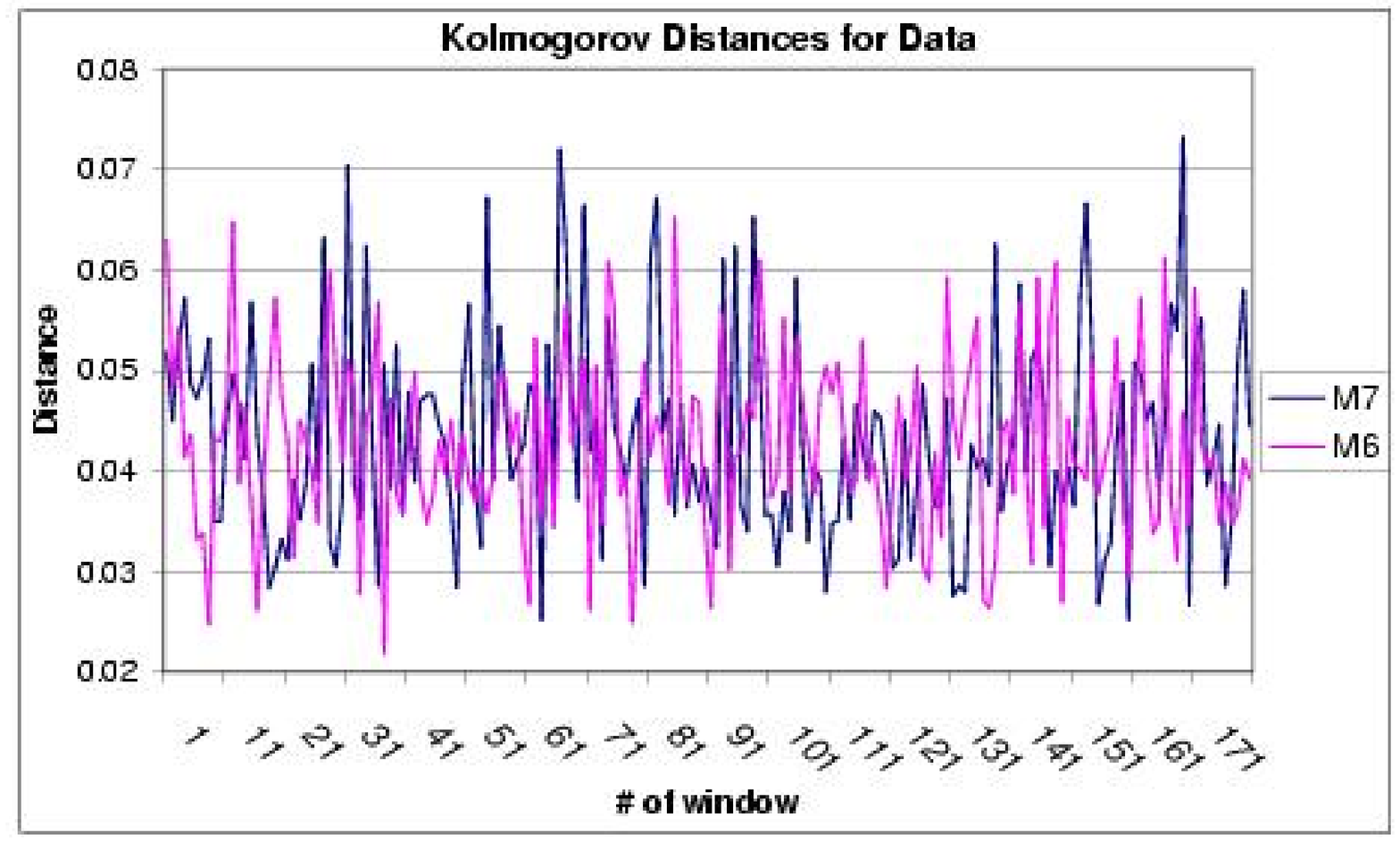}}
 \subfigure[]{\includegraphics[height=150pt, width=200pt]{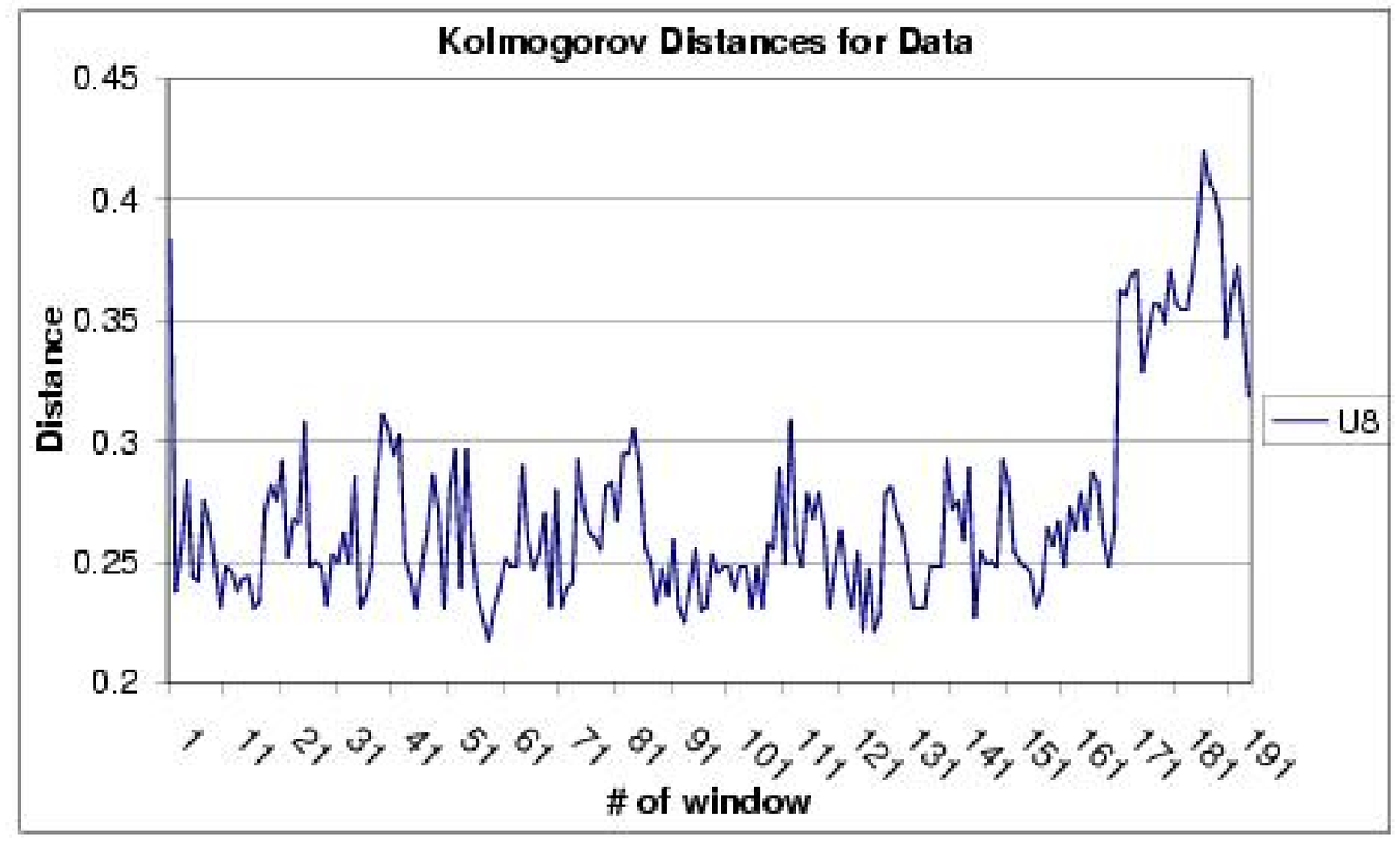}}
 \caption{Kolmogorov Distances of trace marginals -(a), (c), (d)- and Slow Start simulation
          marginals -(b)- from the N(0,1) distribution. The numbers in (b)'s index denote $M$,
          the Slow Start maximum. Marginals have been computed using non-overlapping consecutive
          windows of $2^{9}=512$ bins: the oscillations, and the distances themselves,  
          are too large to be attributed to the finite number of observations used for the 
          computation of the empirical marginals. A time bin of 100ms was used in (a) and of 1ms in 
          (c) and (d). Note the different vertical scales in the different figures.} 
 \label{KolDist}         
\end{figure}
\addtocounter{figure}{-1}
\stepcounter{figure}

\subsection{Is oscillation between Gaussianity and p-Stability a reality?}

\begin{figure}
 \centering
 \subfigure[]{\includegraphics[height=150pt, width=200pt]{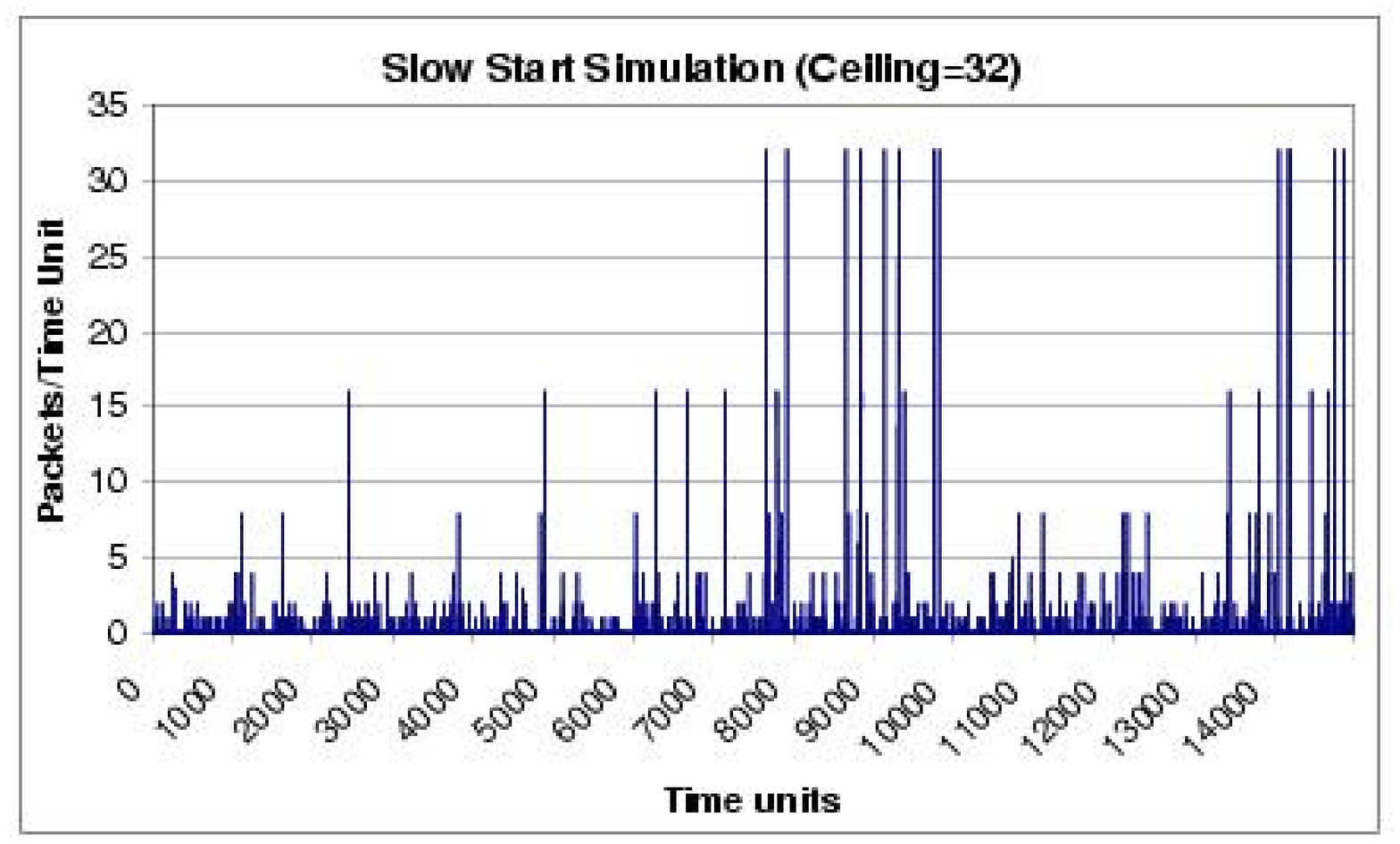}}
 \subfigure[]{\includegraphics[height=150pt, width=200pt]{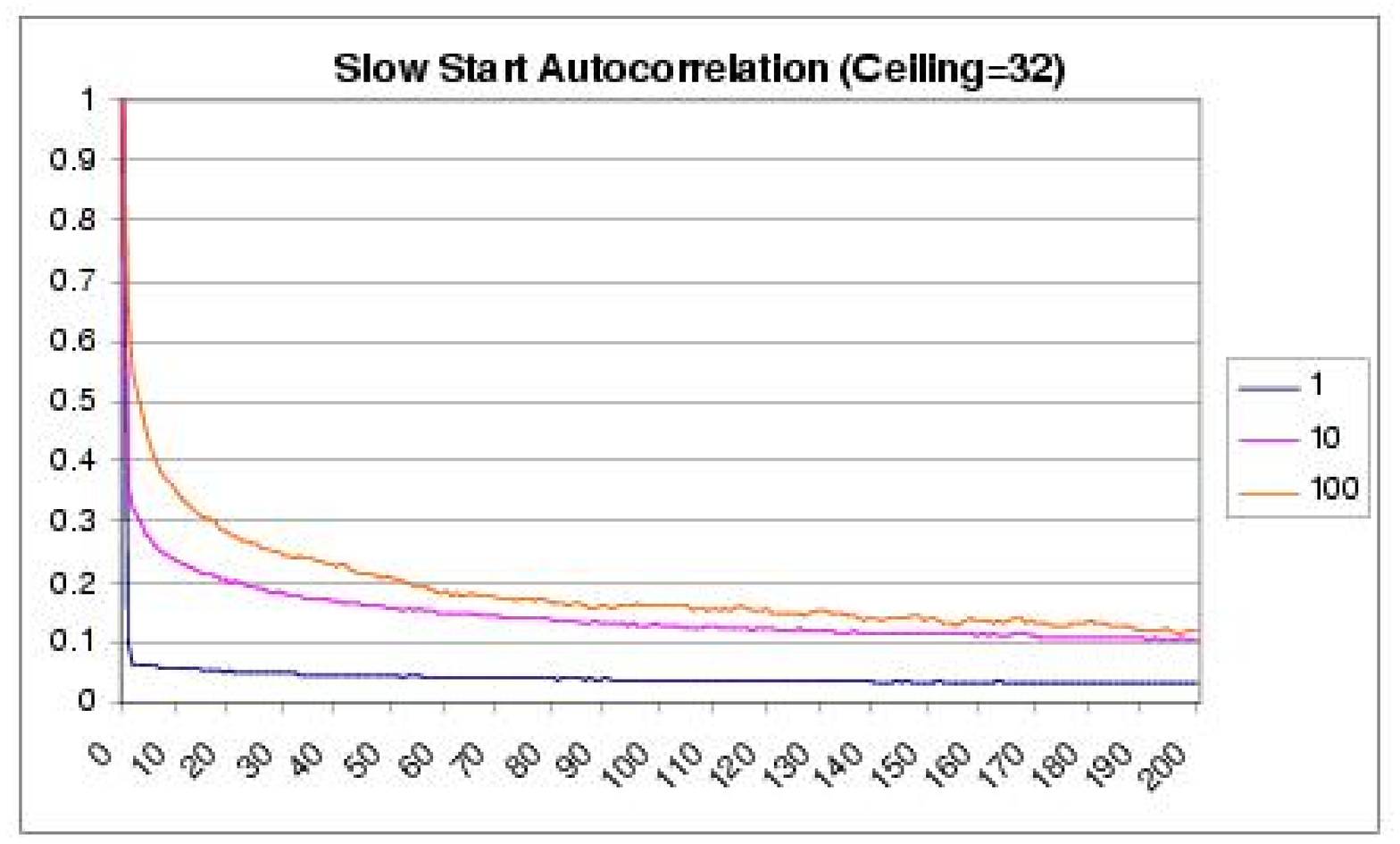}}
 \subfigure[]{\includegraphics[height=150pt, width=200pt]{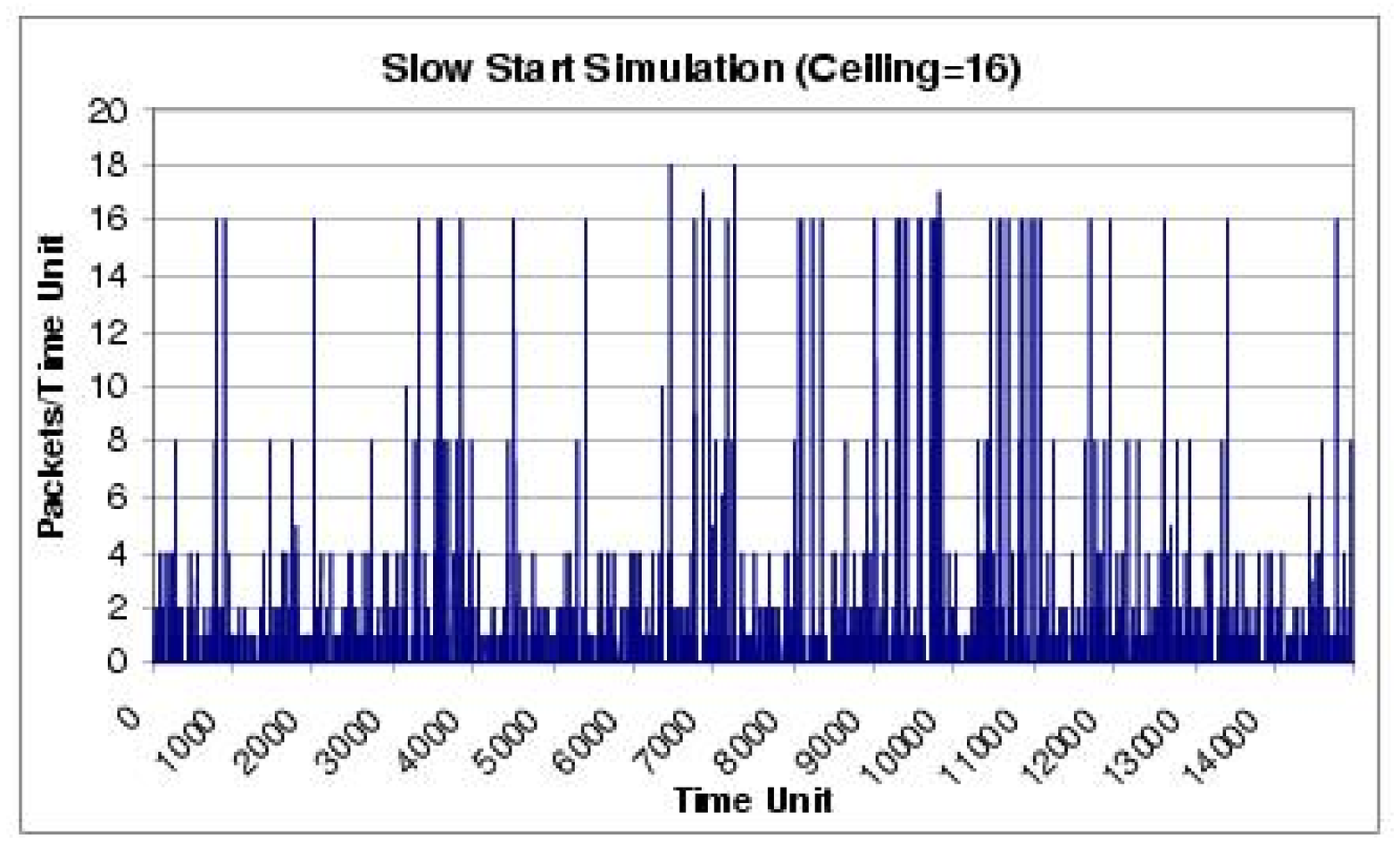}}
 \subfigure[]{\includegraphics[height=150pt, width=200pt]{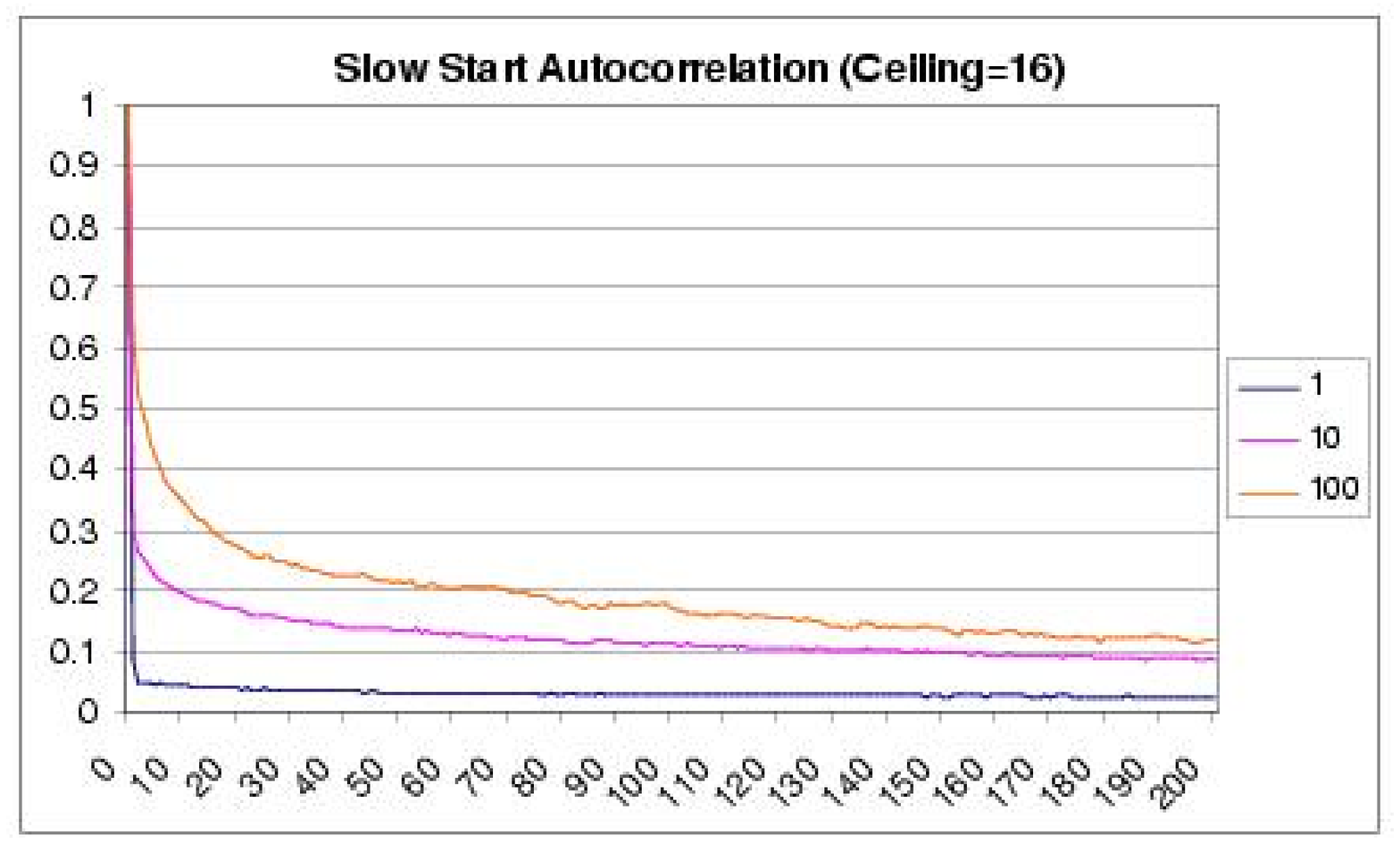}}
 \subfigure[]{\includegraphics[height=150pt, width=200pt]{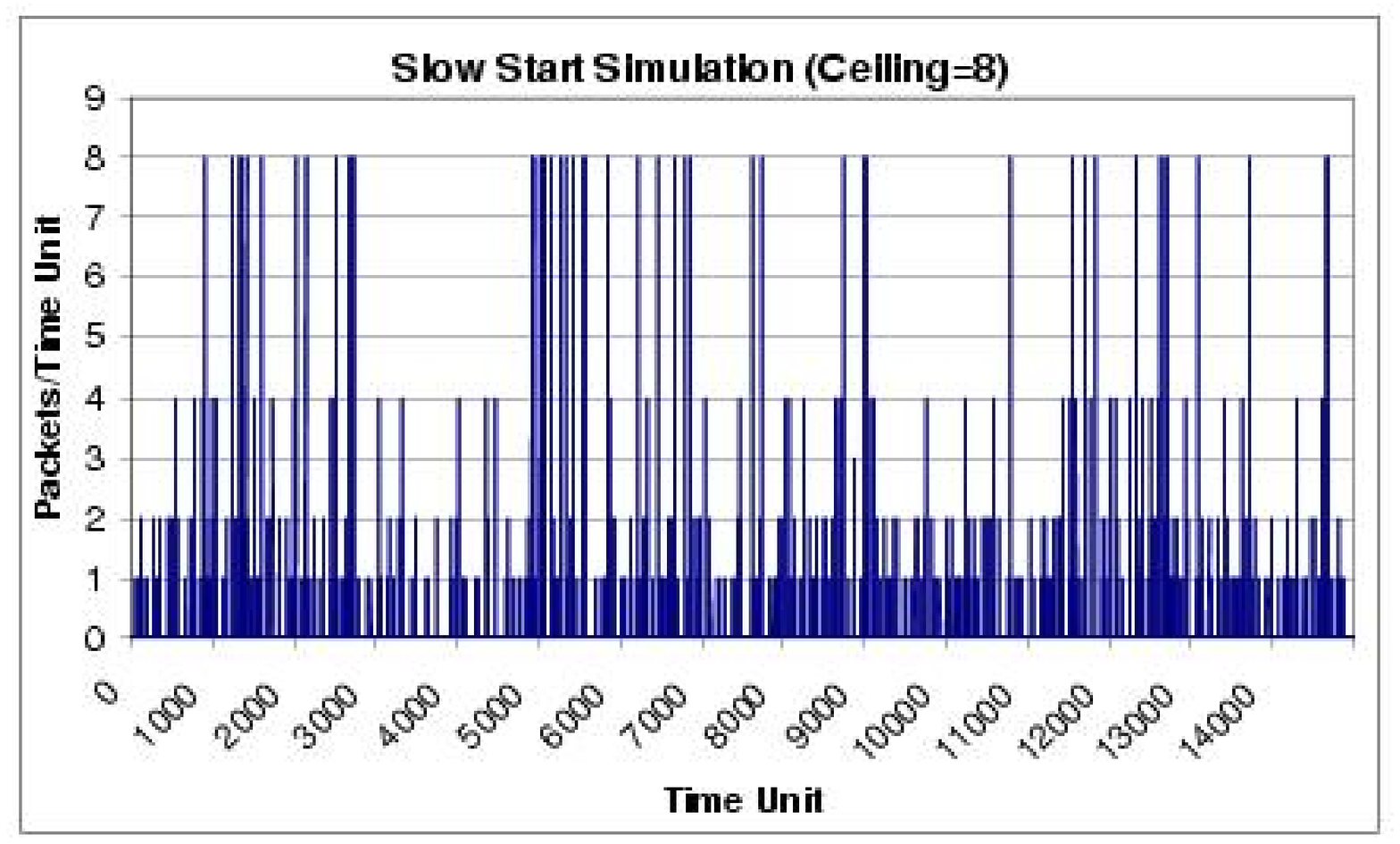}}
 \subfigure[]{\includegraphics[height=150pt, width=200pt]{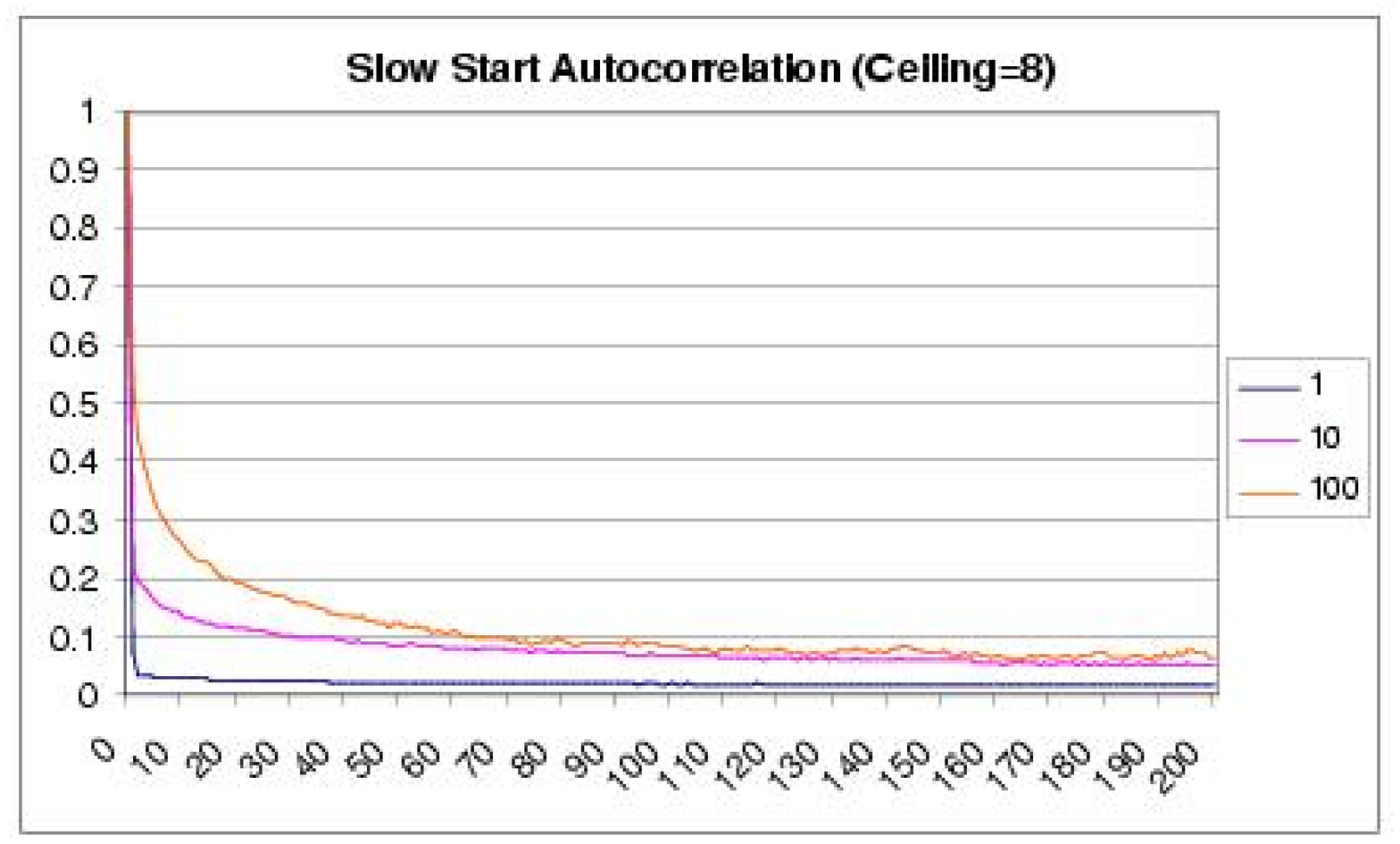}}
 \caption{Excerpts of Slow Start simulation traces with (a) $M=32$, (c) $M=16$, 
          (e) $M=8$, and their corresponding -(b),(d),(f)- autocorrelations.}
 \label{SSSim}         
\end{figure}

Simulations of model C seem to be very successful in resembling the true traffic visually, and in capturing
its long range dependence (see Fig. \ref{SSSim}, \ref{Traffic9497}, and \ref{TrafficM6U8}). To test whether the other prediction of Theorem C, namely the oscillation between Gaussianity and
p-Stability, is indeed seen in real traffic, the following experiment was designed. For simplicity, the description will be based on the 94 trace, when it 
comes to choosing values of parameters, but it can be easily translated into the context of newer traces,  
using the observation that transmission appears to be 100 times faster. 

Start by binning the trace using a bin of $0.1s$. Since under a time bin of $0.1ms$ it is mostly true that at most one packet falls into a particular time bin, the $0.1ms$ bin size can be considered a ``natural'' time bin for the trace; taking a marginal distribution of the trace at this bin size would yield the marginal of the packet size itself. Consequently, the value of every bin is the sum of 1000 ``natural'' bins, i.e. the sum of 1000 identically distributed random variables, and actually distributed as the packet size, if 0 is allowed to be a legitimate packet size. One would hope that this number is sufficiently large to yield a reasonably accurate convergence to the Gaussian limit predicted by the CLT, if convergence to a Gaussian r.v. were to hold, or yield results that are observably at least non-Gaussian, if convergence to a p-Stable limit were to hold. Actual convergence to the p-Stable limit is too much to ask for, because a) a different normalization is required for that, compared to the Gaussian limit, and b) infinite variance is impossible to capture with a finite trace, so that a p-Stable variable cannot be distinguished from a large variance variable, in practice (see Fig. \ref{KolDistTest}(e) and (f)).

In order to examine the trace behavior locally, we divide the trace using consecutive, equisized, non-overlapping windows and take the normalized marginal (corresponding to mean 0 and variance 1) within each window. We then compute the Kolmogorov distance (see (\ref{KolD})) between this marginal and a normal $N(0,1)$ distribution, for each window; graphing this Kolmogorov distance as a function of the bin label illustrates changes in time of the local behavior. As in the case of the windowed Fourier transform, the window size is the result of a trade-off: small windows will yield well localized inaccurate results, whereas large windows will yield poorly localized accurate results. A size of 512 bins seemed a good choice, thus ending up with approximately 140 points in the graph (the trace consists of approximately 72000 bins). The resulting Kolmogorov distance graphs for the data sets are shown in Fig. \ref{KolDist}(a), (c), and (d). In some cases, the Kolmogorov distance is small throughout (M6, M7, L4); in others it oscillates violently. 

In order to calibrate these results, we performed the same computations on realizations of known distributions; the results are shown in Fig. \ref{KolDistTest}(a), (b), and (c). 

In Fig. \ref{KolDistTest}(a) we compare a realization of the normal distribution with the theoretical curve; as expected, we find a very small Kolmogorov distance $d_K$. Based on this, we place the threshold of ``not differing fro Gaussianity'' at $d_K\approx 0.08$.

In Fig. \ref{KolDistTest}(b) one can see the results for the three heavy-tailed distributions of $Y=ax^{-b}$, where $b=0.5,\ 0.7,\ 0.9$ respectively; in these cases the oscillations of the local Kolmogorov distance to Normal have a larger amplitude (the difference between minimum and maximum is about 0.12 in all three cases, which is between 25 and 35\% of the mean Kolmogorov distance), but they all are solidly above 0.2, and move further away from Normal as $b$ increases, as is to be expected. 

Fig. \ref{KolDistTest}(c) shows the local Kolmogorov distance to Normal for four different Weibull distributions, with cumulative distributions $F(x)=1-e^{-ax^b}$. For $b=2$, the graph is in the Gaussian range, as could be expected; as $b$ decreases, one moves away from the Gaussian range, and the oscillations in the local behavior become more marked; nevertheless, their amplitude never exceeds 35\% of the mean. For $b\leq 1$ they remain well segregated from the Gaussian reange at all times. 

We can now evaluate the results for the data sets based on these calibrations. The data sets M6 and M7 look solidly Gaussian, as does the Slow Start simulation with $Max=8$; data 97 and simulations with $Max=16$ dart between Gaussianiy and the region $d_K\geq 0.1$, with a total oscillation amplitude approaching 60\% of their mean; data sets 94, L4 and L5, and simulations with $Max=32$ reach much larger $d_K$, but occasionally reach into Gaussianity, and their total oscillation amplitude exceeds their mean; finally, trace U8 is completely in the ``far from Gaussianity'' region $d_K \geq 0.2$, with again maximum oscillation exceeding its mean. However, what distinguishes these traces is that they reach local Kolmogorov distances typical of distributions with large variance, yet that, unlike e.g. the large variance Weibull distributions, they have very large amplitude that occasionally take them back to or close to Gaussianity. For $Max$ in a ``critical'' range ($Max=16$ for our parameters), our simulations have the same character. This suggests that the data, as well as the simulation, exhibit an oscillation between Gaussian and far-from-Gaussian behavior, consistent with Theorem C.     

\addtocounter{figure}{-1}
\stepcounter{figure}

\subsection{Traffic amount is correlated with marginals locally}

An experiment related to the previous one is the following: compute the Kolmogorov distances of the traces, as before, and also calculate the total traffic within each window. This procedure generates two new sequences. Table \ref{CrossCorr} gives their correlation for our different data sets, for 0.128s-binned traffic. Interestingly, it is non-negligible and negative for all the data sets that are clearly not uniformly Gaussian according to the discussion above (i.e. 94, L4, L5, U8, and 89). 

This is again consistent with Theorem C, by the following argument: Theorem C states that the more numerous the active users are, the more Gaussian the traffic will be, and thus the smaller the Kolmogorov distance. On the other hand, the number of active users is highly and positively correlated with the traffic volume. The combination of the two statements implies that traffic volume and the Kolmogorov distances should be negatively correlated, unless, of course, the traffic is already (almost) Gaussian, in which case no prediction according to the theorems is possible (the data sets close to Gaussianity are traces 97, M6, and M7, whose correlations are indeed either negligible or positive). Note that Gaussian-type traffic is \emph{not} extremely spiky. The negative correlations in Table \ref{CrossCorr} thus illustrate that when more users are active, and when traffic is heavier, it becomes less spiky. This seems counterintuitive; yet it is consistent with Theorem C, and borne out by the data. 

\begin{table}
\begin{center}
{Cross-correlation between Kolmogorov distances and traffic:}\\ \bigskip
\begin{tabular}{||c|c|c|c|c|c|c|c|c||} \hline
\emph{Trace} & 94 & 97 & L4 & L5 & M6 & M7 & U8 & 89\\ \hline
\emph{Correlation} & -0.300 & 0.392 & -0.161 & -0.244 & 0.031 & -0.048 & -0.787 & -0.806\\ \hline
\end{tabular}
\caption{\label{CrossCorr} Cross-correlation between Kolmogorov distances of windowed empirical marginals from $N(0,1)$ and the total traffic within the window, for the real traces.}
\end{center}
\end{table}

Note here that the correlation can vary dramatically with the time
bin used, due to the presence of ``levels'', which will be discussed in Section \ref{modelD}. But, in general, for large time bins, the correlations are large, as Theorem C predicts. 

\subsection{A comment on Slow Start and heavy-tailed session sizes}

Note that the importance of the impact of Slow Start on the traffic is intimately coupled to the fact that the
distribution of $L$ is heavy-tailed. If $L$ is always
very large, the window size will reach $M$ and will stay there. Thus, with
the exception of the first few times, the user will always be sending $M$
packets at a time. This situation is little, if at all, different from model
B. If, on the contrary, $L$ is always very small, then
the \emph{effective} $M$ will be very small, since Slow Start will be
stopping at a very early stage, the session being over.

Consequently, Slow Start can be the star only when most of the sessions are small, 
but occasionally huge ones show up. It
is precisely these latter ones that form spikes, taking advantage of a large
value of $M$. Rephrasing in network terminology, ``spikes
are created by large sessions over fast links'' \cite{SRB1}. Statistical
analysis of the \emph{flows} of e.g. trace 94 reveals that they are indeed very
short, with occasional large ones.

\begin{figure}
 \centering
 \subfigure[]{\includegraphics[height=150pt, width=200pt]{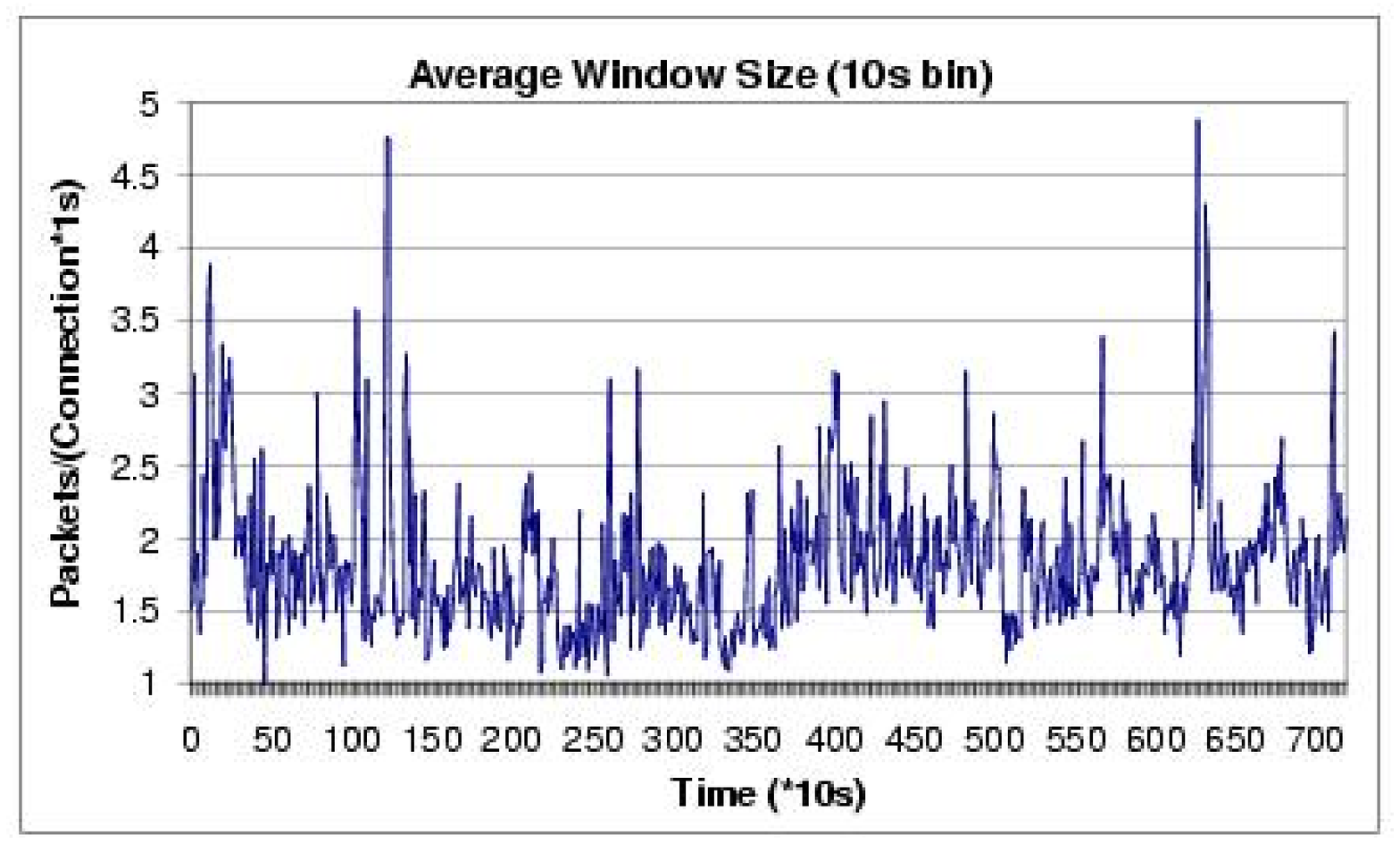}}
 \subfigure[]{\includegraphics[height=150pt, width=200pt]{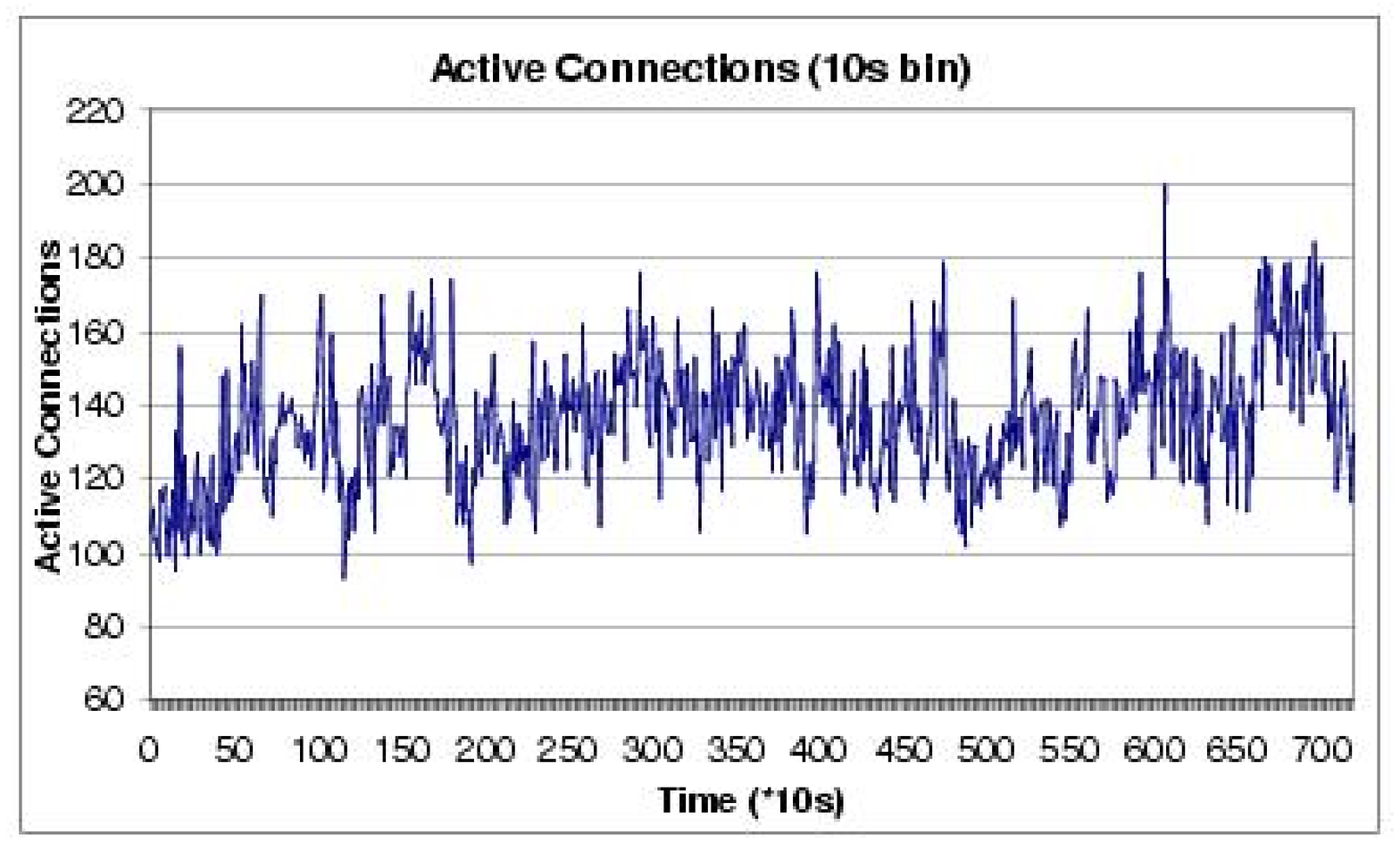}}
 \caption{Average window size (a) and number of active connections (b) for 
          94 trace, using a time bin of 10s.}
 \label{Stat94}         
\end{figure}

The assumption made implicitly in the discussion above is that a spike is the
result of one user, not the aggregate result of many of them. Although this
has already been confirmed \cite{SRB1}, an independent confirmation is
offered here: divide the time axis into bins of constant size, and in each
bin measure the number of active connections (i.e. connections that send in
this interval more than one packet), and the average window size (i.e. the
number of packets sent divided by the number of active connections). If, for
a spiky traffic, it is the latter that is spiky, and not the former, then
this indicates that some users must be sending a lot of packets. Fig. \ref{Stat94} 
shows that this is indeed the case for real traffic. 

\addtocounter{figure}{-1}
\stepcounter{figure}

\subsection{Earlier work on the Gaussianity and p-Stability of the traffic}

The idea of a process that can be either self-similar of $p$-Stable, contingent upon the value of a parameter, already appears in the form of theorems in \cite{MRRS1}, and of observations in real world traces in \cite{SRB1}. More precisely:

\begin{itemize}

	\item \cite{SRB1} analyzed empirically the oscillation of the traffic between Gaussianity and non-Gaussianity and offered a plausible networking explanation. No attempts were made, however, to explain it mathematically. Consequently, the authors were unable to justify the observed negative correlation between amount of traffic and spikiness. 
	
	\item \cite{MRRS1} provides an interesting variant of model A, the ISPM, and considers a process equivalent to $Z^a_n(Tt)$, according to the notation of section \ref{modelA}. It proves then that the limit process, as 
$T\rightarrow\infty$ and $n\rightarrow\infty$ simultaneously, can be either p-Stable or Gaussian, depending on the rates by which $n$ and $T$ tend to $\infty$. In our opinion, though, the limit process does not relate to Internet traffic in an obvious way, and its mathematical predictions, at least as stated in
\cite{MRRS1}, cannot be directly translated in the language of networks, except yielding a result for the traffic at 
coarse scales on very busy links, which is exactly the situation the two limits above describe. Moreover, no attempt
was made in \cite{MRRS1} to corroborate the model's assumptions by empirical evidence, or to compare the
properties of the simulated traffic it produces with the corresponding properties of real traffic. Finally, the limit processes in \cite{MRRS1} are strictly either self-similar or p-Stable; no allusion is made to the possibility of oscillation of the process between these two extremes. 

\end{itemize}

\section{Model D}

\label{modelD}

\subsection{Introduction}

\begin{figure}
 \centering
 \subfigure[]{\includegraphics[height=150pt, width=200pt]{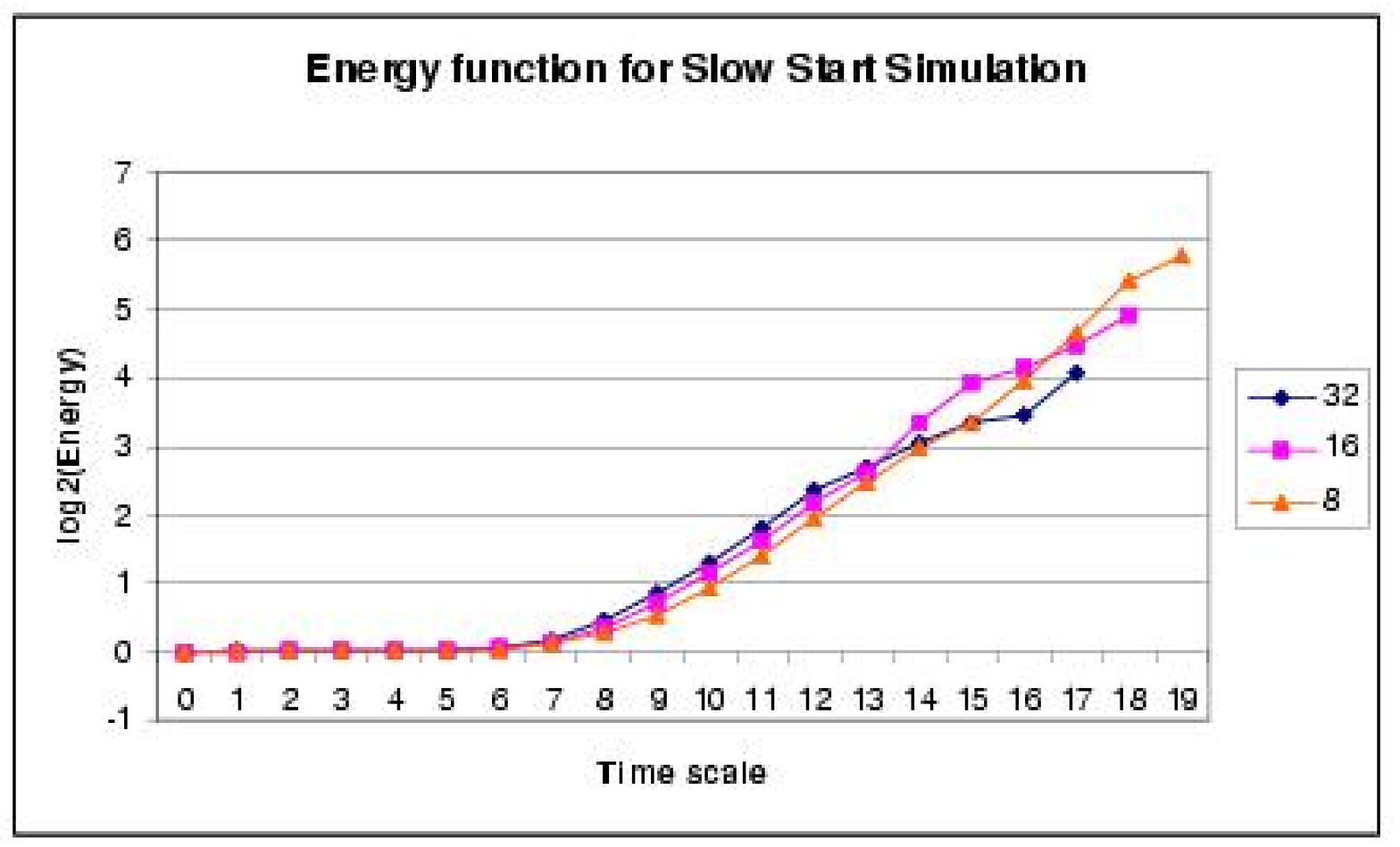}}
 \subfigure[]{\includegraphics[height=150pt, width=200pt]{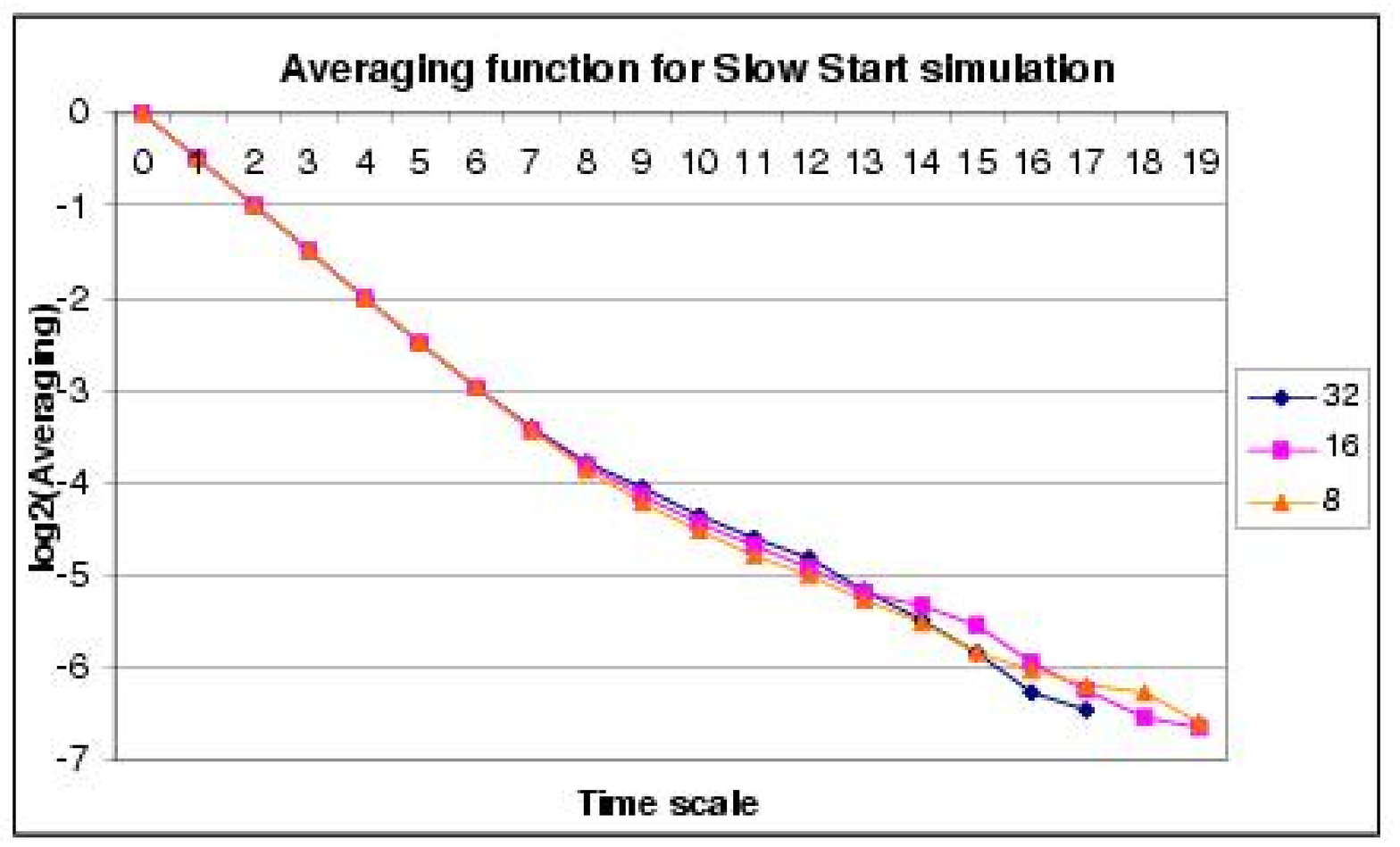}}
 \caption{Energy (a) and Averaging (b) functions for the Slow Start Simulations: although some slight curving is present in these simulations, it is much more modest than what is observed in real traffic (see Fig.~\ref{En}).}
 \label{EnAvSS}
\end{figure}

Models B and C were able to capture and explain crucial features of
Internet traffic, such as its long range dependence, its different behavior in different time scales,
and its spikiness (compare Fig. \ref{SSSim}, \ref{Traffic9497}, \ref{TrafficM6U8}, and \ref{AC}).
Unfortunately, the Energy function of the simulation is still far from
reality (compare Fig. \ref{EnAvSS} and \ref{En}). This suggests that pieces of the puzzle are still missing.

\addtocounter{figure}{-1}
\stepcounter{figure}

\begin{figure}
 \centering
 \subfigure[]{\includegraphics[height=150pt, width=200pt]{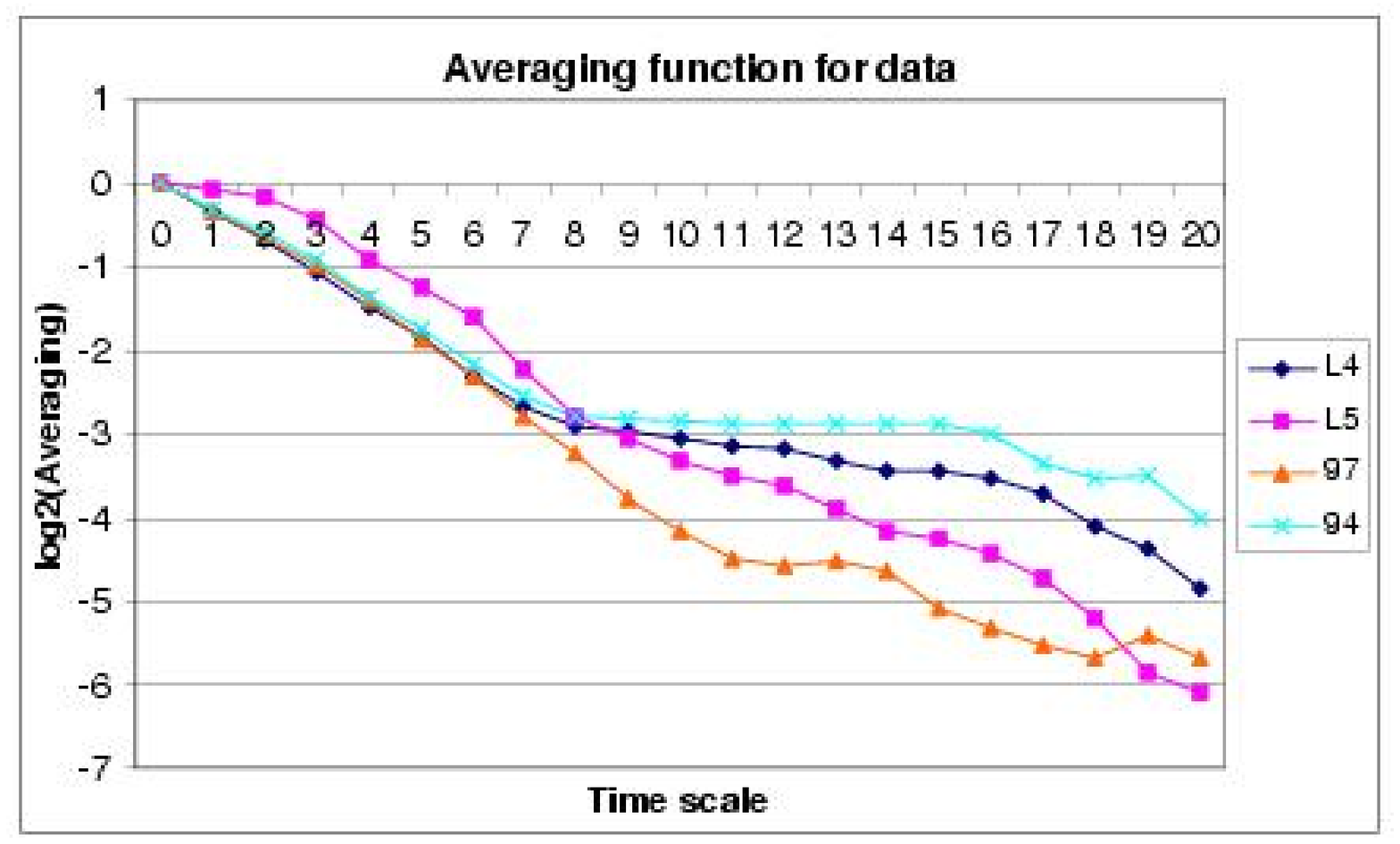}}
 \subfigure[]{\includegraphics[height=150pt, width=200pt]{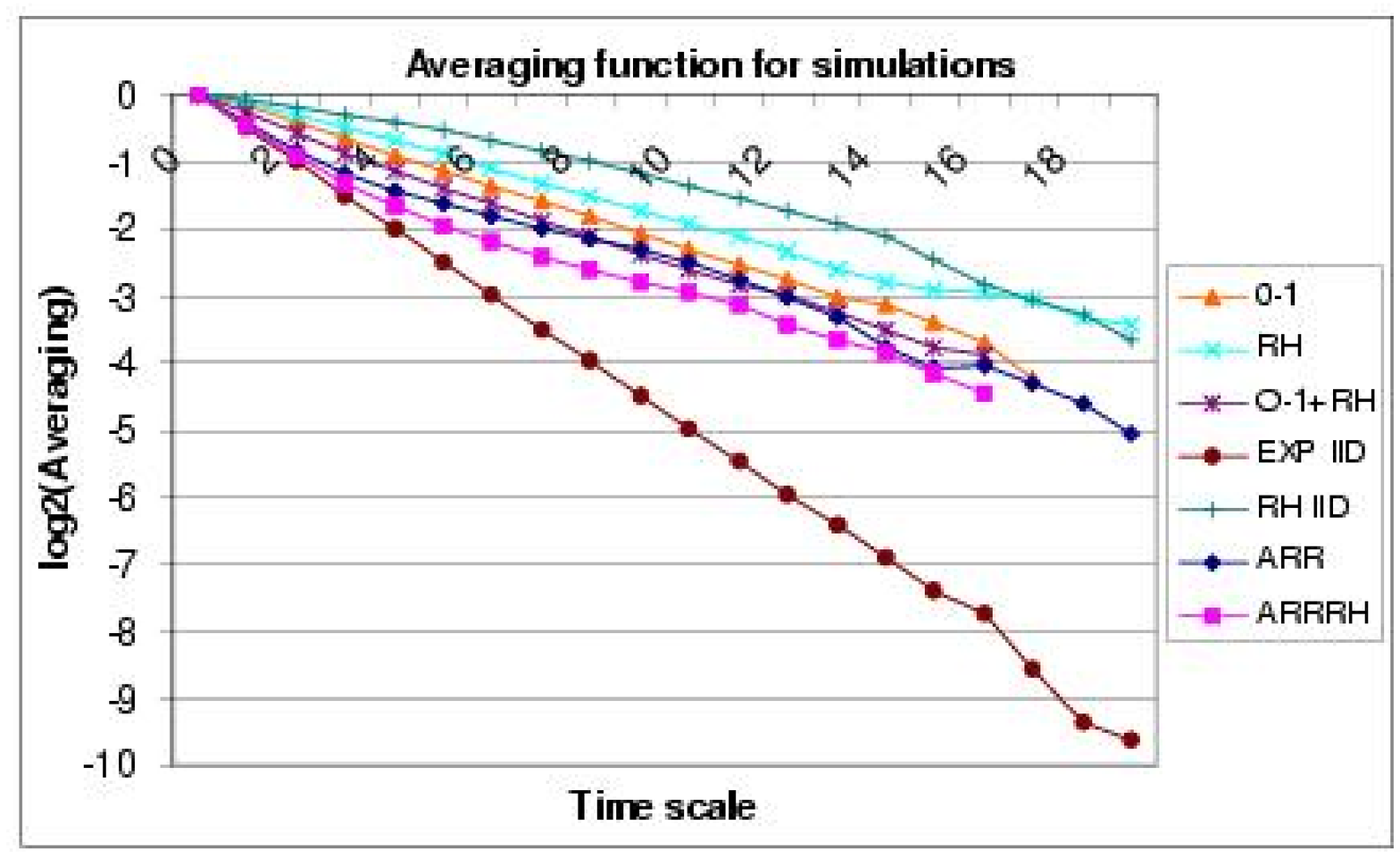}}
 \subfigure[]{\includegraphics[height=150pt, width=200pt]{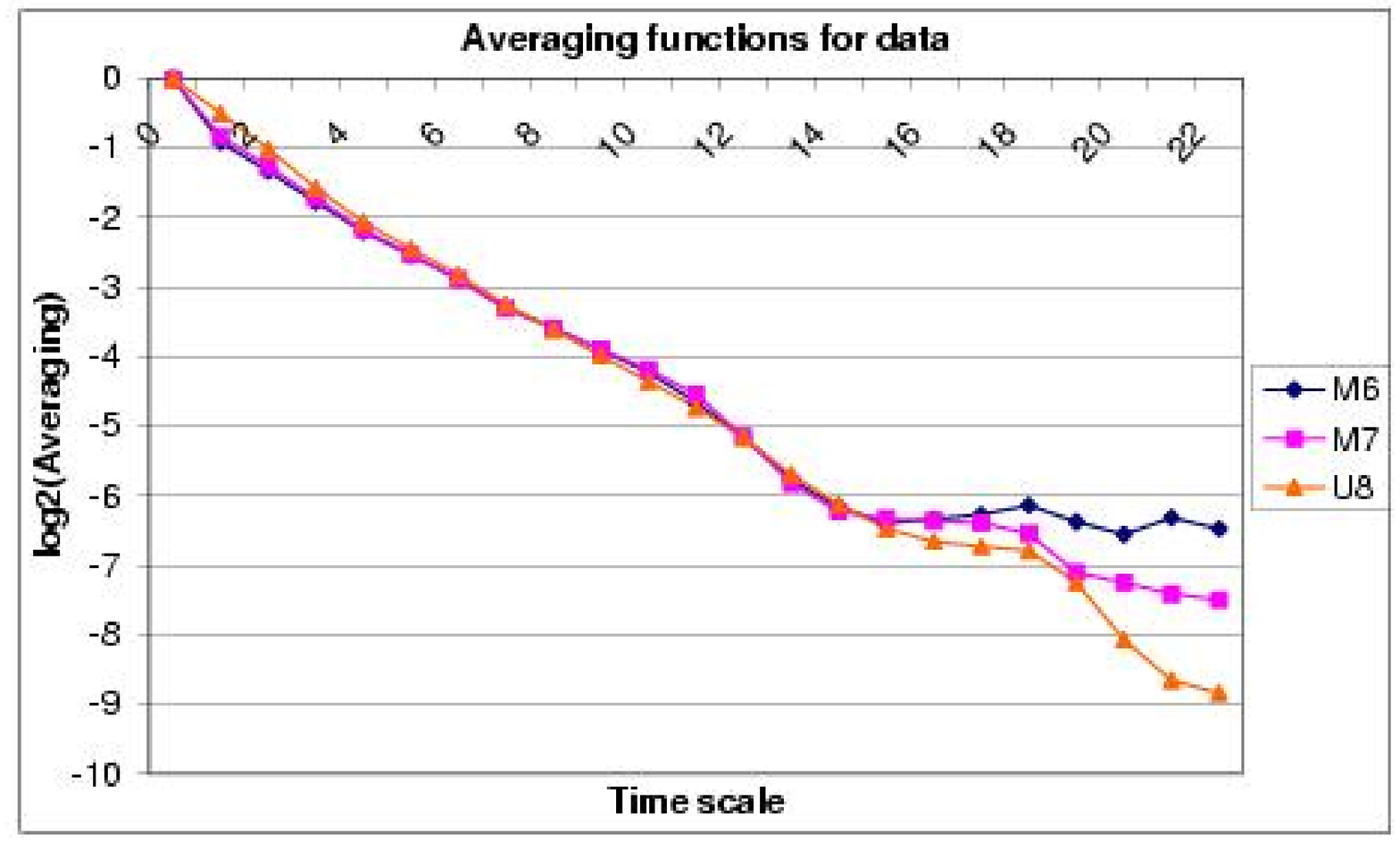}}
 \subfigure[]{\includegraphics[height=150pt, width=200pt]{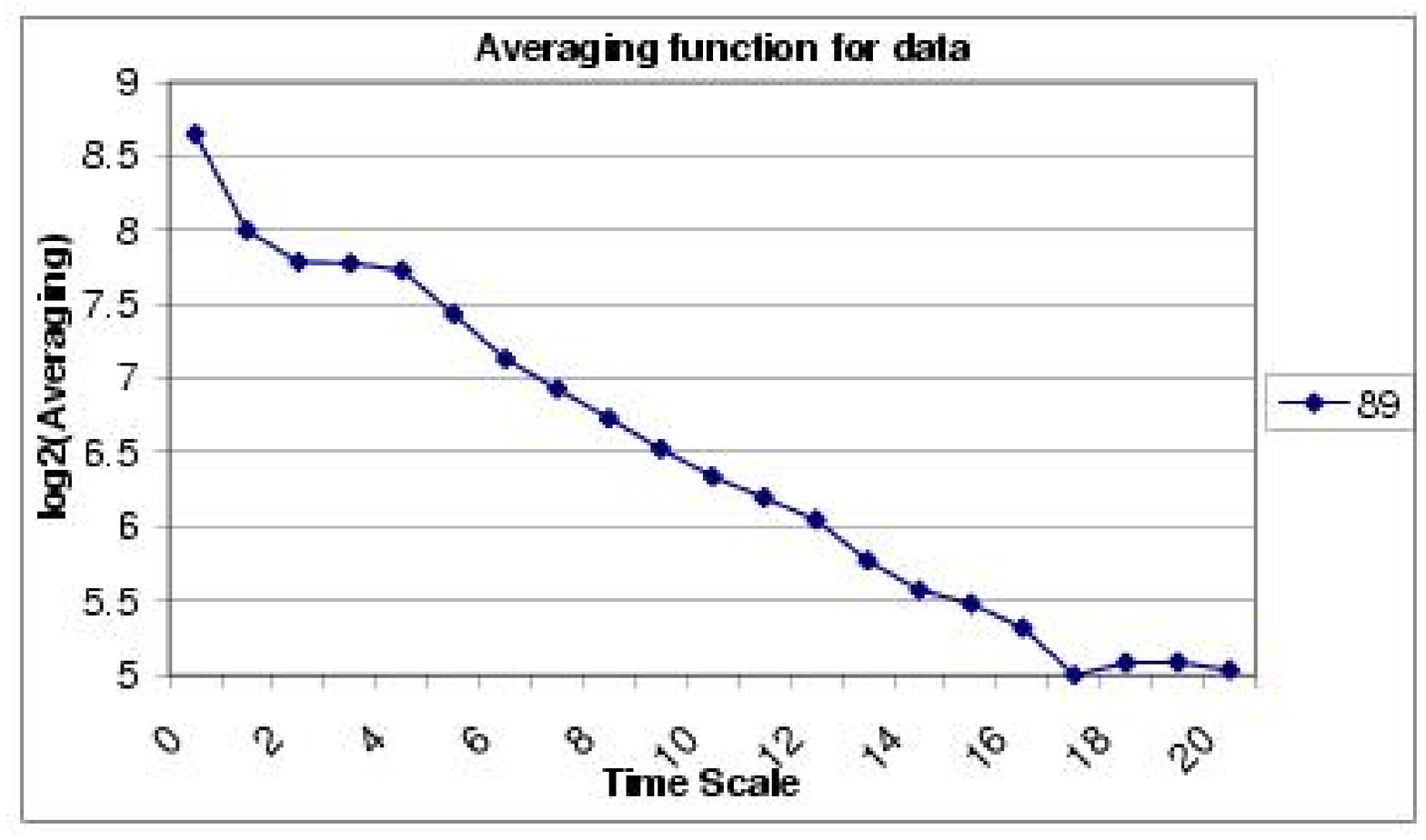}}
 \caption{Averaging functions for data -(a), (c), (d)- and simulation (b) traces.}
 \label{Av}
\end{figure}
\addtocounter{figure}{-1}
\stepcounter{figure}

\subsection{A closer look into Averaging functions}

So far, the simulations failed to introduce a substantial curvature into the
Averaging or Energy function (Fig. \ref{En}, \ref{Av}), although the processes used spanned a 
variety of attributes (heavy/light-tailed marginals, independence/long-range dependence, etc.), 
and even combined them (as in the $\alpha$-$\beta$ traffic model). It was precisely this
characteristic that stimulated research of ``multifractality'', questioned
the accuracy of traffic description through Poisson or self-similar
processes, and ended up demonstrating that the traffic is a fundamentally
different process from these two \cite{SRB1,FGHW1,FGW1}.

The Averaging functions of the data traces (Fig. \ref{Av}) are truly remarkable:
recall that a slope $\alpha $ at time scale $\log_{2}(n)$ indicates a behavior of
type $n^{-\alpha }$ for the difference between two mean value estimators of $n$ terms; 
if $\alpha$ depends on $n$, the process is not
self-similar. This is the case for all of the traces. Moreover, as explained in the heuristic argument in section \ref{EnAvDef}, one expects, if the $X_i$ are asymptotically independent, that the Averaging function will decrease monotonically as a function of $n$, at least for large $n$. This is not the case for our data. For instance, the Averaging function of trace 94 is constant in the region from $n=2^8=256$ till $n=2^{15}=32768$. In other traces the Averaging function even \emph{increases} with $n$: trace 97 exhibits a little ``bump'' at
time scale 13 (Fig. \ref{Av}(a)), and trace M6 exhibits a bigger bump at time scale 18 
(Fig. \ref{Av}(c)); simulations carried out later in this section achieve an even more spectacular result (Fig. \ref{AvLev}).

In our heuristic argument we assumed asymptotical independence of the $X_i$. We now know that the $X_i$ in our traces cannot be considered asymptotically independent, at least relatively to the trace duration, and we interpret the ``deviant behavior'' of the Averaging function as a manifestation of this lack of independence. Why should it show up much more at some scales than in others? This is the question we address in this section. 

\subsection{Averaging functions and multiscale structure of the traffic --- \textit{levels}}

There are essentially two important aspects of Internet traffic that have been overlooked 
throughout the discussion of the previous models:

On the one hand, the exact attributes of time intervals between packet
emissions were hitherto supposed to be of little importance: in the
theorems, they were assumed to follow any distribution of finite variance,
and simulations were performed with an exponential one. The exact
distribution has indeed been an object of study, and it turned out that
the Poisson model for a single user fails \cite{PF1}, and that (bi)Pareto, Weibull,
and empirical (Tcplib \cite{PF1}) distributions come into play \cite{NSSW1}.
The former may violate the finite variance assumption.

On the other hand, traffic seems to comprise many phenomena that are hierarchically structured:
for example, users organize texts in chapters, chapters in paragraphs,
paragraphs in sentences etc., thus one expects a variety of durations of
pauses and activity. This will have an effect on traffic, if one supposes
the user types over the net. Also, during Web browsing, users read, process
the information, and act accordingly. Finally, in the network there is
protocol activity spanning many time scales: congestion control occurs
every fraction of a second, but routing information is exchanged every
half an hour, or so.

It appears then that many time scales are of importance, each for a
different reason, but all because some kind of activity is taking place there. 
The exact values of these time scales, or other details
about them, are still unclear, to the best of our knowledge, although they
are generally accepted to exist \cite{PF1, RRCB1, FGWK1}, as the result of both user behavior and
protocols. This activity in multiple time scales, however, can alter completely the
autocorrelation of the traffic, introducing strong correlation in some time scales and weak
in others. This, in turn, will impact the Averaging function as well, according to (\ref{ACAv}).
Consequently, the curving of the Averaging function may have its roots in this multiscale activity, 
and this seems worth investigating.

\subsection{Evidence for multiscale activity}

A useful and simple tool for the detection of these time scales (\emph{levels}) in
the traffic is the \emph{characteristic function} of a connection\footnote{%
Similar diagrams exist in \cite{PF1}.}. A \emph{connection} is defined to consist of
the entries of the trace corresponding to a particular value of sender host,
receiver host, sender port, and receiver port, and is essentially what the 
$W_{i}$ stand for in Theorems B and C. The characteristic function of a connection is a 
function of time, equal to 1 if that particular time is a time stamp belonging to the connection,
and 0 otherwise.

\begin{figure}
 \centering
 \subfigure[]{\includegraphics[height=150pt, width=200pt]{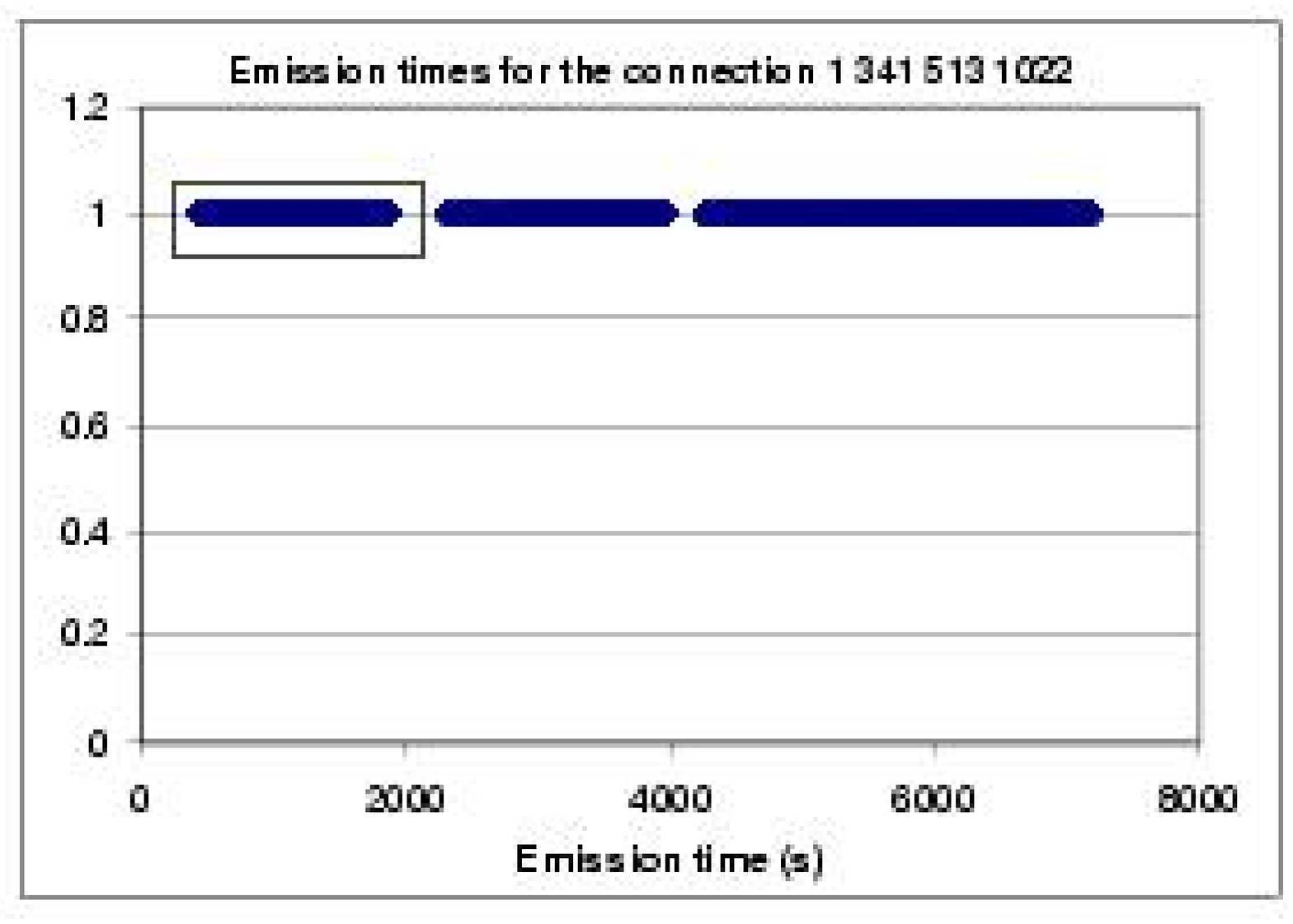}}
 \subfigure[]{\includegraphics[height=150pt, width=200pt]{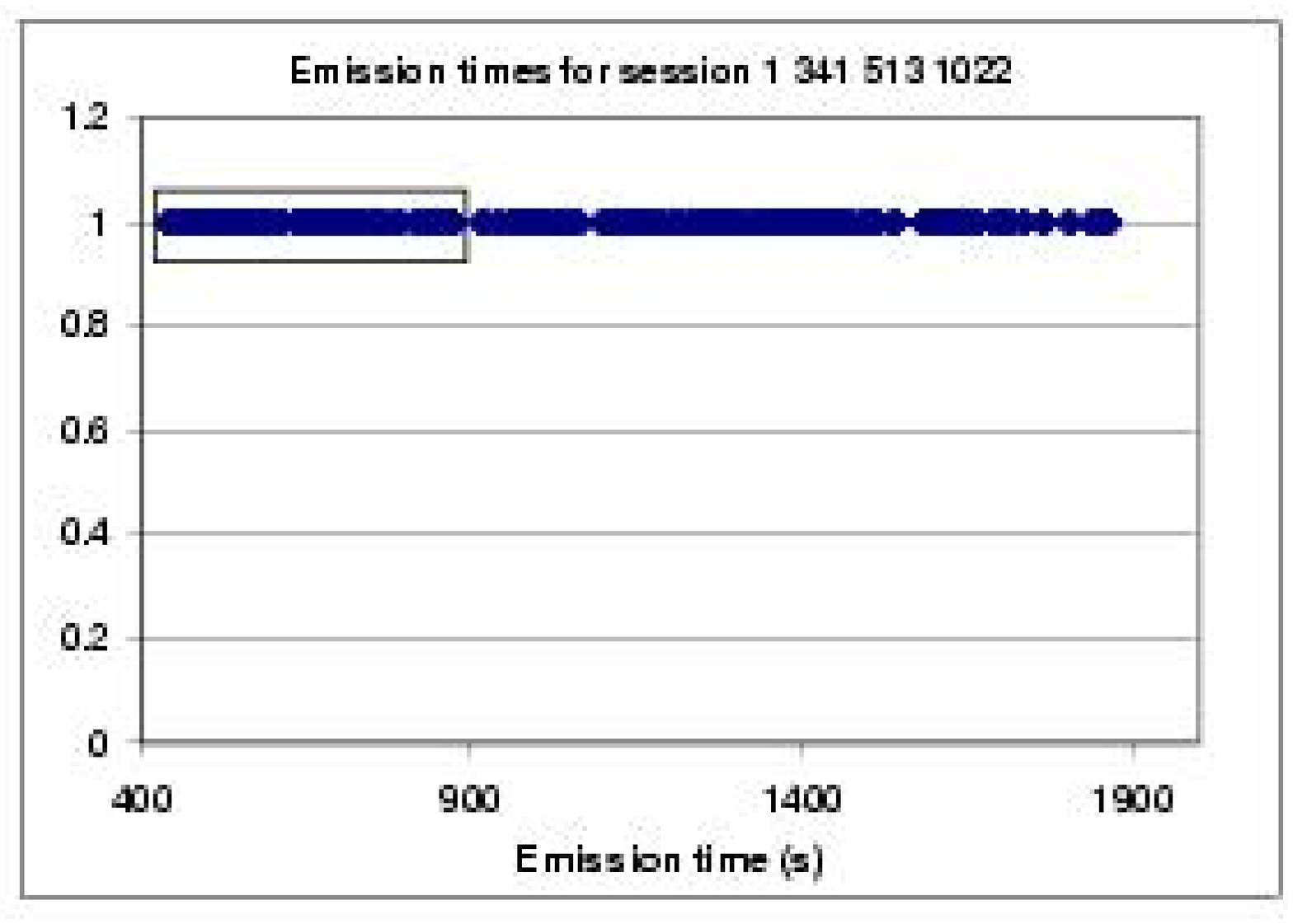}}
 \subfigure[]{\includegraphics[height=150pt, width=200pt]{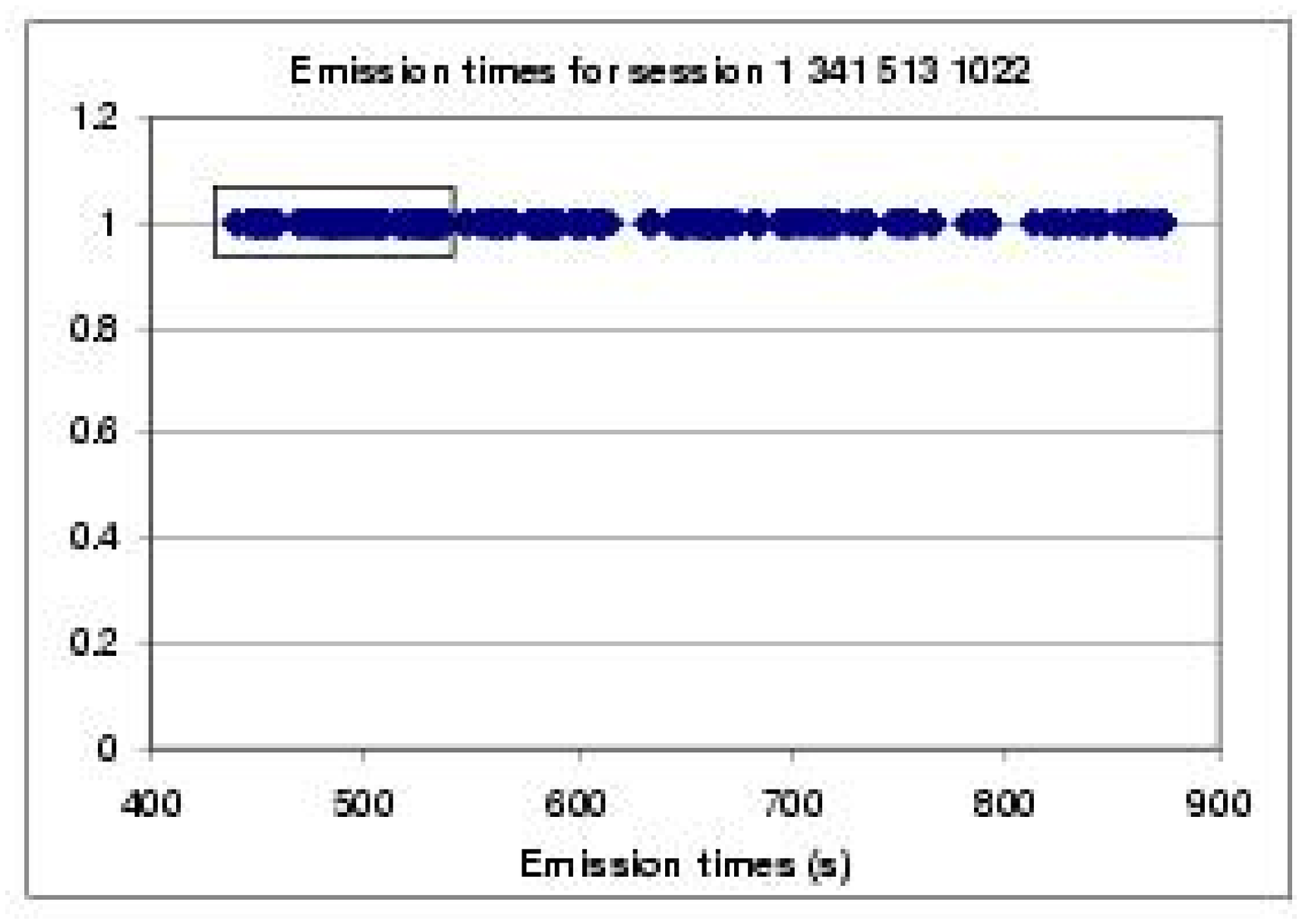}}
 \subfigure[]{\includegraphics[height=150pt, width=200pt]{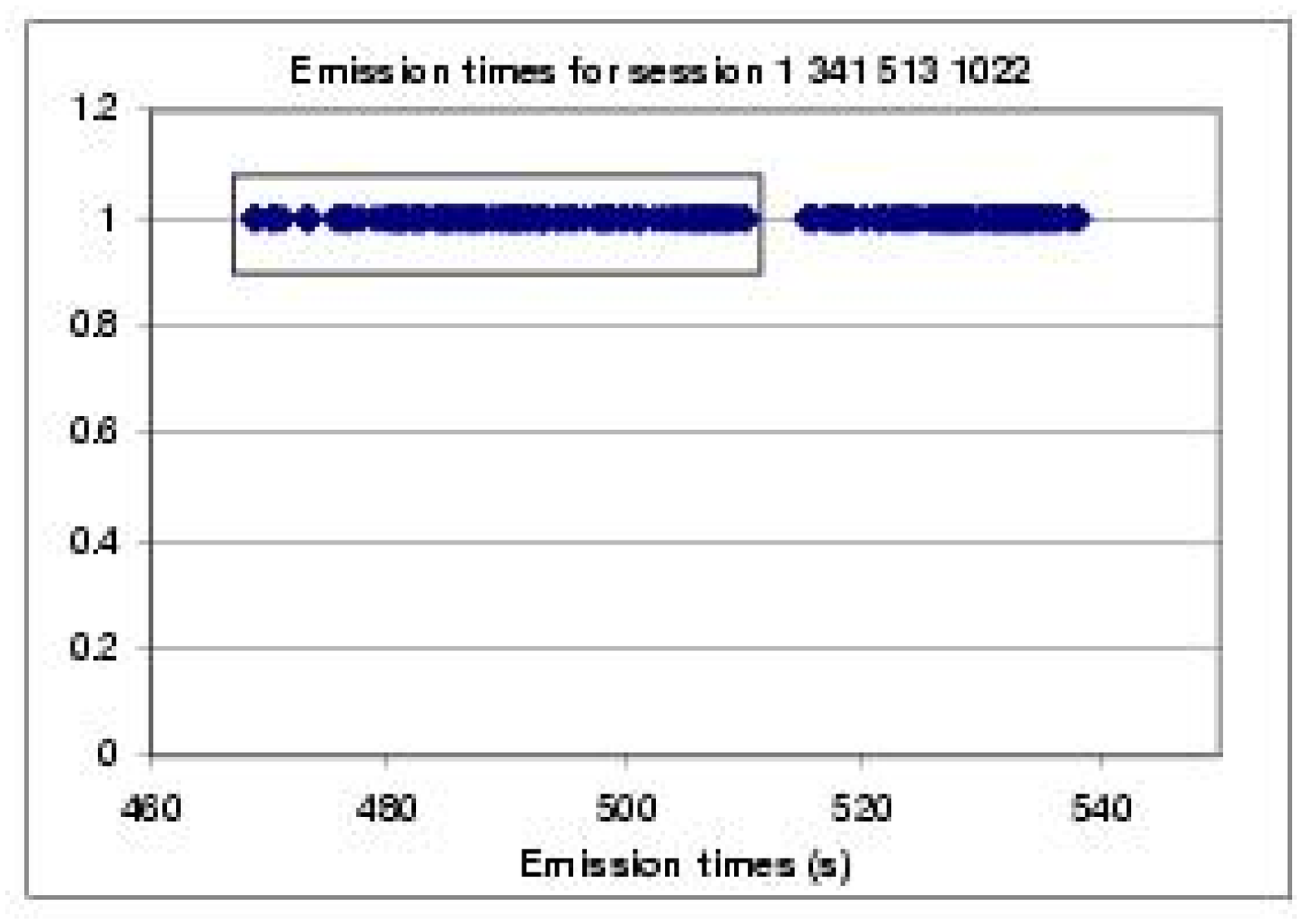}}
 \subfigure[]{\includegraphics[height=150pt, width=200pt]{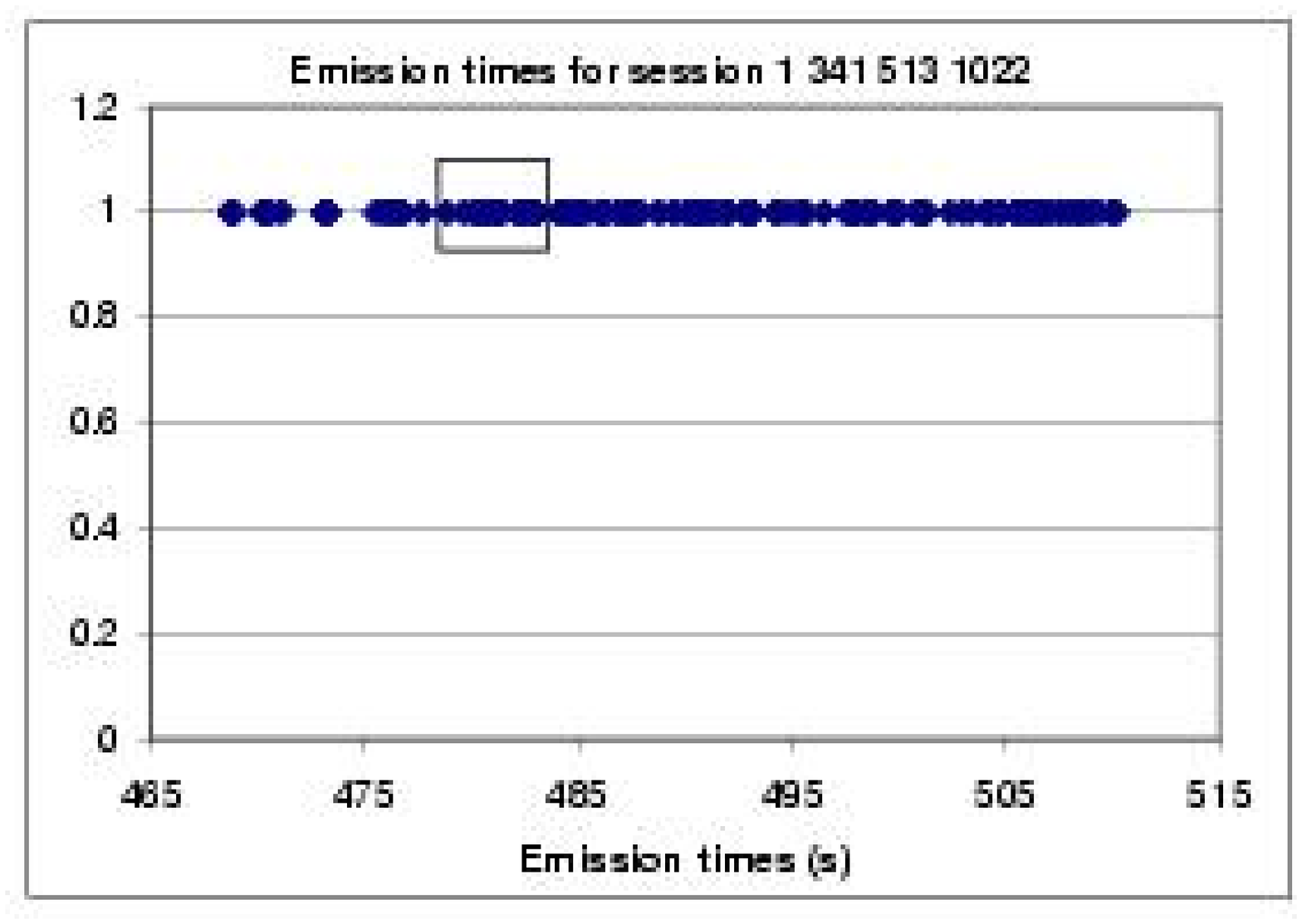}}
 \subfigure[]{\includegraphics[height=150pt, width=200pt]{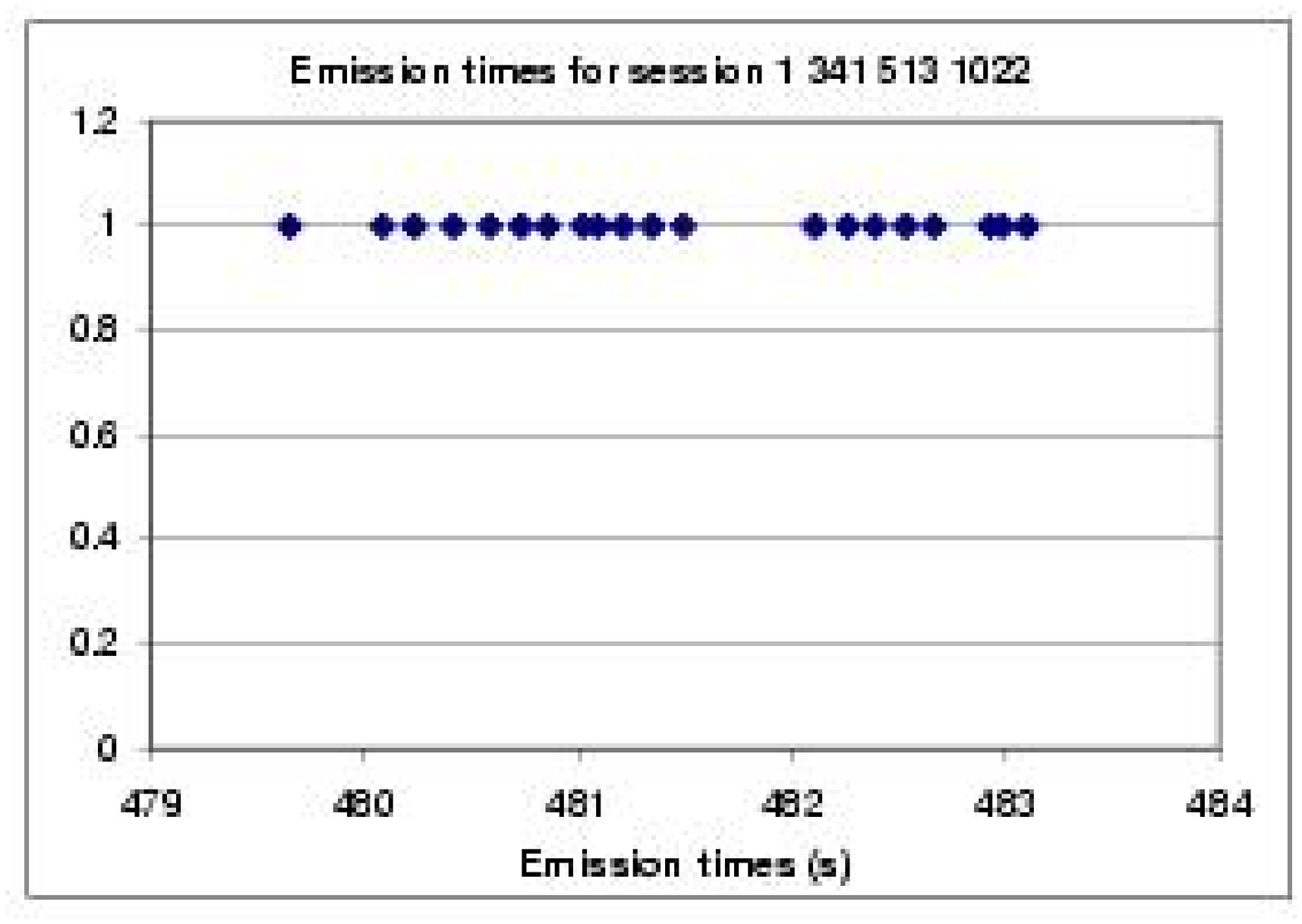}}
 \caption{The characteristic function of connection 1 341 513 1022 from 
          trace 94, in 6 detail levels. Each figure presents a ``zoom-in'' in the boxed segment
          of the previous one, illustrating a fractal-like behavior. Here, 6 levels can be detected: the (approximate) lengths of their $1$-intervals are $2000,\ 500,\ 100,\ 40,\ 5,\ \text{and}\ 0.1s$, whereas the lengths of their $0$-intervals are much smaller in each case: $400,\ 15,\ 10,\ \text{and}\ 4s$. Note that $0$s are not plotted.}
\label{Lev}         
\end{figure}

Characteristic functions seem to have a fractal structure (see Fig. \ref{Lev} and \ref{LevND});
moreover, analysis of the characteristic functions of the 50 largest connections of
trace 94, and connections of trace M6, leads to the conclusion that levels of
connections do not average out in the aggregate traffic, but, on the
contrary, have a strong influence on it. Capturing this effect will be the main goal of this section; but first, 
let us propose a more quantitative level detecting method in individual user sessions, which establishes the presence of multiscale levels more firmly than the anecdotal evidence of Fig. \ref{Lev} and \ref{LevND}. 

\addtocounter{figure}{-1}
\stepcounter{figure}

\begin{figure}[t]
 \centering
 \subfigure[]{\includegraphics[height=150pt, width=200pt]{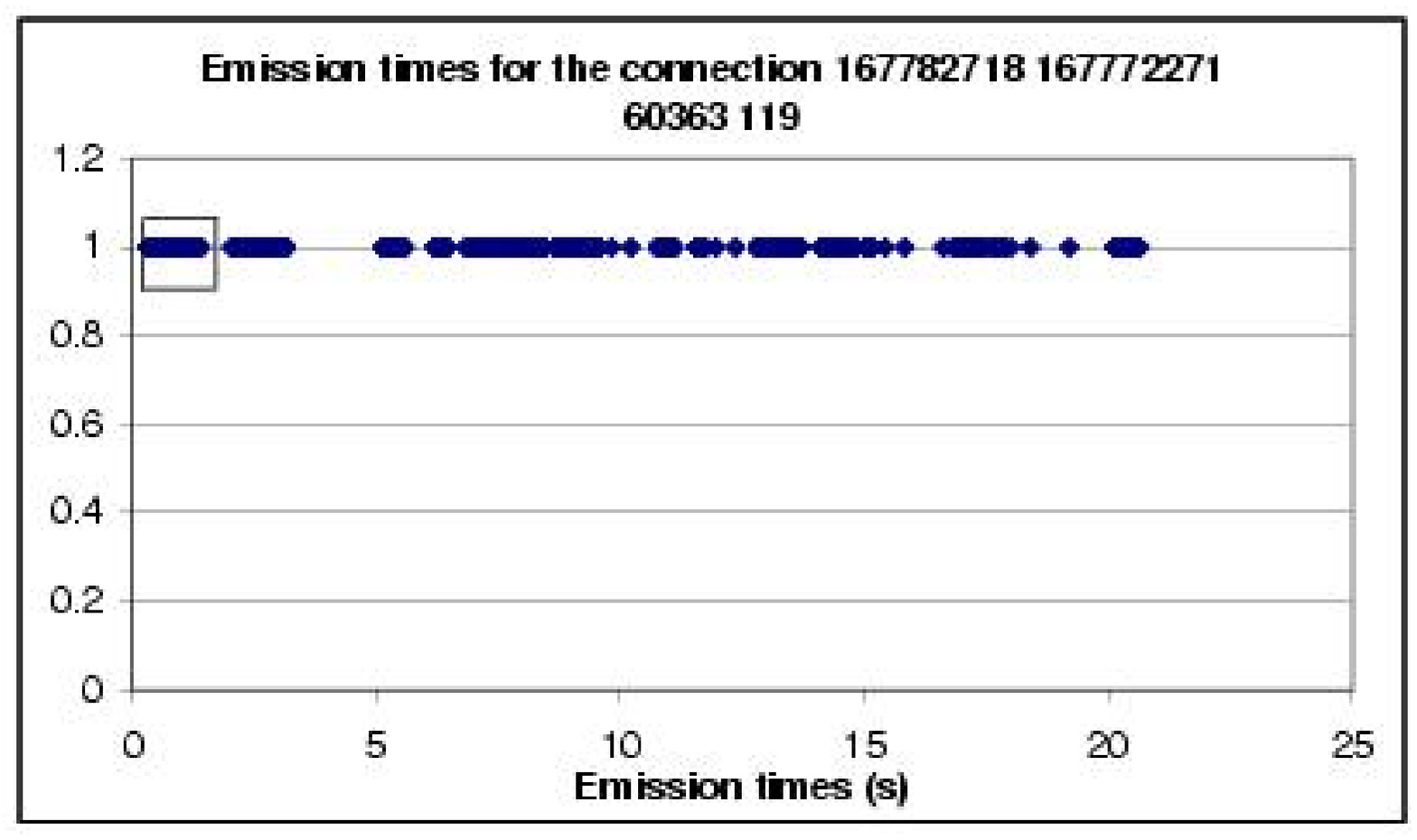}}
 \subfigure[]{\includegraphics[height=150pt, width=200pt]{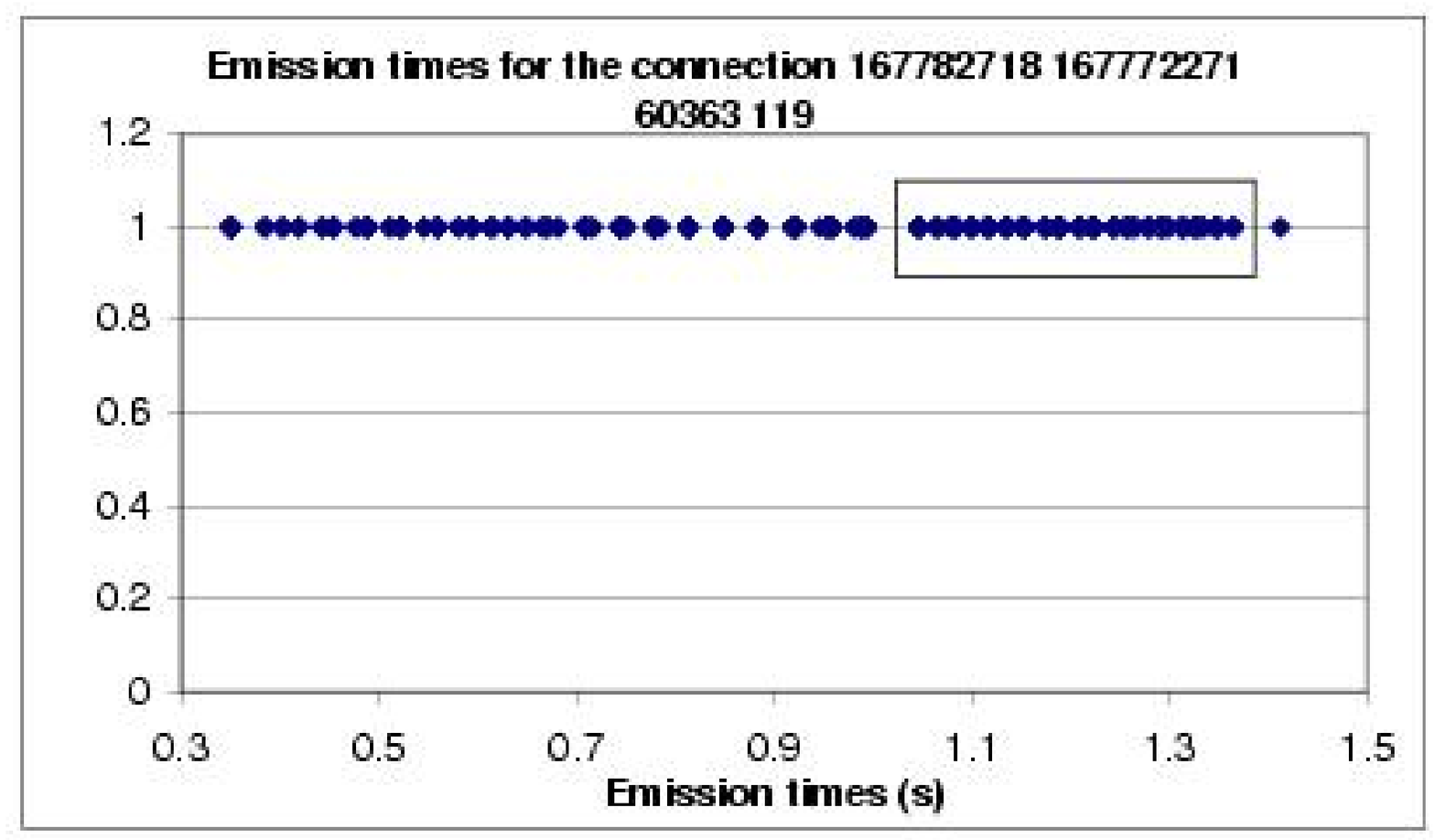}}
 \subfigure[]{\includegraphics[height=150pt, width=200pt]{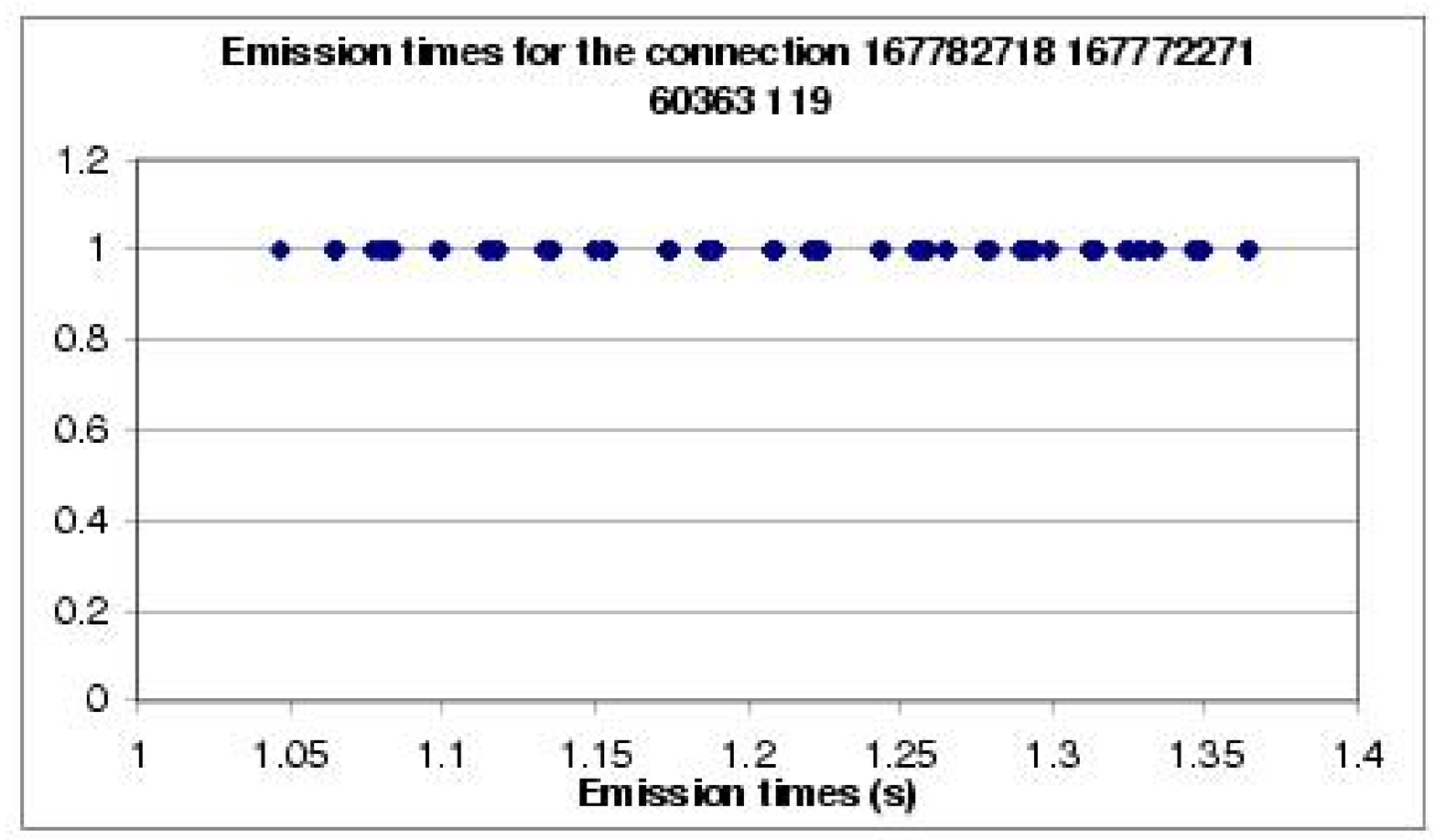}}
 \caption{The characteristic function of connection 167782718 167772271 60363 119 from 
          trace M6, in 3 detail levels. Here, 3 levels can be detected: the ``approximate'' length of their $1$-intervals is $2$, $0.44$, and $0.01$, whereas the corresponding $0$-intervals are much smaller. Notice that
          coarser levels are not present in the data set, due to its short duration. Notice also that $0$s are not plotted.}
\label{LevND}         
\end{figure}
\addtocounter{figure}{-1}
\stepcounter{figure}

\subsection{Detecting levels quantitatively in user sessions: the Interval Detection Algorithm (IDA)} 

\label{quant}

The purpose of this experiment is to provide quantitative information on the levels present in a session. The tool developed herein bears many similarities to Tools 1 and 2 of section \ref{discussion}, providing similar results, but by following a different (and much simpler) method, based on the idea of detecting the lengths of periods of activity and inactivity.

We start by describing the Interval Detection Algorithm; we will come back to its motivation below. 

Take a user session and bin it using an initial time bin, but in a special way: omitting packet sizes, assign a value of 1 to every bin in which a packet emission took place, and 0 otherwise. Subsequently, choose a base $b$, and run on the session the \emph{Interval Detection Algorithm}:

\begin{enumerate}
\item Part 1:

\begin{enumerate}
	\item Record the lengths of all intervals of value 0 (\emph{gaps}) in the data set.
	\item For every $i\geq 0$, find the sum $s(i)$ of the lengths of gaps $G$ that satisfy $[\log_b(|G|)]=i$. 
\end{enumerate}

\item Part 2:

For all $i\geq 0$ such that $i\leq \left\lfloor \log_b(\text{Session Length})\right\rfloor+1$
\begin{enumerate}
	\item Invert the session bin values logically.
	\item Run Part 1.
	\item Store the result as the $i$th column of the array $a(j,i)$ (we will call these results the $i$th stage of the detection of $1$-intervals of levels or simply the $i$th stage of the IDA).
	\item Invert the session bin values logically.
	\item Assign value 1 to all bins of the gaps $G$ that satisfy $[\log_b(|G|)]\leq i$.
	\item Set $w(i)$ to be the session values sum. 
\end{enumerate}

\item Part 3:
Choose $0< \gamma <1$. Then, for all $j\geq 0$
\begin{enumerate}
	\item Compute $M=\max_i(a(j,i))$.
	\item Set $i$ equal to the maximum possible value, and start decreasing it, until $a(j,i)<\gamma M$, but $a(j,i-1)\geq \gamma M$. Set $I=i$.
	\item For all $i$, set $a(j,i)\leftarrow a(j,i)/w(I)$.
	\item Set $s(i)\leftarrow s(i)/w(I)$.
\end{enumerate}

\item Part 4 (Representation of the results):
\begin{enumerate}
	\item 
	Let $im$ be an array. Choose numbers $c_1>c_2$. Then, for all $i\geq 0$:
\begin{enumerate}
	\item Find the maximum of the $i$th column of $a$, $m_1$.
	\item Find the second maximum of the $i$th column of $a$, $m_2$.
	\item If $m_1>c_1m_2$ and $m_2>c_2m_1$ then $\forall j: a(j,i)\leftarrow a(j,i)/m_2$ else $\forall j: a(j,i)\leftarrow \min(a(j,i)/m_2,1)$.
\end{enumerate}

\item Repeat the three steps above for $s$ (\emph{optional}).
\item $\forall i,j$, let $im(j,i)=a(j,i)$.
\item Append to $im$ a column of $0$s (\emph{optional}).
\item Append $s$ to $im$ as a column.
\item \label{GraphIDA} Mesh $im^T$ in greyscale: the higher a value of $im$ is, the darker the point representing it (this will be the \emph{greyscale representation}).
\item Set $v_1$ to be a vector, whose $i$th element is the sum of the elements of the $i$th row of $im$, excepting its last column (the $0$-interval stage). 
\item Set $v_0$ to be the last column of $im$. 
\item Normalize $v_1$ and $v_0$ so that they both have maxixmum value equal to $1$ and plot them (this will be the \emph{sum representation}). 

\end{enumerate}

\end{enumerate}

Within the framework of the greyscale representation, the IDA will indicate that a level for the $0$-intervals exists where $s(\cdot)$ has a high value, whereas a level for the $1$-intervals will exist where local maxima of $a(\cdot,i)$ align for many different $i$. On the other hand, within the framework of the sum representation, levels for the $1$- and the $0$-intervals will be assumed to exist at the local maxima of $v_1$ and $v_0$ respectively. 

Notice that the two representations are not equivalent, since the sum representation does not show the individual stages of the $1$-interval detection. However, as far as level detection is concerned, they can be considered equivalent: if many stage maxima cluster at a particular time scale, then the sum of the stage values at this particular time scale will tend to be rather large, and vice versa. Therefore, the two representations can be used interchangeably for our purposes.  

To explain how the IDA gives us this information, we first describe our way to simulate multi-scale sessions similar to the real sessions appearing in Fig. \ref{Lev} and \ref{LevND}, and see how the IDA allows us to observe their structure. The following is our simulation strategy:  

\begin{enumerate}
\label{LevAl}
	\item Generate a vector of isolated $1$s separated by intervals of $0$s; these intervals represent the RTTs and should be independent and identically distributed. 
	\item Generate a vector of alternating intervals of value 1 and 0. These intervals should also be identically and independently distributed, but also cross-independent. The two generating distributions, though (for intervals of value 1 an 0), need not be the same. 
	\item Repeat the previous step as many times as necessary, by choosing distributions with progressively smaller mean. Simulated traffic with $n$ levels will result, if the step is repeated $n$ times.  
	\item Multiply the vectors obtained above.
\end{enumerate}
This strategy contains the main idea behind levels, that coarse $1$-intervals get recursively subdivided into finer and finer $1$- and $0$-intervals. 

In order then to detect the structure of a user session (whether simulated or real), it will suffice to parse this user session only once in order to find the (unnormalized) histogram of the $0$-interval lengths, which is exactly what Part 1 of the IDA does. The detection of $1$-interval lengths will be trickier, however, since such intervals are not immediately present in the session: the $0$-intervals of finer levels fragment them into small pieces, and their continuity has to be restored before they can be detected. Part 2 takes care of this by computing sequentially the histogram of the $1$-interval lengths, after ``filling in'' larger and larger $0$-intervals.  

Finally, the IDA deals with the normalization of the histograms it produces. Making a traditional histogram by counting how many interval lengths can be logarithmically rounded to a certain time scale, and turning the histogram into a probability density by dividing with the total number of intervals that contributed to it will not do in this case, since the coarser the intervals, the more rare they are, but this does not mean they become less important, as this normalisation would suggest. The normalization must reveal the strength of the presence of a level within the session. One way to do this is to add the lengths of all intervals that contribute to a certain entry of the histogram and divide by the session length, which is equivalent to finishing the IDA in Part 2. But this will bias in favor of the coarse levels a lot, since a level does not really live on the whole session, but on the immediately coarser level only. Thus, the best normalization will be obtained by dividing not by the session length, but by the sum of $1$-interval lengths of the immediately coarser level, on which the level in question resides. An easy way to determine this normalization factor is to observe when the values of the unnormalized histograms corresponding to the time scale in question become negligible: if, after filling the gaps with length approximately equal to $b^j$, the histogram values corresponding to time scale $b^i$ become negligible, the proper normalization is $\forall k, a(i,k)\leftarrow a(i,k)/w(j)$ and $s(i)\leftarrow s(i)/w(j)$, as Part 3 explains. This method is approximate, of course: its main drawback is that the normalization depends entirely on the $1$-intervals. However, it seems to work in practice, as we are about to see.  

The normalization issue discussed above suggests that the location of the maxima may be misleading, if one attempts to detect the $1$-intervals of the levels based on them: indeed, despite the improved normalization discussed above, our simulations show that the maxima of the different stages of the detection of $1$-intervals still span a wide numerical range. This is why we decided to suppress this information in the final presentation of the results of the IDA, and use other information instead, as suggested in Part 4 of the IDA algorithm.   

\begin{figure}
 \centering
 \subfigure[]{\includegraphics[height=150pt, width=200pt]{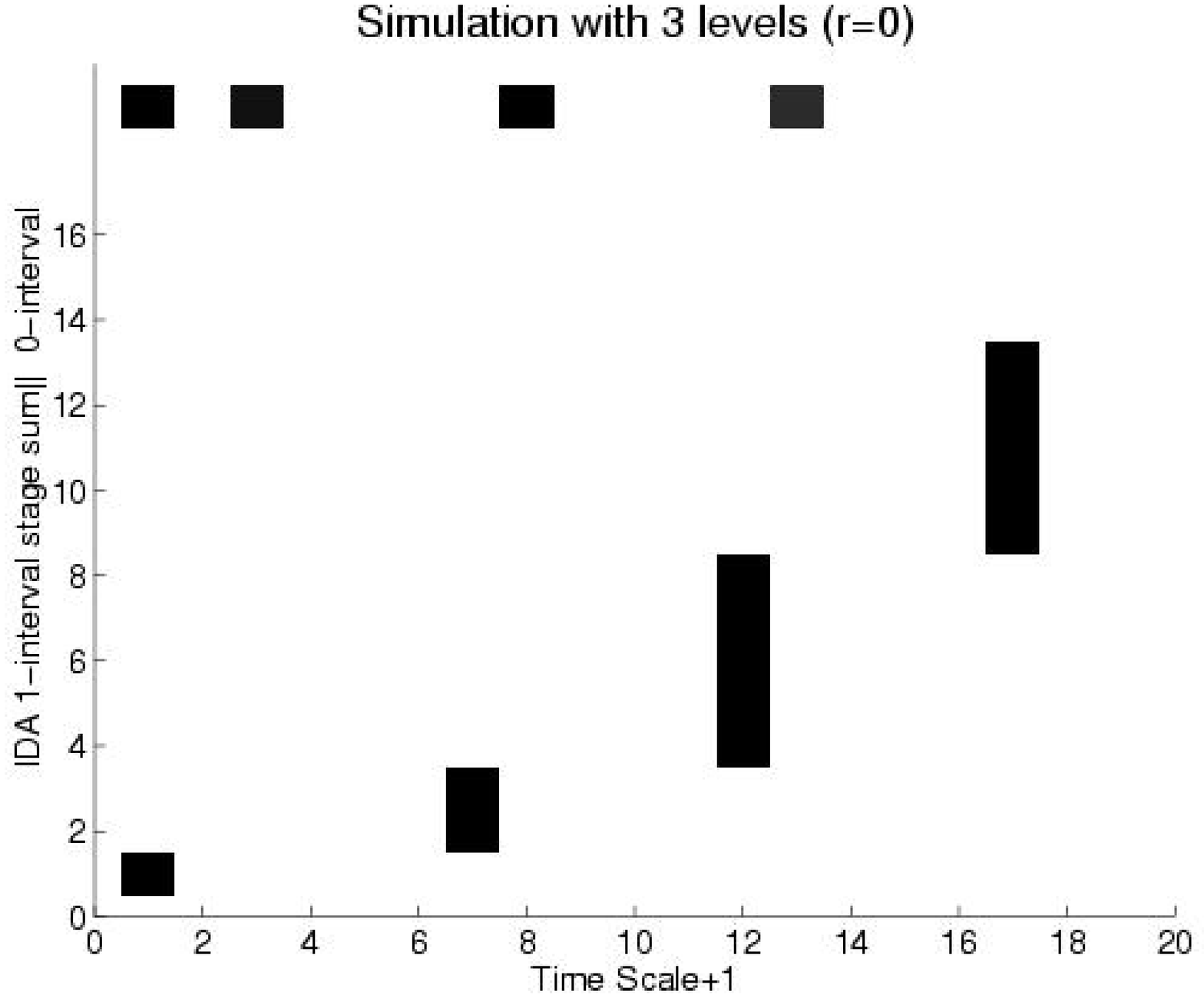}}
 \subfigure[]{\includegraphics[height=150pt, width=200pt]{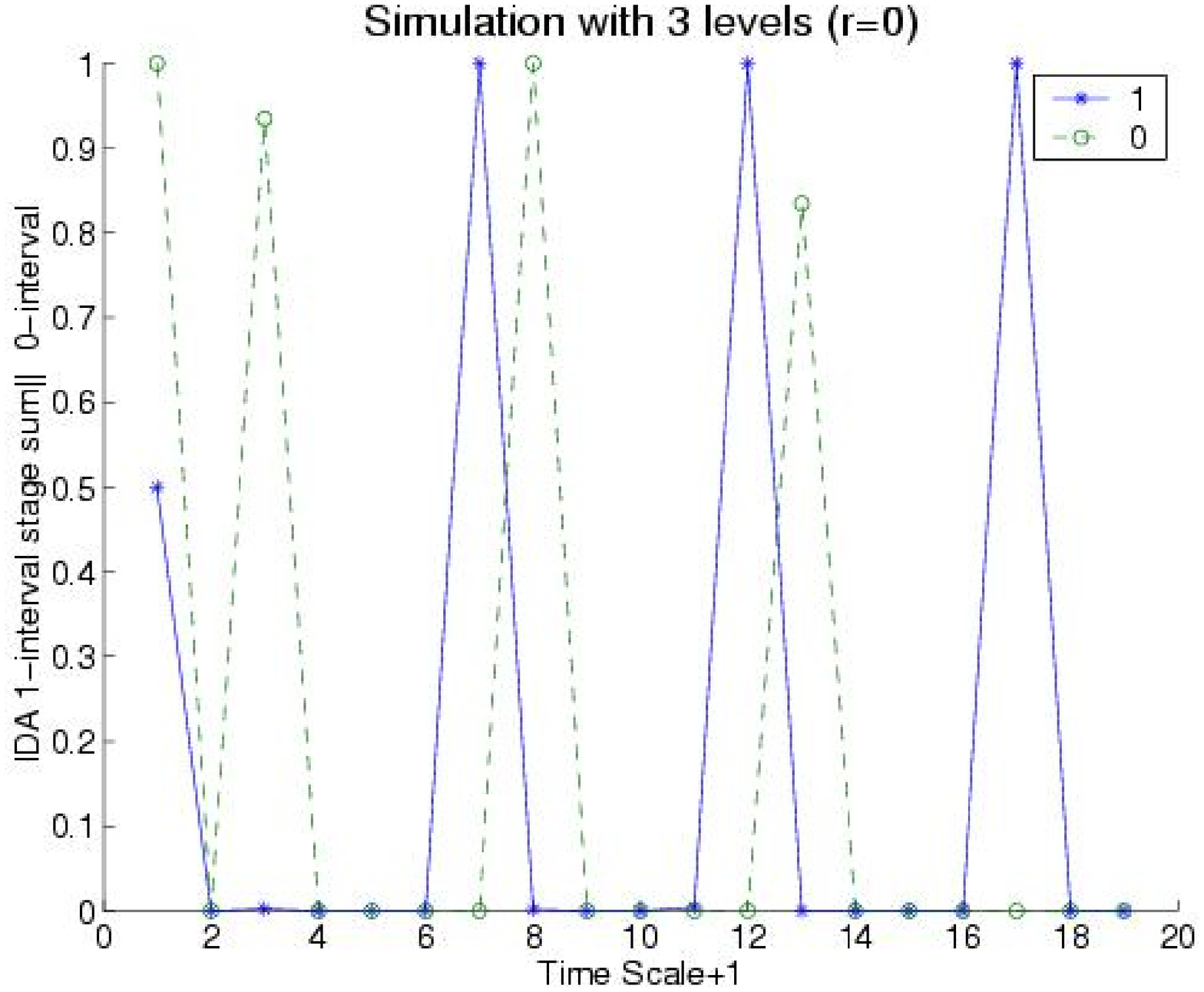}}
 \subfigure[]{\includegraphics[height=150pt, width=200pt]{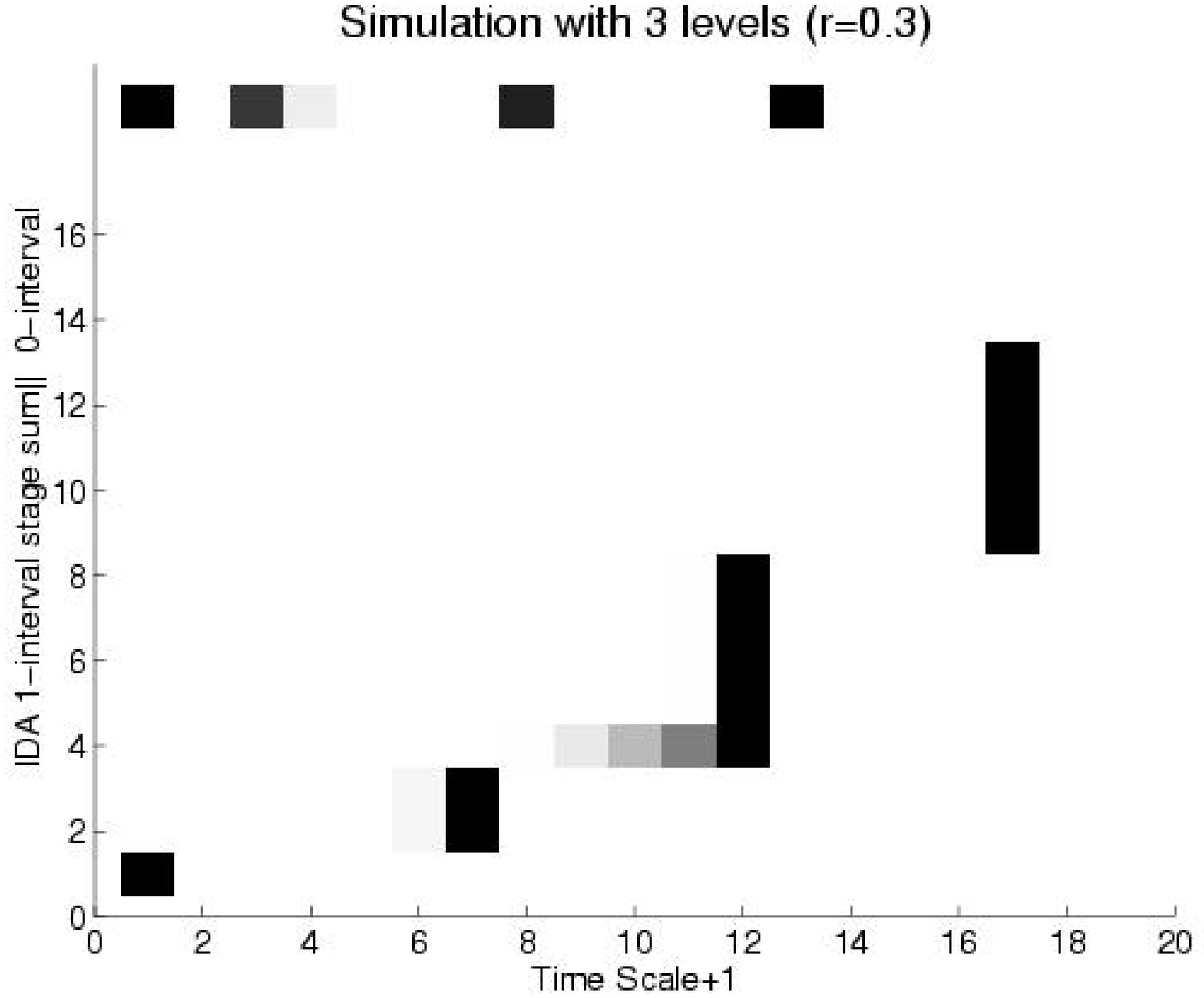}}
 \subfigure[]{\includegraphics[height=150pt, width=200pt]{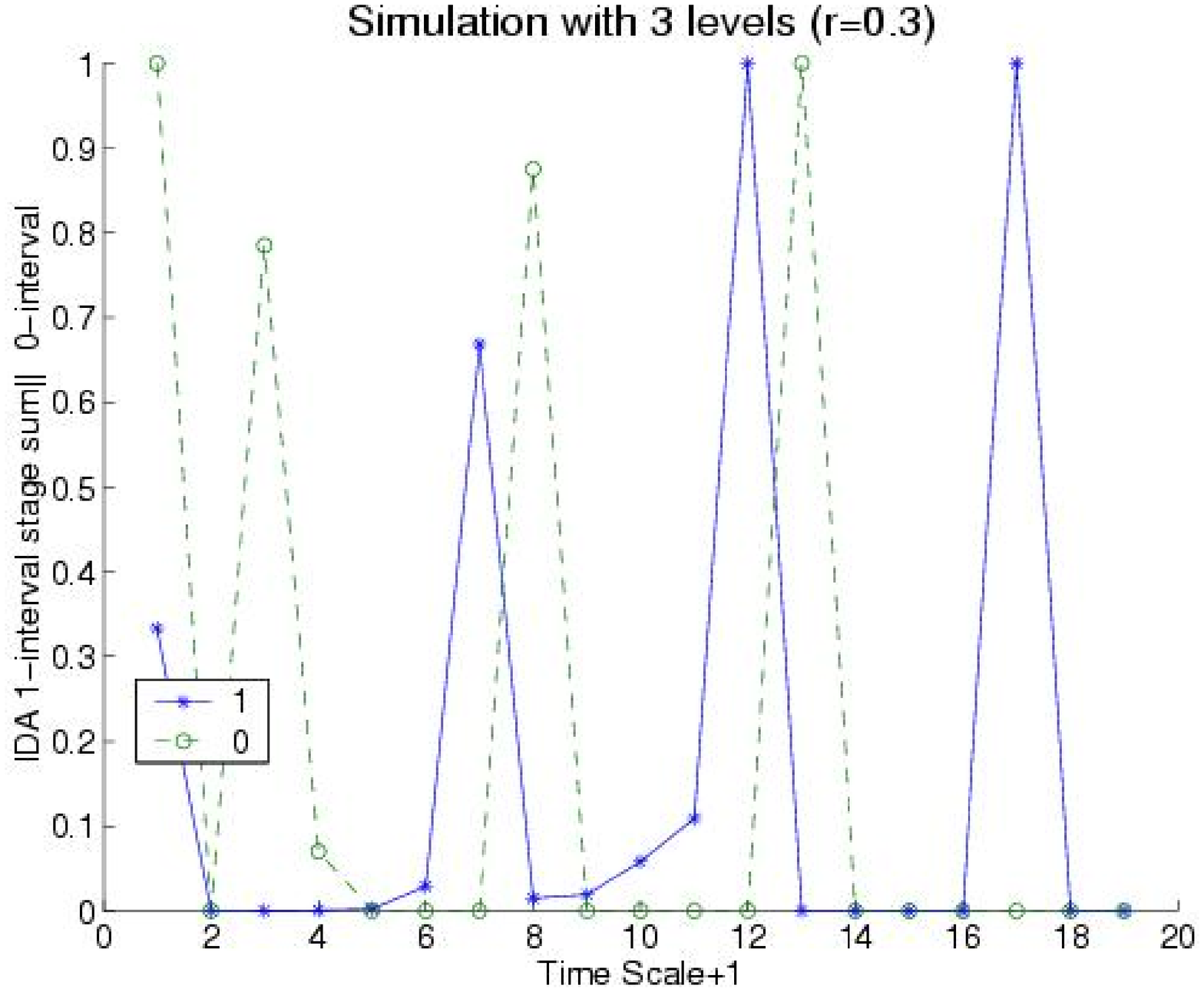}}
 \subfigure[]{\includegraphics[height=150pt, width=200pt]{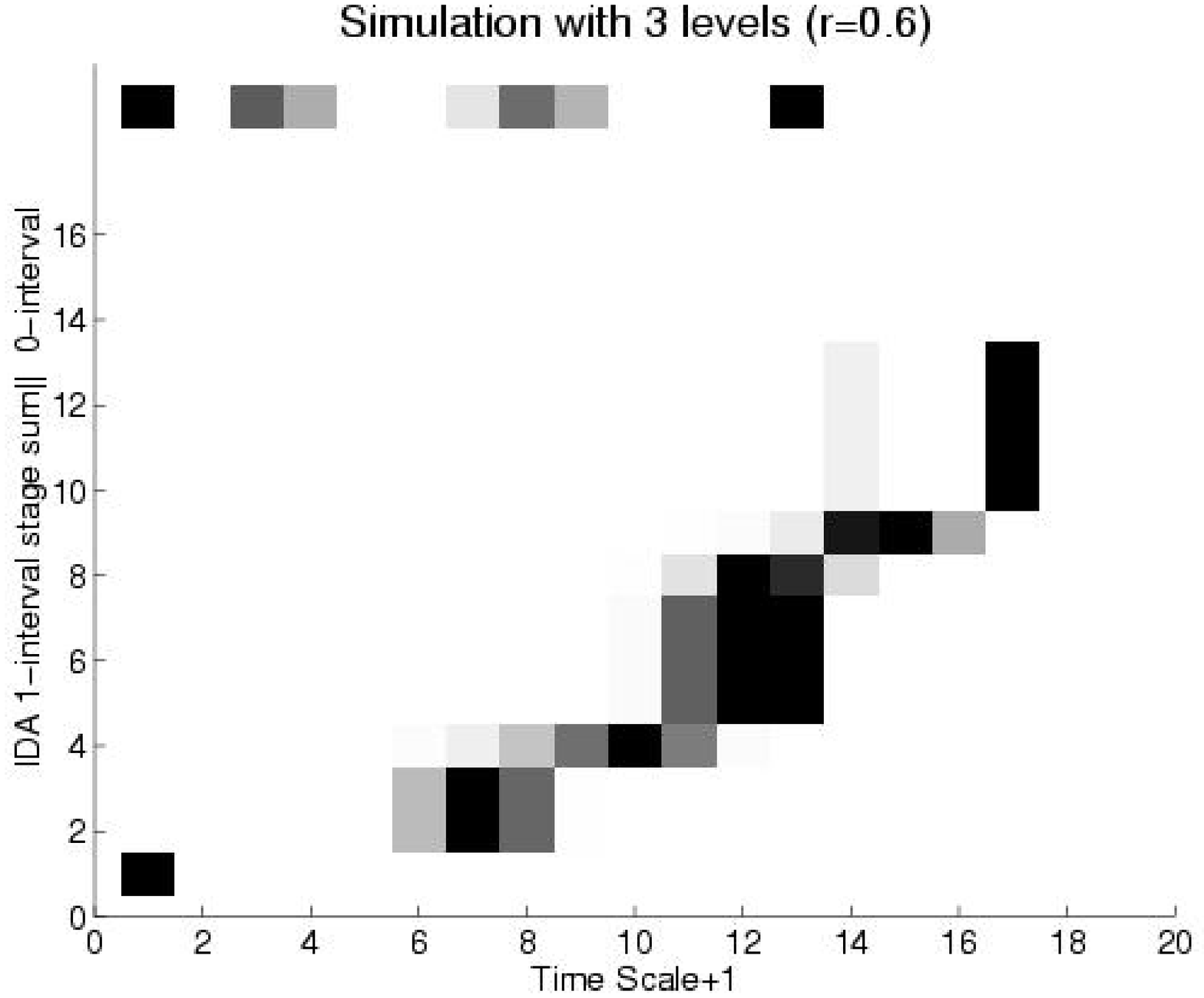}}
 \subfigure[]{\includegraphics[height=150pt, width=200pt]{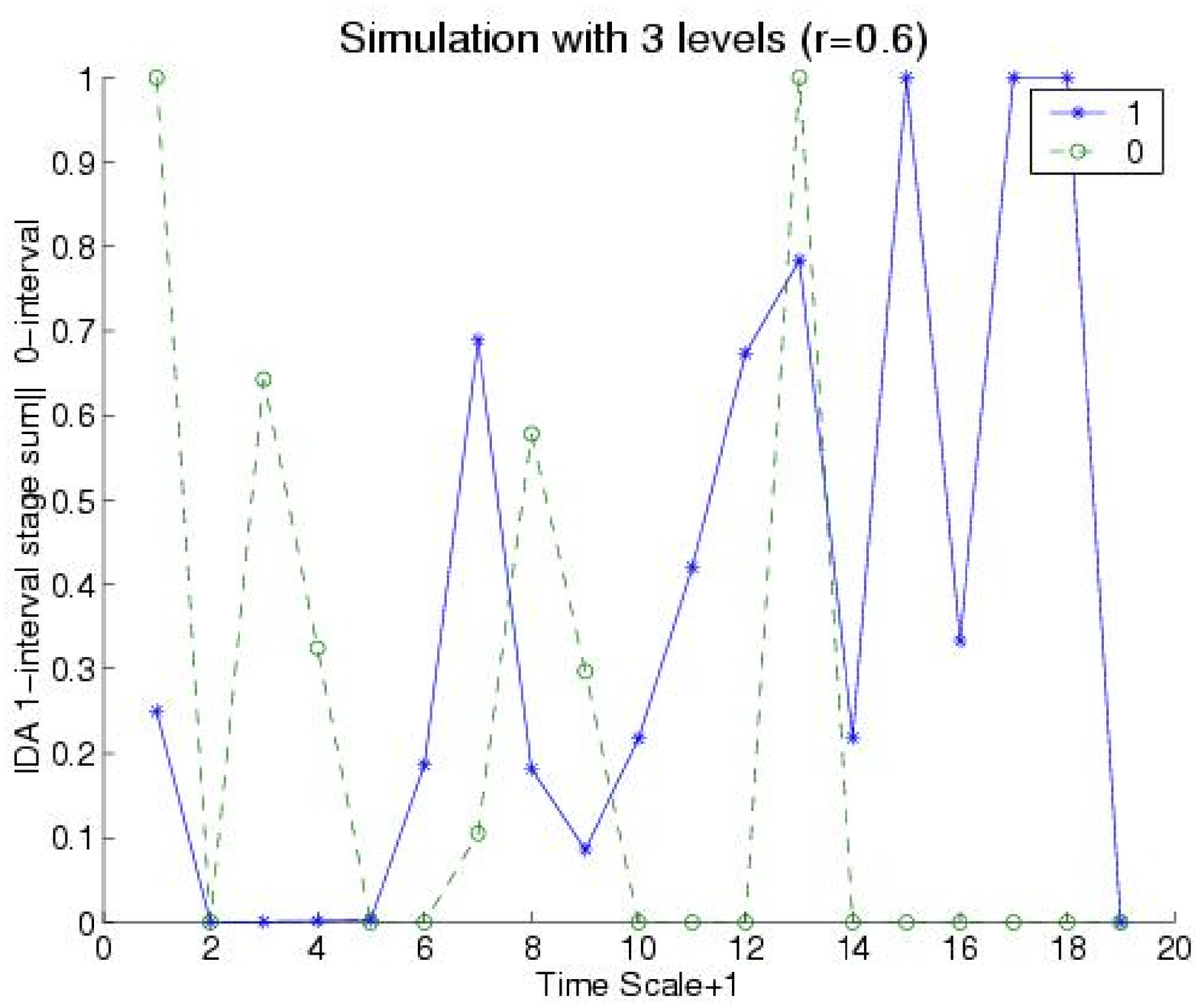}}
 \caption{Detection of $0$- and $1$-intervals of levels in simulated sessions (part I): the $1$-interval mean sizes are at $2^6$, $2^{11}$, and $2^{16}$, and the $0$-interval mean sizes are at $2^2$, $2^7$, and $2^{12}$. The intervals are uniformly distributed around their mean, and the width of variation over the mean is is $0$ in (a), (b), $0.3$ in (c), (d), and $0.6$ in (e), (f) . Notice that after filling the gaps of the coarsest level, the IDA considers the whole data set as an $1$-interval, hence there should be an indication of a level at time scale $18$; this artifact has been removed from all pictures, except (f), where it was left for demonstration purposes.}
\label{LevelReader1}         
\end{figure}

As the algorithm evolves and gaps of the data set get filled, the maxima of the stages will be moving towards larger scales. In the case of simulations with clearly separated and defined levels, this motion will be discontinuous: if the $1$-intervals of a level have approximate lengths of $b^i$ and the ones of the immediately coarser level $b^j$, then for stages immediately before $i$ the maximum will be at $i$, but for stages $i$ up to $j$ it will be at $j$: indeed, filling gaps of lengths between $b^i$ and $b^j$ will be of no effect, because the number of these gaps, by assumption, is not significant. The IDA then produces results such as in Fig. \ref{LevelReader1}(a), (b), where the sudden jumps and the long immobility periods of the maxima are very eloquently depicted (more information about Fig. \ref{LevelReader1} will follow shortly). 

In the case of simulations with diffused levels or of real data sets, though, the maximum will be moving more or less smoothly. When depicting the results of the IDA in the greyscale representation, this motion will appear as a diagonal line of maxima in the figure (see e.g. Fig. \ref{LevelReader2}(a), (c), and \ref{LevelReader2}(a)). However, there will still be time scales where maximum (or nearly maximum) values will seem to persist for more than one stages (such as time scales $12$ and $18$ in Fig. \ref{LevelReader2}(a), (b)), just as in the case of non-diffused levels, indicating the existence of a level there. Thus, the algorithm will not rely upon the maximal values of the stages for the detection of $1$-intervals, but will instead pick the time scales where large values of the stages tend to ``cluster'', as Part 4 of the IDA description states.   

In order to test the performance of the IDA, we applied it on six simulated sessions, generated by the session simulation algorithm presented above; the results are shown in Fig. \ref{LevelReader1} and \ref{LevelReader2}. Each one of the six simulations, which have a behavior very similar to Fig. \ref{Lev} and \ref{LevND}, features 3 levels, whose $1$-intervals have mean size $2^6$, $2^{11}$, and $2^{16}$, and whose $0$-intervals have mean size $2^2$, $2^7$, and $2^{12}$, i.e. they are $16$ times smaller. The interval lengths are uniformly distributed around their mean in Fig. \ref{LevelReader1} and Fig. \ref{LevelReader2}(a), (b), and the ratio of the allowed variation over the mean is $0$ in Fig. \ref{LevelReader1}(a), (b), $0.3$ in Fig. \ref{LevelReader1}(c), (d), $0.6$ in Fig. \ref{LevelReader1}(e), (f), and $0.9$ in Fig. \ref{LevelReader2}(a), (b). The levels in Fig. \ref{LevelReader2}(c), (d) are exponentially distributed around their mean, and levels in Fig. \ref{LevelReader2}(e), (f) follow a Gaussian distribution, whose mean in each case coincides with the level mean, and whose variance is proportional to the mean: $\sigma^2=4\mu$. The IDA picks out the positions of the $0$- and $1$-intervals reasonably accurately, even when the variance becomes large: for scales (horizontal axis) at which $1$-intervals of levels exist, the greyscale column is much ``taller'' than at other scales; the scales at which $0$-intervals of levels exist are picked out by the top greyscale line. In particular, it is always clear, except in the exponential case in Fig. \ref{LevelReader2}(c), (d), that there are 3 distinct levels. 

The Averaging functions corresponding to two of the above simulations (\ref{LevelReader1}(b) and \ref{LevelReader2}(a)) are shown in Fig. \ref{LevelReaderAv}. At this point, we would like to jumb ahead a bit, and ask the reader to notice the clear correspondence between the way the Averaging function curves, and the position of the levels, as captured by the IDA. Indeed, the following sections will demostrate the relation between these two phenomena. 

\addtocounter{figure}{-1}
\stepcounter{figure} 

\begin{figure}
 \centering
 \subfigure[]{\includegraphics[height=150pt, width=200pt]{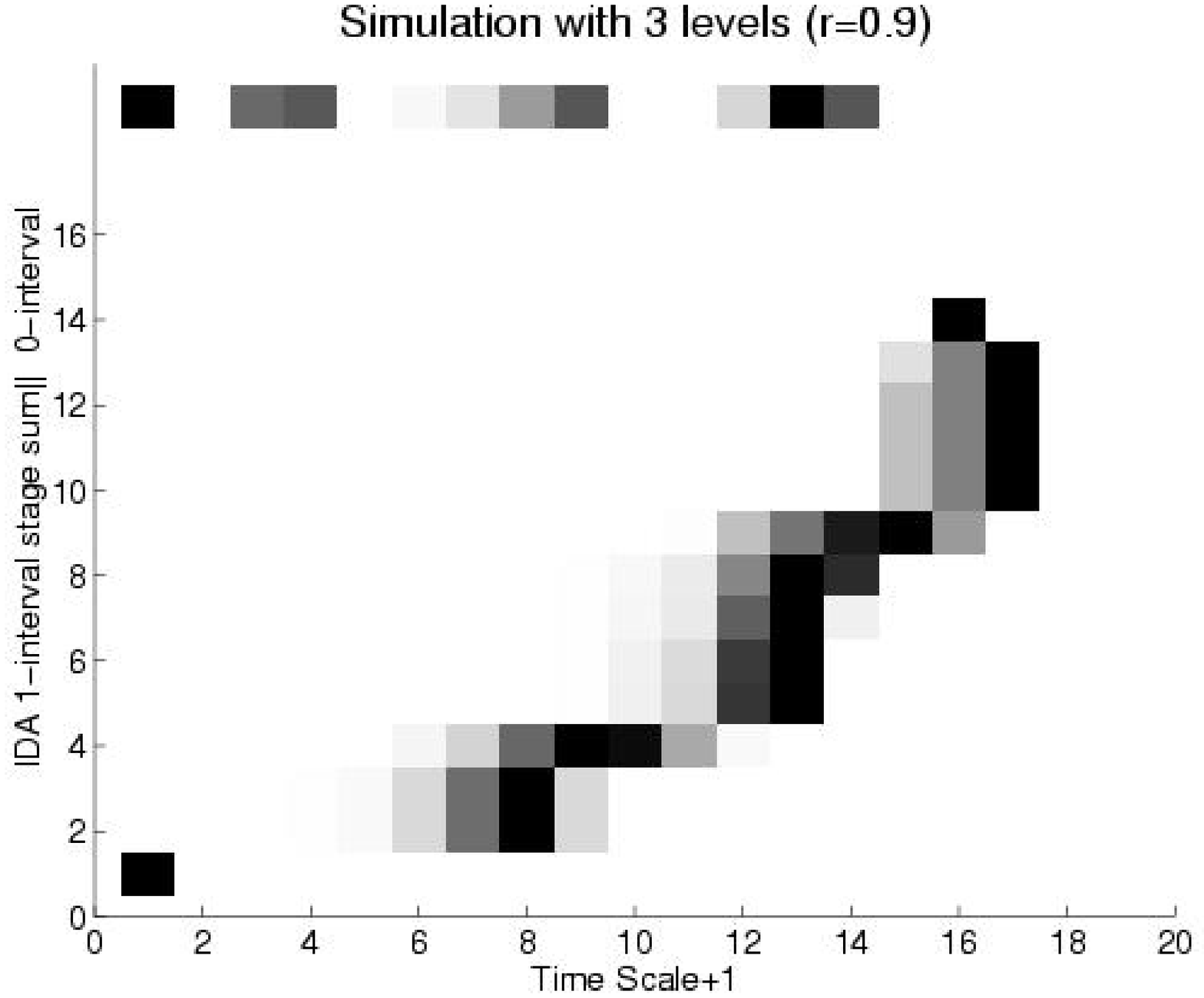}}
 \subfigure[]{\includegraphics[height=150pt, width=200pt]{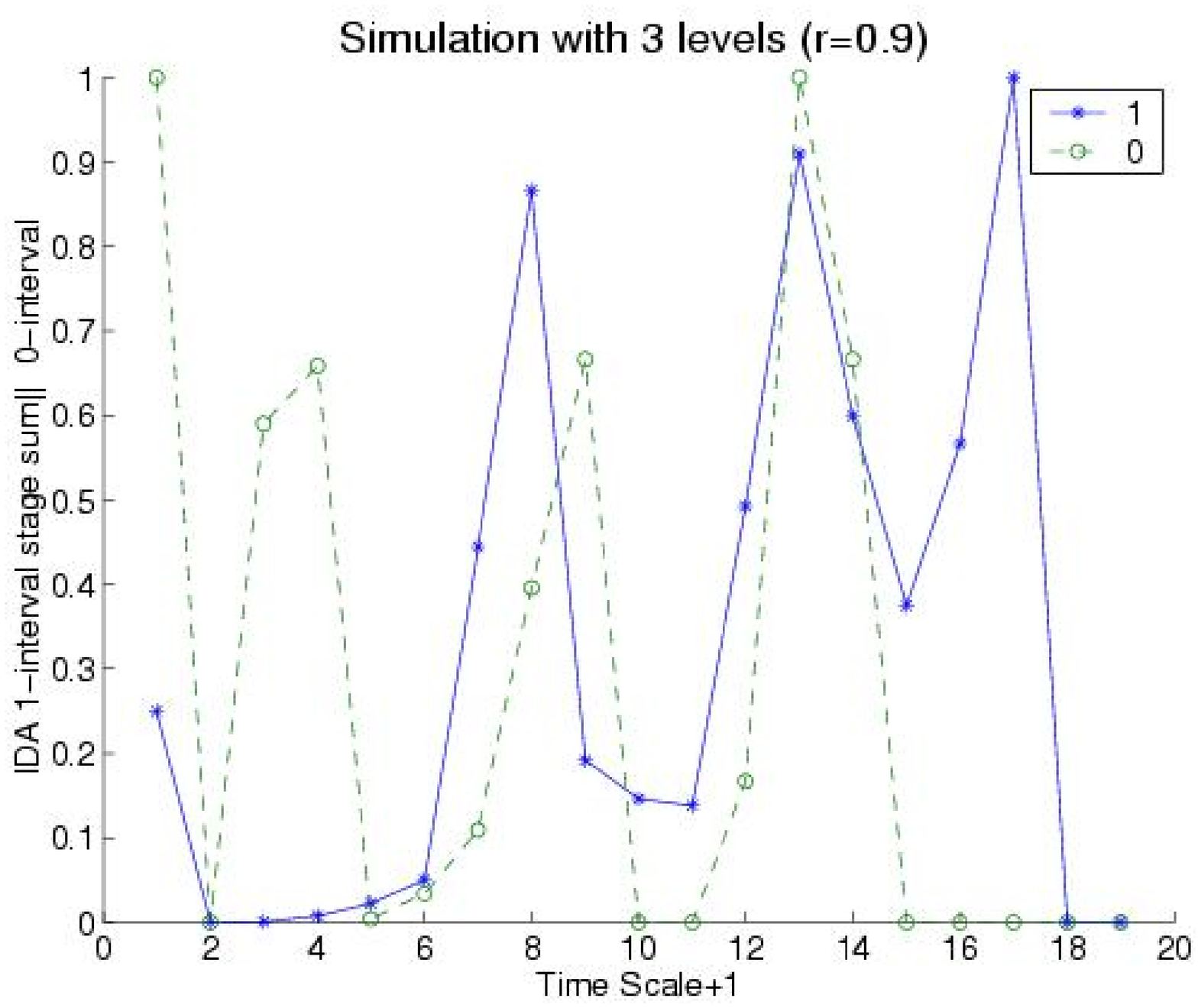}}
 \subfigure[]{\includegraphics[height=150pt, width=200pt]{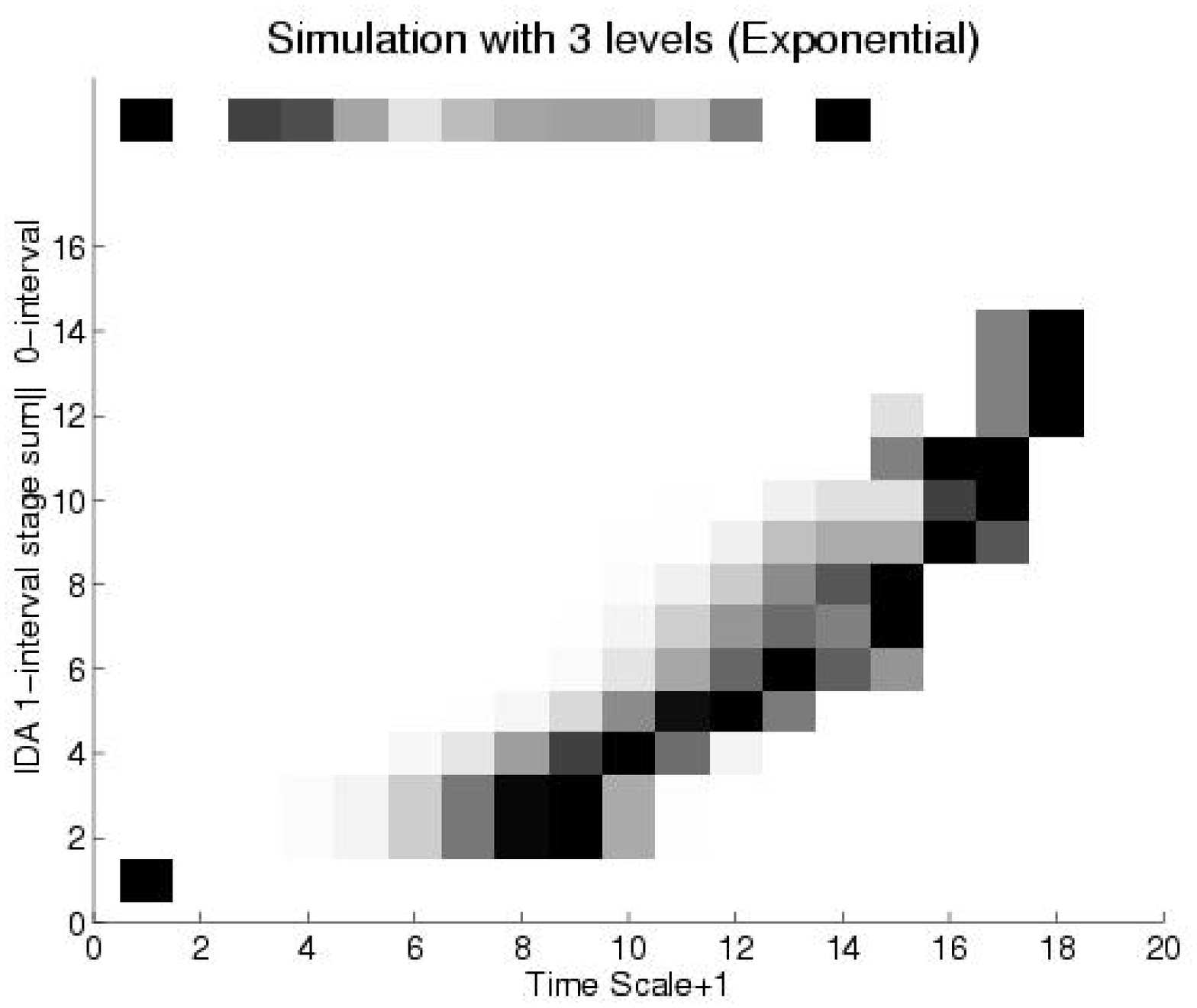}}
 \subfigure[]{\includegraphics[height=150pt, width=200pt]{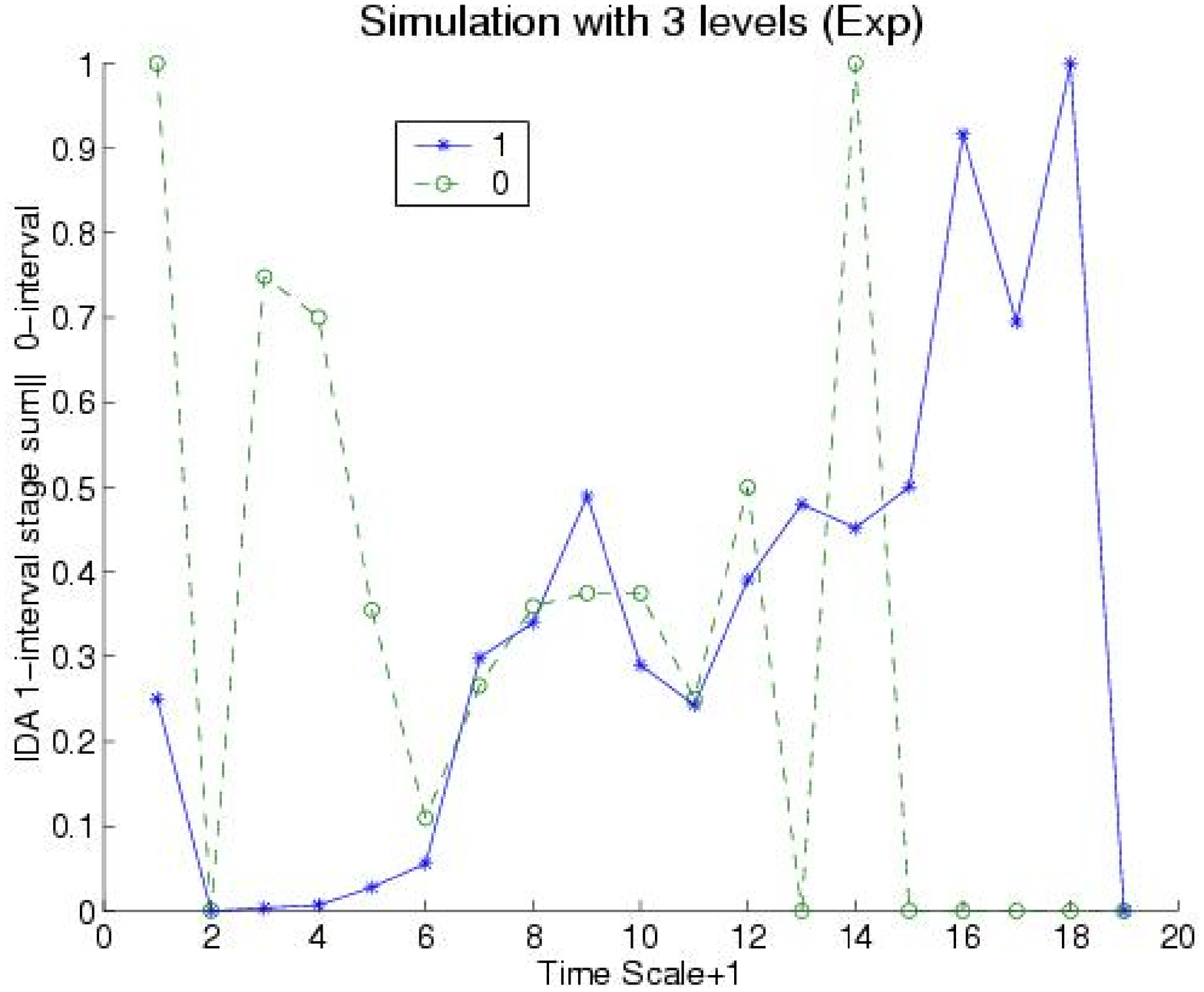}}
 \subfigure[]{\includegraphics[height=150pt, width=200pt]{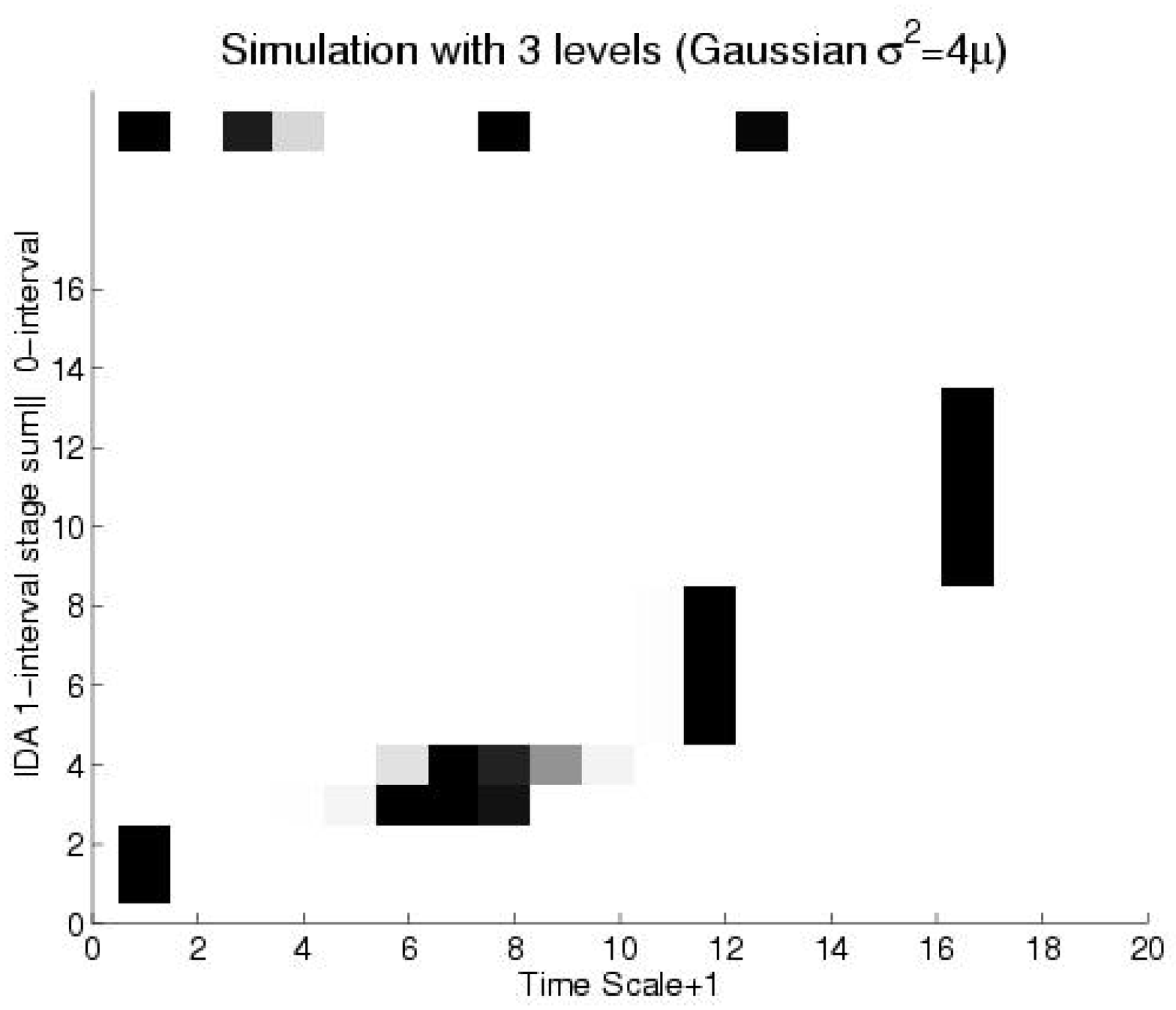}}
 \subfigure[]{\includegraphics[height=150pt, width=200pt]{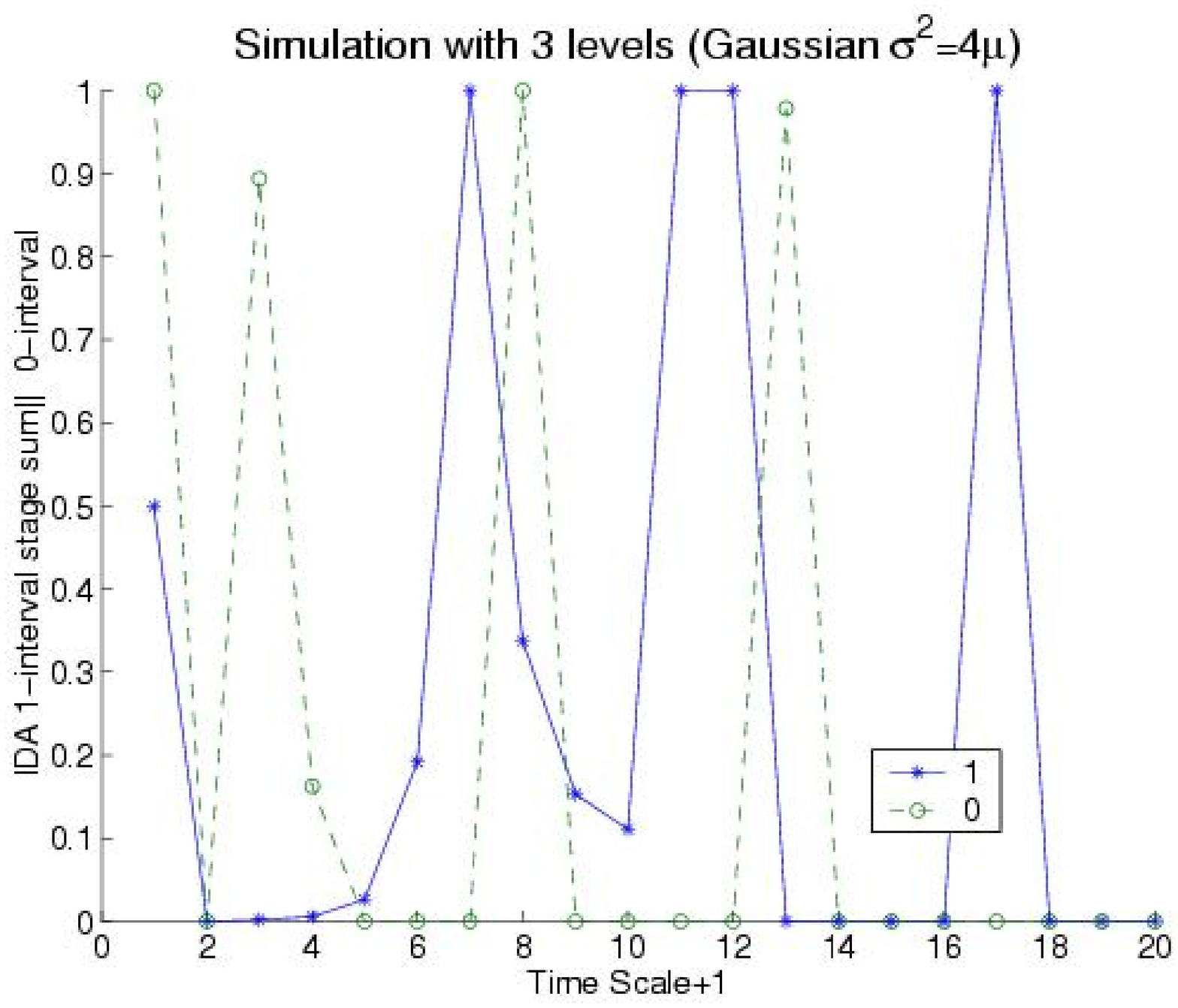}}
 \caption{Detection of $0$- and $1$-intervals of levels in simulated sessions, part II: the $1$-interval mean sizes are at $2^6$, $2^{11}$, and $2^{16}$, and the $0$-interval mean sizes are at $2^2$, $2^7$, and $2^{12}$. In figures (a) and (b) the intervals are uniformly distributed around their mean, and the width of variation over the mean is $0.9$. In (c) and (d) the intervals are exponentially distributed, whereas in (e) and (f) they are gaussianly distributed, the variance for each level being proportinal to the mean: $\sigma^2=4\mu $. Notice that after filling the gaps of the coarsest level, the IDA considers the whole data set as an $1$-interval, hence there should be an indication of a level at time scale $18$; this artifact has been removed from all pictures, except (c) and (d), where it was left for demonstration purposes.}
\label{LevelReader2}         
\end{figure}
\addtocounter{figure}{-1}
\stepcounter{figure} 

\begin{figure}
 \centering
 \includegraphics[height=150pt, width=200pt]{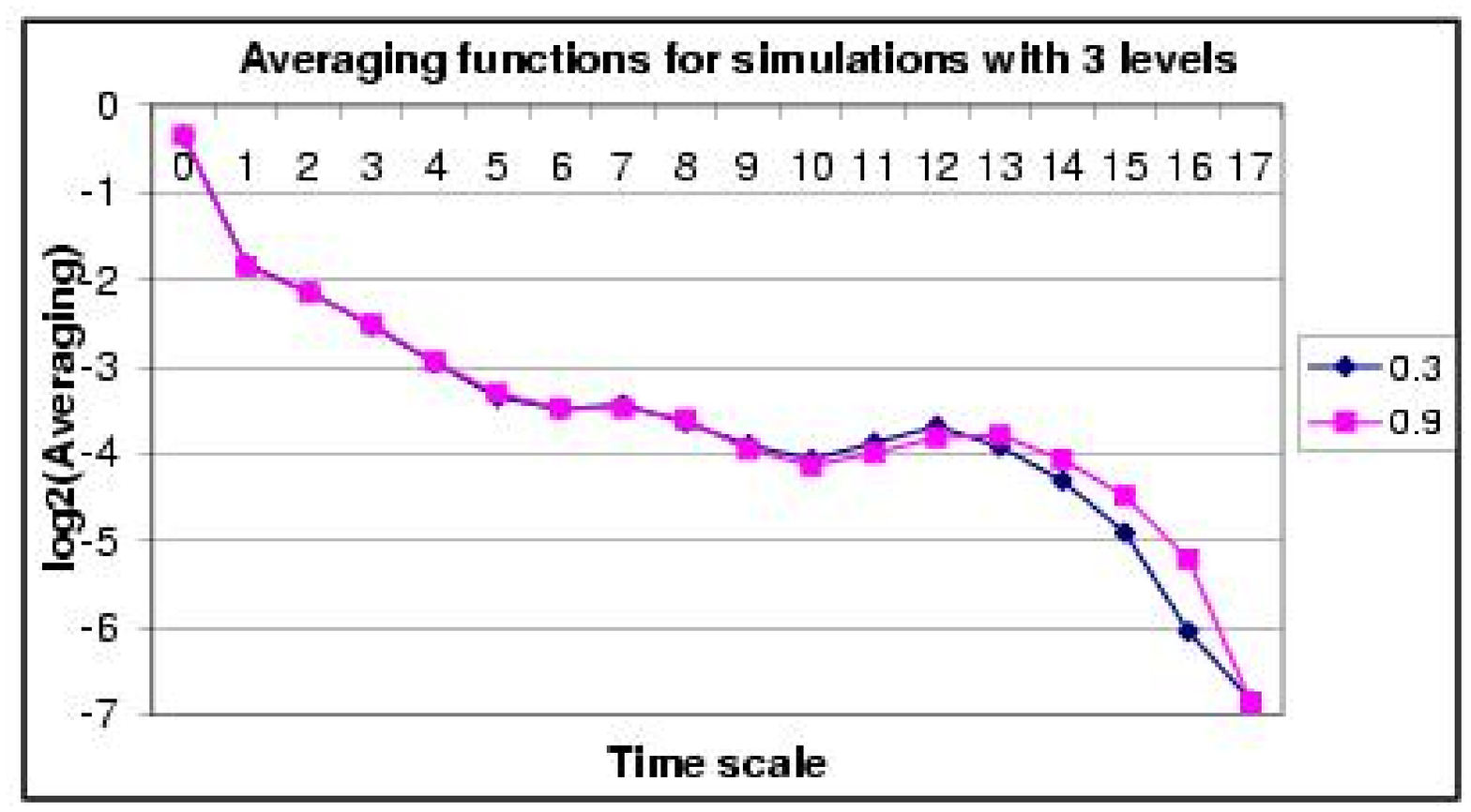}
 \caption{Averaging functions corresponding to the simulations of Fig. \ref{LevelReader1}(b) and \ref{LevelReader2}(a).}
\label{LevelReaderAv}         
\end{figure}
\addtocounter{figure}{-1}
\stepcounter{figure} 

\begin{figure}[t]
 \centering
 \begin{tabular}{cc}
 \subfigure[]{\includegraphics[height=150pt, width=200pt]{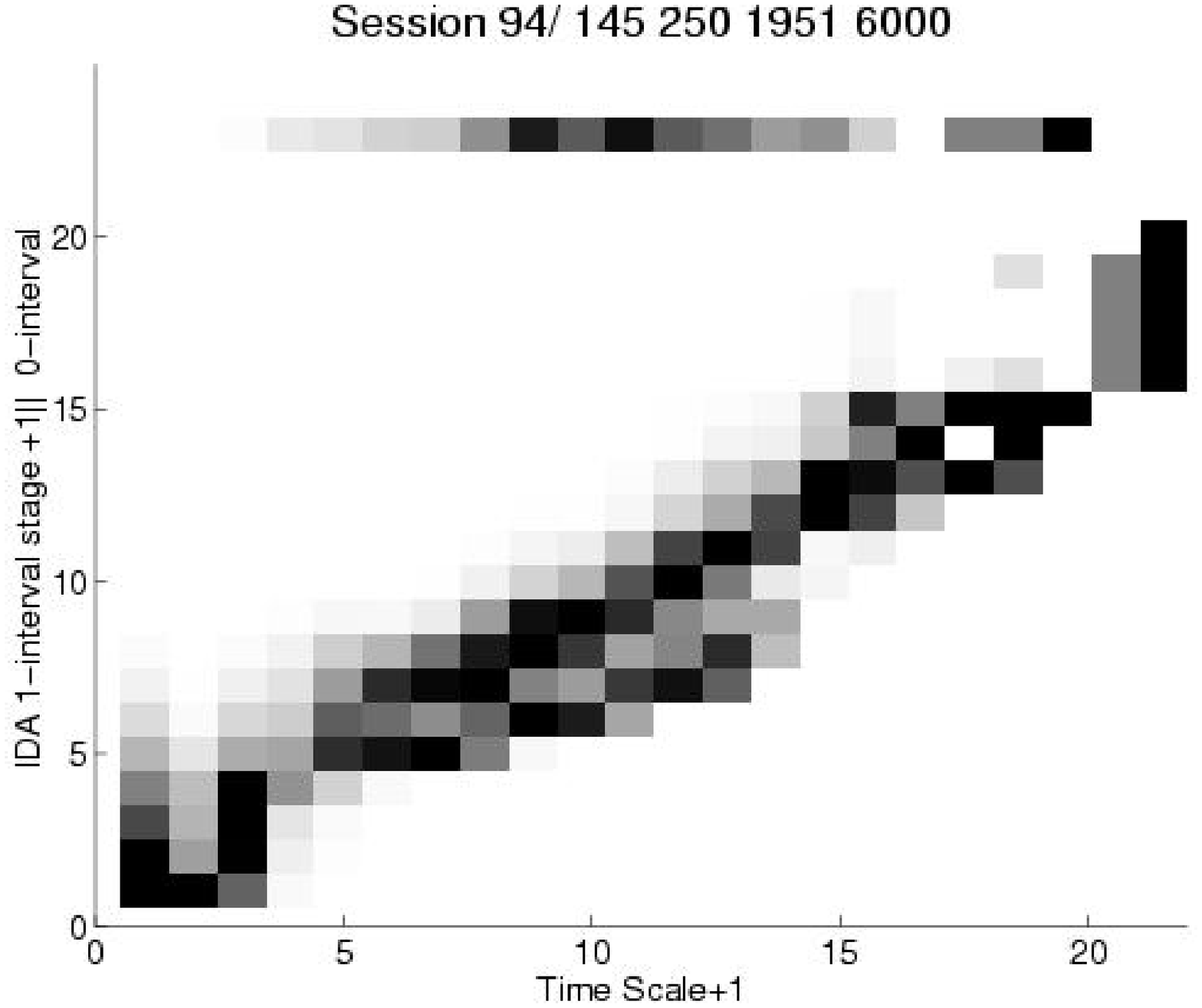}} &
 \subfigure[]{\includegraphics[height=150pt, width=200pt]{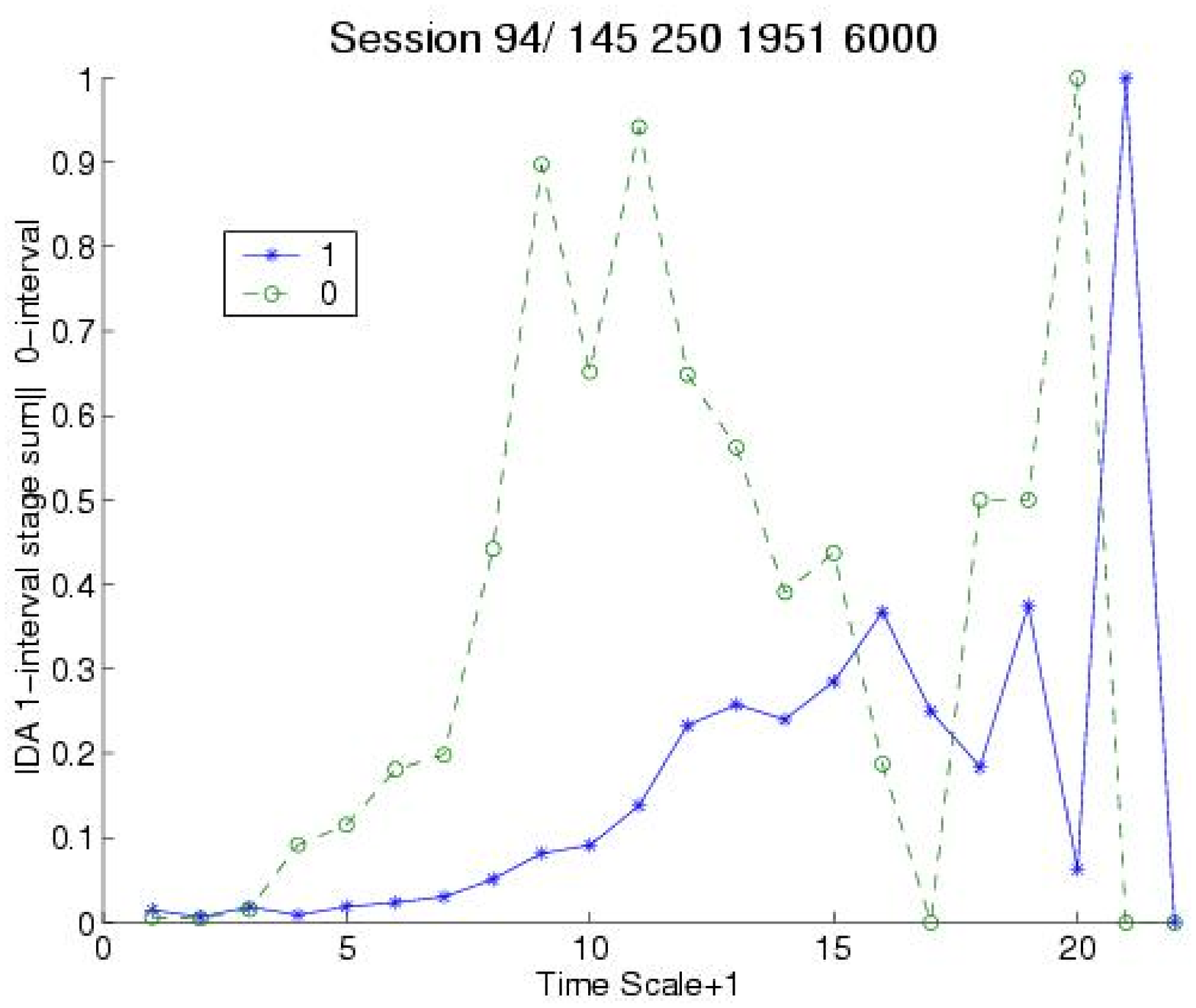}}\\
 & \subfigure[]{\includegraphics[height=150pt, width=200pt]{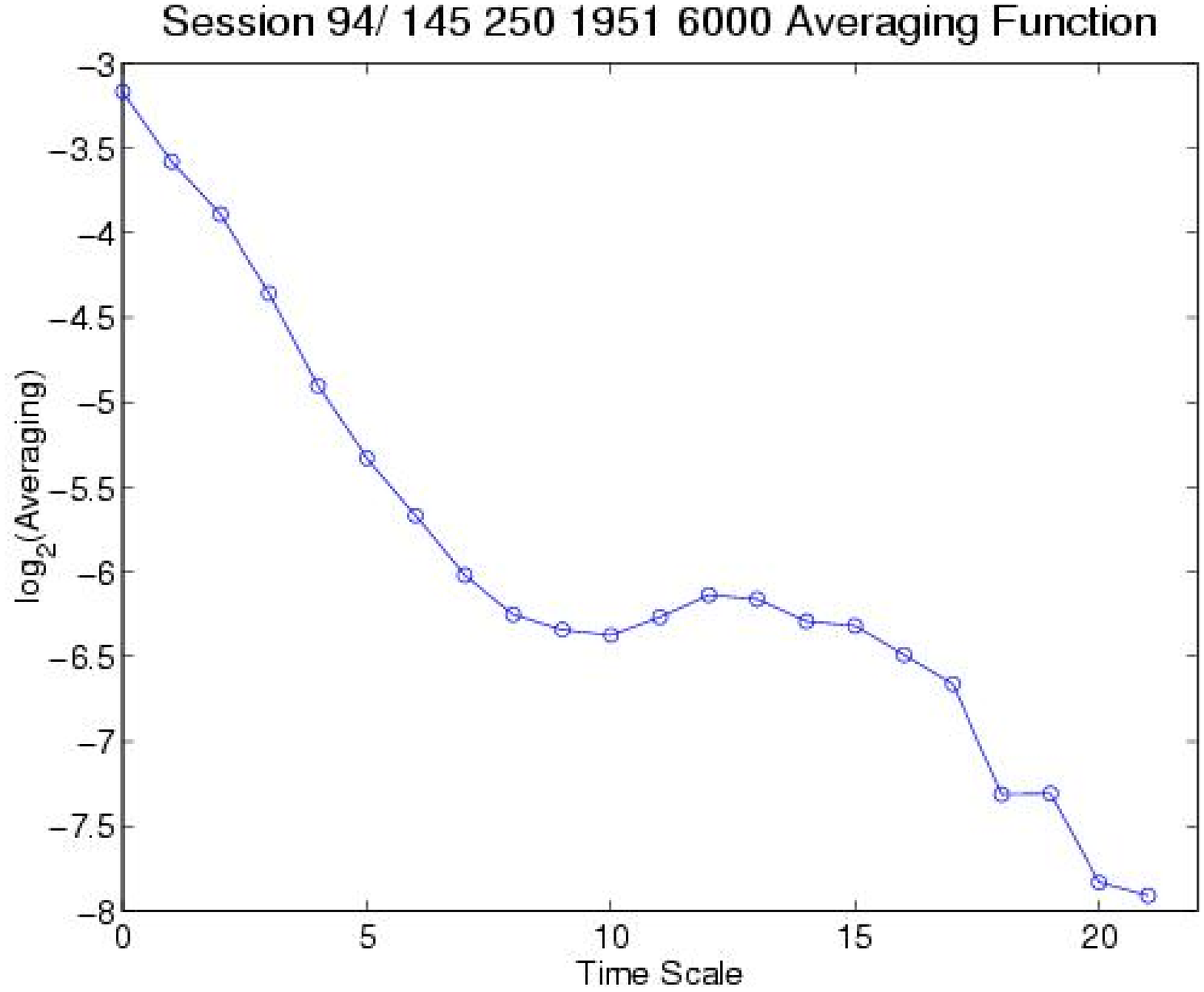}}
 \end{tabular}
 \caption{Detection of $0$- and $1$-intervals of levels in the longest session of trace 94: (a) and (b) are the greyscale and sum representation of the IDA results, respectively; (c) shows the Averaging function, which is included here in order to facilitate comparison (see section \ref{curvexplan}).}
\label{LevelReaderReal1}         
\end{figure}

Having thus assured ourselves of the efficiency of the IDA on simulated sessions, we use it in Fig. \ref{LevelReaderReal1} and  \ref{LevelReaderReal2} on two individual sessions taken from the data sets. In each case, the IDA identifies time scales that stand out. In particular,  regarding the longest session of trace 94 (Fig. \ref{LevelReaderReal1}), the IDA predicts that the $1$-intervals of levels have lengths between $2^{10}ms$ and $2^{20}ms$, whereas the $0$-interval lengths lie between $2^4ms-2^14ms\ \text{and}\ 2^{17}ms-2^{19}ms$, which means that this session appears to have levels all over the place. Jumbing ahead once more, let us note that the Averaging function seems to suggest the same thing (see the discussion in sections \ref{curvexplan} and \ref{metcomp} for the connection between Averaging function curving and levels), as the position of the ``bumps'' of the Averaging function (Fig. \ref{LevelReaderReal1}(c)) corresponds exactly to the position of the levels predicted by the IDA.

In two occasions above, a simulated and a real session, the IDA and the Averaging function where found to yield compatible results about where the $1$- and $0$-intervals of the session levels lie. Is then this always the case? Unfortunately, another real session, the one that was also used in Fig. \ref{Lev}, answers the question in the negative. 

For this session, the IDA detects $0$-intervals of levels in two areas (see Fig. \ref{LevelReaderReal2}(a), (b)): at $2^{18} ms - 2^{19}ms\approx 250s-500s$, and at $2^8ms-2^{14}ms\approx 0.25s-16s$. In Fig. \ref{Lev}, we had identified by mere observation levels at $400,\ 15,\ 10,\ \text{and}\ 4s$ which fits very well. On the other hand, the IDA detects $1$-intervals at $2^{15}ms \approx 30s$ and at $2^{18}ms-2^{21}ms\approx 250s-2000s$, whereas in Fig. \ref{Lev} we observe $1$-intervals with lengths of $1000s-2000s,\ 500s,\ \text{and}\ 40s$. Hence, visual inspection and the IDA yield almost identical results. But how does the IDA compare to the Averaging function (Fig. \ref{LevelReaderReal2}(c))? Although the IDA detects a very prominent level with $0$-intervals around $2^{13}ms$, the Averaging function seems to ignore it completely. Although this appears to be a potential failure for our model at first sight, it is not. We will come back to that issue in section \ref{uneq}, where we will discuss explicitly levels with unequal $1$- and $0$-intervals. 

As a general remark, compared to the simulations of Fig. \ref{LevelReader1}, the behavior of real traces seems to be more complex: actually, the percentage is almost nowhere equal to 0, which should mean that levels, even weak ones, exist in almost all time scales.  
   
\addtocounter{figure}{-1}
\stepcounter{figure} 

\begin{figure}[t]
 \centering
 \begin{tabular}{cc}
 \subfigure[]{\includegraphics[height=150pt, width=200pt]{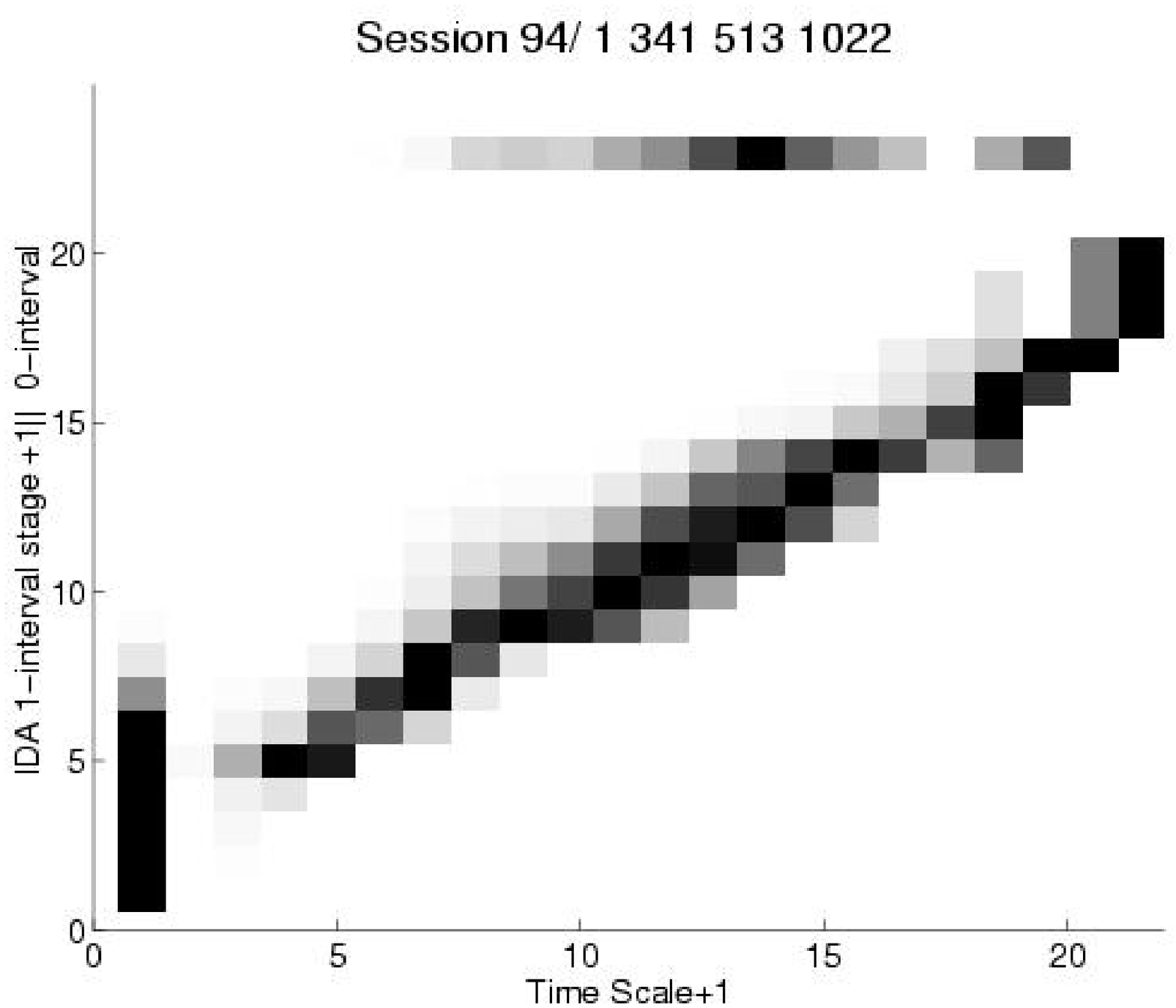}} &
 \subfigure[]{\includegraphics[height=150pt, width=200pt]{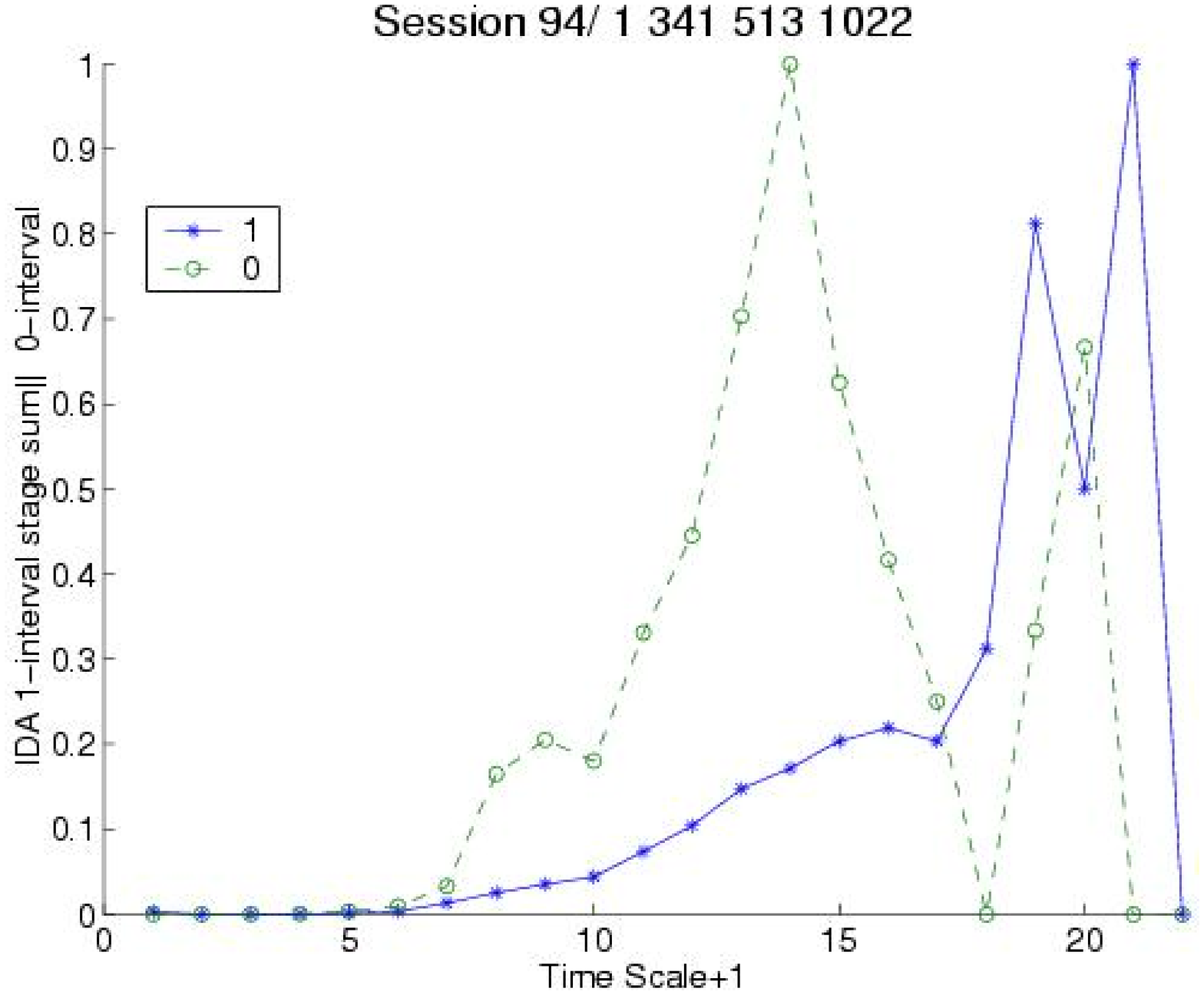}}\\
 & \subfigure[]{\includegraphics[height=150pt, width=200pt]{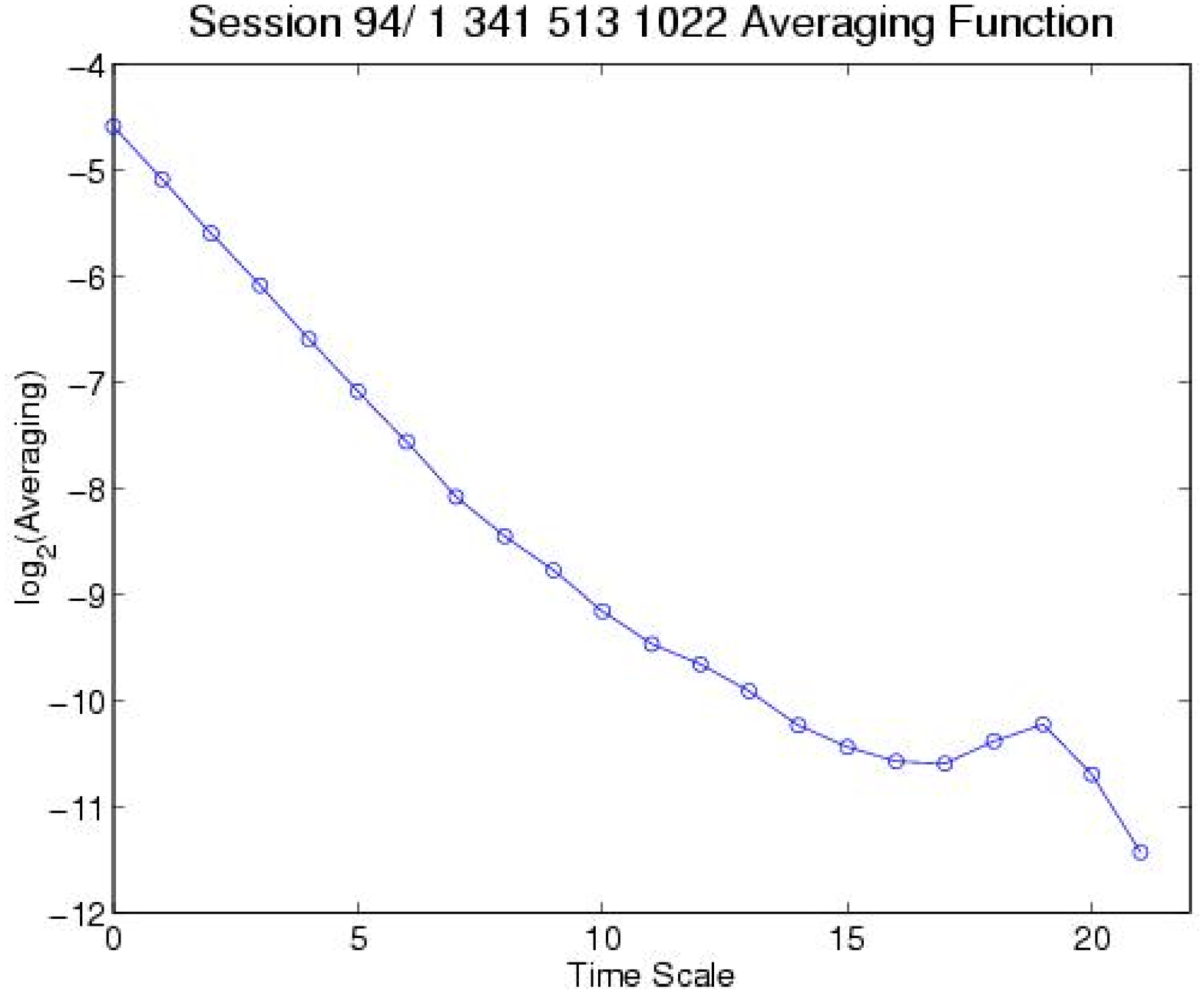}}
 \end{tabular}
 \caption{Detection of $0$- and $1$-intervals of levels in the session of trace 94 used in Fig. \ref{Lev}: (a) and (b) are the greyscale and sum representation of the IDA results, respectively; (c) shows the Averaging function, which is included here in order to facilitate comparison. According to the results of section \ref{curvexplan}, one might think that here the IDA and the Averaging function disagree; we will come back to that in sections \ref{uneq} and \ref{discussion}.}
\label{LevelReaderReal2}         
\end{figure}
\addtocounter{figure}{-1}
\stepcounter{figure} 

\begin{figure}[t]
 \centering
 \begin{tabular}{cc}
 \subfigure[]{\includegraphics[height=150pt, width=200pt]{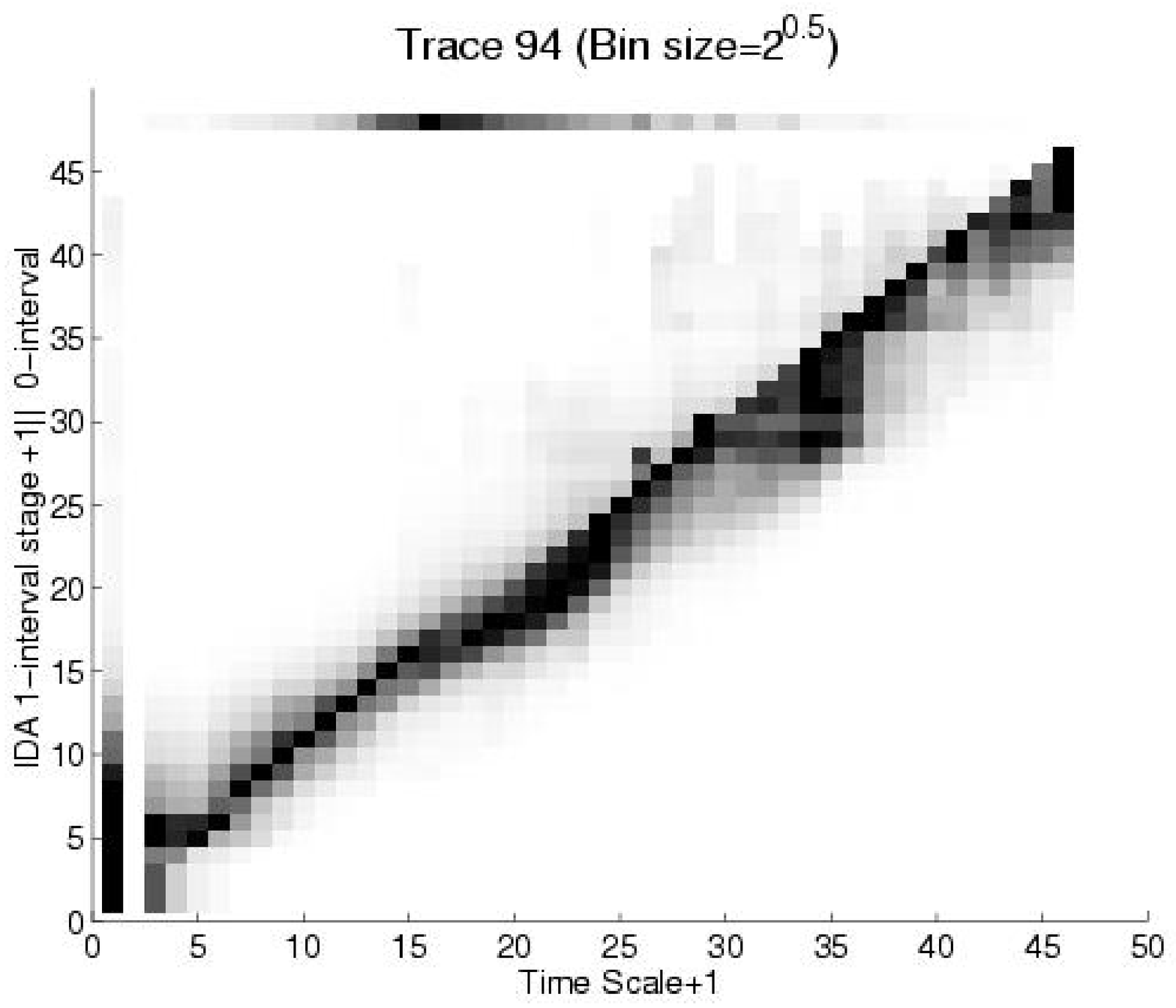}}&
 \subfigure[]{\includegraphics[height=150pt, width=200pt]{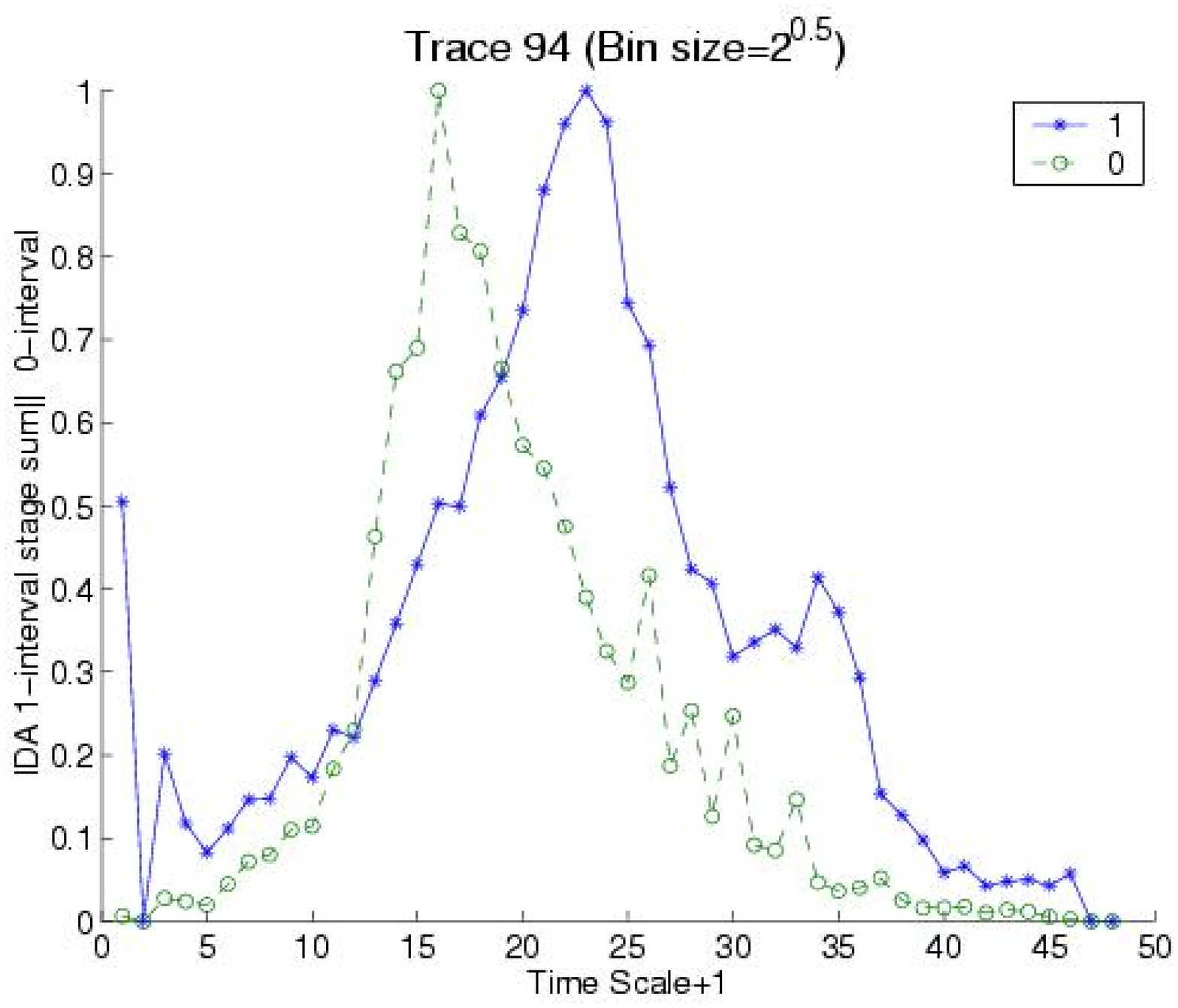}}\\
 & \subfigure[]{\includegraphics[height=150pt, width=200pt]{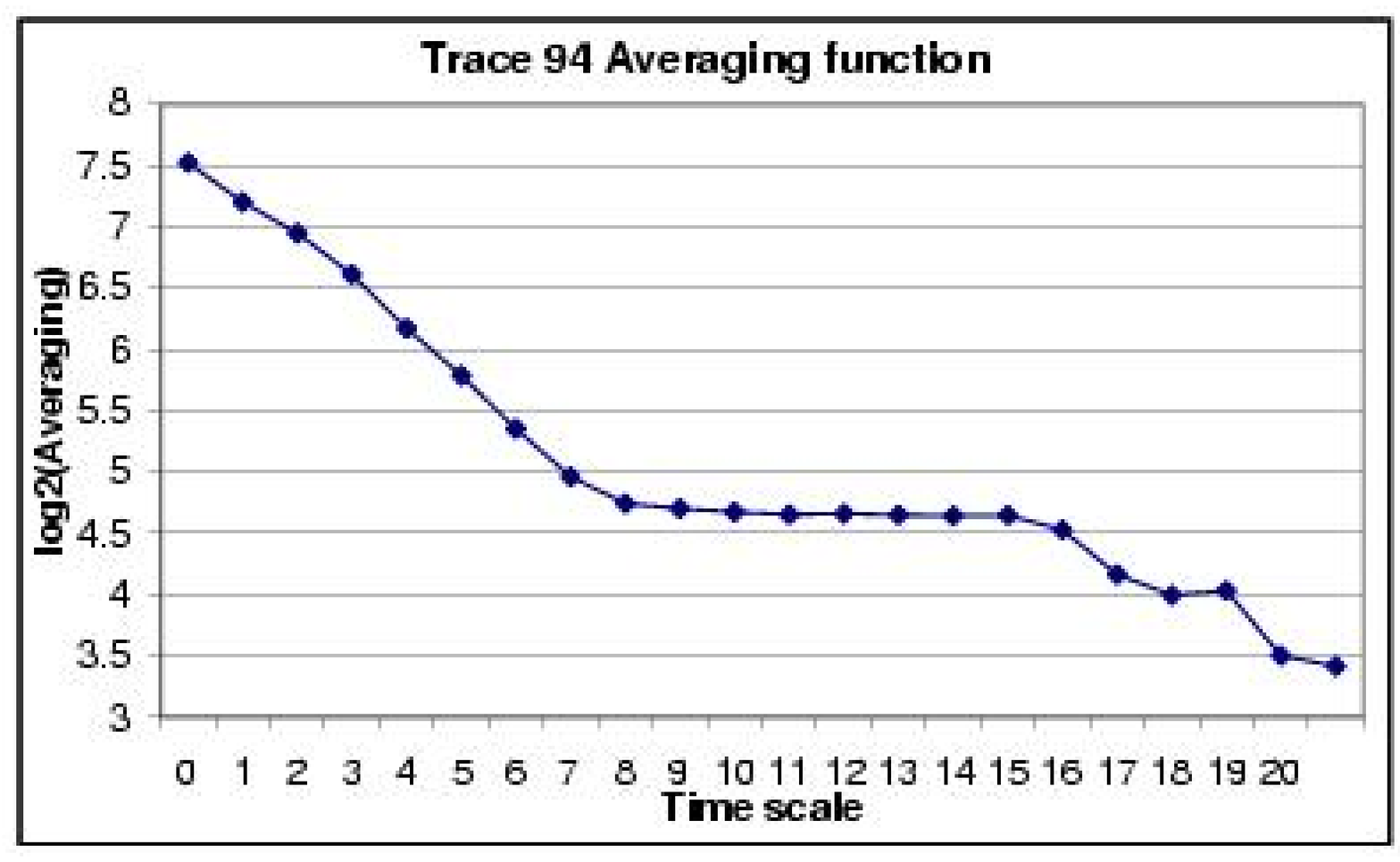}}
 \end{tabular}
 \caption{Aggregate results of the detection of $1$- and $0$-intervals of levels in all trace 94 sessions are shown in (a) and (b): the fact that the algorithm detects levels in this figure suggests that individual sessions have levels in common, which combine in a constructive, not destructive, way. Therefore, it is legitimate to assign levels to entire traces. Note here that the bins increase geometrically, as usual, but with a step of $\sqrt{2}$ instead of $2$. The Averaging function of trace 94 is shown in (c).}
\label{LevDet}         
\end{figure}

Having established the effectiveness of the IDA on individual sessions, simulated or observed, let us now turn our attention to aggregate information: we apply the algorithm to \emph{each} session in trace 94, and superpose the individual results. The result is shown in Fig. \ref{LevDet} (note here that bin sizes increase geometrically by $\sqrt{2}$ rather than by $2$), which shows clearly that some levels persist in the whole trace: their $1$-interval lengths lie between $2^9ms$ and $2^{13}ms$, and also between $2^{15}ms\ \text{and}\ 2^{18}ms$, whereas their $0$-interval lengths lie between $2^7ms$ and $2^{15}ms$. Comparison with the Averaging function of trace 94 (see Fig. \ref{LevDet}(c)), which shows an extended period of activity between $2^7ms$ and $2^{16}ms$, and also a little ``bump'' at $2^{18}ms$, demonstrates that the two methods are compatible, yielding similar results.     

\addtocounter{figure}{-1}
\stepcounter{figure} 

The fact that the IDA detects levels in Fig. \ref{LevDet} suggests that the majority of the connections of trace 94 share some levels, which appear then in the aggregate results, and that, consequently, it is possible to attribute levels to an entire data set, thus extending the notion of levels. The excellent agreement of the two methods (IDA and Averaging function --- see also section \ref{metcomp}) corroborates strongly the legitimacy of the extension of the notion of levels to entire traces, and also implies that the much simpler Averaging function can be used for their detection, instead of the IDA. Indeed, on the same computer, the IDA needed two days to process trace 94 entirely, whereas the Averaging function only takes a few seconds.

\subsection{Theorem D: statement and proof}

Theorem D below implies that the covariance structure of total traffic is
strongly dependent on its levels. 

\begin{theorem}
Let $\displaystyle A^{\Delta }(W)=\Delta
^{-1}\left(\bold{Var}\left(\int_{0}^{\Delta }W(t)dt\right)-\bold{Cov}\left(\int_{0}^{\Delta
}W(t)dt,\int_{\Delta }^{2\Delta }W(t)dt\right)\right)^{1/2}$ and let $%
\widetilde{W}_{i}^{d}(t)$ be the same as $W_{i}^{b}(t)$ but
with $O_{f}^{b}(i)=0\ a.s.$ and $O_{n}^{b}(i)=\infty\ 
a.s.$ Also, take $\widehat{W}_{i}^{d}(t)$ to be the same as $%
W_{i}^{a}(t)$ (notice that $\widehat{W}^{d}$ represents a continuous traffic (\emph{no spikes traffic}), while $\widetilde{W}^{d}$ represents a discrete traffic, where RTT is present (\emph{spikes traffic})). 

\begin{list}{\alph{temp})}{\usecounter{temp}}
\item 
$
\displaystyle
A^{\Delta }\left(\widetilde{W}_{1}^{d}\right)\approx \Delta ^{-1/2}\rightarrow 0\ 
as\ \Delta \rightarrow \infty
$

\item 
$\displaystyle
A^{\Delta }\left(\widehat{W}_{1}^{d}\right)\rightarrow 0\ as\ \Delta \rightarrow 0\ 
or\ \Delta \rightarrow \infty
$
and there exists a constant $K>0$ such that  
\begin{equation*}
A^{\Delta }\left(\widehat{W}_{1}^{d}\right)\geq K\ for\ \Delta
=\bold{E}(O_{n}^{b}+O_{f}^{b})/2
\end{equation*}
\end{list}
\end{theorem}

\bigskip

\begin{proof}

\textit{For a)}. Assume that the $R^b_{i,j}$ are exponentially
distributed with $\bold{E}(R)=\gamma >0.$ Let $Q_{k,\Delta }=\int_{(k-1)\Delta
}^{k\Delta }W_{1}(t)dt$, then clearly $\{Q_{k,\Delta }\}_{k=1}^{\infty }$ are
i.i.d.\ Poisson with parameter $\lambda =\Delta /\gamma .$ Thus 
\begin{equation*}
A^{\Delta }=\Delta ^{-1}\sqrt{\bold{Var}\left(Q_{1,\Delta }\right)-\bold{Cov}\left(Q_{1,\Delta
},Q_{2,\Delta }\right)}=\gamma ^{-1/2}\Delta ^{-1/2}\rightarrow 0\text{, as }\Delta
\rightarrow \infty .
\end{equation*}
The ``exponentiallity'' assumption on $R^b$ makes the proof easier, since it implies $%
\bold{Cov}(Q_{1,\Delta },Q_{2,\Delta })=0$ and it simplifies the computation of $%
\bold{Var}(Q_{1,\Delta })$. However, it can be shown that under much more
general assumptions on the $R$, it is still true that $\left|\bold{Cov}(Q_{1,\Delta },Q_{2,\Delta
})\right|\leq \lambda _{1}\bold{Var}(Q_{1,\Delta })$ for some $\lambda _{1}<1$ and $%
\bold{Var}(Q_{1,\Delta })=\lambda _{2}\Delta $ for some $\lambda _{2}>0.$

\bigskip

\textit{For b).} The first step is to establish that $A^{\Delta }\rightarrow 0$ as $\Delta
\rightarrow \infty $. \cite{TWS1} showed that 
\begin{equation*}
\bold{Var}(Q_{1,\Delta })\underset{\text{large }\Delta }{\approx }\Delta ^{2/p}%
\text{ and }\bold{Cov}(Q_{1,\Delta },Q_{2,\Delta })\geq 0
\end{equation*}
and these two results imply that $A^{\Delta }\underset{\text{large }\Delta }{\approx }%
\Delta ^{1/p-1}\rightarrow 0$, as $\Delta \rightarrow \infty $.

\bigskip 

For the case $A^{\Delta }\rightarrow 0$, as $\Delta \rightarrow 0$,
we observe that, for $\Delta $ small enough, $\displaystyle Q_{k,\Delta }=\int_{(k-1)\Delta
}^{k\Delta }W_{1}(t)dt$ has the following structure 
\begin{equation}
Q_{k,\Delta }=\left\{ 
\begin{array}{rl}
0 & \text{with probability }p_{1}(\Delta ) \\ 
u & \text{with probability }c\Delta \\ 
\Delta & \text{with probability}\ p_{2}(\Delta )
\end{array}
\right.  \label{SmallDelta}
\end{equation}
where $u$ is a r.v. uniformly distributed on the interval $[0,\Delta ]$, $%
c=2/(\bold{E}(O_{n})+\bold{E}(O_{f}))$, $p_{1}(\Delta )=(\bold{E}(O_{f})-\Delta )/(\bold{E}(O_{n})+\bold{E}(O_{f}))$, and $p_{2}(\Delta )=(\bold{E}(O_{n})-\Delta )/\bold{E}(EO_{n})+\bold{E}(O_{f}))$. In order to simplify the computation,
assume that $\bold{E}(O_{n})=\bold{E}(O_{f})=1$. Then, 
\begin{equation*}
\bold{E}(Q_{1,\Delta })=\Delta /2,\ \bold{E}((Q_{1,\Delta })^{2})=\Delta^{2}/2-\Delta ^{3}/3.
\end{equation*}
It is clear, using the very same argument as in (\ref{SmallDelta}), that, for
small $\Delta $ and $u$ a uniform r.v. on interval $[0,\Delta ]$, 
\begin{equation*}
Q_{1,\Delta}Q_{2,\Delta }=\left\{
\begin{array}{rl}
0 & \text{with probability}\ 1/2 \\
\Delta & \text{with probability}\ 1/2-2\Delta \\
u\Delta & \text{with probability}\ 2\Delta 
\end{array}
\right.
\end{equation*}

Thus 
\begin{equation*}
\bold{E}(Q_{1,\Delta }Q_{2,\Delta })=\Delta ^{2}/2-2\Delta ^{3}+2\Delta ^{5}/3
\end{equation*}
Putting this together 
\begin{equation*}
A^{\Delta }=\Delta ^{-1}\sqrt{\bold{E}(Q_{1,\Delta })^{2}-\bold{E}(Q_{1,\Delta},Q_{2,\Delta })}
\end{equation*}
\begin{equation*}
\underset{\text{small }\Delta }{\approx }\Delta ^{-1}\sqrt{5\Delta
^{3}/3-2\Delta ^{5}/3}\underset{\text{small }\Delta }{\approx }\Delta
^{1/2}\rightarrow 0\text{ as }\Delta \rightarrow 0.
\end{equation*}
Finally, the last statement of part $b)$ has to be established. Here, the
assumptions will be simplified considerably. Let $O_{n},$ $O_{f}$ be uniformly distributed on
the interval $[1,2]$ and let $\Delta =1$. Then, trivially 
\begin{equation*}
Q_{1,\Delta }=\left\{ 
\begin{array}{rl}
0 & \text{ with probability }1/8 \\ 
u & \text{with probability }3/4 \\ 
1 & \text{with probability }1/8
\end{array}
\right.
\end{equation*}
where $u$ is a uniformly distributed r.v. on the interval $[0,1]$. Also, since $%
O_{n}\leq \dot{2}$ a.s., $Q_{1,\Delta }Q_{2,\Delta }=0$ a.s.
Therefore:
\begin{equation*}
A^{\Delta }=\Delta ^{-1}\sqrt{\bold{E}(Q_{1,\Delta })^{2}-\bold{E}(Q_{1,\Delta},Q_{2,\Delta })}=3/8
\end{equation*}
Clearly, $O_{n}$ and $O_{f}$ have been chosen in order to simplify the
computation. The extension to more general distributions for $O_{n}$ and $%
O_{f}$ is much more involved, but nevertheless straightforward. Namely, one
can show that under much more general assumptions it is still true $%
|\bold{Cov}(Q_{1,\Delta },Q_{2,\Delta })|\leq \lambda _{1}\bold{Var}(Q_{1,\Delta })$ for
some $\lambda _{1}<1$ and that $\bold{Var}(Q_{1,\Delta })=\lambda _{2}\Delta ^{2}$ for
some $\lambda _{2}>0$ (provided that $\Delta \approx \bold{E}(O_{n})$ ).

\end{proof}

\subsection{Explanation of the curvature of the Averaging/Energy function}

\label{curvexplan}

The following argument explains the link between Theorem D and the curving of the Averaging/Energy function. Namely, assume there is a level at time scale $j_0$, and no other levels nearby. It will be assumed that $2^{j_0}\gg RTT$. Then:

\begin{itemize}
\item For $j< j_{0}$ and $2^j \gg RTT$, the ON-periods appear continuous, and can be modeled by no spikes traffic
$\widehat{W}^{d}$. Consequently, part b) suggests that the Averaging function increases.

\item For $j> j_{0}$, the ON-periods again appear continuous, so part b) suggests that the Averaging function 
decreases, as $j\rightarrow\infty$.
\end{itemize}

\begin{figure}
 \centering
 \subfigure[]{\includegraphics[height=150pt, width=200pt]{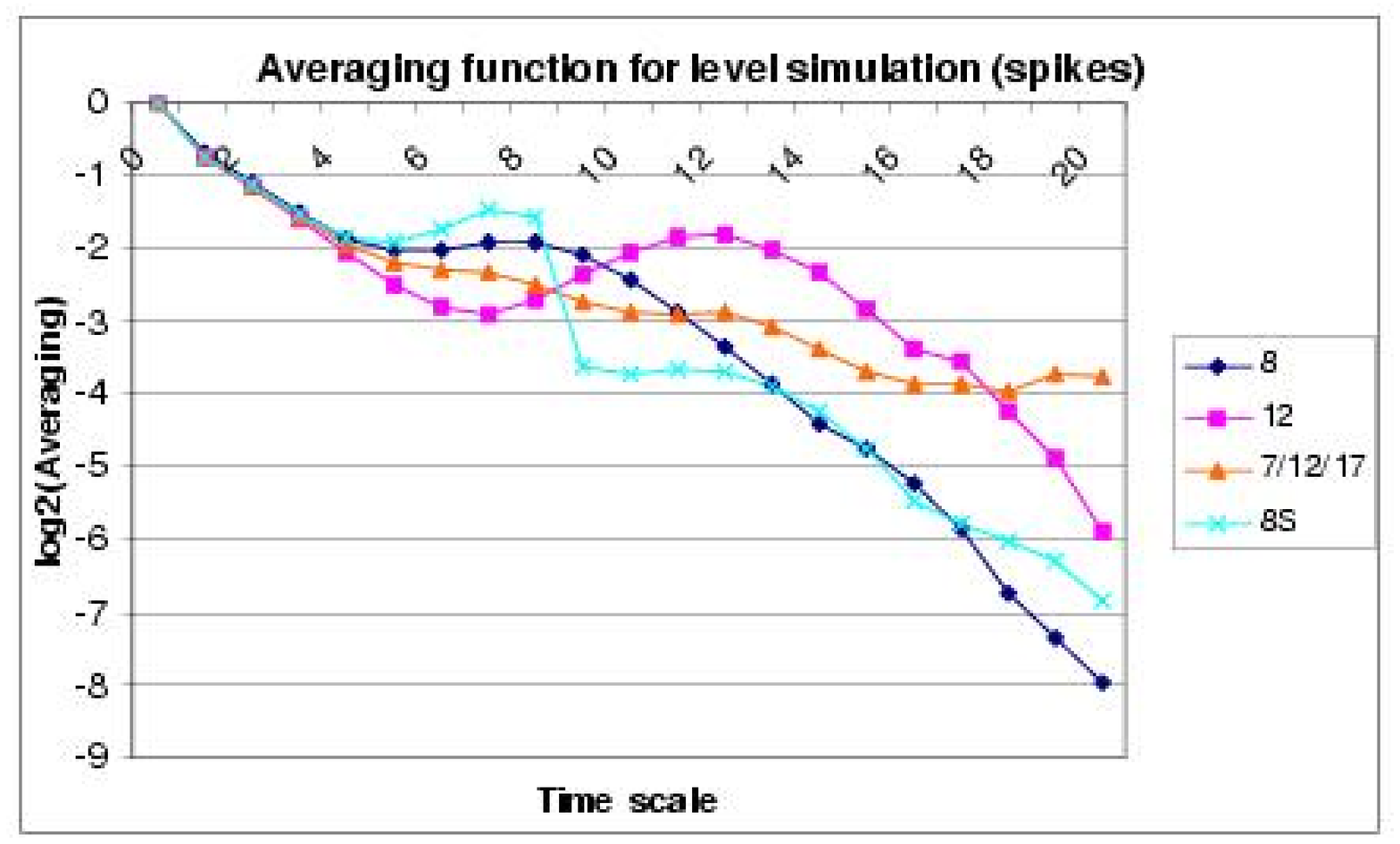}}
 \subfigure[]{\includegraphics[height=150pt, width=200pt]{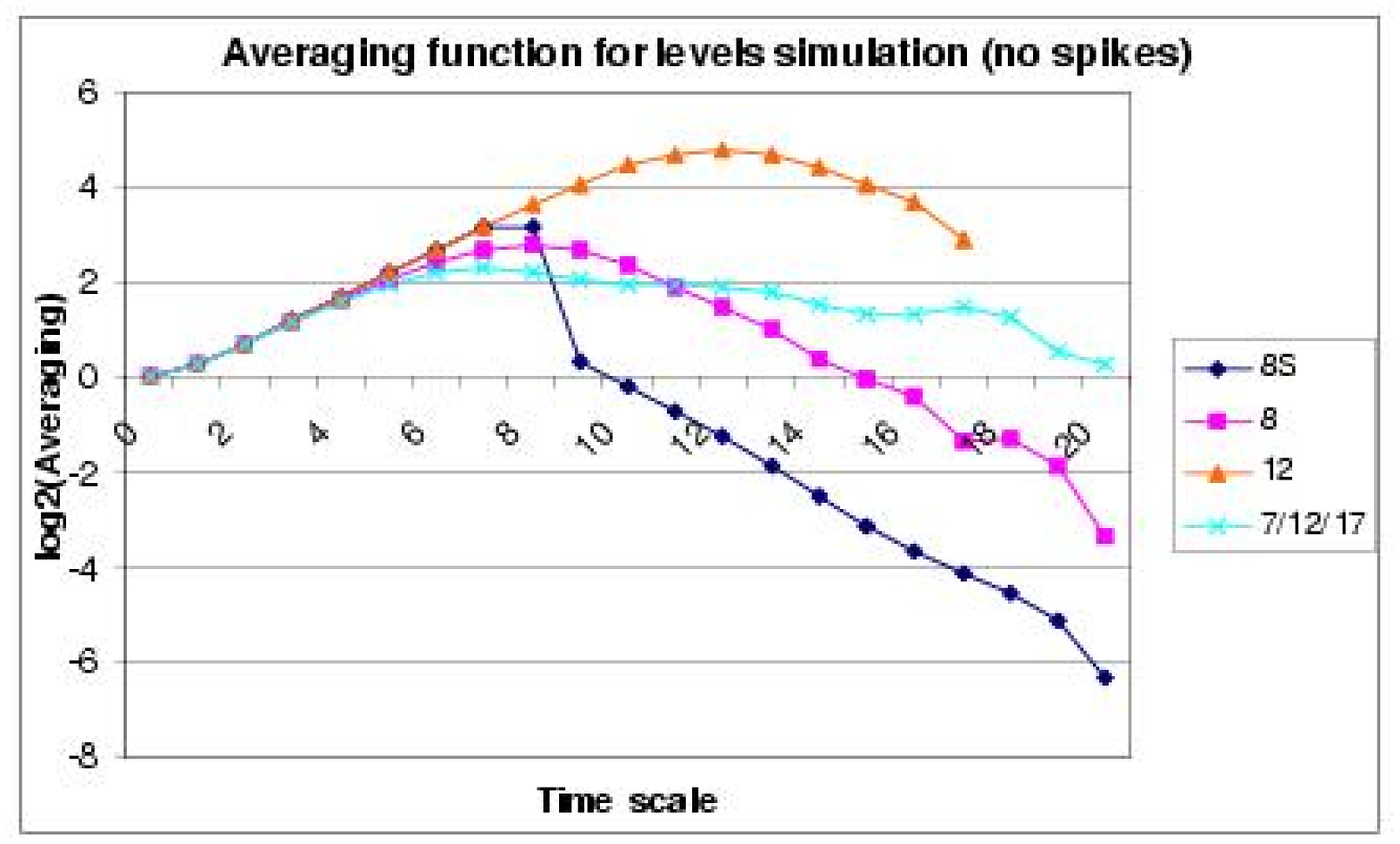}}
 \caption{Average functions of some model D simulations: 8 has a level at $2^{8}$,
          12 at $2^{12}$ and 7/12/17 three levels at $2^{7},2^{12},2^{17}$. ON and OFF 
          intervals of the levels are exponentially distributed, the mean being the level. In 8S, 
          ON/OFF intervals are uniformly distributed within $\pm10\% $ around the mean. In (a), RTT is 
          present: the ON intervals of the finest level consist of spikes of height 1, separated by exponentially 
          distributed intervals with mean 4. In (b), RTT is not present, and the ON intervals of the finest level are 
          continuous.} 
 \label{AvLev}         
\end{figure}

Because the Averaging function is strictly positive, the points above suggest that it must have a 
maximum around the level $j_0$. This has indeed been confirmed by simulations (Fig. \ref{AvLev}(b)). 

However, the Averaging function of real traces initially decays, unlike the simulations in Fig. \ref{AvLev}(b). Part a)
of Theorem D explains this behavior as well: 

\begin{itemize}

\item For $j_{0}< j$ and $2^j\approx RTT$, the ON-periods appear discrete and can be  modeled by spikes traffic
$\widetilde{W}^{d}$. Part a) takes over and ensures White noise-like behavior, i.e. $\log _{2}A_{j}\approx -j/2$.

\end{itemize}

In case more than one levels are present, the argument above can be repeated, the role of RTT being now played by
the size of the previous (finer) level. Therefore, the Averaging function should initially decrease, then change
slope around the level, and have a ``bump'' around every level after the
first. Simulations indeed verify this statement (Fig. \ref{AvLev}(a)). The exact form
of these ``bumps'' has to do with the ``sharpness'' of the level, i.e. with
the variance of the durations of ON and OFF intervals. Very sharp levels,
i.e. with small variance, seem to lead, at the time scale of the level, to
the phenomenon of \emph{superconvergence}\footnote{ 
In real traces, we observed this behavior in the fine and middle scales of traces M6, M7, U8, and 89.  
See also 8S in Fig. \ref{AvLev}.}, i.e. convergence to the mean
value faster than what the CLT predicts, which manifests itself through a
slope less than $-0.5$ (Fig. \ref{AvLev}).

We would like to emphasize that we have introduced levels in our model consistent with our wish for ``local'' models (in the sense that the sessions are constructed by putting together appropriate building blocks, rather than the ``global'' approach of e.g. Random Cascades.) Our inspiration for the model was taken from direct inspection of observed traces; it does however differ from our earlier step in that we cannot point to a particular feature in user behavior or network protocol that could cause them. The scales of the levels, their sharpness, and their possible evolution in time may well be of interest to network engineers. For this reason we propose tools, in section \ref{discussion}, that are faster than the level-detector intoduced above. Note also that heavy-tailed distributions were \emph{not} needed in this construction: 
long range dependence and levels thus appear to be completely separate characteristics.

\addtocounter{figure}{-1}
\stepcounter{figure}

Now that we have understood the role played by levels, we are in position to
explain the shape of the Energy and Averaging functions in the case of Slow
Start simulations (Fig. \ref{SSSim}), where the two functions are actually piecewise linear,
not linear. The change of slope occurs between scales 7 and 9. But the mean
value of the connection size was 5, and the mean RTT was 100, which implies
that the mean total duration of the session was $500\approx2^{9}$, and
this is exactly the level of the traffic. 

\begin{figure}
 \centering
 \subfigure[]{\includegraphics[height=150pt, width=200pt]{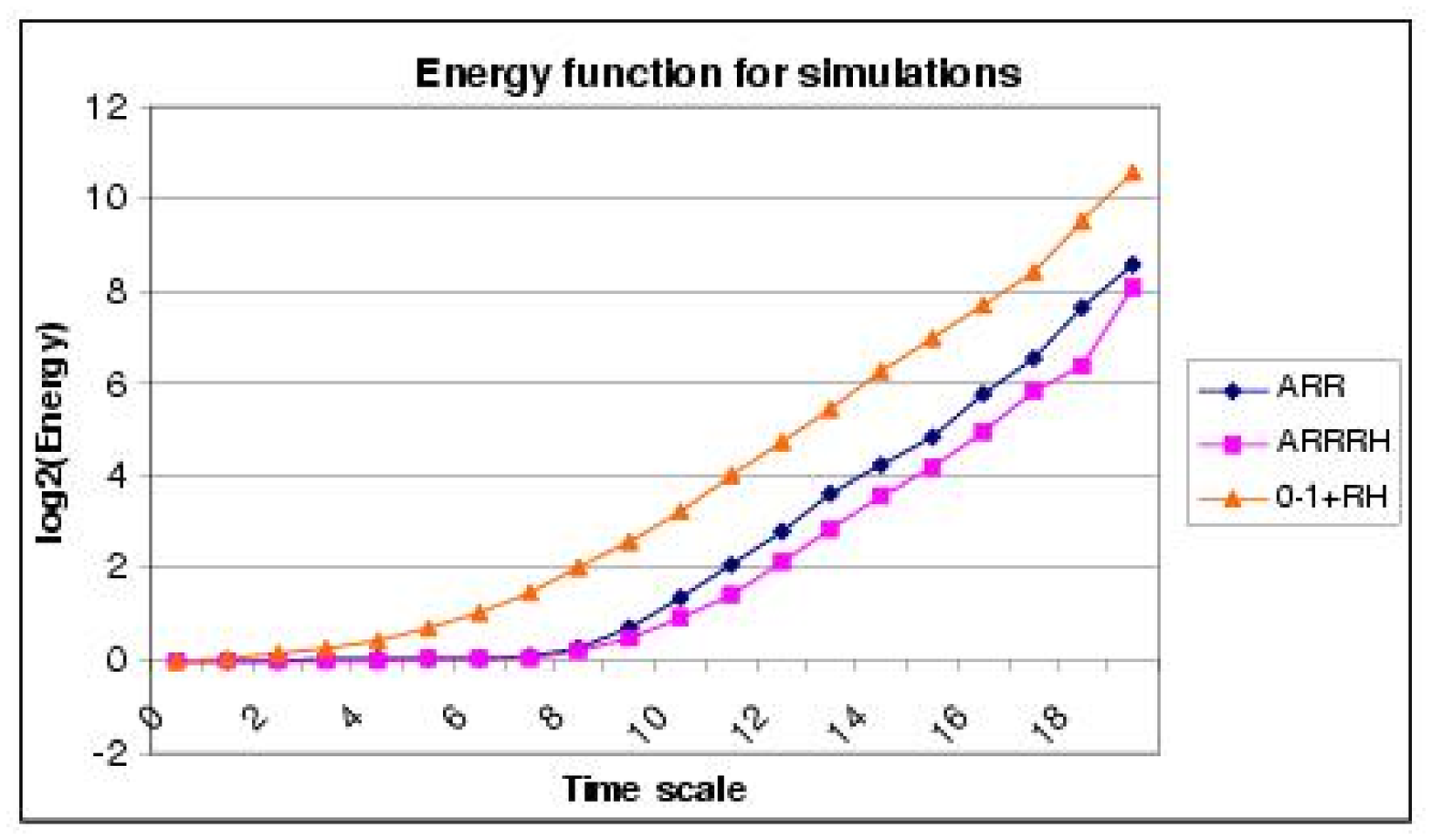}}
 \subfigure[]{\includegraphics[height=150pt, width=200pt]{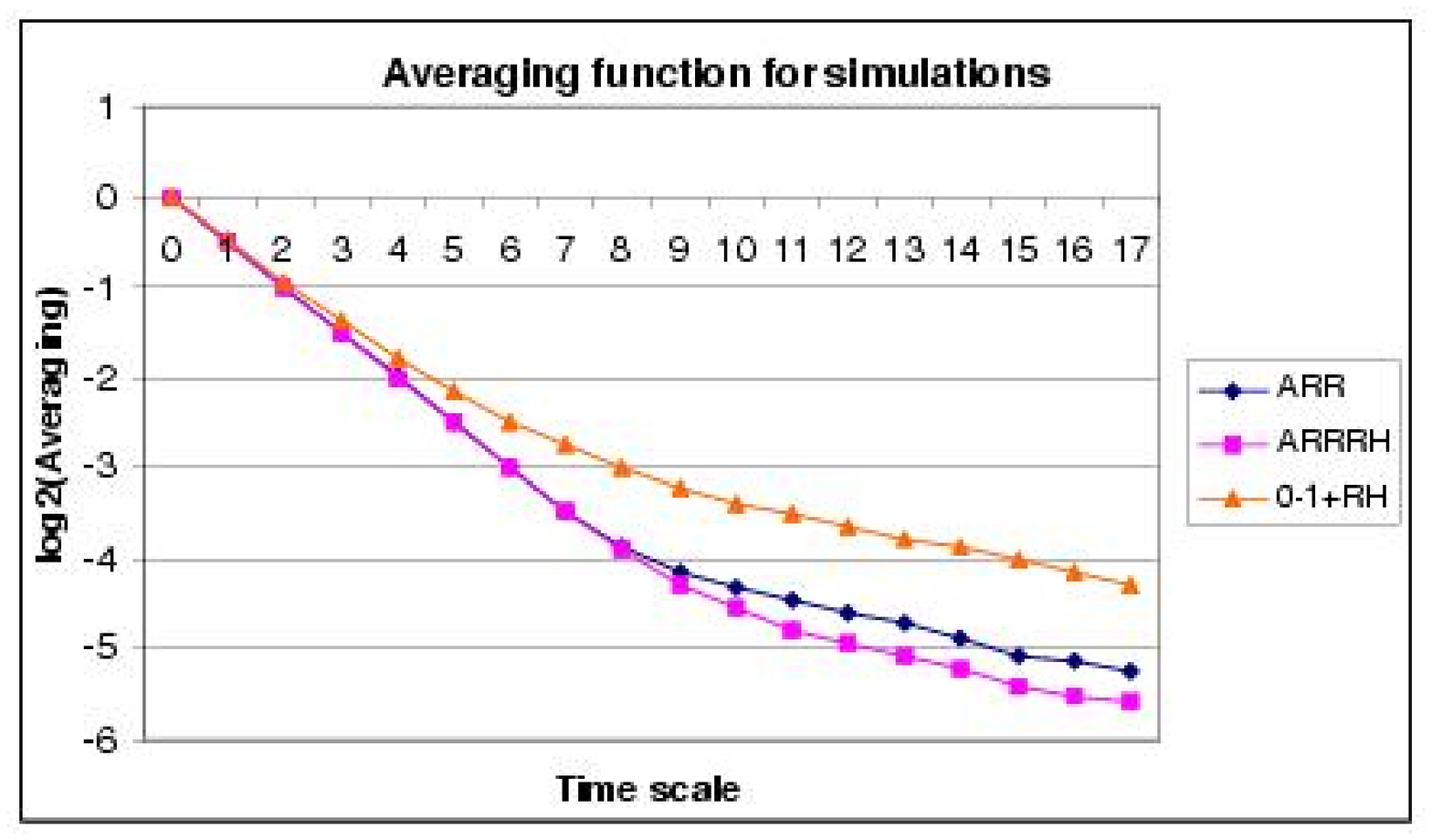}}
 \caption{Some more simulations: curvature is present, since there is one single
          change of slope, the smoothness of which varies. Hence, there is one
          level in these traces.}
 \label{EnAvCur}         
\end{figure}

One could argue that if
the simulations of Fig. \ref{Av} and \ref{En} were repeated with different parameters, some
sort of curving might be achieved. This is true, to some extent (Fig. \ref{EnAvCur}), but 
the curvature achieved is slight and does not resemble the one of real traces.

Finally, let us reexamine Fig. \ref{Av}, armed with this qualitative understanding of levels. Trace 94 seems to have 
several levels between the time scales $8$ and $15$, and maybe around $20$. A comparison with our discussion of Fig. \ref{LevDet} shows a good correspondence, although Fig. \ref{LevDet} gives a more detailed picture. Trace 97 has probably just one strong level at $13$; trace L4 has one level at $8$ and one at $15$; trace L5 seems to have levels all over the place, but mainly one at some small scale (between $0$ and $4$), one at $8$ and one at $15$; traces M6, M7, and U8 appear to have levels at $12$, $16$, and $19$, and M6 maybe has one more at $21$; finally, trace 89 has levels at $2$, $13$, $16$ and $19$. In the next session, tools will be developed to quantify such guesswork.

\addtocounter{figure}{-1}
\stepcounter{figure}

\subsection{Levels with unequal $0$- and $1$-intervals}

\label{uneq}

\begin{figure}
 \centering
 \subfigure[]{\includegraphics[height=150pt, width=200pt]{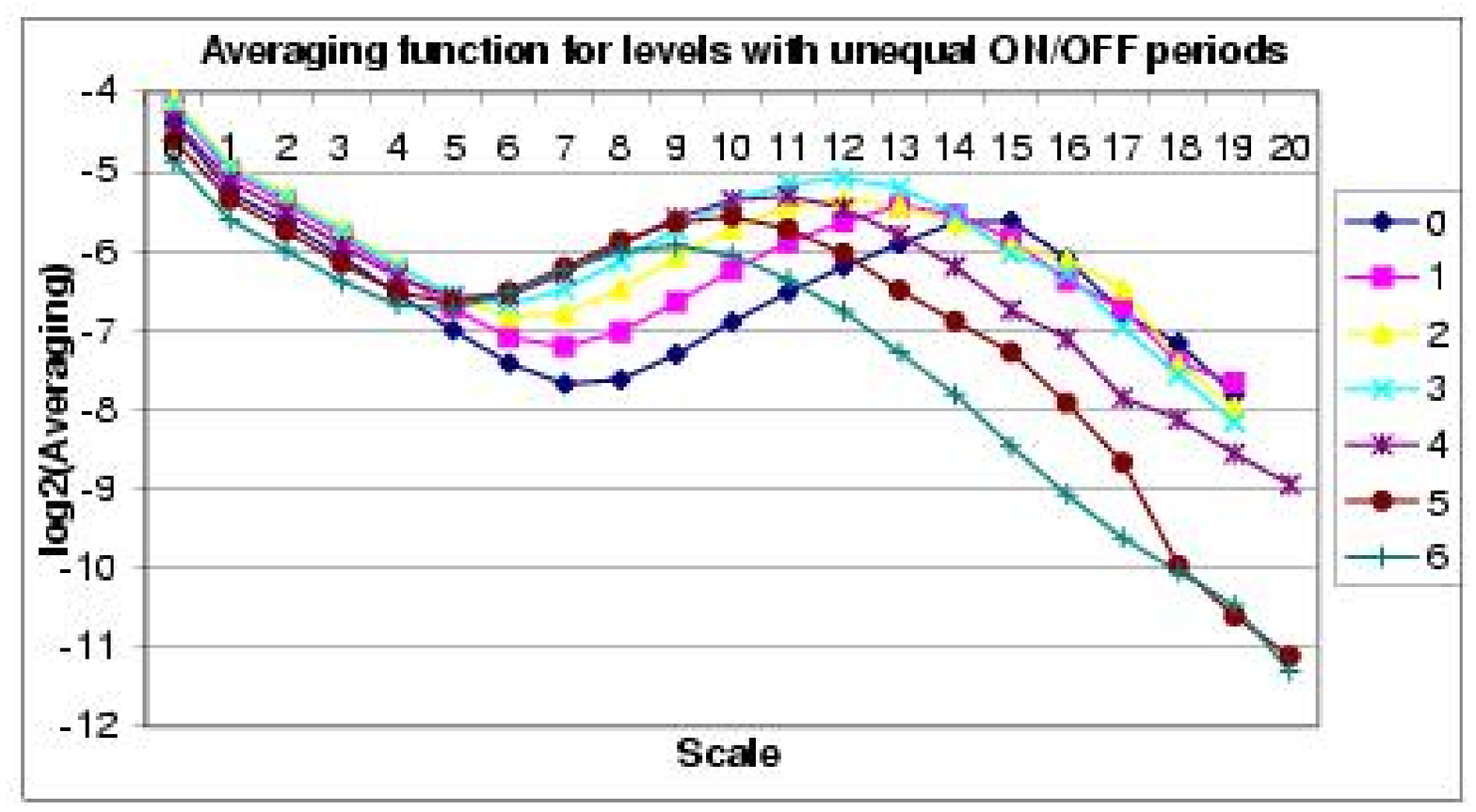}}
 \subfigure[]{\includegraphics[height=150pt, width=200pt]{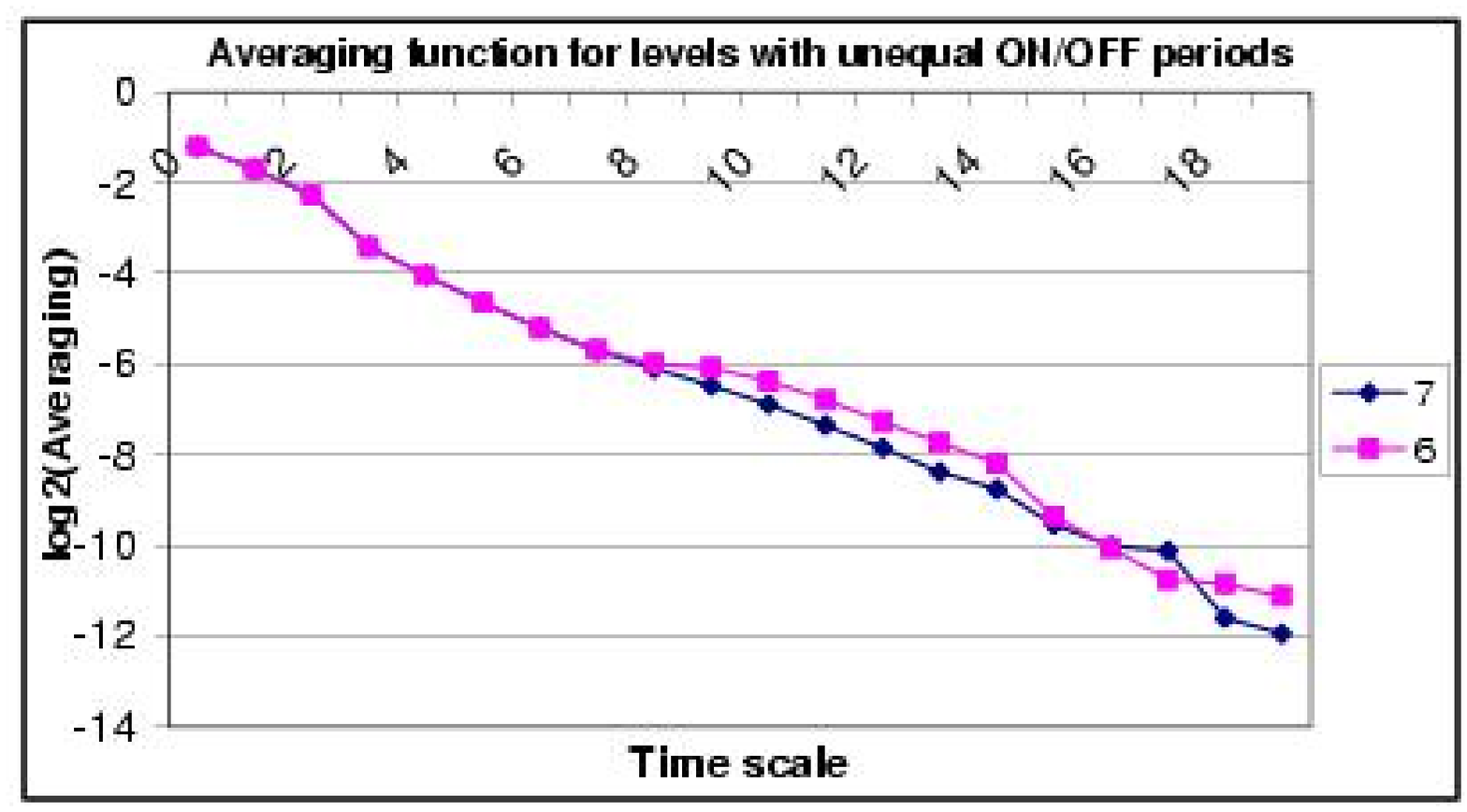}}
 \caption{Averaging function of a simulated session with only one level whose $1$-intervals live in time scale $15$, but whose $0$-intervals live in time scale $15-i$, where, for each curve, $i$ can be read from the legend. In (a), as a rule, the local maximum of the bump ``locks'' on $15-i$, the scale of the $0$-intervals. In (b), using a more diffused interval length distribution, the bump is negligible.}
 \label{UneqLev}         
\end{figure}

Although Theorem D does not require the $0$- and $1$-intervals to follow the same distributions, the discussion in section \ref{curvexplan} considered exclusively equally distributed $0$- and $1$-intervals, in order to enhance clarity. Unfortunately, as we have seen so far (for example, see Fig. \ref{Lev} and \ref{LevND}, or Fig. \ref{LevelReaderReal1} and \ref{LevelReaderReal2}), different distributions seem to be the rule rather than the exception, with the $1$-intervals being typically much larger than the $0$-intervals, so we need to determine how the Averaging function will respond to such a situation. Fig. \ref{UneqLev}(a) shows the Averaging function of a simulated session containing just one level whose $1$-intervals have lengths around $2^{15}$, but whose $0$-intervals have average lengths that are $1,\ 2,\ 4,\ldots$ times smaller. Observe that the Averaging function appears to ``lock'' on the $0$-intervals, the existence of the larger $1$-intervals being only alluded to by the size of the ``bumps'', which are more prolonged than usual. Observe also (in Fig. \ref{UneqLev}(b)) that, if the lengths of the $0$-intervals are much smaller than the lengths of the $1$-intervals, and their distributions sufficiently diffused, it is possible for the bump to become negligible, and thus for the Averaging function to give the impression, by visual inspection, that no level exists.

An excellent example of a real session where the above observations apply is the session of Fig. \ref{LevelReaderReal2} in section \ref{quant} (remember the last remarks made therein): it appears there is a quite pronounced level with $0$-intervals around $2^{13}ms$, but comparison with Fig. \ref{LevelReaderReal2}(b) shows that the Averaging function does not seem to ``feel'' it. Does this suggest a potential failure of the model? Hardly so, in view of the discussion of the previous paragraph. Indeed, according to Fig. \ref{Lev}, the $0$-intervals of the session levels are much smaller than the corresponding $1$-intervals, and therefore it is highly probable that the Averaging function responds to their presence very ``moderately'', just as is the case in Fig. \ref{UneqLev}(b). But even such a moderate response can still be detected by means of suitably designed tools, if not by eye alone, such as the ones presented in section \ref{discussion}.

\addtocounter{figure}{-1}
\stepcounter{figure}

\subsection{Averaging function vs. IDA}

\label{metcomp}

Up to this point, two completely independent methods have been developed for the detection of levels: the Averaging function and the IDA. Naturally, two questions need to be answered now:
\begin{enumerate}
	\item Are the results they give compatible?
	\item Which one is more efficient to use? 
\end{enumerate}

The discussion in sections \ref{quant} and \ref{curvexplan} implies that in principle the two methods yield equivalent results. In turn, Fig. \ref{LevelReader2} corroborates this statement: the session on the left seems to have very diffused levels that only slightly influence the curving of the Averaging function. However, the level at $19$ has a strong impact on the Averaging function. The session on the right, on the contrary, seems to have sharp levels, as the correspondence of the two methods is almost 1-1. Indeed, the extended region of activity between the scales $5$ and $14$ corresponds to a region of intensive curving of the Averaging function. The fact that the latter extends in higher scales also can be attributed to the different distributions of $0$- and $1$-intervals, as explained in section \ref{uneq}. Notice that the coarse levels above scale $17$ impact the Averaging function strongly, as was the case for the previous session, too. 

Regarding efficiency, the Averaging function is far easier and faster to compute than to run the IDA on all of the trace sessions: the Averaging function can be computed in seconds, whereas the IDA typically requires days. The advantage of the IDA, however, is that it gives more detail, as it can distinguish the distributions of $0$- and $1$-intervals. Because of its computational speed, only the Averaging function will be used in the sequel.  

\section{The complete model}

It was mentioned earlier that spikiness and levels are two phenomena highly uncorrelated, and that, as such, they can be studied independently. Following this principle, spikiness was studied in section \ref{modelC} and levels in section \ref{modelD}, each ignoring the existence of the other. A successful simulation, though, should contain both, so our model would not be complete without an algorithm that combines spikiness and levels in a simulated trace. 

The independence of the two features suggests that this combination can be done trivially; just use the algorithm of page \pageref{LevAl}, but interject the following step between the original Step 1 and Step 2:

\begin{enumerate}
	\item[1-2.] 
	
	Set the Slow Start maximum at M, and consider two indices, $i$ and $j$. Set $j\leftarrow 1$. Then:
\begin{enumerate}
	\item Generate a heavy tailed integer random variable $N$.
	\item Move $i$ to the next $1$ value of the vector (generated by Step 1 at page \pageref{LevAl}). If end of vector found, then stop.
	\item Set the vector value at $i$ equal to $j$, set $N\leftarrow N-j$, set $j\leftarrow \min(2*j,M)$. If $N\leq 0$, go to Step (a) and set $j\leftarrow 1$.
	\item Go to Step (b).
\end{enumerate}
\end{enumerate}

All this step does is to supply the level algorithm with a Slow Start background. The shortcoming of this construction, though, is that the simulated Slow Start sessions in the final simulated trace do not necessarily correspond to sessions that could have possibly been generated in practice; indeed, the later steps of the algorithm may ``chop'' the beginning or the middle of a Slow Start session, and obviously such behavior cannot exist in practice. So, at this point, our construction deviates slightly from the actual network mechanisms, but this does not affect the results at all. 

It is possible, however, to change the algorithm slightly and remove these artifacts. Just go back to the algorithm of page \pageref{LevAl}, and add the following step at the end this time:

\begin{enumerate}
	\item[5.] 
	
	Set the Slow Start maximum at M, and consider two indices, $i$ and $j$. Set $j\leftarrow 1$. Choose a number of consecutive levels, starting with the finest one, and label them as ``RTT levels''. Then:
\begin{enumerate}
	\item Generate a heavy tailed integer random variable $N$.
	\item Move $i$ to the next $1$ value of the vector (prdoduced at the end of the algorithm of page \pageref{LevAl}). If end of vector found, then stop.
	\item If the gap between the previous and the current value of $i$ is an RTT or belongs to an RTT level, generate a new heavy tailed integer random variable $N$ and set $j\leftarrow 1$. 
	\item Set the vector value at $i$ equal to $j$, set $N\leftarrow N-j$, set $j\leftarrow \min(2*j,M)$. If $N\leq 0$, go to Step (a) and set $j\leftarrow 1$.
	\item Go to Step (b).
\end{enumerate}
\end{enumerate}

Here, the ``RTT level'' $0$-intervals play the role of additional RTTs that a session may encounter. This additional step in the algorithm, then, essentially identifies whether a gap should be considered a RTT or a $0$-interval: in the former case, the Slow Start session continues over it normally; in the latter, the current Slow Start session ends, and a new one starts. 

\section{Discussion and Applications}

\label{discussion}

\textbf{Correlation of user activity}. It has been systematically assumed in the previous models that the user
streams $W_{i}(t)$ are independent. By adopting this assumption, one neglects effect of network topology on the traffic. Congestion is taken into account only to the extent that it randomizes the
RTT and maybe decreases the maximum $M$ of Slow Start, but this clearly does not
introduce cross-correlation among the users. Let us look a little closer at this seemingly drastic limitation.

The possible causes for interdependence between users are \emph{path sharing}
and \emph{server sharing}. The former means that the data belonging to different users flow along network paths that share some links, or that exhibit correlation otherwise (ISP policies etc.). The latter means that two or more users access the same server, and perhaps the same information, i.e. that the end of the two paths is the same.

As regards to path sharing, the size of the Internet makes it very
improbable today that, in a WAN setting, a major number of users will share
a significant number of links; traffic typically goes through many hops (8 to
15), so that, if correlation builds up somewhere, it will vanish later, as user
streams split again, assuming that delays on different links are
independent. Of course, this will not be the case if the majority of users 
share the first link, but this will mean that the users belong to the same
LAN (perhaps Ethernet, FDDI ring etc. \cite{PD1}), which violates the WAN
setting assumption. 

As regards to server sharing, the variety of interests of Internet users
suggests that this will hardly ever be a problem. Extreme scenarios, such as users
massively accessing news sites in order to get news about an
emergency, are clearly rare, and no attempt was made here to model them.

As a final comment, the combination of models C and D produces traffic
statistically indistinguishable from the real one, from the point of view of the statistical tools
hitherto presented: the Averaging function, the marginal distributions,
and the autocorrelations yield identical results for simulations and real traffic. The fact that this was achieved
without modeling any user dependence suggests that its effect is,
in reality, less important than it is believed to be, or at least not detectable by these measurements. 

\bigskip

\noindent \textbf{Burstiness and levels}: In our analysis, spikiness of the traffic and
curvature of the Enegry/Averaging function are two separate issues, due to different and independent
mechanisms (geometric increase in transmission rate versus fragmentation of ON-intervals). Heavy-tailed 
distributions have no effect on the latter, but influence directly the former.

\bigskip

Tools that could detect spikiness and levels could prove useful. Some tools, based on 
the theorems proved earlier, are given below:

\bigskip

\noindent \textbf{Level-detecting tools}: Since levels leave a trace on the Averaging
function, the Averaging function can be used to detect the level structure in a
particular data set. This would allow for monitoring of user behavior
without resorting to a connection-wise study of the trace, which is a very inefficient technique: for trace 94, the Averaging function was about 20000 times faster than computing the function in Fig. \ref{LevDet}, on the same computer.  We propose the following two algorithms:

Given the binary logarithms of the Averaging function: $\log _{2}(A\left(
0\right) ),...,\log _{2}(A\left( n-1\right) )$, we compute slopes $%
S_{i}=\log _{2}(A\left( i\right) )-\log _{2}(A\left( i-1\right) )$, $%
i=1,...,n-1$ and slope changes $\left|S_{i+1}-S_{i}\right|,\ i=1,...,n-2$. These last
quantities are equivalent to second derivatives, which, as is well known,
measure the local curvature: as we argued in our discussion of Theorem D, the Averaging function has ``bumps'' 
around its levels, thus any curvature indicates the presence of levels.  
There is, though, one more point: it is not
really important how much adjacent slopes differ, but how much they differ
with respect to their absolute magnitudes. It seems then more appropriate to
consider: 
\begin{equation*}
Sc_{i}=\frac{\left|S_{i+1}-S_{i}\right|}{\max \left( \min
(\left|S_{i}\right|,\left|S_{i+1}\right|),\varepsilon \right)},\ i=1,...,n-2
\end{equation*}
where $\varepsilon $
is a small positive number, which will safeguard against division by 0, or
an extremely small number (we used $\varepsilon=0.01$).

\bigskip

\noindent\textbf{Tool 1: \emph{(Levels detector)}} 
\begin{equation*}
L(i)=Sc_{i},i=1,...,n-2
\end{equation*}
\noindent The local maxima are considered to be levels.

\bigskip

Since White Noise leads to an Averaging function whose binary logarithm is
just a straight line of slope $-0.5$ (see corollary B), a much simpler idea would be to declare
that a level exists where the slope becomes too large, as a signed quantity (ideally, it should become positive). 
How large it should be exactly is,
in view of the lack of an accurate method to determine this, a matter of taste. We think that, due to the ``noise'' 
induced by the finite number of
measurements available, it is very hard to distinguish a decay of type $\left(\cdot\right)^{-0.1}$ from no decay at all, 
thus we set a threshold at $-0.1$.

\bigskip

\noindent \textbf{Tool 2: \emph{(Detecting flat regions)}} 
\begin{equation*}
F(i)=1_{S_{i}>-0.1},i=1,...,n-1
\end{equation*}
\noindent Levels are considered to exist within flat (or concave) regions where $F(i)=1$.

\begin{figure}[t]
 \centering
 \subfigure[]{\includegraphics[height=150pt, width=200pt]{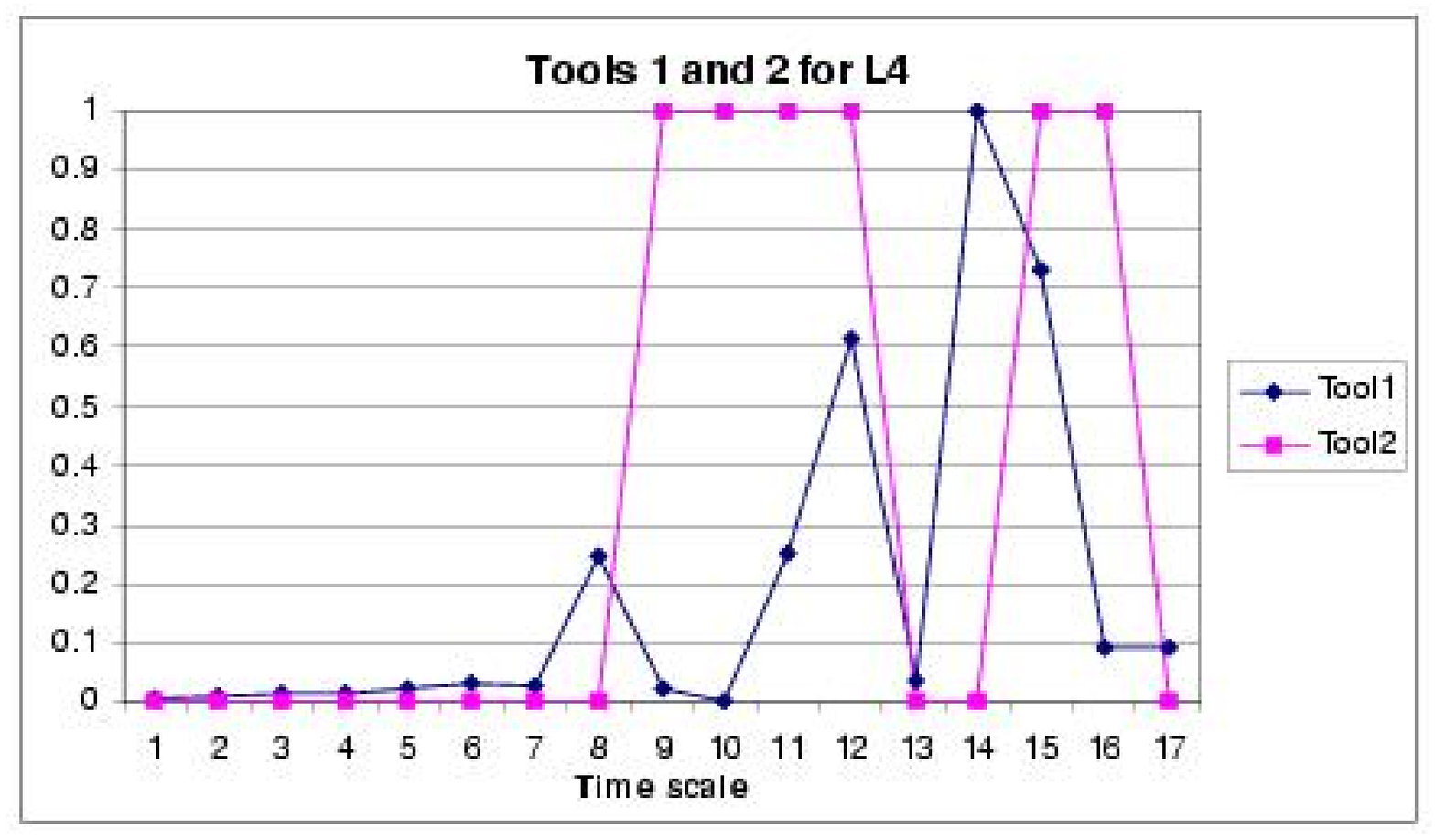}}
 \subfigure[]{\includegraphics[height=150pt, width=200pt]{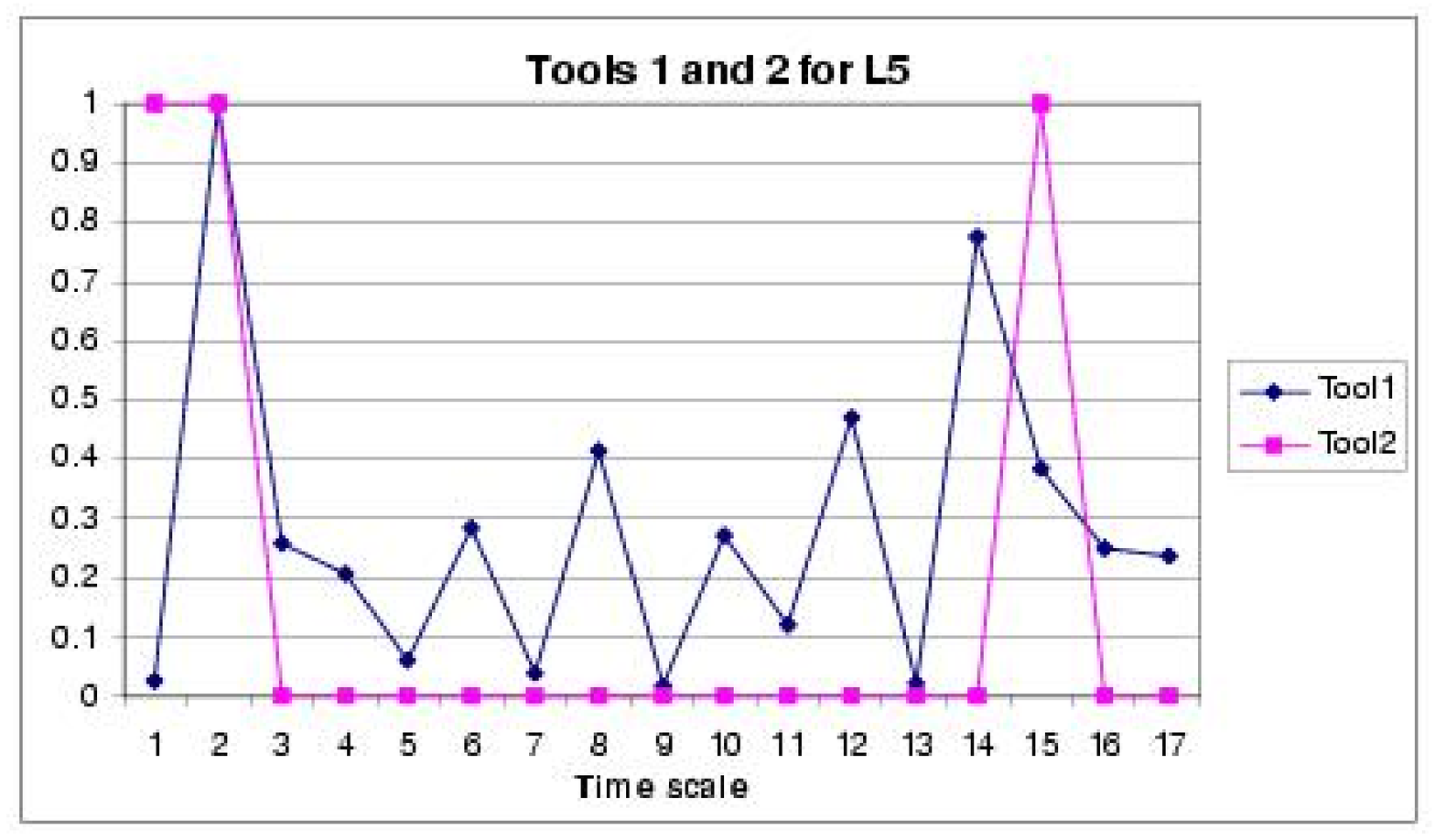}}
 \subfigure[]{\includegraphics[height=150pt, width=200pt]{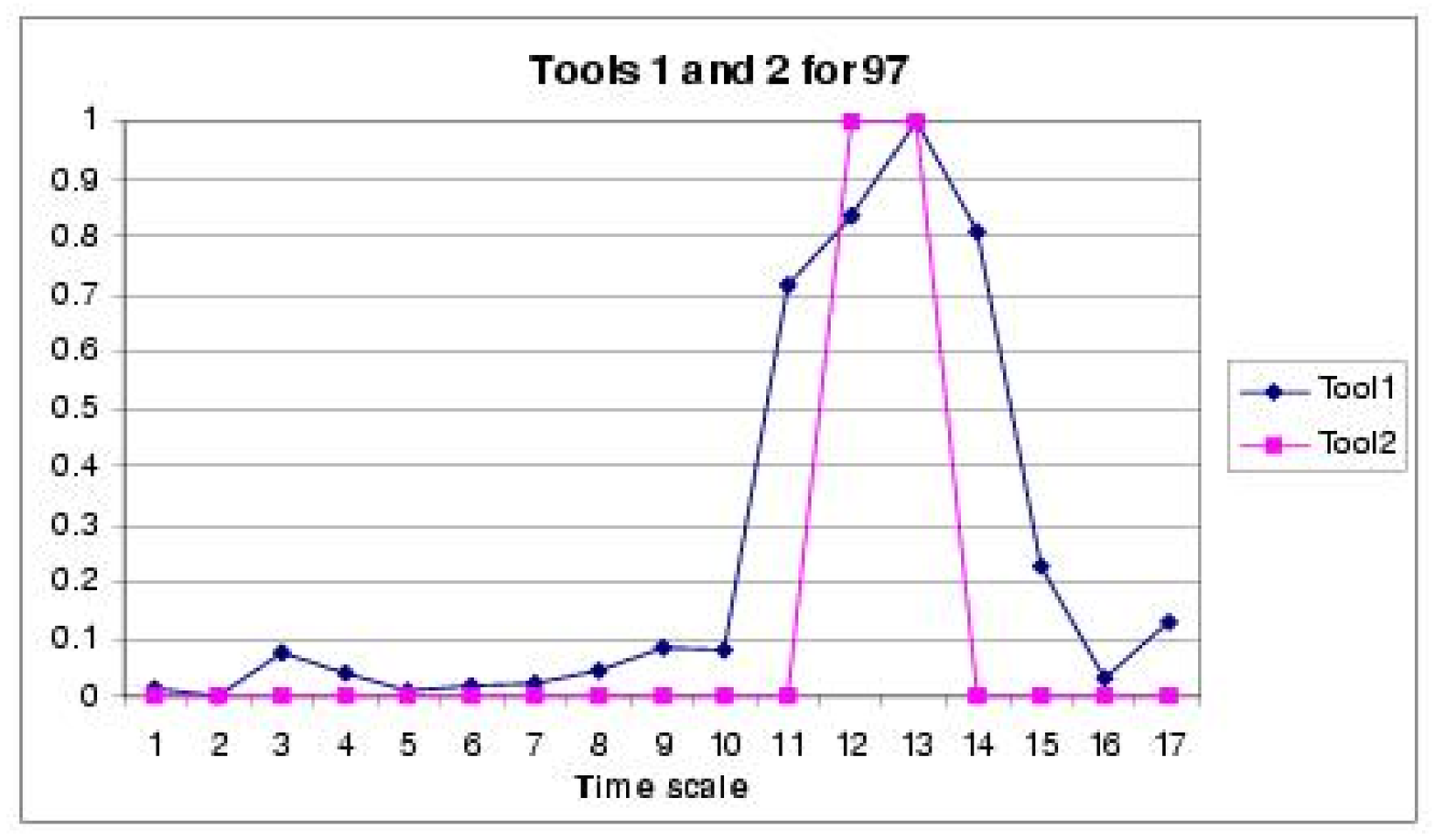}}
 \subfigure[]{\includegraphics[height=150pt, width=200pt]{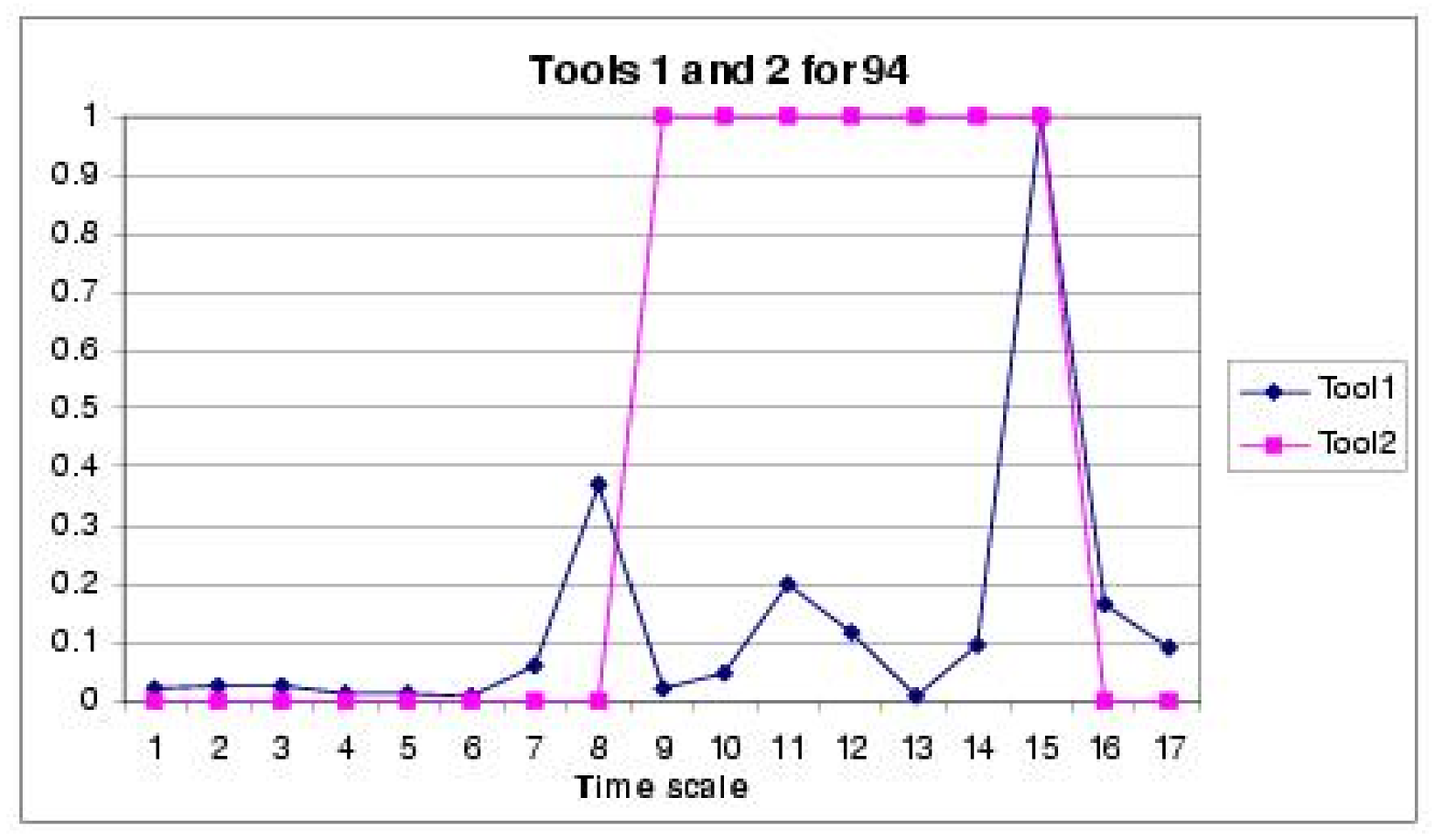}}
 \caption{Tools 1 and 2 applied to real traces. $\epsilon=0.01$ was used.}
 \label{ToolData}
\end{figure}

The results of Tools 1 and 2 are shown in Fig. \ref{ToolData}, \ref{ToolDataNew}  and \ref{ToolLev} 
for the real traces and some of the simulations, respectively. 

\bigskip

\noindent \textbf{Spikiness calibration tools}: Spikiness varies widely, as was
explained before, according to the protocol specifications, the amount of
the transported data etc. For example, trace 94 seems much ``wilder'' than
trace 97 (Fig. \ref{Traffic9497}). Since spikiness is intimately related to the marginal
distribution, the Kolmogorov distance between the trace distribution and the
Gaussian curve seems to be a good measure of it. But a knowledge of the
degree of persistence of spikes through different time scales is also
desirable, therefore the Kolmogorov distance should be computed for each of
them.

A possible extension would be to add information about how much the degree
of spikiness varies within a particular time scale; this may be of
particular importance in view of the oscillation of the process between
Gaussianity and p-Stability (recall the discussion of model C). Thus,
``windowed'' marginals can be taken: the total trace, in a certain time
scale, will be divided into adjacent parts (sub-traces) of consecutive
points, the marginal of each and the corresponding Kolmogorov distances will
be computed, and finally their mean value will be taken, in order to combine
all this information into a single number-index.

In order to combine into a single number the results from every time scale, we compute a
(weighted) average. There is considerable liberty here regarding the choice of
the weights: Tool 3 below utilizes an equiweighted sum, whereas Tool 4 uses weighted averages,  
inspired by the discussion in section \ref{modelC}.

\bigskip

\noindent \textbf{Tool 3: \emph{(Calibrating Deviation from Gaussian Process)}}

Given data $X_{i}$: for $i=1,...,N=2^{n}$ and their coarser versions $%
X_{i}^{j}=X_{(i-1)2^{j}+1}+...+X_{i2^{j}}$, for $i=1,...,N_{j}=N2^{-j}$ and $ 
j=0,...,k-1$ , compute the empirical cumulative distribution $\displaystyle 
F_{j}(t)=N_{j}^{-1}\sum_{i=1}^{N_{j}}1_{\widehat{X}_{i}^{j}\leq t}$, 
where $\widehat{X}_{i}^{j}$ stands for the centered and normalized version of $%
X_{i}^{j}$ (i.e. $\widehat{X}_{i}^{j}=(X_{i}^{j}-\overline{X_{i}^{j}}
)/Std(X_{i}^{j})$), and then $D_{j}=\underset{t\in \mathbb{R}}{\sup }%
|F_{j}(t)-\Phi (t)|$, where $\Phi (t)$ is the standard normal cumulative
distribution. The index of burstiness is: 
\begin{equation*}
D=k^{-1}\sum_{j=0}^{k-1}D_{j}
\end{equation*}
\noindent For Gaussian traffic, $D=0$; the larger $D$ is, the more the traffic 
deviates from Gaussianity.

Different $k$ can be considered, but we suggest $k\leq n-10,$ since for
larger $k$ sample size would not be sufficient, and thus the deviation of
the empirical distribution from the true one would be no more negligible.

\bigskip

\textbf{Tool 4: \emph{(Calibration of Burstiness of the Traffic)}}

Given data $X_{i}$ and $X_{i}^{j}$ as before and constants $k$ and $s$ (we
suggest $s\geq 9$ and $k\leq n-s-7$) let $N_{j,s}=N2^{-j-s}$.
Compute $\displaystyle F_{j,l}(t)=2^{-s}\sum_{i=(l-1)2^{s}+1}^{l2^{s}}1_{\widehat{X}%
_{i}^{j}\leq t}$ for $j=0,...,k-1$ and $l=1,...,N_{j,s}$ (where $\widehat{X}%
_{i}^{j}$ stands for centered and studentized version of $X_{i}^{j}$ )$.$ In
other words, $F_{j,l}(t)$ is the empirical cumulative distribution of the $l$'th
window of the $j$'th time scale. Next, compute $D_{j,l}=\underset{t\in 
\mathbb{R}}{\sup }|F_{j,l}(t)-\Phi (t)|$, $\displaystyle T_{j,l}=
\sum_{i=(l-1)2^{s}+1}^{l2^{s}}\widehat{X}_{i}^{j},$ and $\overline{D}\displaystyle 
_{j}=N_{j,s}^{-1}\sum_{l=1}^{N_{j,s}}D_{j,l}$, which represent the
Kolmogorov distance of each window, the total traffic of each window and the
mean Kolmogorov distance for the time scale $j$, respectively. Next, for a
fixed $j$, compute the cross-correlation between the sequences $%
\{D_{j,l}\}_{l=1}^{N_{j,s}}$ and $\{T_{j,l}\}_{l=1}^{N_{j,s}}$: 
\begin{equation*}
C_{j}=\bold{Corr}(\{D_{j,l}\}_{l=1}^{N_{j,s}}and\{T_{j,l}\}_{l=1}^{N_{j,s}})
\end{equation*}

The result is computed as: 
\begin{equation*}
O=\frac{\sum_{j=0}^{k-1}C_{j}\overline{D}_{j}}{\sum_{j=0}^{k-1}\overline{D}_{j}}
\end{equation*}
\noindent The larger $O$ is, the burstier is the traffic.

Note here that the main quantity considered is the correlation and that Kolmogorov distances are used as weights. The 
results will in general be negative (see discussion in section \ref{modelC}). 

\addtocounter{figure}{-1}
\stepcounter{figure}

\begin{figure}[t]
 \centering
 \subfigure[]{\includegraphics[height=150pt, width=200pt]{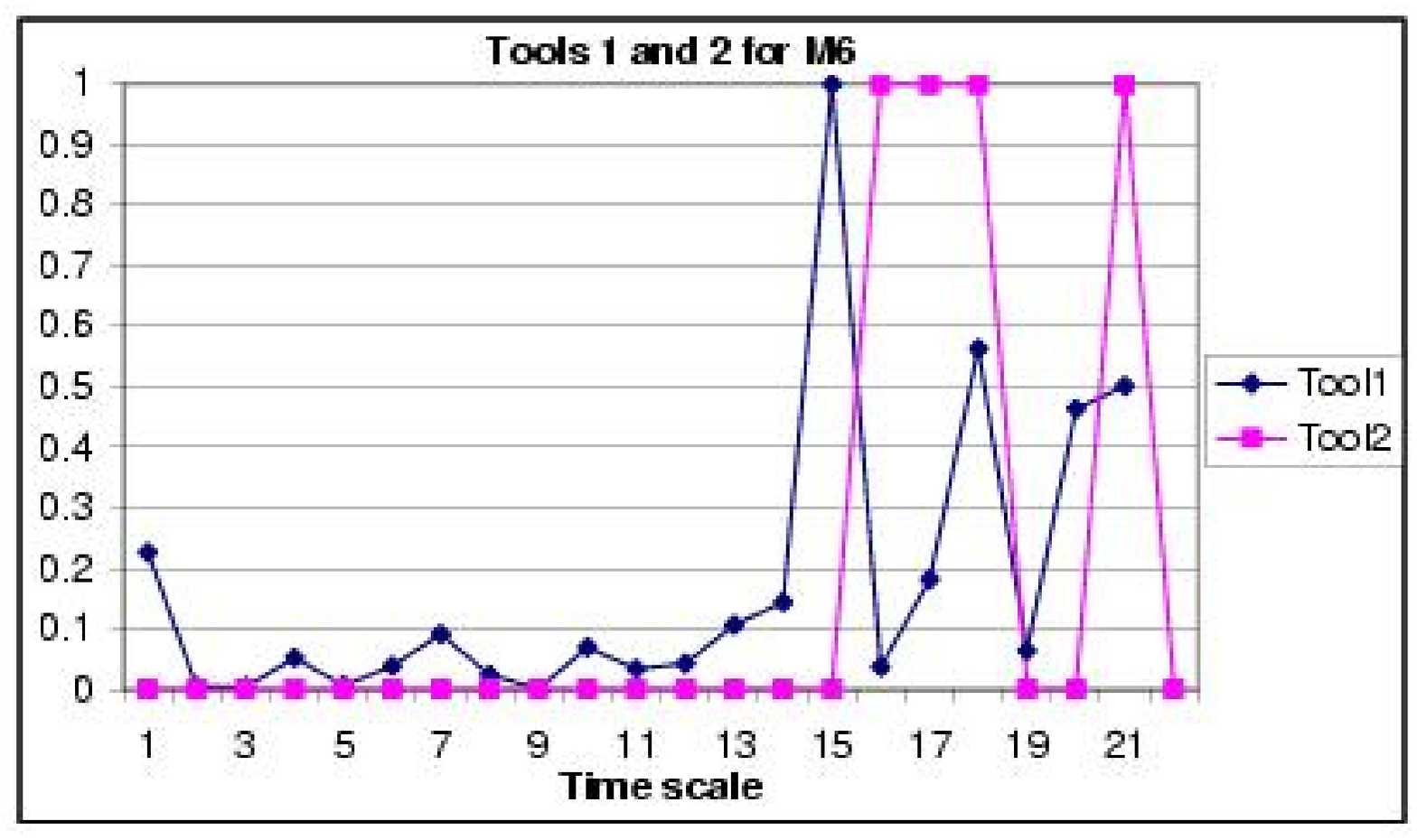}}
 \subfigure[]{\includegraphics[height=150pt, width=200pt]{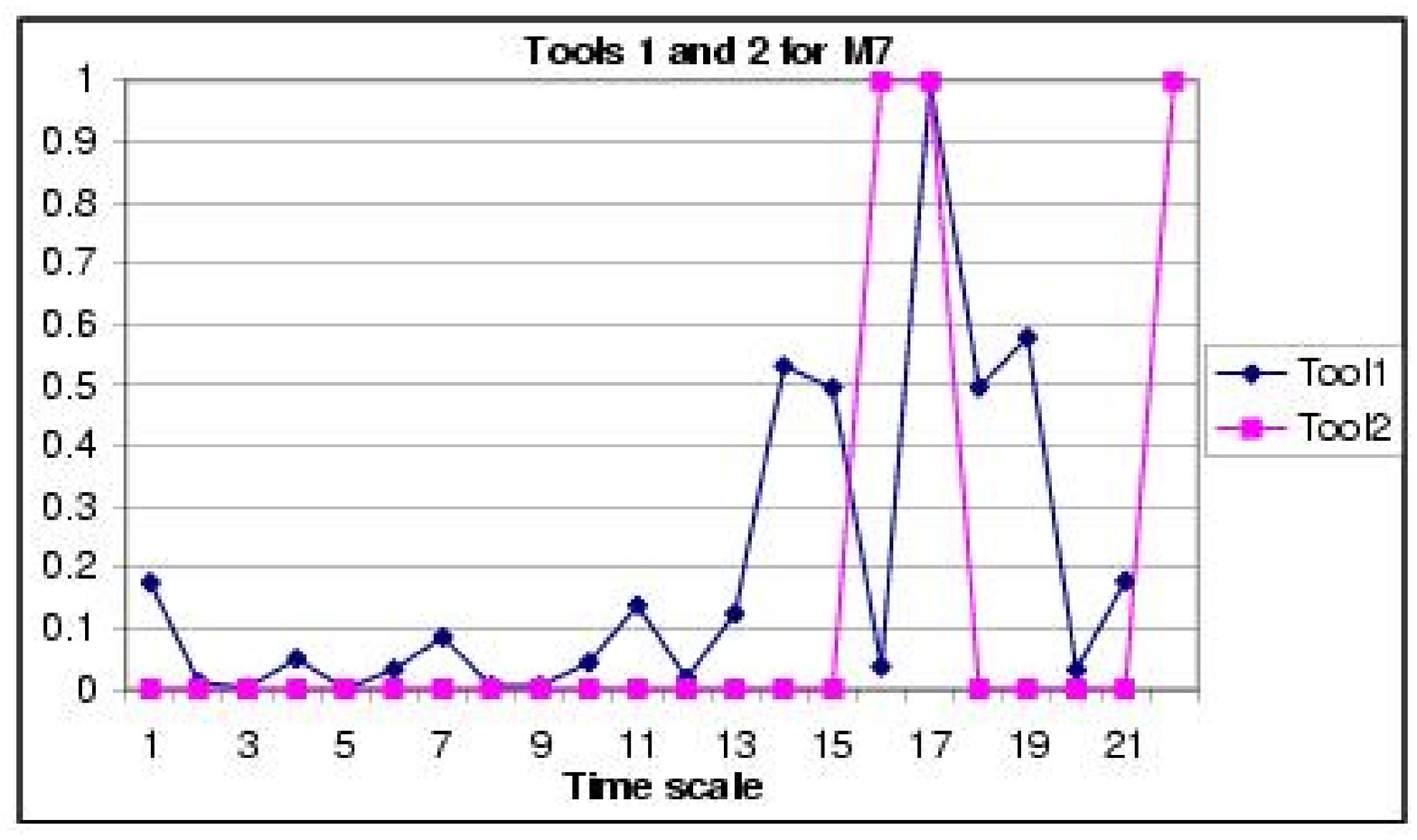}}
 \subfigure[]{\includegraphics[height=150pt, width=200pt]{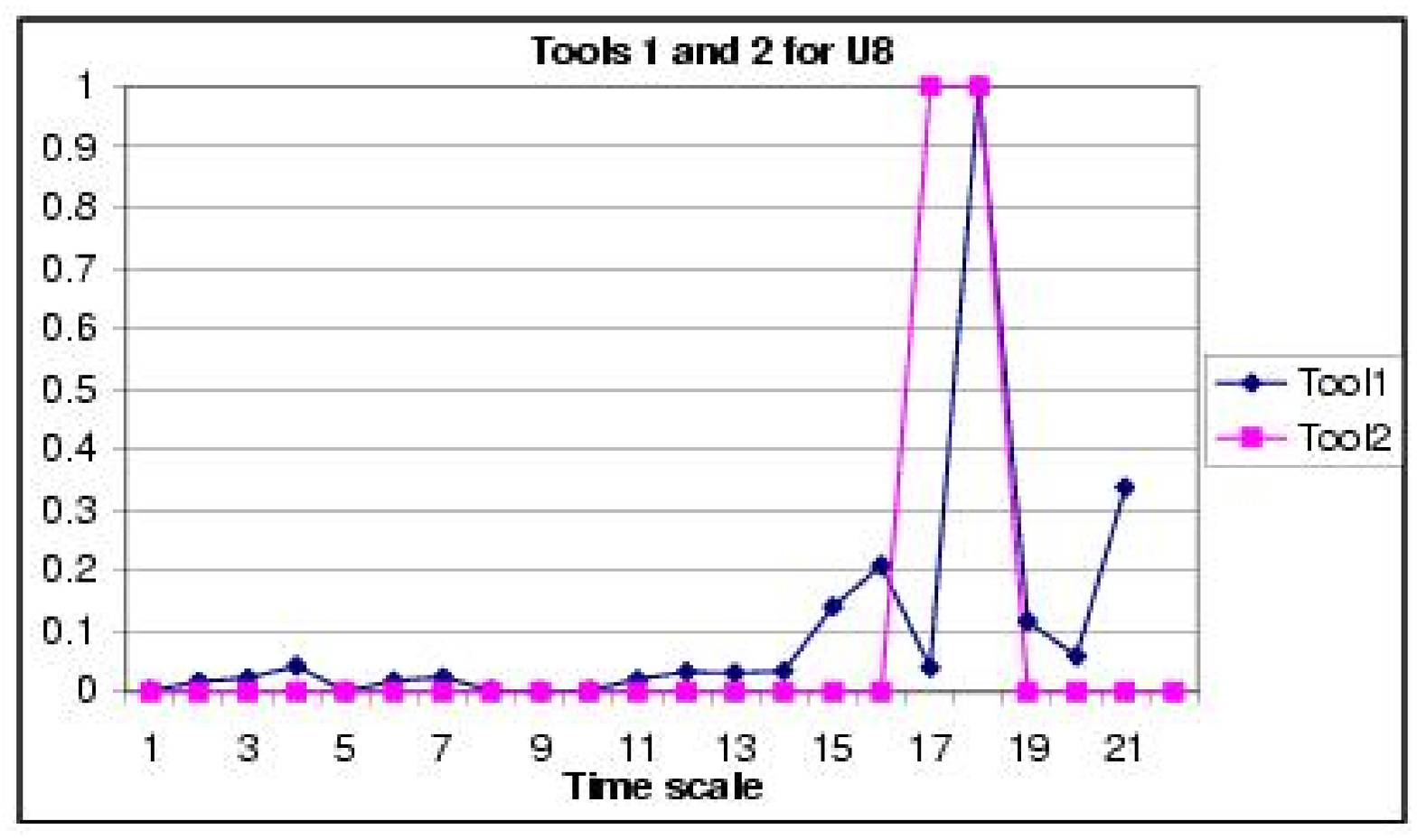}}
 \subfigure[]{\includegraphics[height=150pt, width=200pt]{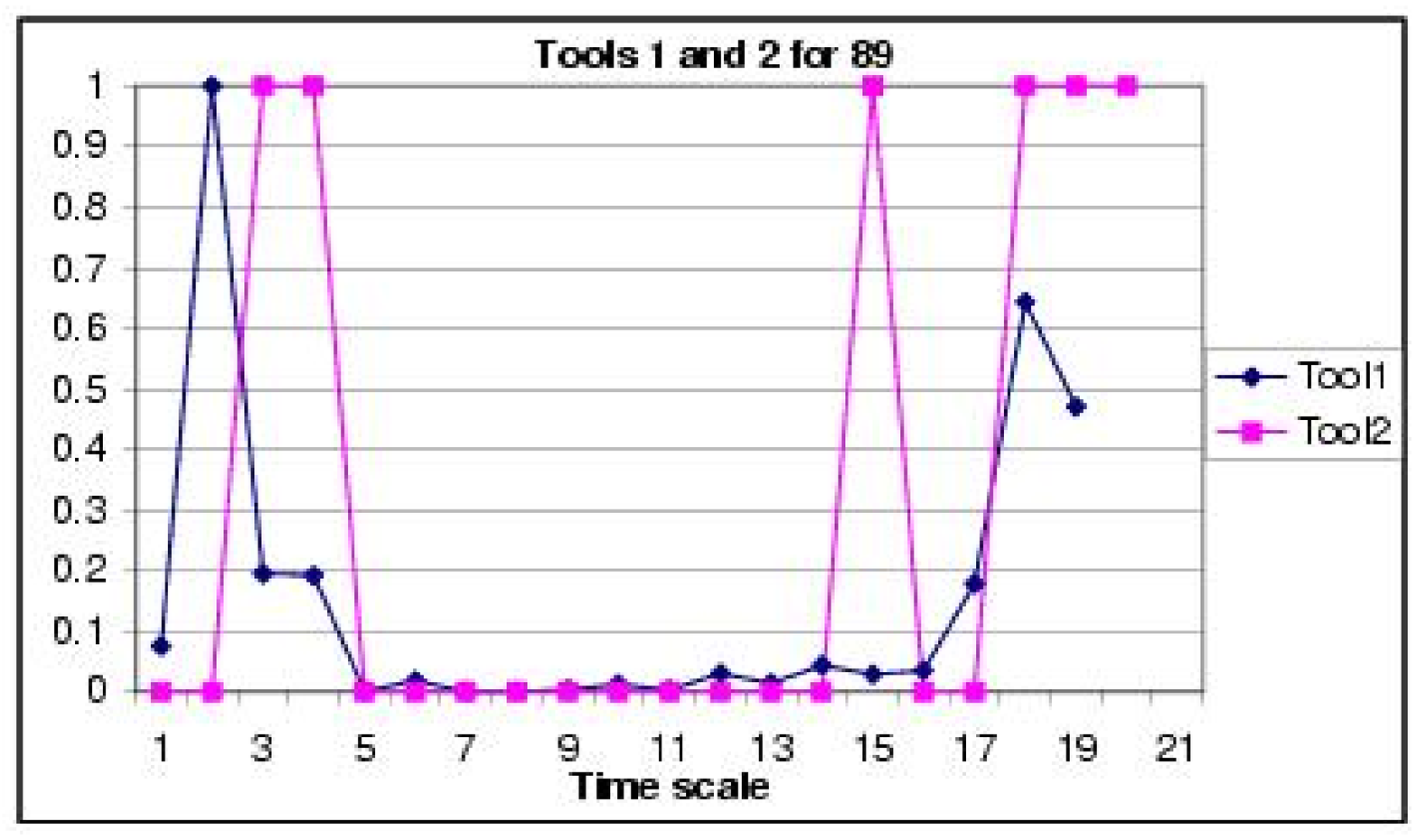}}
 \caption{Tools 1 and 2 applied to real traces. $\epsilon=0.01$ was used.}
 \label{ToolDataNew}
\end{figure}
\addtocounter{figure}{-1}
\stepcounter{figure}

Tools 1 and 2 seem to reveal more or less the same information, as to where the
levels are, but Tool 1 is a bit more detailed than Tool 2, in showing how
``prominent'' levels are (Fig. \ref{ToolData} and \ref{ToolDataNew}). The accuracy of the tools was tested by
their results on simulated data (Fig. \ref{ToolLev}): they nicely match the curvature of the
Averaging function for real traces, and they perform reasonably well with the simulations, where 
levels were artificially induced. As for tools 3 and 4, they reveal that trace 97 is the least spiky trace.
This is in agreement with visual inspection (Fig. \ref{Traffic9497} and \ref{TrafficM6U8}).

Before closing this section, let us submit Tool 2 to a final (and hard) test: can it detect correctly the levels of the trace 94 session presented in Fig. \ref{Lev} and \ref{LevelReaderReal2}? We saw earlier in Fig. \ref{LevelReaderReal2} in section \ref{quant} that there seems to be a very prominent level with $0$-intervals around $2^{13}ms$ (actually extending from $2^{10}ms$ to $2^{16}ms$), which nonetheless does not seem to affect the Averaging function much. Later, in section \ref{UneqLev}, we explained why this is the case. But is it then still possible for Tool 2 to detect the ``traces'' of this level on the Averaging function, however weak these may be? Actually it can, as Fig. \ref{ussest2} demonstrates. Tool 2 detects levels near $2^7ms$, $2^{11}ms$, and $2^{16}ms-2^{19}ms$, which are very close to the ones detected by the IDA in Fig. \ref{LevelReaderReal2}: in particular, the one detected at $2^{11}ms$ corresponds to the IDA-detected level in question.    

\clearpage

\begin{figure}[t]
 \centering
 \subfigure[]{\includegraphics[height=150pt, width=200pt]{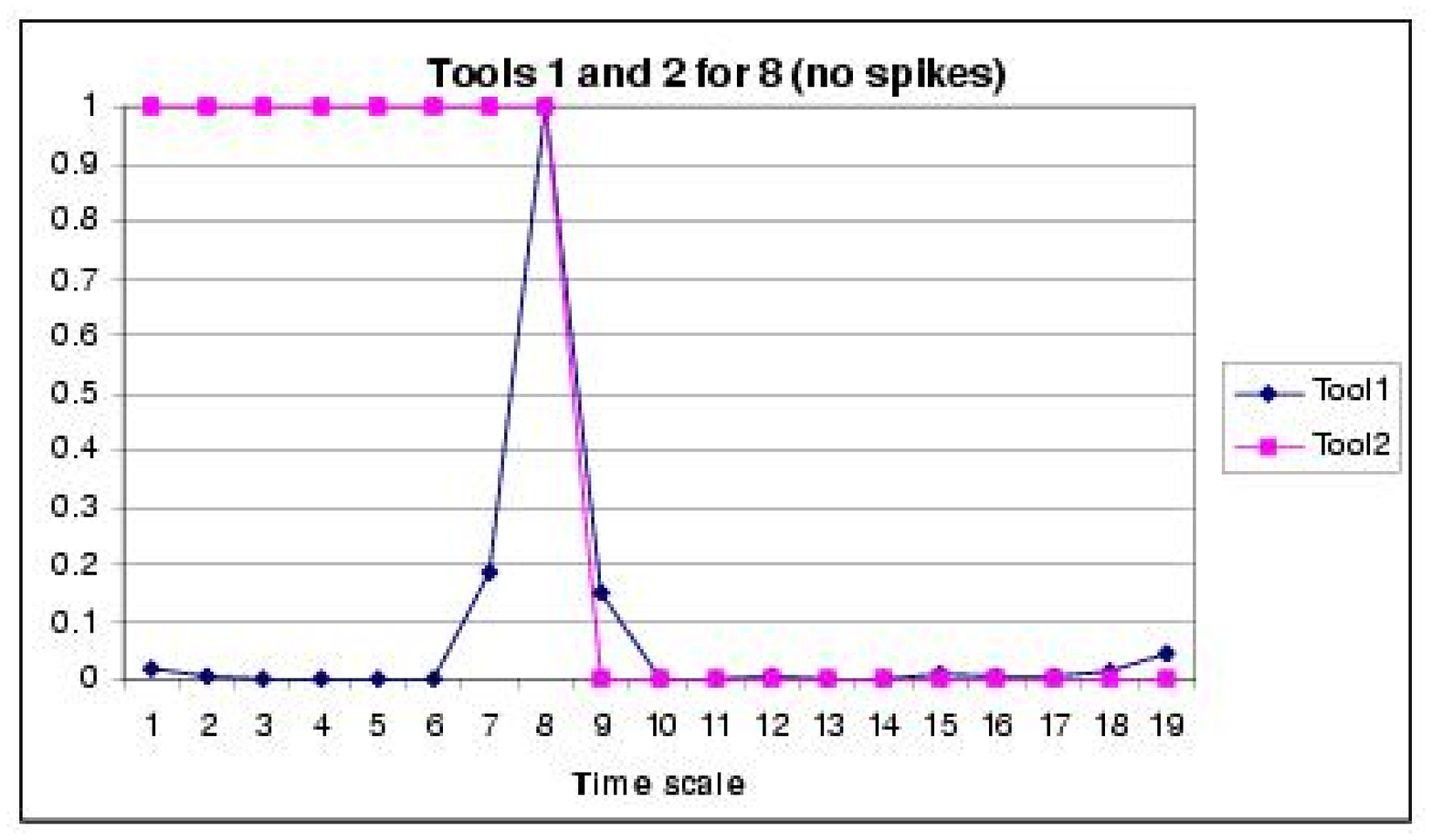}}
 \subfigure[]{\includegraphics[height=150pt, width=200pt]{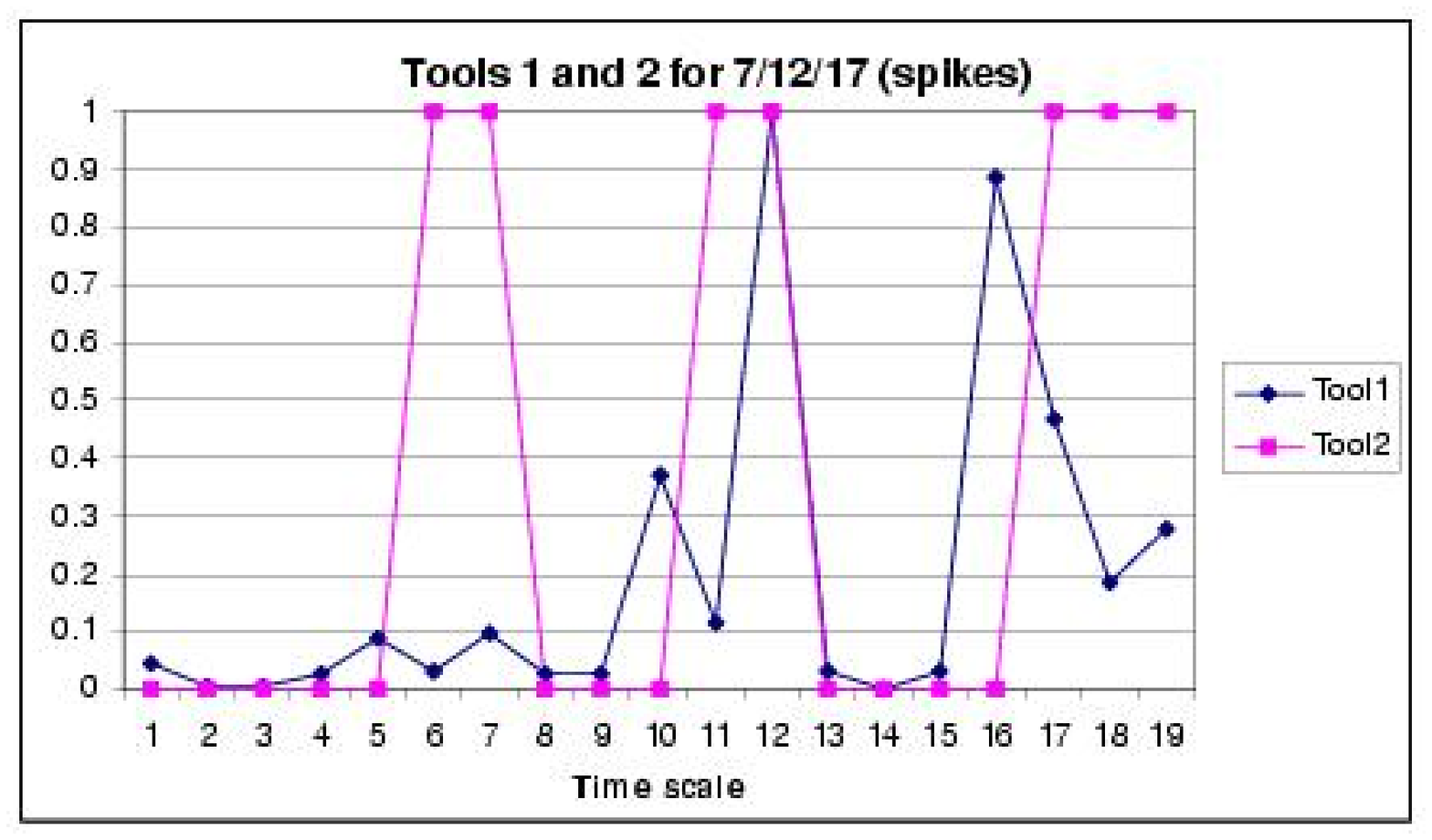}}
 \caption{Tools 1 and 2 applied to two model D simulations: 8 (No spikes) 
          and 7/12/17 (Spikes)}
 \label{ToolLev}         
\end{figure}
\addtocounter{figure}{-1}
\stepcounter{figure}

\begin{figure}[t]
 \centering
 \includegraphics[height=150pt, width=200pt]{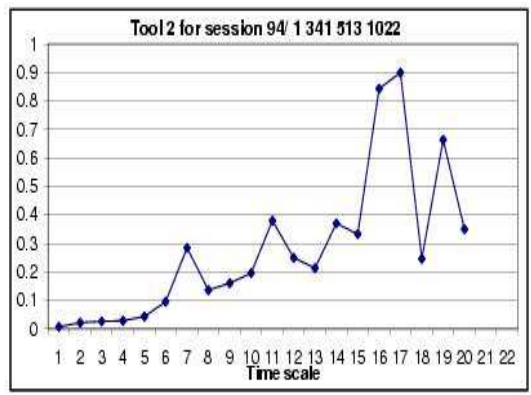}
 \caption{Tool 2 applied on session 94/ 1 341 514 1022: it clearly identifies levels near $2^7ms$, $2^{11}ms$, and $2^{16}ms-2^{19}ms$}.
 \label{ussest2}         
\end{figure}
\addtocounter{figure}{-1}
\stepcounter{figure}

\begin{table}
\begin{center} 
{Results of Tools 3 and 4 for real traces:} \\ \bigskip
\begin{tabular}{||c|c|c|c|c|c|c|c|c||} \hline
\emph{Trace} & 94 & 97 & L4 & L5 & M6 & M7 & U8 & 89 \\ \hline
\emph{Tool 3} & 0.274 & 0.084 & 0.272 & 0.289 & 0.016 & 0.017 & 0.145 & 0.125\\ \hline
\emph{Tool 4} & -0.479 & 0.101 & -0.541 & -0.425 & -0.084 & -0.123 & -0.340 & -0.797\\ \hline
\end{tabular}
\end{center}
\caption{\label{Tool34}The results of Tools 3 and 4 when applied to the real traces: notice that both of them indicate
         that trace 97 is very close to Gaussianity.}
\end{table}

\section{Acknowledgements}

The authors would like to thank Dr.\ Walter Willinger and Dr.\ Andy Ogielski of AT\&T for their comments, ideas, and support, Anja Feldmann of the Computer Science Department, University of Saarbr\"ucken, for the procurence of the data set 97, and Prof.\ Ingrid Daubechies of the Department of Mathematics, Princeton University, for contributing to the shaping of the key ideas, and improving decisively this treatise both \ae sthetically and scientifically, working as hard as the authors themselves did. K.\ Drakakis would also like to thank the Lilian Boudouris Foundation for the scholarship it granted him. Finally, the authors would like to thank NSF, AFOSR, ONR, DARPA, and AT\&T for their financial support.


\bigskip

\begin{thebibliography}{10}

\bibitem[1]{PF1}  V. Paxson, S. Floyd. \textit{Wide-Area Traffic: The
Failure of Poisson Modelling}\ \ IEEE/ACM Transactions on Networking, Vol. 3 No. 3, pp. 226-244, June 1995

\bibitem[2]{ENW1}  A. Erramilli, O. Narayan, W. Willinger. \textit{%
Experimental Queueing Analysis with Long-Range Dependent Packet Traffic}\ \ IEEE 1996

\bibitem[3]{RRCB1}  V. Ribeiro, R. Riedi, M. Crouse, R. Baraniuk. \textit{%
Simulation of non-Gaussian Long-Range-Dependent Traffic using Wavelets}\ \ ACM 1999 

\bibitem[4]{WTSW1}  W. Willinger, M. Taqqu, R. Sherman, D. Wilson. \textit{%
Self-Similarity Through High-Variability: Statistical Analysis of Ethernet
LAN Traffic at the Source Level}\ \ IEEE/ACM Transactions on Networking, Vol. 5, No.1,
February 1997

\bibitem[5]{TTW1}  M. Taqqu, V. Teverovsky, W. Willinger. \textit{Is network
traffic self-similar or multifractal?}\ \ Fractals. \textbf{5} (1997) 63-73

\bibitem[6]{ZPS1}  Y. Zhang, V. Paxson, S. Shenker. \textit{The
Stationarity of Internet Path Properties: Routing, Loss, and Throughput}\ \ 
ACIRI Technical Report, May 2000

\bibitem[7]{RE1}  B. Ryu, A. Elwalid. \textit{The importance of Long-Range
Dependence of VBR Video Traffic in ATM\ Traffic Engineering:\ Myths and
Realities}\ \ Proc. ACM SIGCOMM '96, Stanford University, CA, Aug. 1996

\bibitem[8]{NSSW1}  C. Nuzman, I. Saniee, W. Sweldens, A. Weiss. \textit{A
compound Model for TCP Connection Arrivals}\ \ February 2000 (work in progress)

\bibitem[9]{CCLS1}  J. Cao, W. Cleveland, D. Lin, D. Sun. \textit{On the
nonstationarity of Internet Traffic}\ \ Performance Evaluation Review: Proc. ACM Sigmetrics 2001, 29:102-112

\bibitem[10]{LT1}  J. Levy, M. Taqqu. \textit{Renewal-reward processes with
heavy-tailed interrenewal times and heavy-tailed rewards}\ \ Bernoulli. \textbf{6} (2000) 23-44

\bibitem[11]{PT1}  V. Pipiras, M. Taqqu. \textit{The limit of a
renewal-reward process with heavy-tailed rewards is not a linear fractional
stable motion}\ \ Bernoulli. \textbf{6} (2000) 607-614

\bibitem[12]{VB1}  A. Veres, M. Boda. \textit{The Chaotic Nature of TCP
Congestion Control}\ \ Proc. IEEE INFOCOM 2000, Tel Aviv, March 2000 

\bibitem[13]{CL1}  M. Crovella, C. Lindermann. \textit{Internet Performance
Modeling:\ The State of the Art at the Turn of the Century}\ \ Performance Evaluation, March 2000

\bibitem[14]{FP1}  V. Paxson, S. Floyd. \textit{Why We Don't Know How To
Simulate The Internet}\ \ Proceedings of the 1997 Winter Simulation Conference, December 1997

\bibitem[15]{FGW1}  A. Feldmann, A. Gilbert, W. Willinger. \textit{%
Investigating the multifractal nature of Internet WAN traffic}\ \ Computer Communication Review, Vol.\ 28. No.\ 4, 
Proceedings of the ACM/SIGCOMM'98, September 1998, Vancouver, Canada pp. 12--55, 1998 

\bibitem[16]{FGHW1}  A. Feldmann, A. Gilbert, P. Huang, W. Willinger. 
\textit{Dynamics of IP traffic: A study of the role of variability and the
impact of control}\ \ Proceedings of the ACM/SIGCOMM'99, August 29--September 1, 1999, Cambridge, MA

\bibitem[17]{PKC1}  K. Park, G. Kim, M. Crovella. \textit{On the
relationship between file sizes, transport protocols, and self-similar
network traffic}\ \ ICNP '96

\bibitem[18]{FGWK1}  A. Feldmann, A. Gilbert, W. Willinger, T. Kurz. \textit{%
The changing nature of network traffic: Scaling phenomena}\ \ Computer Communication Review, April 1998

\bibitem[19]{LTWW1}  W. Leland, M.S. Taqqu, W. Willinger, D. Wilson. \textit{%
On the self-similar nature of Ethernet traffic (extended version)}\ \ ACM/SIGCOMM'93. 
Computer Communication Review, \textbf{23} (1993), 183-193

\bibitem[20]{RLV1}  R. Riedi, J. L\'evy-V\'ehel. \textit{Multifractal 
Properties of TCP Traffic: a Numerical Study}\ \ INRIA research report 3129, February 1997

\bibitem[21]{R1}  R. Riedi. \textit{Introduction to Multifractals}\ \ Long range dependence: theory and applications, 
eds. Doukhan, Oppenheim and Taqqu (to appear 2001) 

\bibitem[22]{R2}  R. Riedi. \textit{Multifractal Processes}\ \ Long range dependence: theory and applications, 
eds. Doukhan, Oppenheim and Taqqu (2000)

\bibitem[23]{SRB1}  S. Sarvotham, R. Riedi, R. Baraniuk. \textit{%
Network Traffic Analysis and Modeling at the Connection Level}\ \ 
Technical Report, ECE Dept., Rice University, June 30, 2001 

\bibitem[24]{TWS1}  M. Taqqu, W. Willinger, R. Sherman. \textit{Proof of a
Fundamental Result in Self-Similar Traffic Modelling}\ \ Computer Communication Review, \textbf{27} (1997) 5-23

\bibitem[25]{PD1}  L. Peterson, B. Davie. \textit{Computer Networks:\ A
System Spproach (Second Edition)}\ \ Morgan--Kaufmann, 2000

\bibitem[26]{Gine}  A. Araujo, E. Gin\'e. \textit{The central Limit theorem for
Real and Banach Valued Random Variables}\ \ Wiley, New York (1980)

\bibitem[27]{Pollard}  D. Pollard \textit{Convergence of Stochastic Processes}\ \ Springer,
New York (1984)

\bibitem[28]{Tallagrand}  M. Ledoux, M. Talagrand. \textit{Probability in
Banach Spaces}\ \ Springer--Verlag (1991)

\bibitem[29]{MRRS1}  T. Mikosch, S.I. Resnick, H. Rootz\'en, A. Stegeman. (2001) \textit{Is network traffic approximated by stable Levy motion or fractional Brownian motion?} \ \ Ann. Appl. Probab. (2001), to appear. 


\end{thebibliography}
\end{document}